\documentclass[DIV=12]{scrartcl}

\usepackage{amsmath}
 \numberwithin{equation}{section}
\usepackage{amssymb}
\usepackage{amsthm}

\theoremstyle{definition}
 
\theoremstyle{remark}
 
\usepackage[main=english,german,french]{babel}
\usepackage[noisbn,noissn,fixlanguage]{babelbib}%
 \setbtxfallbacklanguage{english}
 
\usepackage{cite}%
\usepackage[T1]{fontenc}
\usepackage{hyperref}
\usepackage{lmodern}
\usepackage{mathscinet}
\usepackage{mathtools}
 \mathtoolsset{mathic=true}
\usepackage{microtype}
\usepackage{pmat}
\usepackage{url}

\usepackage{183arxiv}

\title{A Schur--Nevanlinna type algorithm for the truncated matricial Hausdorff moment problem}
\author{Bernd Fritzsche \and Bernd Kirstein \and Conrad M\"adler}
\date{}

\begin{document}
\maketitle

\begin{abstract}
 The main goal of this paper is to achieve a parametrization of the solution set of the truncated matricial \tHausdorff{} moment problem in the non-degenerate and degenerate situation.
 We treat the even and the odd cases simultaneously.
 Our approach is based on \tSchur{} analysis methods.
 More precisely, we use two interrelated versions of \tSchur{}-type algorithms, namely an algebraic one and a function-theoretic one.
 The algebraic version, worked out in our former paper \zita{MR4051874}, is an algorithm which is applied to finite or infinite sequences of complex matrices.
 The construction and discussion of the function-theoretic version is a central theme of this paper.
 This leads us to a complete description via \tStieltjes{} transform of the solution set of the moment problem under consideration.
 Furthermore, we discuss special solutions in detail.
\end{abstract}

\begin{description}
 \item[Mathematics Subject Classification (2010)] 44A60, 47A57, 30E05.
 \item[Keywords:] Truncated matricial Hausdorff moment problem, Schur--Nevanlinna type algorithm, parametrization of the solution set via Stieltjes transform.
\end{description}

\section{Introduction}

 The main aim of this paper is to work out an algorithm of \tSchur{}--\tNevanlinna{} type for functions, which leads via \tStieltjes{} transforms to a full description of the solution set of the matricial \tHausdorff{} moment problem in the general case.
 In order to realize this goal, we use our former investigations in \zitas{MR3775449,MR3979701,MR4051874} on matricial \tHausdorff{} moment sequences, where such sequences were studied from the point of view of \tSchur{} analysis and an algebraic version of a corresponding \tSchur{}--\tNevanlinna{} type algorithm was worked out.
 The synthesis of these two \tSchur{}--\tNevanlinna{} type algorithms will provide us finally the desired result.

 This strategy was already used by the authors to study the matricial \tHamburger{} moment problem (see \zitas{MR3014199,MR3380267}) and the matricial \tStieltjes{} moment problem (see \zitas{MR3611479,MR3611471,MR3765778}).
 
 It will turn out that the case of the matricial \tHausdorff{} moment problem is much more difficult in comparison with the matricial problems named after \tHamburger{} and \tStieltjes{}.
 This phenomenon could be already observed in the discussion of the so-called \tnd{} case, where the \tHamburger{} problem (see \zitas{MR703593,MR1018213}) and the \tStieltjes{} problem (see \zitas{MR645305,MR752057,MR686076,MR2053150}) were studied considerably earlier than the \tHausdorff{} problem (see \zitas{MR2222521,MR2342899}). 
 The main reason for the greater complexity of the \tHausdorff{} moment problem is caused by the fact that the localization of the measure in a prescribed compact interval of the real axis requires to satisfy simultaneously more conditions.
 This implies that the possible moment sequences have a more complicated structure (see \zitas{MR3775449,MR3979701,MR4051874}).
 
 Continuing the work done in \zitas{CR01,MR2342899,MR2222521} A.~E.~Choque-Rivero \zitas{MR3745463,MR3704346,MR3377989,CR17} investigated further aspects of the \tnd{} truncated matricial \tHausdorff{} moment problem.
 As in \zitas{CR01,MR2342899,MR2222521} he distinguished between the case of an odd or even number of prescribed matricial moments.
 The approach used in \zitas{CR01,MR2342899,MR2222521} is based on V.~P.~Potapov's method of Fundamental Matrix Inequalities.
 In the scalar case \(q=1\) the classical method used in \zita{MR0044591} or \zitaa{MR0458081}{\cch{4}, \S~7} is based on the application of orthogonal polynomials.
 A.~E.~Choque-Rivero \zita{MR3377989} obtained a matrix generalization of the Krein--Nudelman representation of the resolvent matrix by using four families of orthogonal matrix polynomials on the interval \(\ab\).
 From the point of view of V.~P.~Potapov the \tSchur{} algorithm is interpreted as a multiplicative decomposition of a \(J\)\nobreakdash-contractive matrix function into simplest elementary factors.
 Such multiplicative decompositions were constructed for the resolvent matrix for the moment problem under study by A.~E.~Choque-Rivero in \zitas{MR3704346,CR17}.
 
 An important feature of this paper is to achieve a simultaneous treatment of the even and odd truncated matricial \tHausdorff{} moment problems in the general case.
 Our strategy is based on the application of \tSchur{} analysis methods.
 
 This paper is organized as follows.
 In \rsec{S.MP}, we state some general facts on matricial power moment problems on \tBorel{} subsets of the real axis.
 In \rsec{S1506}, we summarize essential insights about the structure of matricial \(\ab\)\nobreakdash-\tHausdorff{} moment sequences, which were mostly obtained in our former papers \zitas{MR2570113,MR645305,MR752057}.
 A key observation is that we do not treat the original matricial moment problem, but an equivalent problem in the class \(\RFqab \) of holomorphic matrix functions.
 In \rsec{F.s2.RFab}, we summarize some facts on the class \(\RFqab\), which are needed in the sequel.
 In \rsec{S1723}, we formulate a problem for functions of \(\RFqab\) which is equivalent to the original matricial moment problem.
 This equivalence is caused by \tSTF{}.
 In \zita{MR3979701}, we constructed a \tSchur{}--\tNevanlinna{} type algorithm for \tnnH{} \tqqa{measures} on \(\ab\) by translating the \tSchur{}--\tNevanlinna{} type algorithm for matricial \(\ab\)\nobreakdash-\tHausdorff{} moment sequences into the language of measures.
 In \rsec{S1323}, we translate now this algorithm via \tSTF{} into the class \(\RFqab\).
 On this way our main goal is to achieve a description of all solutions of the truncated matricial \(\ab\)\nobreakdash-\tHausdorff{} moment problem via a linear fractional transformation of matrices.
 This requires to find the generating matrix function of this transformation and the corresponding domain.
 In the first step, we concentrate on that domain.
 Remember that in \zita{MR2222521} we already studied the problem under consideration in the \tnd{} case, however by use of Potapov's method of fundamental matrix inequalities.
 Doing this, we were led to a particular class \(\PRFabq\) of ordered pairs of \tqqa{matrix}-valued functions which are meromorphic in \(\Cab\) (see \zitaa{MR2222521}{\cdefn{5.2}{}}).
 In \rsec{F.s2.PRF}, we summarize some important facts about the class \(\PRFabq\) and present an example of a remarkable element of this class (see \rexam{E1410}).
 The experiences from \zitas{MR2222521,MR2342899} teach us that it is necessary to introduce an equivalence relation within \(\PRFabq\).
 In order to take into account possible degeneracies of the moment problem, we have to single out an appropriate subclass of \(\PRFabq\), which is adapted to the prescribed matricial moments.
 Furthermore, we have to ensure that the construction of this subclass stands in harmony with the above mentioned equivalence relation in \(\PRFabq\).
 The just mentioned two themes are treated in \rsec{F.S.PM}. 
 The goal of the following considerations is to prepare basic instruments for the version of our \tSchur{}--\tNevanlinna{} type algorithm for functions.
 This algorithm should stand in correspondence with the \tSchur{}--\tNevanlinna{} type algorithm for matricial \tFnnd{} sequences, which was developed in \zita{MR4051874}.
 A remarkable feature of this version is that this algorithm contains two elements of different nature.
 More precisely, the first step of the algorithm differs from the remaining steps.
 There occur (equivalence classes) of ordered pairs of matrix-valued functions in the first step, whereas the further steps require only matrix-valued functions. 
 In \rsecss{S1229}{S1230}, we study the corresponding transformations and its inverses for the first and remaining steps, \tresp{}
 The transformations are defined by using \tMoore{}--\tPenrose{} inverses of matrices.
 It will turn out that under special conditions, which will be indeed satisfied in the case of interest for us, these transformations can be rewritten as usual linear fractional transformations of matrices the generating matrix-valued functions of which are quadratic \taaa{2q}{2q}{matrix} polynomials. 
 Having a closer look at the considerations in \rsecss{S1229}{S1230} one can observe that the basic tools used there are of rather algebraic nature.
 In \rsecss{S1531}{S1141}, we study the elementary steps of the forward and backward algorithm in more detail, \tresp{}
 Namely, we turn our attention to the concrete classes of meromorphic matrix-valued functions occurring there.
 Moreover, we demonstrate that these elementary steps of the algorithm for functions are concordant with the elementary steps of the algebraic algorithm applied to the matricial moment sequences. 
 In \rsec{S1245}, we check (see \rthm{F.T.Phi}) that the iteration of the elementary steps leads to a parametrization of the set of \tSTF{s} of all solutions of the original matricial moment problem. 
 In \rsec{F.S.glt}, we translate \rthm{F.T.Phi} into the language of linear fractional transformations of matrices and obtain our main results \rthmss{F.P.PVR}{F.P.RFabmred}.
 
 In \rsec{S0948}, via \tSTF{} we determine all those solutions \(\sigma\) of the moment problem associated with a sequence \(\seqs{m}\) for which the sequence \(\seqs{m+1}\), where \(\su{m+1}\) is the \(\rk{m+1}\)\nobreakdash-th power moment of \(\sigma\), is \tabHd{}.
 
 The main object of study in \rsec{S1014} is that solution of the moment problem associated with a sequence \(\seqs{m}\) which corresponds to the central extension of \(\seqs{m}\).
 We determine the position of the \tSTF{} of this solution within the general parametrization obtained in \rthm{F.P.PVR}.
 
 In several appendices we summarize some needed facts on various topics such as particular aspects of matrix theory, \tnnH{} measures and corresponding integration theory, \tStieltjes{}  of \tnnH{} measures on the real line, ordered pairs of matrices corresponding to linear relations, linear fractional transformations of matrices, holomorphic matrix-valued functions.

\section{Matricial moment problems on Borel subsets of the real axis}\label{S.MP}
 In this section, we are going to formulate a certain class of matricial power moment problems.
 Before doing this, we have to introduce some terminology.
 We denote by \(\Z\) the set of all integers.
 Let \(\N\defeq\setaca{n\in\Z}{n\geq1}\).
 Furthermore, we write \(\R\) for the set of all real numbers and \(\C\) for the set of all complex numbers.
 In the whole paper, \(p\) and \(q\) are arbitrarily fixed integers from \(\N\).
 We write \(\Cpq\) for the set of all complex \tpqa{matrices} and \(\Cp\) is short for \(\Coo{p}{1}\).
 When using \(m,n,r,s,\dotsc\) instead of \(p,q\) in this context, we always assume that these are integers from \(\N\).
 We write \(A^\ad\) for the conjugate transpose of a complex \tpqa{matrix} \(A\).
 Denote by \(\Cggq\) (\tresp{}\ \(\Cgq\)) the set of all complex \tnnrpH{} \tqqa{matrices}.
 If \(\mathcal{X}\) is a subset of \(\Coo{q}{r}\) and if \(A\in\Cpq\), then let \(A\mathcal{X}\defeq\setaca{AX}{X\in\mathcal{X}}\).

 Let \((\mathcal{X},\mathfrak{X})\) be a measurable space.
 Each countably additive mapping whose domain is \(\mathfrak{X}\) and whose values belong to \(\Cggq\) is called a \tnnH{} \tqqa{measure} on \((\mathcal{X}, \mathfrak{X})\).
 For the integration theory with respect to \tnnH{} measures, we refer to Kats~\zita{MR0080280} and Rosenberg~\zita{MR0163346}.
 For the convenience of the reader, a summary concerning this matter, sufficient for our purposes, is given in \rapp{B.se.IT}.

 Let \(\BsAR\) (\tresp{}\ \(\BsAC \)) be the \(\sigma\)\nobreakdash-algebra of all Borel subsets of \(\R\) (\tresp{}\ \(\C \)).
 In the whole paper, \(\Omega\) stands for a non-empty set belonging to \(\BsAR\).
 Let \(\BsAO \) be the \(\sigma\)\nobreakdash-algebra of all Borel subsets of \(\Omega\) and let \(\MggqO \) be the set of all \tnnH{} \tqqa{measures} on \((\Omega,\BsAO )\).
 Observe that \(\Mggoa{1}{\Omega}\) coincides with the set of ordinary measures on \((\Omega,\BsAu{\Omega})\) with values in \([0,\infp)\).

 Let \(\NO\defeq\setaca{m\in\Z}{m\geq0}\).
 Throughout this paper, \(\kappa\) is either an integer from \(\NO\) or \(\infi\).
 In the latter case, we have \(2\kappa\defeq\infi\).
 Given \(\upsilon,\omega\in\R\cup\set{-\infty,\infp}\), we set \(\mn{\upsilon}{\omega}\defeq\setaca{k\in\Z}{\upsilon\leq k\leq\omega}\).
 Let \(\Mggoua{q}{\kappa}{\Omega}\) be the set of all \(\mu\in\MggqO \) such that for each \(j\in\mn{0}{\kappa}\) the power function \(x\mapsto x^j\) defined on \(\Omega\) is integrable with respect to \(\mu\).
 If \(\mu\in\MggquO{\kappa}\), then, for all \(j\in\mn{0}{\kappa}\), the matrix
\beql{I.G.mom}
 \mpm{\mu}{j}
 \defeq\int_\Omega x^j\mu\rk{\dif x}
\eeq
 is called (power) \emph{moment of \(\mu\)} of order \(j\).
 Obviously, we have \(\Mggoa{q}{\Omega}=\Mggoua{q}{0}{\Omega}\subseteq \MggquO{\ell}\subseteq \MggquO{\ell+1}\subseteq\MggquO{\infi}\) for every choice of \(\ell\in\NO \) and, furthermore, \(\mpm{\mu}{0}=\mu\rk{\Omega}\) for all \(\mu\in\Mggoa{q}{\Omega}\).
 If \(\Omega\) is bounded, then one can easily see that \(\MggqO =\MggquO{\infi}\).
 We now state the general form of the moment problem lying in the background of our considerations:

\begin{Problem}[\mprob{\Omega}{\kappa}{=}]
 Given a sequence \(\seqska\) of complex \tqqa{matrices}, parametrize the set \(\MggqOsg{\kappa}\) of all \(\sigma\in\Mggoua{q}{\kappa}{\Omega}\) satisfying \(\mpm{\sigma}{j}=\su{j}\) for all \(j\in\mn{0}{\kappa}\).
\end{Problem}

 In the whole paper, let \(\ug\) and \(\obg\) be two arbitrarily given real numbers satisfying \(\ug<\obg\) and let \(\ba\defeq\obg-\ug\).
 In what follows, we mainly consider the case that \(\Omega\) is the compact interval \(\ab\) of the real axis \(\R\).
 As mentioned above, we have \(\MggqF=\MggquF{\infi}\).
 
 Since each solution of~\mprob{\ab}{\kappa}{=} generates in a natural way solutions of~\mprob{\rhl}{\kappa}{=},~\mprob{\lhl}{\kappa}{=}, and~\mprob{\R}{\kappa}{=}, we will also use results concerning the treatment of these moment problems.

 \section{Matricial \(\ab\)-Hausdorff moment sequences}\label{S1506}

 In this section, we recall a collection of results on the matricial Hausdorff moment problem and corresponding moment sequences of \tnnH{} \tqqa{measures} on the interval \(\ab\), which are mostly taken from~\zitas{MR3775449,MR3979701,MR4051874}.
 To state a solvability criterion, we introduce the relevant class of sequences of complex matrices.

\bnotal{N.HKG}
 Let \(\seqska \) be a sequence of complex \tpqa{matrices}.
 Then let the \tbHms{} \(\Hu{n}\), \(\Ku{n}\), and \(\Gu{n}\) be given by \(\Hu{n}\defeq\matauuuo{s_{j+k}}{j}{k}{0}{n}\) for all \(n\in\NO\) with \(2n\leq\kappa\), by \(\Ku{n}\defeq\matauuuo{s_{j+k+1}}{j}{k}{0}{n}\) for all \(n\in\NO\) with \(2n+1\leq\kappa\), and by \(\Gu{n}\defeq\matauuuo{s_{j+k+2}}{j}{k}{0}{n}\) for all \(n\in\NO\) with \(2n+2\leq\kappa\), \tresp{}
\enota

 To emphasize that a certain (block) matrix \(X\) is built from a sequence \(\seqska\), we sometimes write \(X^{\ok{s}}\) for \(X\).

\bnotal{F.N.sa}
 Assume \(\kappa\geq1\) and let \(\seqska\) be a sequence of complex \tpqa{matrices}.
 Then let the sequences \(\seqsa{\kappa-1}\) and \(\seqsb{\kappa-1}\) be given by
\begin{align*}
 \sau{j}& \defeq-\ug\su{j}+\su{j+1}&
&\text{and}&
 \sub{j}& \defeq\obg\su{j}-\su{j+1},
\end{align*}
 \tresp{}
 Furthermore, if \(\kappa\geq2\), then let the sequence \(\seqsab{\kappa-2}\) be given by
\[
 \sab{j}
 \defeq-\ug\obg\su{j}+(\ug+\obg)\su{j+1}-\su{j+2}.
\]
 Similarly to our more algebraic considerations \zitas{MR3775449,MR3979701,MR4051874} it will turn out that a characteristic feature of~\mprob{\ab}{\kappa}{=} is to analyze and organize the interplay between the four matrix sequences \(\seqska\), \(\seqsa{\kappa-1}\), \(\seqsb{\kappa-1}\) and \(\seqsab{\kappa-2}\) under the view of corresponding \tbHms{} generated by them.
 For each matrix \(X_k=X_k^{\ok{s}}\) built from the sequence \(\seqska\), we denote (if possible) by \(X_\aur{k}\defeq X_k^{\ok{\saus}}\), by \(X_\lub{k}\defeq X_k^{\ok{\subs}}\), and by \(X_\aub{k}\defeq X_k^{\ok{\sabs}}\) the corresponding matrix built from the sequences \(\seqsa{\kappa-1}\), \(\seqsb{\kappa-1}\), and \(\seqsab{\kappa-2}\) instead of \(\seqska\), \tresp{}
\enota

 In view of \rnota{N.HKG}, we get in particular
\begin{align*}%
 \Hau{n}&=-\ug\Hu{n}+\Ku{n}&
 &\text{and}&
 \Hub{n}&=\obg\Hu{n}-\Ku{n}
\end{align*}
 for all \(n\in\NO\) with \(2n+1\leq\kappa\) and
\[%
 \Hab{n}
 =-\ug\obg\Hu{n}+(\ug+\obg)\Ku{n}-\Gu{n}
\]
 for all \(n\in\NO\) with \(2n+2\leq\kappa\).
 In the classical case \(\ug=0\) and \(\obg=1\), we have furthermore \(\sau{j}=\su{j+1}\) and \(\sub{j}=\su{j}-\su{j+1}\) for all \(j\in\mn{0}{\kappa-1}\) and \(\sab{j}=\su{j+1}-\su{j+2}\) for all \(j\in\mn{0}{\kappa-2}\).

\breml{F.R.012} 
 Let \(\seqska\) be a sequence of complex \tpqa{matrices}.
 Then \(\ba\su{j}=\sau{j}+\sub{j}\) and \(\ba\su{j+1}=\obg\sau{j}+\ug\sub{j}\) for all \(j\in\mn{0}{\kappa-1}\).
 Furthermore, \(\ba\su{j+2}=\obg^2\sau{j}+\ug^2\sub{j}-\ba\sab{j}\) for all \(j\in\mn{0}{\kappa-2}\).
\erem

\bdefnl{D0745}
 Let \(\Fggqu{0}\) (\tresp{}\ \(\Fgqu{0}\)) be the set of all sequences \(\seqs{0}\) of complex \tqqa{matrices}, for which the \tbHm{} \(\Hu{0}\) is \tnnrpH{}, \tie{}, for which \(\su{0}\in\Cggq\) (\tresp{}\ \(\su{0}\in\Cgq\)) holds true.
 For each \(n\in\N\), denote by \(\Fggqu{2n}\) (\tresp{}\ \(\Fgqu{2n}\)) the set of all sequences \(\seqs{2n}\) of complex \tqqa{matrices}, for which the \tbHms{} \(\Hu{n}\) and \(\Hab{n-1}\) are both \tnnrpH{}.
 For each \(n\in\NO\), denote by \(\Fggqu{2n+1}\) (\tresp{}\ \(\Fgqu{2n+1}\)) the set of all sequences \(\seqs{2n+1}\) of complex \tqqa{matrices} for which the \tbHms{} \(\Hau{n}\) and \(\Hub{n}\) are both \tnnrpH{}.
 Furthermore, denote by \(\Fggqinf\) (\tresp{}\ \(\Fgqinf\)) the set of all sequences \(\seqsinf \) of complex \tqqa{matrices} satisfying \(\seqs{m}\in\Fggqu{m}\) (\tresp{}\ \(\seqs{m}\in\Fgqu{m}\)) for all \(m\in\NO\).
 The sequences belonging to \(\Fggqu{0}\), \(\Fggqu{2n}\), \(\Fggqu{2n+1}\), or \(\Fggqinf\) (\tresp{}\  \(\Fgqu{0}\), \(\Fgqu{2n}\), \(\Fgqu{2n+1}\), or \(\Fgqinf\)) are said to be \emph{\tFnnd} (\tresp{}\ \emph{\tFpd}).
\edefn

 (Note that in~\zita{MR3775449}, the sequences belonging to \(\Fggqka \) were called \emph{\(\ab\)\nobreakdash-Hausdorff non-negative definite}.) 
 A necessary and sufficient condition for the solvability of \mprob{\ab}{\kappa}{=} is the following:

\bthmnl{\tcf{}~\zitaa{MR2222521}{\cthm{1.3}{127}} and~\zitaa{MR2342899}{\cthm{1.3}{106}}}{I.P.ab}
 Let \(\seqska\) be a sequence of complex \tqqa{matrices}.
 Then \(\MggqFksg\neq\emptyset\) if and only if \(\seqska\in\Fggqka\).
\ethm

 Since \(\Omega=\ab\) is bounded, one can easily see that \(\MggqF=\MggquF{\infi}\), \tie{}, the power moment \(\mpm{\sigma}{j}\) defined by \eqref{I.G.mom} exists for all \(j\in\NO\).
 If \(\sigma\in\MggqF\), then we call \(\seqmpm{\sigma}\) given by \eqref{I.G.mom} the \emph{\tfpmfa{\(\sigma\)}}.

 Given the complete sequence of prescribed power moments \(\seqsinf\), the moment problem on the compact interval \(\Omega=\ab\) differs from the moment problems on the unbounded sets \(\Omega=\rhl\) and \(\Omega=\R\) in having necessarily a unique solution, assumed that a solution exists:

\bpropl{I.P.ab8}
 If \(\seqsinf\in\Fggqinf\), then the set \(\MggqFsg{\infi}\) consists of exactly one element.
\eprop

 \rprop{I.P.ab8} is a well-known result, which can be proved, in view of \rthm{I.P.ab}, using the corresponding result in the scalar case \(q=1\) (see, \zita{MR1544592} or~\zitaa{MR0184042}{\cthm{2.6.4}{74}}).

 We can summarize \rprop{I.P.ab8} and \rthm{I.P.ab} for \(\kappa=\infi\):

\bpropl{I.P.ab8Fgg}
 The mapping \(\Xi_\ab\colon\MggqF\to\Fggqinf\) given by \(\sigma\mapsto\seqmpm{\sigma}\) is well defined and bijective.
\eprop

 For each \(n\in\NO\), denote by \(\Hggq{2n}\) the set of all sequences \(\seqs{2n}\) of complex \tqqa{matrices}, for which the corresponding \tbHm{} \(\Hu{n}\) is \tnnH{}. 
 Furthermore, for each \(\iota\in\NOinf\) and each non-empty set \(\mathcal{X}\), denote by \(\seqset{\iota}{\mathcal{X}}\) the set of all sequences \(\seq{X_j}{j}{0}{\iota}\) from \(\mathcal{X}\).
 Obviously, \(\Fggqu{0}\) coincides with the set of all sequences \((s_j)_{j=0}^0\) with \(\su{0}\in\Cggq\).
 Furthermore, we have
\beql{Fgg2n}
 \Fggqu{2n}
 =\setaca*{\seqs{2n}\in\Hggq{2n}}{\seq{\sab{j}}{j}{0}{2(n-1)}\in\Hggq{2(n-1)}}
\eeq
 for all \(n\in\N\) and
\beql{Fgg2n+1}
 \Fggqu{2n+1}
 =\setaca*{\seqs{2n+1}\in\seqset{2n+1}{\Cqq}}{\set*{\seqsa{2n},\seqsb{2n}}\subseteq\Hggq{2n}}
\eeq
 for all \(n\in\NO\). 
 Note that the following \rpropss{ab.R0933}{F.L.sabF}, which are proved in a purely algebraic way in~\zita{MR3775449}, can also be obtained immediately from \rthm{I.P.ab}.
 
\bpropnl{\tcf{}~\zitaa{MR3775449}{\cprop{7.7(a)}{18}}}{ab.R0933}
 If \(\seqska  \in \Fggqka \), then \(\seqs{m} \in \Fggqu{m}\) for all \(m\in \mn{0}{\kappa}\).
\eprop

 In view of \rprop{ab.R0933}, the definition of the class \(\Fggqinf\) seems to be natural.

\bpropnl{\zitaa{MR3775449}{\cprop{9.1}{27}}}{F.L.sabF}
 Let \(\seqska\in\Fggqka\).
 If \(\kappa\geq1\), then \(\set{\seqsa{\kappa-1},\seqsb{\kappa-1}}\subseteq\Fggqu{\kappa-1}\).
 If \(\kappa\geq2\), then furthermore \(\seqsab{\kappa-2}\in\Fggqu{\kappa-2}\).
\eprop

 We write \(\ran{A}\defeq\setaca{Ax}{x\in\Cq}\) and \(\nul{A}\defeq\setaca{x\in\Cq}{Ax=\Ouu{p}{1}}\) for the column space and the null space of a complex \tpqa{matrix} \(A\), \tresp{}
 Denote by \(\Dpqka\) the set of all sequences \(\seqska\) of complex \tpqa{matrices} satisfying
 \(\bigcup_{j=0}^\kappa\ran{s_j}\subseteq\ran{s_0}\) and \(\nul{s_0}\subseteq\bigcap_{j=0}^\kappa\nul{s_j}\).
 
\bremnl{\zitaa{MR4051874}{\cprop{7.11}{161}}}{F.R.Fgg<D} 
 Let \(\seqska\in\Fggqka\).
 Then \(\seqska\in\Dqqu{\kappa}\).
 If \(\kappa\geq1\), furthermore \(\set{\seqsa{\kappa-1},\seqsb{\kappa-1}}\subseteq\Dqqu{\kappa-1}\).
 If \(\kappa\geq2\), moreover \(\seqsab{\kappa-2}\in\Dqqu{\kappa-2}\).
\erem

 The set \(\CHq\defeq\setaca{M\in\Cqq}{M^\ad=M}\) of \tH{} matrices from \(\Cqq\) is a partially ordered \(\R\)\nobreakdash-vector space with positive cone \(\Cggq\).
 For two complex \tqqa{matrices} \(A\) and \(B\), we write \(A\lleq B\) or \(B\lgeq A\) if \(A,B\in\CHq\) and \(B-A\in\Cggq\) are fulfilled.
 This partial order \(\lleq\) on the set of \tH{} matrices is sometimes called \emph{L\"owner semi-ordering}.
 
\blemnl{\zitaa{MR3979701}{\clem{5.7}{2138}}}{F.R.Fgg-s}
 Let \(\seqska\in\Fggqka\).
 Then \(\su{j}\in\CHq\) for all \(j\in\mn{0}{\kappa}\) and \(\su{2k}\in\Cggq\) for all \(k\in\NO\) with \(2k\leq\kappa\).
 Furthermore, \(\ug\su{2k}\lleq\su{2k+1}\lleq\obg\su{2k}\) for all \(k\in\NO\) with \(2k+1\leq\kappa\).
\elem

 Let \(\Opq\) be the zero matrix from \(\Cpq\) and let \(\Iq\defeq\matauuo{\Kronu{jk}}{j,k}{1}{q}\) be the identity matrix from \(\Cqq\), where \(\Kronu{jk}\) is the Kronecker delta.
 Sometimes, if the size is clear from the context, we will omit the indices and write \(\NM\) and \(\EM\), \tresp{}
 Taking into account \rrem{A.R.rA<rB}, we obtain from \rlem{F.R.Fgg-s}:

\breml{F.R.Fgg-r}
 If \(\seqska\in\Fggqka\), then \(\ran{\sau{2k}}\cup\ran{\sub{2k}}\subseteq\ran{\su{2k}}\) and \(\nul{\su{2k}}\subseteq\nul{\sau{2k}}\cap\nul{\sub{2k}}\) for all \(k\in\NO\) with \(2k\leq\kappa-1\).
\erem

 Finite sequences from \(\Fggqu{m}\) can always be extended to sequences from \(\Fggqu{\ell}\) for all \(\ell\in\minf{m+1}\), which is due to the fact that a \tnnH{} measure on the bounded set \(\ab\) possesses power moments of all non-negative orders.
 One of the main results in~\zita{MR3775449} states that the possible one-step extensions \(\su{m+1}\in\Cqq\) of a sequence \(\seqs{m}\) to an \tFnnd{} sequence \(\seqs{m+1}\) fill out a matricial interval.
 In order to give an exact description of this interval, we are now going to introduce several matrices and recall their role in the corresponding extension problem for \tFnnd{} sequences, studied in~\zita{MR3775449}.

 Given an arbitrary \(n\in\N\) and arbitrary rectangular complex matrices \(A_1,A_2,\dotsc,A_n\), we write \(\col\seq{A_j}{j}{1}{n}=\col\rk{A_1,A_2,\dotsc,A_n}\) (\tresp{}, \(\row\seq{A_j}{j}{1}{n}\defeq\mat{A_1,A_2,\dotsc,A_n}\)) for the block column (\tresp{}, block row) built from the matrices \(A_1,A_2,\dotsc,A_n\) if their numbers of columns (\tresp{}, rows) are all equal.

\bnotal{N.yz}
 Let \(\seqska\) be a sequence of complex \tpqa{matrices}.
 For all \(\ell,m\in\NO\) with \(\ell\leq m\leq\kappa\), then let \(\yuu{\ell}{m}\defeq\col\seq{s_j}{j}{\ell}{m}\) and \(\zuu{\ell}{m}\defeq\row\seq{s_j}{j}{\ell}{m}\).
\enota

 The \tbHm{} \(\Hu{n}\) admits the following \tbr{s}:

\breml{H.R.Hblock}
 If \(\kappa\geq2\) and if \(\seqska\) is a sequence of complex \tpqa{matrices}, then \(\Hu{n}=\tmat{\Hu{n-1} & \yuu{n}{2n-1} \\\zuu{n}{2n-1} & s_{2n}}\) and \(\Hu{n}=\tmat{s_0&\zuu{1}{n}  \\\yuu{1}{n} & \Gu{n-1}}\) for all \(n\in\N\) with \(2n\leq\kappa\).
\erem

 In this paper, the Moore--Penrose inverse of a complex matrix plays an important role.
 For each matrix \(A\in\Cpq\), there exists a unique matrix \(X\in\Cqp\), satisfying the four equations
\begin{align}\label{mpi}
 AXA&=A,&
 XAX&=X,&
 \rk{AX}^\ad&=AX,&
&\text{and}&
 \rk{XA}^\ad&=XA
\end{align}
 (see \teg{}~\zitaa{MR1152328}{\cprop{1.1.1}{14}}).
 This matrix \(X\) is called the \emph{Moore--Penrose inverse of \(A\)} and is denoted by \(A^\mpi\).
 Concerning the concept of Moore--Penrose inverse we refer to~\zita{MR0338013},~\zitaa{MR1105324}{\cch{1}}, and~\zitaa{MR1987382}{\cch{1}}.
 For our purposes, it is convenient to apply~\zitaa{MR1152328}{\csec{1.1}}.
 
 If \(\tmat{A & B\\ C & D}\) is the \tbr{} of a complex \taaa{(p+q)}{(r+s)}{matrix} \(M\) with \taaa{p}{r}{block} \(A\), then the matrix
\beql{E/A}
 M\schca A
 \defeq D-CA^\mpi B
\eeq
 is called the \emph{Schur complement of \(A\) in \(M\)}.
 Concerning a variety of applications of this concept in a lot of areas of mathematics, we refer to~\zitas{MR2160825}.

 In the paper, various kinds of concrete \tSchur{} complements in block matrices will play an essential role.
 By virtue of \rrem{H.R.Hblock}, we use in the sequel the following notation:

\bnotal{H.N.L}
 If \(\seqska\) is a sequence of complex \tpqa{matrices}, then let \(\Lu{0}\defeq\Hu{0}\) and let \(\Lu{n}\defeq\Hu{n}\schca \Hu{n-1}\) for all \(n\in\N\) with \(2n\leq\kappa\).
\enota

 We write \(\rank A\) for the rank of a complex matrix \(A\) and \(\det B\) for the determinant of a square complex matrix \(B\).
 
\bremnl{\tcf{}~\zitaa{MR3014201}{\clem{3.5}{223}}}{R1533} 
 Let \(n\in\NO\) and let \(\seqs{2n}\in\Hggq{2n}\).
 Then \(\rank\Hu{n}=\sum_{k=0}^n\rank\Lu{k}\) and \(\det\Hu{n}=\prod_{k=0}^n\det\Lu{k}\).
\erem

\bnotal{N.Trip}
 Let \(\seqska\) be a sequence of complex \tpqa{matrices}.
 Then let \(\Tripu{0}\defeq\Opq\) and \(\Tripu{n}\defeq\zuu{n}{2n-1}\Hu{n-1}^\mpi\yuu{n}{2n-1}\) for all \(n\in\N\) with \(2n-1\leq\kappa\).
\enota

\bdefnl{D0750}
 If \(\seqska\) is a sequence of complex \tpqa{matrices}, then (using \rnota{F.N.sa}) the sequences \(\seq{\umg{j}}{j}{0}{\kappa}\) and \(\seq{\omg{j}}{j}{0}{\kappa}\) given by \(  \umg{2k}\defeq\ug\su{2k}+\Tripa{k}\) and \(\omg{2k}\defeq\obg\su{2k}-\Tripb{k}\) for all \(k\in\NO\) with \(2k\leq\kappa\) and by \(\umg{2k+1}\defeq\Tripu{k+1}\) and \(\omg{2k+1}\defeq-\ug\obg\su{2k}+(\ug+\obg)\su{2k+1}-\Tripab{k}\) for all \(k\in\NO\) with \(2k+1\leq\kappa\) are called the \emph{\tflep{\(\seqska\)}} and the \emph{\tfrep{\(\seqska\)}}, \tresp{}
\edefn

 By virtue of \rnota{N.Trip}, we have in particular
\begin{align}\label{F.G.uo01}
 \umg{0}&=\ug\su{0},&
 \omg{0}&=\obg\su{0},&
 \umg{1}&=\su{1}\su{0}^\mpi\su{1},&
&\text{and}&
 \omg{1}&=-\ug\obg\su{0}+(\ug+\obg)\su{1}.
\end{align}

 Using \rlem{F.R.Fgg-s} and \rrem{A.R.A++*}, we easily obtain:

\breml{ab.L0911}
 If \(\seqska\in\Fggqka\), then \(\set{\umg{j},\omg{j}}\subseteq\CHq\) for all \(j\in\mn{0}{\kappa}\).
\erem

 Observe that for arbitrarily given \tH{} \tqqa{matrices} \(A\) and \(B\), the (closed) matricial interval
\beql{matint}
 \matint{A}{B}
 \defeq\setaca{X\in\CHq}{A\lleq X\lleq B}
\eeq
 is non-empty if and only if \(A\lleq B\).
 
\bthmnl{\zitaa{MR3775449}{\cthm{11.2(a)}{44}}}{165.T112}
 If \(m\in\NO\) and \(\seqs{m}\in\Fggqu{m}\), then the matricial interval \(\matint{\umg{m}}{\omg{m}}\) is non-empty and coincides with the set of all complex \tqqa{matrices} \(\su{m+1}\) for which \(\seqs{m+1}\) belongs to \(\Fggqu{m+1}\).
\ethm

\bdefnl{D1861}
 If \(\seqska \) is a sequence of complex \tpqa{matrices}, then we call \(\seq{\dia{j}}{j}{0}{\kappa}\) given by \(\dia{j}\defeq\omg{j}-\umg{j}\) the \emph{\tfdfa{\(\seqska\)}}.
\edefn

 By virtue of \eqref{F.G.uo01}, we have in particular
\begin{align}\label{F.G.d01}
 \dia{0}&=\ba\su{0}&
&\text{and}&
 \dia{1}&=-\ug\obg\su{0}+(\ug+\obg)\su{1}-\su{1}\su{0}^\mpi\su{1}.
\end{align}

\breml{F.R.diatr}%
 Let \(\seqska\) be a sequence of complex \tpqa{matrices} with \tfdf{} \(\seqdiaka \).
 For each \(k\in\mn{0}{\kappa}\), the matrix \(\dia{k}\) is built  from the matrices \(\su{0},\su{1},\dotsc,\su{k}\).
 In particular, for each \(m\in\mn{0}{\kappa}\), the \tfdfa{\(\seqs{m}\)} coincides with \(\seqdia{m}\).
\erem

\bremnl{\zitaa{MR4051874}{\crem{7.26}{163}}}{ab.L0907} 
 Suppose \(\kappa\geq1\).
 If \(\seqska\in\Dpqkappa\), then \(\dia{1}=\sau{0}\su{0}^\mpi\sub{0}\) and \(\dia{1}=\sub{0}\su{0}^\mpi\sau{0}\).
\erem

\bdefnl{D2972}
 Let \(\seqska\) be a sequence of complex \tpqa{matrices}.
 Then the sequence \(\seq{\usc{j}}{j}{0}{\kappa}\) given by \(\usc{0}\defeq\su{0}\) and by \(\usc{j}\defeq\su{j}-\umg{j-1}\) is called the \emph{\tflsc{\(\seqska\)}}.
 Furthermore, if \(\kappa\geq1\), then the sequence \(\seq{\osc{j}}{j}{1}{\kappa}\) given by \(\osc{j}\defeq\omg{j-1}-\su{j}\) is called the \emph{\tfusc{\(\seqska\)}}.
\edefn

 Because of \eqref{F.G.uo01}, we have in particular
\begin{align}\label{F.G.AB1B2}
 \usc{1}&=\sau{0},&
 \osc{1}&=\sub{0},&
&\text{and}&
 \osc{2}&=\sab{0}.
\end{align}

\breml{F.R.ABL}
 Let \(\seqska\) be a sequence of complex \tpqa{matrices}.
 Then \(\usc{2n}=\Lu{n}\) for all \(n\in\NO\) with \(2n\leq\kappa\) and \(\usc{2n+1}=\Lau{n}\) for all \(n\in\NO\) with \(2n+1\leq\kappa\).
 In particular, if \(n\geq1\), then \(\usc{2n}\) is the Schur complement of \(\Hu{n-1}\) in \(\Hu{n}\) and \(\usc{2n+1}\) is the Schur complement of \(\Hau{n-1}\) in \(\Hau{n}\). 
 Furthermore, \(\osc{2n+1}=\Lub{n}\) for all \(n\in\NO\) with \(2n+1\leq\kappa\) and \(\osc{2n+2}=\Lab{n}\) for all \(n\in\NO\) with \(2n+2\leq\kappa\).
 In particular, if \(n\geq1\), then \(\osc{2n+1}\) is the Schur complement of \(\Hub{n-1}\) in \(\Hub{n}\) and \(\osc{2n+2}\) is the Schur complement of \(\Hab{n-1}\) in \(\Hab{n}\).
\erem

 If \(A\) and \(B\) are two complex \tpqa{matrices}, then the matrix
\beql{ps}
 A\ps B
 \defeq A(A+B)^\mpi B
\eeq
 is called the \emph{parallel sum of \(A\) and \(B\)}.

\bpropnl{\zitaa{MR3775449}{\cthm{10.14}{195}}}{ab.P1422}
 If \(\seqska\in\Fggqka\), then \(\dia{0}=\ba\usc{0}\) and furthermore \(\dia{k}=\ba(\usc{k}\ps\osc{k})\) and \(\dia{k}=\ba(\osc{k}\ps\usc{k})\) for all \(k\in\mn{1}{\kappa}\).
\eprop

\bpropnl{\zitaa{MR3775449}{\cprop{10.15(a)}{37}}}{ab.C0929}
 If \(\seqska\in\Fggqka\), then \(\dia{j}\in\Cggq\) for all \(j\in\mn{0}{\kappa}\).
\eprop

\bpropnl{\zitaa{MR3775449}{\cprop{10.18}{38}}}{ab.R13371422}
 Let \(\seqska\in\Fggqka\).
 Then \(\ran{\dia{0}}=\ran{\usc{0}}\) and \(\nul{\dia{0}}=\nul{\usc{0}}\).
 Furthermore, \(\ran{\dia{j}}=\ran{\usc{j}}\cap\ran{\osc{j}}\) and \(\nul{\dia{j}}=\nul{\usc{j}}+\nul{\osc{j}}\) for all \(j\in\mn{1}{\kappa}\), and \(\ran{\dia{j}}=\ran{\usc{j+1}}+\ran{\osc{j+1}}\) and \(\nul{\dia{j}}=\nul{\usc{j+1}}\cap\nul{\osc{j+1}}\) for all \(j\in\mn{0}{\kappa-1}\).
\eprop

 The ranks of the matrices considered in \rprop{ab.R13371422} are connected by means of the well-known formula for the dimension of the sum of two arbitrary finite-dimensional linear subspaces:

\breml{A.R.dim+}
 If \(\mathcal{U}_1\) and \(\mathcal{U}_2\) are finite-dimensional linear subspaces of some vector space, then \(\dim(\mathcal{U}_1+\mathcal{U}_2)=\dim\mathcal{U}_1+\dim\mathcal{U}_2-\dim(\mathcal{U}_1\cap\mathcal{U}_2)\).
\erem

\bcorl{ab.R0938}
 Let \(\seqska\in\Fggqka\).
 Then \(\rank \dia{0}=\rank\usc{0}\) and
\(
 \rank\dia{j-1}+\rank\dia{j}
 =\rank\usc{j}+\rank\osc{j}
\)
 for all \(j\in\mn{1}{\kappa}\).
\ecor
\bproof
 From \rprop{ab.R13371422} we obtain \(\rank \dia{0}=\rank\usc{0}\) and, for all \(j\in\mn{1}{\kappa}\), furthermore \(\rank\dia{j-1}=\dim(\ran{\usc{j}}+\ran{\osc{j}})\) and \(\rank\dia{j}=\dim(\ran{\usc{j}}\cap\ran{\osc{j}})\).
 The application of \rrem{A.R.dim+} to the linear subspaces \(\ran{\usc{j}}\) and \(\ran{\osc{j}}\) of the finite-dimensional vector space \(\Cq\) yields then \(\rank\dia{j-1}=\rank\usc{j}+\rank\osc{j}-\rank\dia{j}\).
\eproof

 Using \rcor{ab.R0938}, we are able to derive certain relations between the ranks of the matrices \(\dia{j}\) and the ranks of the underlying \tbHms{}:

\bleml{ab.R0929} 
 Let \(\seqska\in\Fggqka\).
 Then \(\rank\dia{0}=\rank\Hu{0}\).
 Furthermore,
 \(
  \sum_{\ell=0}^{2n+1}\rank\dia{\ell}
  =\rank\Hau{n}+\rank\Hub{n}
 \)
 for all \(n\in\NO\) with \(2n+1\leq\kappa\) and
 \(
  \sum_{\ell=0}^{2n}\rank\dia{\ell}
  =\rank\Hu{n}+\rank\Hab{n-1}
 \)
 for all \(n\in\N\) with \(2n\leq\kappa\).
\elem
\bproof
 Because of \(\Hu{0}=\su{0}=\Lu{0}=\usc{0}\) and \rcor{ab.R0938}, we have \(\rank\dia{0}=\rank\Hu{0}\).
 Now consider an arbitrary \(n\in\NO\) with \(2n+1\leq\kappa\).
 From \rprop{ab.R0933} and \eqref{Fgg2n+1} we see that the sequences \(\seqsa{2n}\) and \(\seqsb{2n}\) both belong to \(\Hggq{2n}\).
 Thus, we apply \rrem{R1533} to obtain \(\rank\Hau{n}=\sum_{k=0}^n\rank\Lau{k}\) and \(\rank\Hub{n}=\sum_{k=0}^n\rank\Lub{k}\).
 Using \rcor{ab.R0938} and \rrem{F.R.ABL}, we get then
 \[\begin{split}
  \sum_{\ell=0}^{2n+1}\rank\dia{\ell}
  &=\sum_{k=0}^n\rk{\rank\dia{2k}+\rank\dia{2k+1}}
  =\sum_{k=0}^n\rk{\rank\usc{2k+1}+\rank\osc{2k+1}}\\
  &=\sum_{k=0}^n\rank\usc{2k+1}+\sum_{k=0}^n\rank\osc{2k+1}
  =\rank\Hau{n}+\rank\Hub{n}.
 \end{split}\]
 Now consider an arbitrary \(n\in\N\) with \(2n\leq\kappa\).
 From \rprop{ab.R0933} and \eqref{Fgg2n} we infer \(\seqs{2n}\in\Hggq{2n}\) and \(\seqsab{2(n-1)}\in\Hggq{2(n-1)}\).
 Thus, \rrem{R1533} yields \(\rank\Hu{n}=\sum_{k=0}^n\rank\Lu{k}\) and \(\rank\Hab{n-1}=\sum_{k=0}^{n-1}\rank\Lab{k}\).
 Using \rcor{ab.R0938} and \rrem{F.R.ABL}, we get then
\[\begin{split}
  \sum_{\ell=0}^{2n}\rank\dia{\ell}
  &=\rank\dia{0}+\sum_{m=1}^n\rk{\rank\dia{2m-1}+\rank\dia{2m}}
  =\rank\usc{0}+\sum_{m=1}^n\rk{\rank\usc{2m}+\rank\osc{2m}}\\
  &=\sum_{k=0}^n\rank\usc{2k}+\sum_{k=0}^{n-1}\rank\osc{2k+2}
  =\rank\Hu{n}+\rank\Hab{n-1}.\qedhere
 \end{split}\]
\eproof

\bpropnl{\zitaa{MR3775449}{\ccor{10.21}{202}}}{ab.C1343}
 Let \(\seqska\in\Fggqka\) and assume \(\kappa\geq1\).
 For all \(j\in\mn{1}{\kappa}\), then \(\dia{j}=\ba\usc{j}\dia{j-1}^\mpi\osc{j}\) and \(\dia{j}=\ba\osc{j}\dia{j-1}^\mpi\usc{j}\).
\eprop

\bcorl{F.C.detdia} 
 If \(\seqska\in\Fggqka\), then \(\det\dia{0}=\ba^q\det\usc{0}\) and, for all \(j\in\mn{1}{\kappa}\), furthermore
 \beql{F.C.detdia.B1}
  \det\dia{j-1}\det\dia{j}
  =\ba^q\det\usc{j}\det\osc{j}.
 \eeq
\ecor
\bproof
 Because of \rprop{ab.P1422} we have \(\det\dia{0}=\ba^q\det\usc{0}\).
 Now assume \(\kappa\geq1\) and let \(j\in\mn{1}{\kappa}\).
 First we consider the case that \(\det\dia{j-1}=0\).
 From \rprop{ab.R13371422} we can infer \(\nul{\dia{j-1}}\subseteq\nul{\usc{j}}\).
 Consequently, \(\det\usc{j}=0\) follows.
 Hence, \eqref{F.C.detdia.B1} is fulfilled.
 Now we consider the case \(\det\dia{j-1}\neq0\).
 In view of \rrem{A.R.A-1}, then \eqref{F.C.detdia.B1} is a consequence of \rprop{ab.C1343}.
\eproof

\bleml{F.L.dHdia} 
 Let \(\seqska\in\Fggqka\).
 Then \(\det\dia{0}=\ba^q\det\Hu{0}\).
 Furthermore,
 \(
  \prod_{\ell=0}^{2n+1}\det\dia{\ell}
  =\ba^{(n+1)q}\det\rk{\Hau{n}}\det\rk{\Hub{n}}
 \)
 for all \(n\in\NO\) with \(2n+1\leq\kappa\) and
 \(
  \prod_{\ell=0}^{2n}\det\dia{\ell}
  =\ba^{(n+1)q}\det\rk{\Hu{n}}\det\rk{\Hab{n-1}}
 \)
 for all \(n\in\N\) with \(2n\leq\kappa\).
\elem
\bproof
 Because of \(\Hu{0}=\su{0}=\Lu{0}=\usc{0}\) and \rcor{F.C.detdia} we have \(\det\dia{0}=\ba^q\det\Hu{0}\).
 Now consider an arbitrary \(n\in\NO\) with \(2n+1\leq\kappa\).
 With the same reasoning as in \rlem{ab.R0929}, we can infer from \rrem{R1533} then \(\det\Hau{n}=\prod_{k=0}^n\det\Lau{k}\) and \(\det\Hub{n}=\prod_{k=0}^n\det\Lub{k}\).
 Using \rcor{F.C.detdia} and \rrem{F.R.ABL}, we get then
 \[\begin{split}
  \prod_{\ell=0}^{2n+1}\det\dia{\ell}
  &=\prod_{k=0}^n\rk{\det\dia{2k} \det\dia{2k+1}}
  =\prod_{k=0}^n\rk{\ba^q\det\usc{2k+1} \det\osc{2k+1}}\\
  &=\ba^{(n+1)q}\prod_{k=0}^n\det\usc{2k+1} \prod_{k=0}^n\det\osc{2k+1}
  =\ba^{(n+1)q}\det\Hau{n} \det\Hub{n}.
 \end{split}\]
 Now consider an arbitrary \(n\in\N\) with \(2n\leq\kappa\).
 With the same reasoning as in \rlem{ab.R0929}, we can conclude from \rrem{R1533} analogously \(\det\Hu{n}=\prod_{k=0}^n\det\Lu{k}\) and \(\det\Hab{n-1}=\prod_{k=0}^{n-1}\det\Lab{k}\).
 Using \rcor{F.C.detdia} and \rrem{F.R.ABL}, we obtain then
\[\begin{split}
  \prod_{\ell=0}^{2n}\det\dia{\ell}
  &=\det\dia{0} \prod_{m=1}^n\rk{\det\dia{2m-1} \det\dia{2m}}
  =\ba^q\det\usc{0} \prod_{m=1}^n\rk{\ba^q\det\usc{2m} \det\osc{2m}}\\
  &=\ba^{(n+1)q}\prod_{k=0}^n\det\usc{2k}\prod_{k=0}^{n-1}\det\osc{2k+2}
  =\ba^{(n+1)q}\det\Hu{n} \det\Hab{n-1}.\qedhere
 \end{split}\]
\eproof

 Now we state a consequence of \rthm{165.T112}.

\bcornl{\tcf{}~\zitaa{MR3979701}{\ccor{5.25}{2141}}}{ab.R1011}
 Let \(m\in\NO\), let \(\seqs{m}\in\Fggqu{m}\), let \(\lambda\in[0,1]\), and let \(\su{m+1}\defeq\umg{m}+\lambda\dia{m}\).
 Then, the sequence \(\seqs{m+1}\) belongs to \(\Fggqu{m+1}\).
 Furthermore, \(\usc{m+1}=\lambda\dia{m}\), \(\osc{m+1}=(1-\lambda)\dia{m}\), and \(\dia{m+1}=\ba\lambda(1-\lambda)\dia{m}\).
\ecor

 In~\zitaa{MR3979701}{\cdefn{6.1}{2142}}, we subsumed the Schur complements mentioned in \rrem{F.R.ABL} to a parameter sequence:

\bdefnl{D0752}
 Let \(\seqska\) be a sequence of complex \tpqa{matrices}.
 Let the sequence \(\fpseqka\) be given by \(\fpu{0}\defeq\usc{0}\), by \(\fpu{4k+1}\defeq\usc{2k+1}\) and \(\fpu{4k+2}\defeq\osc{2k+1}\) for all \(k\in\NO\) with \(2k+1\leq\kappa\), and by \(\fpu{4k+3}\defeq\osc{2k+2}\) and \(\fpu{4k+4}\defeq\usc{2k+2}\) for all \(k\in\NO\) with \(2k+2\leq\kappa\). 
 Then we call \(\fpseqka\) the \emph{\tfpfa{\(\seqska \)}}.
\edefn

 In view of \eqref{F.G.AB1B2} and \eqref{F.G.uo01}, we have in particular
\begin{align}
 \fpu{0}&=\su{0},&
 \fpu{1}&=\sau{0}=\su{1}-\ug\su{0},&
&\text{and}&
 \fpu{2}&=\sub{0}=\obg\su{0}-\su{1}.\label{F.G.f012}
\end{align}
 The \tFnnd{ness} as well as rank constellations among the \tnnH{} \tbHms{} \(\Hu{n}\), \(\Hau{n}\), \(\Hub{n}\), and \(\Hab{n}\) can be characterized in terms of \tfp{s} (\tcf~\zitaa{MR3979701}{\cpropss{6.13}{2145}{6.14}{2145}}).

\breml{F.R.fpftr}
 Let \(\seqska \) be a sequence of complex \tpqa{matrices} with \tfpf{} \(\fpseqka\).
 Then \(\fpu{0}=\su{0}\).
 Furthermore, for each \(k\in\mn{1}{\kappa}\), the matrices \(\fpu{2k-1}\) and \(\fpu{2k}\) are built from the matrices \(\su{0},\su{1},\dotsc,\su{k}\).
 In particular, for each \(m\in\mn{0}{\kappa}\), the \tfpfa{\(\seqs{m}\)} coincides with \(\seq{\fpu{j}}{j}{0}{2m}\).
\erem

\bpropnl{\zitaa{MR3979701}{\cprop{6.14}{2145}}}{F.P.FggFP}%
 Let \(\seqska \) be a sequence of complex \tqqa{matrices}.
 Then \(\seqska \in\Fggqka\) if and only if \(\fpu{j}\in\Cggq\) for all \(j\in\mn{0}{2\kappa}\).
\eprop

\bremnl{\zitaa{MR3979701}{\crem{6.16}{2146}}}{F.R.f2n-1}%
 Let \(\seqska\) be a sequence of complex \tpqa{matrices}.
 For all \(k\in\mn{1}{\kappa}\), then \(\fpu{2k-1}=\dia{k-1}-\fpu{2k}\).
\erem

 To single out all sequences \(\fpseqka\) of complex \tqqa{matrices} which indeed occur as \tfp{s} of sequences \(\seqska\in\Fggqka\), we introduced in~\zitaa{MR3979701}{\cnota{6.19}{2146}} the following class:

\bnotal{F.N.cs}
 For each \(\eta\in[0,\infp)\), denote by \(\cs{q}{\kappa}{\eta}\) the set of all sequences \(\seq{f_j}{j}{0}{2\kappa}\) of \tnnH{} \tqqa{matrices} satisfying, in the case \(\kappa\geq1\), the equations \(\eta f_{0}=f_{1}+f_{2}\) and \(\eta\rk{f_{2k-1}\ps f_{2k}}=f_{2k+1}+f_{2k+2}\) for all \(k\in\mn{1}{\kappa-1}\).
\enota

\bthmnl{\tcf{}~\zitaa{MR3979701}{\cthm{6.20}{2147}}}{F.T.FggFP}
 The mapping \(\Gamma_{\ug,\obg}\colon\Fggqka\to\csqkad\) given by \(\seqska\mapsto\fpseqka\) is well defined and bijective.
\ethm

 For each matrix \(A\in\Cggq\), there exists a uniquely determined matrix \(Q\in\Cggq\) with \(Q^2=A\) called the \emph{\tnnH{} square root} \(Q=A^\varsqrt \) of \(A\).
 To uncover relations between the \tfp{s} \(\fpseqka\) and to obtain a parametrization of the set \(\Fggqka\), we introduced in~\zitaa{MR3979701}{\cdefn{6.21}{2147} and \cnota{6.28}{2150}} another parameter sequence \(\seqciaka\) and a corresponding class \(\esqkad\) of sequences of complex matrices.
 (Observe that these constructions are well defined due to \rprop{ab.C0929} and \rrem{A.R.XAX}.)

\bdefnl{D0754}
 Let \(\seqska\in\Fggqka\) with \tfpf{} \(\fpseqka\) and \tfdf{} \(\seqdiaka\).
 Then we call \(\seqciaka\) given by \(\cia{0}\defeq\fpu{0}\) and by \(\cia{j}\defeq\rk{\dia{j-1}^\varsqrt}^\mpi\fpu{2j}\rk{\dia{j-1}^\varsqrt}^\mpi\) 
 for each \(j\in\mn{1}{\kappa}\) the \emph{\tfcfa{\(\seqska\)}}.
\edefn

\blemnl{\tcf{}~\zitaa{MR3979701}{\cprop{6.27}{2149}}}{F.L.cia-fp}
 If \(\seqska\in\Fggqka\), then \(\fpu{2j}=\dia{j-1}^\varsqrt \cia{j}\dia{j-1}^\varsqrt\) for all \(j\in\mn{1}{\kappa}\).
\elem

 With the Euclidean scalar product \(\ipE{\cdot}{\cdot}\colon\x{\Cq}\to\C\) given by \(\ipE{x}{y}\defeq y^\ad x\), which is \(\C\)\nobreakdash-linear in its first argument, the \(\C\)\nobreakdash-vector space \(\Cq\) becomes a unitary space.
 Let \(\mathcal{U}\) be an arbitrary non-empty subset of \(\Cq\).
 The orthogonal complement \(\mathcal{U}^\orth\defeq\setaca{v\in\Cq}{\ipE{v}{u}=0\text{ for all }u\in\mathcal{U}}\) of \(\mathcal{U}\) is a linear subspace of the unitary space \(\Cq\).
 If \(\mathcal{U}\) is a linear subspace itself, the unitary space \(\Cq\) is the orthogonal sum of \(\mathcal{U}\) and \(\mathcal{U}^\orth\).
 In this case, we write \(\OPu{\mathcal{U}}\) for the transformation matrix corresponding to the orthogonal projection onto \(\mathcal{U}\) with respect to the standard basis of \(\Cq\), \tie{}, \(\OPu{\mathcal{U}}\) is the uniquely determined matrix \(P\in\Cqq\) satisfying the three conditions \(P^2=P\), \(P^\ad=P\), and \(\ran{P}=\mathcal{U}\).

\bnotal{F.N.01seq}
 For each \(\eta\in[0,\infp)\), let \(\es{q}{\kappa}{\eta}\) be the set of all sequences \(\seq{e_k}{k}{0}{\kappa}\) from \(\Cggq\) which fulfill the following condition:
 If \(\kappa\geq1\), then \(e_k\lleq\OPu{\ran{d_{k-1}}}\) for all \(k\in\mn{1}{\kappa}\), where the sequence \(\seq{d_k}{k}{0}{\kappa}\) is recursively given by \(d_0\defeq\eta e_0\) and
\[
 d_k
 \defeq\eta d_{k-1}^\varsqrt  e_k^\varsqrt (\OPu{\ran{d_{k-1}}}-e_k) e_k^\varsqrt  d_{k-1}^\varsqrt.
\]
\enota

 Regarding \rthm{165.T112}, in the case \(q=1\) (\tcf{}~\zitaa{MR1468473}{\csec{1.3}}), the (classical) canonical moments \(p_1,p_2,p_3,\dotsc\) of a point in the moment space corresponding to a probability measure \(\mu\) on \(\ab=[0,1]\) are given in our notation by %
\begin{align*}%
 p_k
 &=\frac{\su{k}-\umg{k-1}}{\omg{k-1}-\umg{k-1}}
 =\frac{\usc{k}}{\dia{k-1}},&k&\in\N,
\end{align*}
 where the sequence \(\seqsinf\) of power moments \(\su{j}\defeq\int_{[0,1]}x^j\mu\rk{\dif x}\) associated with \(\mu\) is \txyFnnd{0}{1} with \(\su{0}=1\).
 Observe that the \tfuc{0}{1}s \(\seqciainf\) of \(\seqsinf\) are connected to the canonical moments via
\begin{align}\label{F.G.epq}
 p_1&=1-\cia{1},&
 p_2&=\cia{2},&
 p_3&=1-\cia{3},&
 p_4&=\cia{4},&
 p_5&=1-\cia{5},&
 &\dotsc
\end{align}
 The quantities \(q_k=1-p_k\) occur in the classical framework as well (see, \teg{}~\zitaa{MR1468473}{\csec{1.3}}).
 In the general case \(q\in\N\) we have the following:

\bthmnl{\zitaa{MR3979701}{\cthm{6.30}{2153}}}{F.T.Fggcia}
 The mapping \(\Sigma_{\ug,\obg}\colon\Fggqka\to\esqkad\) given by \(\seqska\mapsto\seqciaka\) is well defined and bijective.
\ethm

\bpropnl{\tcf{}~\zitaa{MR3979701}{\cprop{6.32}{2153}}}{F.P.ed}%
 Let \(\seqska\in\Fggqka\).
 Then \(\cia{0}\lgeq\Oqq\) and \(\cia{j}\in\matint{\Oqq}{\OPu{\ran{\dia{j-1}}}}\) for all \(j\in\mn{1}{\kappa}\).
 Furthermore, \(\dia{0}=\ba\cia{0}\) and \(\dia{j}=\ba\dia{j-1}^\varsqrt  \cia{j}^\varsqrt (\OPu{\ran{\dia{j-1}}}-\cia{j}) \cia{j}^\varsqrt \dia{j-1}^\varsqrt \) for all \(j\in\mn{1}{\kappa}\).
\eprop

\bexanl{\tcf{}~\zitaa{MR4051874}{\cexa{7.35}{166}}}{F.E.e=1/2}%
 Let \(\lambda\in(0,1)\), let \(B\in\Cggq\), and let \(P\defeq\OPu{\ran{B}}\).
 Then \(\seqciainf\) given by \(\cia{0}\defeq B\) and by \(\cia{j}\defeq\lambda P\) for all \(j\in\N\) is the \tfcfa{a} sequence \(\seqsinf\in\Fggqinf\).
\eexa 
 
 We continue by recalling the construction of a certain transformation for sequences of matrices.
 This transformation was introduced in~\zita{MR4051874} and constitutes the elementary step of a Schur type algorithm in the class of \tFnnd{} sequences, \tie{}, Hausdorff moment sequences:

\bdefnl{ab.N1137b}
 Let \(\seqska\) be a sequence of complex \tpqa{matrices}.
 Further let \(\sub{-1}\defeq-\su{0}\) and, in the case \(\kappa\geq1\), let \(\seqsb{\kappa-1}\) be given by \rnota{F.N.sa}.
 Then we call the sequence \(\seqspbka\) given by \(\spb{j}\defeq\sub{j-1}\) the \emph{\tbmodv{\seqska}}.
\edefn

 In particular, if \(\obg=0\), then \(\seqspbka\) coincides with the sequence \(\seq{-\su{j}}{j}{0}{\kappa}\).
 For an arbitrary \(\obg\in\R\), the sequence \(\seqska\) is reconstructible from \(\seqspbka\) as well.

 Let \(\seqska\) and \(\seqt{\kappa}\) be sequences of complex \tpqa{} and \taaa{q}{r}{matrices}, \tresp{}
 As usual, the \emph{\tCP{}} \(\seqa{x}{\kappa}\) of \(\seqska\) and \(\seqt{\kappa}\) is given by \(x_j\defeq\sum_{\ell=0}^js_\ell t_{j-\ell}\).

\bdefnl{D1419}
 Let \(\seqska \) be a sequence of complex \tpqa{matrices}.
 Then we call the sequence \(\seq{s_j^\rez}{j}{0}{\kappa}\) defined by \(s_0^\rez\defeq s_0^\mpi\) and, for all \(j\in\mn{1}{\kappa}\), recursively by
\(
 s_j^\rez
 \defeq-s_0^\mpi\sum_{\ell=0}^{j-1}s_{j-\ell}s_\ell^\rez
\)
 the \emph{\trFa{\(\seq{s_j}{j}{0}{\kappa}\)}}.
\edefn

\breml{M.R.reztr}
 Let \(\seqska\) be a sequence of complex \tpqa{matrices} with \trF{} \(\seqa{r}{\kappa}\).
 For each \(k\in\mn{0}{\kappa}\), then the matrix \(r_k\) is built from the matrices \(s_0,s_1,\dotsc,s_k\).
 In particular, for each \(m\in\mn{0}{\kappa}\), the \trFa{\(\seqs{m}\)} coincides with \(\seqa{r}{m}\).
\erem

 Using the \tCP{} and the \trFa{the \tbmodv{\seqsa{\kappa-1}}}, we introduce now a transformation of sequences of complex matrices:

\bdefnl{ab.N0940}
 Suppose \(\kappa\geq1\).
 Let \(\seqska\) be a sequence of complex \tpqa{matrices}.
 Denote by \(\seqapb{\kappa-1}\) the \tbmodv{\seqsa{\kappa-1}} and by \(\seq{x_j}{j}{0}{\kappa-1}\) the \tCPa{\(\seqsb{\kappa-1}\)}{\(\seq{\apb{j}^\rez}{j}{0}{\kappa-1}\)}.
 Then we call the sequence \(\seqt{\kappa-1}\) given by \(\tu{j}\defeq-\sau{0}\su{0}^\mpi x_j\sau{0}\) the \emph{\tFTv{\(\seqska\)}}.
\edefn

 Since, in the classical case that \(\ug=0\) and \(\obg=1\), the sequence \(\seqsa{\kappa-1}\) coincides with the shifted sequence \(\seq{\su{j+1}}{j}{0}{\kappa-1}\), the \tFuuT{0}{1} is given by \(\tu{j}=-\su{1}\su{0}^\mpi x_j\su{1}\) with the \tCP{} \(\seq{x_j}{j}{0}{\kappa-1}\) of \(\seqsb{\kappa-1}\) and \(\seq{\apb{j}^\rez}{j}{0}{\kappa-1}\), where the sequence \(\seqsb{\kappa-1}\) is given by \(\sub{j}=\su{j}-\su{j+1}\) and the sequence \(\seqapb{\kappa-1}\) is given by \(\apb{0}=-\su{1}\) and by \(\apb{j}=\su{j}-\su{j+1}\) for \(j\in\mn{1}{\kappa-1}\).

\breml{F.R.FTtr}
 Assume \(\kappa\geq1\).
 Let \(\seqska\) be a sequence of complex \tpqa{matrices} with \tFT{} \(\seqt{\kappa-1}\).
 Then one can see from \rrem{M.R.reztr} that, for each \(k\in\mn{0}{\kappa-1}\), the matrix \(t_k\) is built from the matrices \(\su{0},\su{1},\dotsc,\su{k+1}\).
 In particular, for all \(m\in\mn{1}{\kappa}\), the \tFTv{\(\seqs{m}\)} coincides with \(\seqt{m-1}\).
\erem

\blemnl{\zitaa{MR4051874}{\clem{8.29}{175}}}{ab.P1728a} 
 Suppose \(\kappa\geq1\).
 Let \(\seqska\in\Dpqka\).
 Denote by \(\seq{t_j}{j}{0}{\kappa-1}\) the \tFTv{\(\seqska\)}.
 Then \(t_0=\dia{1}\), where \(\dia{1}\) is given by \rdefn{D1861}.
\elem

 Using the parallel sum given via \eqref{ps}, the effect caused by \tFT{ation} on the \tfp{s} can be completely described:
 
\bcornl{\zitaa{MR4051874}{\ccor{9.10}{200}}}{F.C.Salg1}
 Assume \(\kappa\geq1\) and let \(\seqska\in\Fggqka\) with \tFT{} \(\seqt{\kappa-1}\).
 Denote by \(\fgpseq{2(\kappa-1)}\) the \tfpfa{\(\seqt{\kappa-1}\)}.
 Then \(\fgpu{0}=\ba\rk{\fpu{1}\ps\fpu{2}}\) and \(\fgpu{j}=\ba\fpu{j+2}\) for all \(j\in\mn{1}{2(\kappa-1)}\).
\ecor

 We are now going to iterate the \tFT{ation}:

\bdefnl{ab.N1020}
 Let \(\seqska\) be a sequence of complex \tpqa{matrices}.
 Let the sequence \(\seq{\su{j}^\FTa{0}}{j}{0}{\kappa}\) be given by \(\su{j}^\FTa{0}\defeq\su{j}\).
 If \(\kappa\geq1\), then, for all \(k\in\mn{1}{\kappa}\), let the sequence \(\seq{\su{j}^\FTa{k}}{j}{0}{\kappa-k}\) be recursively defined to be the \tFT{} of the sequence \(\seq{\su{j}^\FTa{k-1}}{j}{0}{\kappa-(k-1)}\).
 For all \(k\in\mn{0}{\kappa}\), then we call the sequence \(\seq{\su{j}^\FTa{k}}{j}{0}{\kappa-k}\) the \emph{\tnFTv{k}{\(\seqska \)}}.
\edefn

\breml{F.R.FT1}%
 Suppose \(\kappa\geq1\).
 Let \(\seqska\) be a sequence of complex \tpqa{matrices}.
 Then \(\seq{\su{j}^\FT}{j}{0}{\kappa-1}\) is exactly the \tFTv{\(\seqska\)} from \rdefn{ab.N1020}.
\erem

\breml{F.R.kFTtr} 
 Let \(k\in\mn{0}{\kappa}\) and let \(\seqska \) be a sequence of complex \tpqa{matrices} with \tnFT{k} \(\seq{u_j}{j}{0}{\kappa-k}\).
 In view of \rrem{F.R.FTtr}, we see that, for each \(\ell\in\mn{0}{\kappa-k}\), the matrix \(u_{\ell}\) is built only from the matrices \(\su{0},\su{1},\dotsc,\su{\ell+k}\).
 In particular, for each \(m\in\mn{k}{\kappa}\), the \tnFTv{k}{\(\seqs{m}\)} coincides with \(\seq{u_j}{j}{0}{m-k}\).
\erem

\bpropnl{\zitaa{MR4051874}{\cthm{9.4}{195}}}{ab.P1030}
 If \(\seqska\in\Fggqka\), then \(\seq{\su{j}^\FTa{k}}{j}{0}{\kappa-k}\in\Fggqu{\kappa-k}\) for all \(k\in\mn{0}{\kappa}\).
\eprop

\bpropnl{\zitaa{MR4051874}{\cprop{9.8}{199}}}{F.P.FPFT} 
 Let \(\seqska\in\Fggqka\) with \tfpf{} \(\fpseqka\).
 Then \(\fpu{0}=\su{0}\) and furthermore
\(
 \fpu{4k+1}=\ba^{-2k}\sau{0}^\FTa{2k}\)
 and
\(
 \fpu{4k+2}=\ba^{-2k}\sub{0}^\FTa{2k}
\)
 for all \(k\in\NO\) with \(2k+1\leq\kappa\) and
\(
 \fpu{4k+3}=\ba^{-(2k+1)}\sau{0}^\FTa{2k+1}
\)
 and
\(
 \fpu{4k+4}=\ba^{-(2k+1)}\sub{0}^\FTa{2k+1}
\)
 for all \(k\in\NO\) with \(2k+2\leq\kappa\).
\eprop

 In the following result, we express the \tfdf{} introduced in \rdefn{D1861} by the \tFT{s} of an \tFnnd{} sequence:
 
\bcornl{\zitaa{MR4051874}{\ccor{9.9}{200}}}{F.C.diaFT}
 If \(\seqska\in \Fggqka \), then \(\dia{j}=\ba^{-(j-1)}\su{0}^\FTa{j}\) for all \(j\in\mn{0}{\kappa}\).
\ecor

\bpropnl{\zitaa{MR4051874}{\cprop{9.11}{200}}}{F.P.diaalg} 
 Let \(k\in\mn{0}{\kappa}\) and let \(\seqska\in\Fggqka \) with \tfdf{} \(\seqdiaka\) and \tnFT{k} \(\seq{\su{j}^\FTa{k}}{j}{0}{\kappa-k}\).
 Then \(\seq{\ba^k\dia{k+j}}{j}{0}{\kappa-k}\) coincides with the \tfdfa{\(\seq{\su{j}^\FTa{k}}{j}{0}{\kappa-k}\)}.
\eprop

\bthmnl{\zitaa{MR4051874}{\cthm{9.14}{201}}}{F.T.SalgIP}%
 Let \(\seqska\in\Fggqu{\kappa}\) with \tfdf{} \(\seqdiaka\) and \tfcf{} \(\seqciaka \).
 Let \(k\in\mn{0}{\kappa}\).
 Then the \tnFT{k} \(\seq{\su{j}^\FTa{k}}{j}{0}{\kappa-k}\) of \(\seqska\) belongs to \(\Fggqu{\kappa-k}\) and the \tfcf{} \(\seq{\mathfrak{p}_{j}}{j}{0}{\kappa-k}\) of \(\seq{\su{j}^\FTa{k}}{j}{0}{\kappa-k}\) fulfills \(\mathfrak{p}_{0}=\ba^{k-1}\dia{k}\) and \(\mathfrak{p}_{j}=\cia{k+j}\) for all \(j\in\mn{1}{\kappa-k}\).
\ethm

 Now we look for a characterization of the fixed points of the \tFTion{}.

\bcorl{C1358}%
 Suppose \(\kappa\geq1\).
 Let \(\seqska\in\Fggqka\) with \tfcf{} \(\seqciaka\) and \tFT{} \(\seqt{\kappa-1}\).
 Then \(\tu{0}=\su{0}\) if and only if \(\cia{1}-\cia{1}^2=\ba^{-2}\OPu{\ran{\su{0}}}\).
\ecor
\bproof
 \rprop{F.P.ed} yields \(\dia{0}=\ba\cia{0}\) and \(\dia{1}=\ba\dia{0}^\varsqrt\cia{1}^\varsqrt\rk{\OPu{\ran{\dia{0}}}-\cia{1}}\cia{1}^\varsqrt\dia{0}^\varsqrt\).
 Because of \rlem{F.R.Fgg-s}, the matrix \(\su{0}\) is \tnnH{}.
 According to \eqref{F.G.d01}, we have \(\dia{0}=\ba\su{0}\).
 In view of \(\ba>0\), then \(\dia{0}^\varsqrt=\ba^\varsqrt\su{0}^\varsqrt\) as well as \(\ran{\dia{0}}=\ran{\su{0}}\) and \(\nul{\dia{0}}=\nul{\su{0}}\) follow.
 Consequently, we can conclude \(\dia{1}=\ba^2\su{0}^\varsqrt\cia{1}^\varsqrt\rk{\OPu{\ran{\su{0}}}-\cia{1}}\cia{1}^\varsqrt\su{0}^\varsqrt\).
 Using \rrem{A.R.A++*}, we can infer from \rdefn{D0754} furthermore \(\ran{\cia{1}}\subseteq\ran{\rk{\dia{0}^\varsqrt}^\mpi}=\ran{\rk{\dia{0}^\varsqrt}^\ad}=\ran{\dia{0}^\varsqrt}=\ran{\dia{0}}=\ran{\su{0}}\) and, similarly, \(\nul{\su{0}}\subseteq\nul{\cia{1}}\).
 In particular, \(\ran{\cia{1}^\varsqrt}=\ran{\cia{1}}\subseteq\ran{\su{0}}\) follows, implying \(\cia{1}^\varsqrt\OPu{\ran{\su{0}}}=\cia{1}^\varsqrt\).
 Hence, \(\dia{1}=\ba^2\su{0}^\varsqrt\rk{\cia{1}-\cia{1}^2}\su{0}^\varsqrt\).
 Therefore, \(\ran{\dia{1}}\subseteq\ran{\su{0}^\varsqrt}\) and \(\nul{\su{0}^\varsqrt}\subseteq\nul{\dia{1}}\).
 In view of \rrem{F.R.Fgg<D}, the application of \rlem{ab.P1728a} yields \(t_0=\dia{1}\).
 Hence, it remains to show that \(\ba^2\su{0}^\varsqrt\rk{\cia{1}-\cia{1}^2}\su{0}^\varsqrt=\su{0}\) is equivalent to \(\cia{1}-\cia{1}^2=\ba^{-2}\OPu{\ran{\su{0}}}\).
  
 First assume that \(\ba^2\su{0}^\varsqrt\rk{\cia{1}-\cia{1}^2}\su{0}^\varsqrt=\su{0}\) is fulfilled.
 Taking into account \rrem{A.R.A++*}, we see \(\ran{\cia{1}}\subseteq\ran{\su{0}}=\ran{\su{0}^\varsqrt}=\ran{\rk{\su{0}^\varsqrt}^\ad}=\ran{\rk{\su{0}^\varsqrt}^\mpi}\) and similarly \(\nul{\rk{\su{0}^\varsqrt}^\mpi}\subseteq\nul{\cia{1}}\).
 Thus, we can infer from \rremsss{R.AA+B.A}{R.AA+B.B}{A.R.A++*} that \(\ba^2\rk{\cia{1}-\cia{1}^2}=\rk{\su{0}^\varsqrt}^\mpi\su{0}\rk{\su{0}^\varsqrt}^\mpi\).
 Because of \(\ran{\su{0}}=\ran{\su{0}^\varsqrt}=\ran{\rk{\su{0}^\varsqrt}^\ad}\), the application of \rrem{ab.R1052} yields \(\rk{\su{0}^\varsqrt}^\mpi\su{0}\rk{\su{0}^\varsqrt}^\mpi=\OPu{\ran{\rk{\su{0}^\varsqrt}^\ad}}\OPu{\ran{\su{0}^\varsqrt}}=\OPu{\ran{\su{0}}}^2=\OPu{\ran{\su{0}}}\).
 Consequently, \(\cia{1}-\cia{1}^2=\ba^{-2}\OPu{\ran{\su{0}}}\).
 
 Conversely, assume that \(\cia{1}-\cia{1}^2=\ba^{-2}\OPu{\ran{\su{0}}}\) holds true.
 In view of \(\ran{\su{0}^\varsqrt}=\ran{\su{0}}\), we have \(\su{0}^\varsqrt\OPu{\ran{\su{0}}}=\su{0}^\varsqrt\).
 Thus, we obtain \(\ba^2\su{0}^\varsqrt\rk{\cia{1}-\cia{1}^2}\su{0}^\varsqrt=\su{0}^\varsqrt\OPu{\ran{\su{0}}}\su{0}^\varsqrt=\su{0}\).
\eproof

\bcorl{R1124}%
 Let \(\seqsinf\in\Fggqinf\) with \tfcf{} \(\seqciainf\) and \tFT{} \(\seqtinf\).
 Then the following statements are equivalent:
\baeqi{0}
 \il{R1124.i} \(\seqtinf\) coincides with \(\seqsinf\), \tie{}\ the sequence \(\seqsinf\) is a fixed point of the \tFTion{}.
 \il{R1124.ii} \(\cia{1}-\cia{1}^2=\ba^{-2}\OPu{\ran{\su{0}}}\) and \(\cia{j}=\cia{1}\) for all \(j\in\N\).
\eaeqi
\ecor
\bproof
 In view of \rrem{F.R.FT1}, we see from \rthm{F.T.SalgIP} that \(\seqtinf\) belongs to \(\Fggqinf\) and that the \tfcf{} \(\seq{\mathfrak{p}_{j}}{j}{0}{\infi}\) of \(\seqtinf\) fulfills \(\mathfrak{p}_{0}=\dia{1}\) and \(\mathfrak{p}_{j}=\cia{j+1}\) for all \(j\in\N\).
 Because of \rthm{F.T.Fggcia}, statement~\ref{R1124.i} holds if and only if \(\mathfrak{p}_{j}=\cia{j}\) for all \(j\in\NO\).
 Consequently,~\ref{R1124.i} is fulfilled if and only if \(\dia{1}=\cia{0}\) and \(\cia{j+1}=\cia{j}\) for all \(j\in\N\).
 Hence, it remains to show that \(\dia{1}=\cia{0}\) is equivalent to \(\cia{1}-\cia{1}^2=\ba^{-2}\OPu{\ran{\su{0}}}\).
 In view of \rrem{F.R.Fgg<D}, the application of \rlem{ab.P1728a} yields \(t_0=\dia{1}\).
 By virtue of \rdefn{D0754} and \eqref{F.G.f012}, we have \(\cia{0}=\fpu{0}=\su{0}\).
 The application of \rcor{C1358} completes the proof.
\eproof

 Now we draw special attention to the scalar case \(q=1\).

\bexal{E1749}%
 Let \(\Phi\colon\Fgguuuu{1}{\infi}{\ug}{\obg}\to\Fgguuuu{1}{\infi}{\ug}{\obg}\) be defined by \(\seqsinf\mapsto\seqtinf\), where \(\seqtinf\) is the \tFTv{\(\seqsinf\)}.
 Then:
\benui
 \il{E1749.a} If \(\ba<2\), then \(\Phi\) has one single fixed point \(\seq{\su{0;j}}{j}{0}{\infi}\) given by \(\su{0;j}=0\).
 \il{E1749.b} Suppose \(\ba=2\).
 For each \(M\in[0,\infp)\) there exists a unique fixed point \(\seq{\su{M;j}}{j}{0}{\infi}\) of \(\Phi\) with \(\su{M;0}=M\).
 If \(M=0\), then \(\seq{\su{M;j}}{j}{0}{\infi}\) is given by \(\su{M;j}=0\).
 If \(M>0\), then \(\seq{\su{M;j}}{j}{0}{\infi}\) corresponds to the \tfcf{} \(\seq{\cia{M;j}}{j}{0}{\infi}\) given by \(\cia{M;0}=M\) and by \(\cia{M;j}=\frac{1}{2}\) for \(j\in\N\).
 \il{E1749.c} Suppose \(\ba>2\).
 Then \(\Phi\) has exactly one fixed point \(\seq{\su{0;j}}{j}{0}{\infi}\) with \(\su{0;0}=0\), namely \(\seq{\su{0;j}}{j}{0}{\infi}\) given by \(\su{0;j}=0\).
 For each \(M\in(0,\infp)\), furthermore \(\Phi\) has exactly two fixed points \(\seq{\su{M;j}^\pm}{j}{0}{\infi}\) with \(\su{M;0}^\pm=M\), namely \(\seq{\su{M;j}^\pm}{j}{0}{\infi}\) corresponding to the \tfcf{s} \(\seq{\cia{M;j}^\pm}{j}{0}{\infi}\) given by \(\cia{M;0}^\pm=M\) and by \(\cia{M;j}^\pm=\frac{1}{2}\pm\frac{1}{2\ba}\sqrt{\ba^2-4}\) for \(j\in\N\).
\eenui
\eexa
\bproof
 \eqref{E1749.a} Assume \(\ba<2\).
 It is readily checked that the sequence \(\seq{\su{0;j}}{j}{0}{\infi}\) defined by \(\su{0;j}\defeq0\) belongs to \(\Fgguuuu{1}{\infi}{\ug}{\obg}\) and that the sequence \(\seq{\cia{0;j}}{j}{0}{\infi}\) defined by \(\cia{0;j}\defeq0\) is the \tfcfa{\(\seq{\su{0;j}}{j}{0}{\infi}\)}.
 In view of \(\OPu{\ran{\su{0;0}}}=0\), we can thus infer from \rcor{R1124} that \(\seq{\su{0;j}}{j}{0}{\infi}\) is a fixed point of \(\Phi\).
 
 Now consider an arbitrary fixed point \(\seqsinf\in\Fgguuuu{1}{\infi}{\ug}{\obg}\) of \(\Phi\).
 If \(\su{0}=0\), then from \rrem{F.R.Fgg<D} we can conclude that \(\su{j}=0\) for all \(j\in\N\), \tie{}, \(\seqsinf\) coincides with \(\seq{\su{0;j}}{j}{0}{\infi}\).
 Consider now the case \(\su{0}\neq0\).
 Then \(\OPu{\ran{\su{0}}}=1\).
 Denote by \(\seqciainf\) the \tfcfa{\(\seqsinf\)}.
 Because of \rcor{R1124}, we have \(\cia{1}-\cia{1}^2=\ba^{-2}\OPu{\ran{\su{0}}}\).
 Consequently, \(0<\ba^{-2}=\ba^{-2}\OPu{\ran{\su{0}}}=\cia{1}-\cia{1}^2=\frac{1}{4}-\rk{\cia{1}-\frac{1}{2}}^2\leq\frac{1}{4}\), contradicting \(\ba<2\).
 Thus, \(\seq{\su{0;j}}{j}{0}{\infi}\) is the only fixed point of \(\Phi\).
 
 \eqref{E1749.b} As above, we see that \(\seq{\su{0;j}}{j}{0}{\infi}\in\Fgguuuu{1}{\infi}{\ug}{\obg}\) is a fixed point of \(\Phi\) with \tfcf{} \(\seq{\cia{0;j}}{j}{0}{\infi}\) given by \(\cia{0;j}=0\).
 Consider an arbitrary \(M\in(0,\infp)\) and let \(\seq{\cia{M;j}}{j}{0}{\infi}\) be defined by \(\cia{M;0}\defeq M\) and by \(\cia{M;j}\defeq\frac{1}{2}\) for \(j\in\N\).
 In view of \(\OPu{\ran{M}}=1\), due to \rexam{F.E.e=1/2} there exists a sequence \(\seq{\su{M;j}}{j}{0}{\infi}\in\Fgguuuu{1}{\infi}{\ug}{\obg}\) with \tfcf{} \(\seq{\cia{M;j}}{j}{0}{\infi}\).
 According to \rthm{F.T.Fggcia}, the sequences \(\seq{\su{0;j}}{j}{0}{\infi}\) and \(\seq{\su{M;j}}{j}{0}{\infi}\) are different.
 By virtue of \rdefn{D0754} and \eqref{F.G.f012}, we have \(M=\cia{M;0}=\su{M;0}\).
 In particular \(\OPu{\ran{\su{M;0}}}=1\).
 Taking additionally into account \(\cia{M;1}=\frac{1}{2}\) and \(\ba=2\), hence \(\cia{M;1}-\cia{M;1}^2=\frac{1}{4}=\ba^{-2}\OPu{\ran{\su{M;0}}}\) follows.
 According of \rcor{R1124}, then \(\seq{\su{M;j}}{j}{0}{\infi}\) is a fixed point of \(\Phi\). 
 
 Now consider an arbitrary fixed point \(\seqsinf\in\Fgguuuu{1}{\infi}{\ug}{\obg}\) of \(\Phi\).
 Because of \rlem{F.R.Fgg-s}, we have \(\su{0}\in[0,\infp)\).
 If \(\su{0}=0\), then, as above, the sequence \(\seqsinf\) coincides with \(\seq{\su{0;j}}{j}{0}{\infi}\).
 Now assume \(\su{0}>0\).
 Then \(\OPu{\ran{\su{0}}}=1\).
 Denote by \(\seqciainf\) the \tfcfa{\(\seqsinf\)}.
 Because of \rcor{R1124}, we have \(\cia{1}-\cia{1}^2=\ba^{-2}\OPu{\ran{\su{0}}}\) and \(\cia{1}=\cia{2}=\cia{3}=\dotsb\)
 Taking additionally into account \(\ba=2\), then \(\rk{\cia{1}-\frac{1}{2}}^2=\cia{1}^2-\cia{1}+\frac{1}{4}=-\ba^{-2}\OPu{\ran{\su{0}}}+\frac{1}{4}=-2^{-2}+\frac{1}{4}=0\) follows, \tie{} \(\cia{1}=\frac{1}{2}\).
 By virtue of \rdefn{D0754} and \eqref{F.G.f012}, we have \(\cia{0}=\fpu{0}=\su{0}\).
 Setting \(M\defeq\su{0}\), then \(M\in(0,\infp)\) and \(\seqciainf\) coincides with \(\seq{\cia{M;j}}{j}{0}{\infi}\) defined as above, implying that \(\seqsinf\) coincides with \(\seq{\su{M;j}}{j}{0}{\infi}\), according to \rthm{F.T.Fggcia}.

 \eqref{E1749.c} As above, we see that \(\seq{\su{0;j}}{j}{0}{\infi}\in\Fgguuuu{1}{\infi}{\ug}{\obg}\) is a fixed point of \(\Phi\) with \tfcf{} \(\seq{\cia{0;j}}{j}{0}{\infi}\) given by \(\cia{0;j}=0\).
 Consider an arbitrary \(M\in(0,\infp)\) and let \(\seq{\cia{M;j}^\pm}{j}{0}{\infi}\) be defined by \(\cia{M;0}^\pm\defeq M\) and by \(\cia{M;j}^\pm\defeq\frac{1}{2}\pm\frac{1}{2\ba}\sqrt{\ba^2-4}\) for all \(j\in\N\).
 In view of \(\ba>2\), we have \(\ba^2-4>0\) and consequently \(0<\cia{M;j}^-<\frac{1}{2}<\cia{M;j}^+<1\) for all \(j\in\N\).
 Regarding \(\OPu{\ran{M}}=1\), hence, due to \rexam{F.E.e=1/2}, there exist sequences \(\seq{\su{M;j}^-}{j}{0}{\infi}\) and \(\seq{\su{M;j}^+}{j}{0}{\infi}\) belonging to \(\Fgguuuu{1}{\infi}{\ug}{\obg}\) with \tfcf{} \(\seq{\cia{M;j}^-}{j}{0}{\infi}\) and \(\seq{\cia{M;j}^+}{j}{0}{\infi}\), \tresp{}
 According to \rthm{F.T.Fggcia}, the sequences \(\seq{\su{0;j}}{j}{0}{\infi}\), \(\seq{\su{M;j}^-}{j}{0}{\infi}\), and \(\seq{\su{M;j}^+}{j}{0}{\infi}\) are pairwise different.
 By virtue of \rdefn{D0754} and \eqref{F.G.f012}, we have \(M=\cia{M;0}^\pm=\su{M;0}^\pm\).
 In particular \(\OPu{\ran{\su{M;0}^\pm}}=1\).
 Taking additionally into account \(\rk{\cia{M;1}^\pm-\frac{1}{2}}^2=\frac{\ba^2-4}{4\ba^2}=\frac{1}{4}-\ba^{-2}\), hence \(\cia{M;1}^\pm-\rk{\cia{M;1}^\pm}^2=\frac{1}{4}-\rk{\cia{M;1}^\pm-\frac{1}{2}}^2=\ba^{-2}\OPu{\ran{\su{M;0}^\pm}}\) follows.
 According of \rcor{R1124}, then the sequences \(\seq{\su{M;j}^-}{j}{0}{\infi}\) and \(\seq{\su{M;j}^+}{j}{0}{\infi}\) are both fixed points of \(\Phi\). 
 
 Now consider an arbitrary fixed point \(\seqsinf\in\Fgguuuu{1}{\infi}{\ug}{\obg}\) of \(\Phi\).
 Because of \rlem{F.R.Fgg-s}, we have \(\su{0}\in[0,\infp)\).
 If \(\su{0}=0\), then, as above, the sequence \(\seqsinf\) coincides with \(\seq{\su{0;j}}{j}{0}{\infi}\).
 Assume \(\su{0}>0\).
 Then \(\OPu{\ran{\su{0}}}=1\).
 Denote by \(\seqciainf\) the \tfcfa{\(\seqsinf\)}.
 Because of \rcor{R1124}, we get \(\cia{1}-\cia{1}^2=\ba^{-2}\OPu{\ran{\su{0}}}\) and \(\cia{1}=\cia{2}=\cia{3}=\dotsb\)
 Then \(\rk{\cia{1}-\frac{1}{2}}^2=\cia{1}^2-\cia{1}+\frac{1}{4}=-\ba^{-2}\OPu{\ran{\su{0}}}+\frac{1}{4}=-\frac{1}{\ba^2}+\frac{1}{4}=\frac{\ba^2-4}{4\ba^2}\) follows, \tie{}\ \(\cia{1}=\frac{1}{2}-\frac{1}{2\ba}\sqrt{\ba^2-4}\) or \(\cia{1}=\frac{1}{2}+\frac{1}{2\ba}\sqrt{\ba^2-4}\).
 By virtue of \rdefn{D0754} and \eqref{F.G.f012}, we have \(\cia{0}=\fpu{0}=\su{0}\).
 Setting \(M\defeq\su{0}\), then \(M\in(0,\infp)\) and \(\seqciainf\) coincides with \(\seq{\cia{M;j}^-}{j}{0}{\infi}\) or with \(\seq{\cia{M;j}^+}{j}{0}{\infi}\) defined above, implying that \(\seqsinf\) coincides with \(\seq{\su{M;j}^-}{j}{0}{\infi}\) or with \(\seq{\su{M;j}^+}{j}{0}{\infi}\), according to \rthm{F.T.Fggcia}.
\eproof

\breml{P1338}
 Let \(m\in\NO\) and let \(\seqs{m}\in\Fggqu{m}\).
 Then the following statements are equivalent:
\baeqi{0}
 \il{P1338.i} \(\seqs{m}\in\Fgqu{m}\).
 \il{P1338.ii} \(\dia{j}\in\Cgq\) for all \(j\in\mn{0}{m}\).
 \il{P1338.iii} \(\det\dia{j}\neq0\) for all \(j\in\mn{0}{m}\).
 \il{P1338.iv} \(\dia{m}\in\Cgq\).
 \il{P1338.v} \(\det\dia{m}\neq0\).
\eaeqi
 Indeed, in view of \rprop{ab.C0929} we see that~\ref{P1338.ii} and~\ref{P1338.iii} \tresp{}~\ref{P1338.iv} and~\ref{P1338.v} are equivalent.
 Furthermore, for each \(j\in\mn{0}{m}\), from \zitaa{MR3775449}{\cprop{10.23}{203}} we obtain \(\dia{j}\lgeq\rk{4/\ba}^{m-j}\dia{m}\).
 Thus,~\ref{P1338.iv} implies~\ref{P1338.ii}.
 The equivalence of~\ref{P1338.iii} and~\ref{P1338.i} is a consequence of \rlem{F.L.dHdia}. 
\erem

\breml{C1127} 
 Let \(\seqsinf\in\Fggqinf\).
 In view of \rrem{P1338}, then the following statements are equivalent:
\baeqi{0}
 \il{C1127.i} \(\seqsinf\in\Fgqinf\).
 \il{C1127.ii} \(\dia{j}\in\Cgq\) for all \(j\in\NO\).
 \il{C1127.iii} \(\det\dia{j}\neq0\) for all \(j\in\NO\).
\eaeqi
\erem

\section{The class \(\RFqab\)}\label{F.s2.RFab}
 In this section, we introduce several classes of matrix-valued functions, holomorphic in the sense explained in \rapp{B.S.hol}.
 We consider the following open half-planes in the complex plane: \(\Hla\defeq\setaca{z\in\C}{\re z<\ug}\), \(\Hrb\defeq\setaca{z\in\C}{\re z>\obg}\), \(\lhp\defeq\setaca{z\in\C}{\im z<0}\), and \(\uhp\defeq\setaca{z\in\C}{\im z>0}\).
 Furthermore, we write \(\re A\defeq\frac{1}{2}\rk{A+A^\ad}\) and \(\im A\defeq\frac{1}{2\iu}\rk{A-A^\ad}\) for the real part and the imaginary part of a complex square matrix \(A\), \tresp{}

\bnotal{H.N.NF}
 Denote by \(\NFq\) the set of all matrix-valued functions \(F\colon\uhp\to\Cqq\), which are holomorphic and satisfy \(\im F(z)\in\Cggq\) for all \(z\in\uhp\).
\enota

 By means of
\[
 G(z)
 \defeq
 \begin{cases}
  F(z)\tincase{\im z>0}\\
  [F(\ko z)]^\ad\tincase{\im z<0}
 \end{cases},
\]
 the matrix-valued functions \(F\) of the class \(\NFq\) can be extended to holomorphic matrix-valued functions \(G\colon\C\setminus\R\to\Cqq\), which satisfy \(\im G(z)/\im z\in\Cggq\) for all \(z\in\C\setminus\R\).
 In the scalar case \(q=1\), such a function is called \emph{\(R\)-function} in~\zitaa{zbMATH03455679}{\cpage{1}}.
 The matrix-valued functions of the class \(\NFq\) are also called \emph{Herglotz functions}, \emph{Nevanlinna functions} or \emph{Pick functions}.
 They admit a well-known integral representation, the scalar version of which can be found, \teg{}, in~\zitaa{zbMATH03455679}{\ceq{S1.1.1}{2}}.

 Using the Euclidean norm \(\normE{x}\defeq\sqrt{x^\ad x}\) on \(\Cq\) corresponding to the Euclidean scalar product, we define the operator norm 
\(
 \normS{A}
 \defeq\max\setaca{\normE{Au}}{u\in\Cq\text{ with }\normE{u}=1}
\)
 on \(\Cpq\) induced by the Euclidean norms on \(\Cq\) and \(\Cp\), which is also called \emph{spectral norm} on \(\Cpq\).

\bnotal{N1409}
 Denote by \(\RFOq\) the set of all \(F\in\NFq\) satisfying the growth condition \(\sup_{y\in[1,\infp)}y\normS{F(\iu y)}<\infp\).
\enota

\bthmnl{\tcf{}~\zitaa{MR2222521}{\cthm{8.7}{167}}}{T8-7}
\benui
 \il{T8-7.a} If \(F\in\RFOq\), then there exists a unique \(\sigma\in\MggqR\) such that
\beql{T8-7.A}
 F(z)
 =\int_\R\frac{1}{t-z}\sigma(\dif t)
\eeq
 holds true for all \(z\in\uhp\).
 \il{T8-7.b} If \(\sigma\in\MggqR\), then \(F\colon\uhp\to\Cqq\) defined via \eqref{T8-7.A} belongs to \(\RFOq\).
\eenui
\ethm

\bdefnl{D1441}
 Let \(F\in\RFOq\).
 Then the unique measure \(\sigma\in\MggqR\) such that \eqref{T8-7.A} holds true for all \(z\in\uhp\) is called the \emph{(matricial) \tSpMav{\(F\)}} and is denoted by \(\sigmaF\).
\edefn

 In certain situations, an upper bound for \(y\normS{F(\iu y)}\) can be obtained, using \rlem{L.AEP}:

\blemnl{\zitaa{MR2222521}{\clem{8.9}{168}}}{CRDFK06.L8-9}
 Let \(M\in\Cqq\) and let \(F\colon\uhp\to\Cqq\) be a holomorphic matrix-valued function such that, for all \(z\in\uhp\), the matrix \(\smat{M&F(z)\\\ek{F(z)}^\ad&\frac{1}{\im z}\im F(z)}\) is \tnnH{}.
 Then \(F\in\RFuq{0}\) with \(\sup_{y\in(0,\infp)}y\normS{F(\iu y)}\leq\normS{M}\) and \(\sigmaFa{\R}\lleq M\).
\elem

 Now we introduce that class of holomorphic matrix-valued functions, which is relevant for the \rmprobm{\ab}{\kappa}{=}.

\bnotal{ab.N1603}
 Denote by \(\RFqab\) the set of all matrix-valued functions \(F\colon\Cab\to\Cqq\) which are holomorphic and satisfy the following conditions:
\bAeqi{0}
 \il{ab.N1603.I} \(\im F(z)\in\Cggq\) for all \(z\in\uhp\).
 \il{ab.N1603.II} \(F(x)\in\Cggq\) for all \(x\in\crhl\) and \(-F(x)\in\Cggq\) for all \(x\in\clhl\).
\eAeqi
\enota

 Since such functions are holomorphic in \(\Cab\) with \tnnH{} imaginary part in \(\uhp\), we can think of \(\RFqab\) as a subclass of \(\NFq\), by virtue of the identity theorem for holomorphic functions:

\breml{F.R.Rab<NF}
 By means of restricting matrix-valued functions of \(\RFqab\) to \(\uhp\), an injective mapping from \(\RFqab\) into \(\NFq\) is given.
\erem

 \rlem{F.L.RabSp} will show that the above mentioned restrictions even belong to \(\RFOq\).
 Since each function from \(\RFqab\) is holomorphic in \(\Cab\) with \tH{} values in \(\cab\), we can use the Schwarz reflection principle to obtain the following relation connecting its values on the open upper and lower half-plane:
 
\breml{ab.R1652}
 If \(F\in\RFqab\), then \(\ek{F(z)}^\ad=F(\ko z)\) for all \(z\in\Cab\).
\erem

 The latter result can also be seen from the following integral representation:
 
\bthmnl{\tcf{}~\zitaa{MR2222521}{\cthm{1.1}{124}}}{ab.T1614}
\benui
 \il{ab.T1614.a} If \(F\in\RFqab\), then there exists a unique \(\rabmf\in\MggqF\) such that
\beql{ab.T1614.B}
 F(z)
 =\int_\ab\frac{1}{t-z}\rabmf(\dif t)
\eeq
 holds true for all \(z\in\Cab\).
 \il{ab.T1614.b} If \(\rabmf\in\MggqF\), then \(F\colon\Cab\to\Cqq\) defined via \eqref{ab.T1614.B} belongs to \(\RFqab\).
\eenui
\ethm

\bdefnl{F.D.RabMass}
 Let \(F\in\RFqab\).
 In view of \rthm{ab.T1614}, let \(\rabmf\) be the uniquely determined measure from \(\MggqF\) such that \eqref{ab.T1614.B} holds true for all \(z\in\Cab\).
 Then \(\rabmf\) is called the \emph{\tRabMav{\(F\)}} and is denoted by \(\rabmF\).
\edefn

 As already mentioned in \rsec{S1506}, the power moments of a measure belonging to \(\MggqF\) exist for each non-negative integer order:

\breml{F.R.RabM8}
 If \(F\in\RFqab\), then the \tRabMa{} \(\rabmF\) of \(F\) belongs to \(\MggqinfF\).
\erem

 In view of \rrem{F.R.Rab<NF}, we infer the following relation to the class \(\RFOq\):

\bleml{F.L.RabSp}
 Let \(F\in\RFqab\) with \tRabMa{} \(\rabmF\) and denote by \(f\) the restriction of \(F\) onto \(\uhp\).
 Then \(f\in\RFOq\) and the \tSpMa{} \(\sigmau{f}\) of \(f\) fulfills \(\sigmaua{f}{\cab}=\Oqq\) and \(\sigmaua{f}{B}=\rabmFa{B}\) for all \(B\in\BsAF\).
\elem
\bproof
 Use \rthmss{ab.T1614}{T8-7}.
\eproof

 By virtue of \rthm{ab.T1614}, the following integral representations are readily checked:

\breml{ab.R0857}
 Let \(F\in\RFqab\).
 For all \(z\in\Cab\), then
\begin{align*}
 \re F(z)&=\int_\ab\frac{t-\re z}{\abs{t-z}^2}\rabmFa{\dif t}&
&\text{and}&
 \im F(z)&=\int_\ab\frac{\im z}{\abs{t-z}^2}\rabmFa{\dif t}.
\end{align*}
\erem

 Now we state a useful characterization of the functions belonging to \(\RFqab\).

\bpropl{F.P.RabH}
 Let \(F\colon\Cab\to\Cqq\) be holomorphic.
 Then \(F\in\RFqab\) if and only if the following two conditions are fulfilled:
\bAeqi{0}
 \il{F.P.RabH.I} \(\frac{1}{\im z}\im F(z)\in\Cggq\) for all \(z\in\C\setminus\R\).
 \il{F.P.RabH.II} \(\re F(w)\in\Cggq\) for all \(w\in\Hla\) and \(-\re F(w)\in\Cggq\) for all \(w\in\Hrb\).
\eAeqi
\eprop
\bproof
 If \(F\in\RFqab\), then~\rstat{F.P.RabH.I} and~\rstat{F.P.RabH.II} are readily seen from \rrem{ab.R0857}.
 
 Conversely, suppose that~\rstat{F.P.RabH.I} and~\rstat{F.P.RabH.II} are fulfilled.
 Due to~\rstat{F.P.RabH.I}, we have \(\im F(z)\in\Cggq\) for all \(z\in\uhp\).
 As in the proof of \rlem{B.L.STOJ}, we can conclude from~\rstat{F.P.RabH.I} that, for all \(x\in\cab\), the equation \(\im F(x)=\Oqq\) holds true, implying \(F(x)=\re F(x)\).
 Taking into account~\rstat{F.P.RabH.II}, we thus have \(F(x)\in\Cggq\) for all \(x\in\crhl\) and \(-F(x)\in\Cggq\) for all \(x\in\clhl\).
 Regarding \rnota{ab.N1603}, hence \(F\in\RFqab\).
\eproof

 As an immediate consequence of the following result, we see that the column space \(\ran{F\rk{z}}\) and the null space \(\nul{F\rk{z}}\) of a matrix-valued function \(F\in\RFqab\) are both independent of the argument \(z\in\Cab\).
 
\bpropl{ab.P1648L1409}
 Let \(F\in\RFqab\).
 Then:
\benui
 \il{ab.P1648L1409.a} \(\ran{F(z)}=\ran{\rabmFa{\ab}}\) and \(\nul{F(z)}=\nul{\rabmFa{\ab}}\) 
 for all \(z\in\Cab\).
 \il{ab.P1648L1409.b} \(\ran{\im F(z)}=\ran{\rabmFa{\ab}}\) and \(\nul{\im F(z)}=\nul{\rabmFa{\ab}}\) 
 for all \(z\in\C\setminus\R\).
 \il{ab.P1648L1409.c} \(\ran{\re F(w)}=\ran{\rabmFa{\ab}}\) and \(\nul{\re F(w)}=\nul{\rabmFa{\ab}}\) 
 for all \(w\in\Hla\cup\Hrb\).
\eenui
\eprop
\bproof
 In view of \rthm{ab.T1614}, this follows from \rlem{B.L.rn} applied with \(\Omega=\ab\).
\eproof

 We recall the definitions of two well-studied classes of matrices.

\bdefnl{D.EPAD}
 Let \(A\) be a complex \tqqa{matrix}.
 Then \(A\) is called \emph{\tEP{matrix}} if \(\ran{A^\ad}=\ran{A}\).
 Furthermore, the matrix \(A\) is said to be \emph{\tAD} if each \(x\in\Cq\) with \(x^\ad Ax=0\) necessarily fulfills \(Ax=\Ouu{q}{1}\).
 Denote by \(\CEPq\) and \(\CADq\) the set of \tEP{matrices} and the set of \tAD{} matrices from \(\Cqq\), \tresp{}
\edefn

\bpropl{P1834}
 If \(F\in\RFqab\), then \(F(z)\in\CADq\) for all \(z\in\Cab\).
\eprop
\bproof
 In view of \rthm{ab.T1614}, this follows from \rlem{L1746} applied with \(\Omega=\ab\).
\eproof
 
 According to \rrem{L.AD<EP}, we have \(\CADq\subseteq\CEPq\).
 Hence, \rprop{P1834} implies that the values of a function \(F\in\RFqab\) fulfill \(\ran{\ek{F(z)}^\ad}=\ran{F\rk{z}}\) for all \(z\in\Cab\), a fact that also can be seen from \rprop{ab.P1648L1409} in combination with \rrem{ab.R1652}.

 By means of \rthm{ab.T1614}, a characterization of \(\RFqab\) in terms of the class \(\NFq\) can be obtained:

\bpropnl{\tcf{}~\zitaa{MR2222521}{\clem{3.6}{133}}}{ab.L1642}
 Let \(F\colon\Cab\to\Cqq\) be holomorphic and let the matrix-valued functions \(g,h\colon\uhp\to\Cqq\) be defined by \(g(z)\defeq\rk{z-\ug}F(z)\) and \(h(z)\defeq\rk{\obg-z}F(z)\), \tresp{}
 Then \(F\in\RFqab\) if and only if \(g\) and \(h\) both belong to \(\NFq\).
\eprop

 In view of \rrem{F.R.RabM8}, we can associate to a given function from the class \(\RFqab\) three auxiliary functions, which are intimately connected to the three sequences of complex matrices introduced in \rnota{F.N.sa} (\tcf{}~\rrem{ab.L0921}):
 
\bnotal{ab.N1537}
 Let \(F\in\RFqab\) with \tRabMa{} \(\rabmF\).
 Then let the functions \(\Ffa,\Fb,\Fc\colon\Cab\to\Cqq\) be defined by
\begin{align*}
 \Fav{z}&\defeq\rk{z-\ug}F(z)+\rabmFA{\ab},&
 \Fbv{z}&\defeq\rk{\obg-z}F(z)-\rabmFA{\ab},
\end{align*}
 and
\[
 \Fcv{z}
 \defeq\rk{\obg-z}\rk{z-\ug}F(z)+\rk{\ug+\obg-z}\rabmFA{\ab}-\int_\ab t\rabmFa{\dif t}.
\]
\enota

\bpropl{ab.L1401}
 Let \(F\in\RFqab\) with \tRabMa{} \(\rabmF\).
 Then \(\Ffa\), \(\Fb\), and \(\Fc\) belong to \(\RFqab\) and their \tRabMa{s} \(\rabmFfa\), \(\rabmFb\), and \(\rabmFc\) fulfill
\begin{align*}
 \rabmFav{B}&=\int_B\rk{t-\ug}\rabmFa{\dif t},&
 \rabmFbv{B}&=\int_B\rk{\obg-t}\rabmFa{\dif t},&
 \rabmFcv{B}&=\int_B\rk{\obg-t}\rk{t-\ug}\rabmFa{\dif t}
\end{align*}
 for all \(B\in\BsAF\).
\eprop
\bproof
 Because of \rrem{F.R.RabM8}, the integrals \(\int_\ab t\rabmFa{\dif t}\) and \(\int_\ab t^2\rabmFa{\dif t}\) exist.
 Since \(t-\ug>0\) and \(\obg-t>0\) hold true for all \(t\in\ab\), we can thus conclude that \(\int_B\rk{t-\ug}\rabmFa{\dif t}\), \(\int_B\rk{\obg-t}\rabmFa{\dif t}\), and \(\int_B\rk{\obg-t}\rk{t-\ug}\rabmFa{\dif t}\) are \tnnH{} matrices for all \(B\in\BsAF\).
 Consequently, \(\rabmFfa\), \(\rabmFb\), and \(\rabmFc\) belong to \(\MggqF\).
 Consider now an arbitrary \(z\in\Cab\).
 In view of \rthm{ab.T1614}, we have
\[
 \rk{z-\ug}F(z)
 =\int_\ab\frac{z-\ug}{t-z}\rabmFa{\dif t}
 =\int_\ab\rk*{\frac{t-\ug}{t-z}-1}\rabmFa{\dif t}
 =\int_\ab\frac{t-\ug}{t-z}\rabmFa{\dif t}-\rabmFA{\ab}
\]
 and similarly \(\rk{\obg-z}F(z)=\int_\ab\rk{t-z}^\inv\rk{\obg-t}\rabmFa{\dif t}+\rabmFa{\ab}\).
 Due to \rthm{ab.T1614}, then \(\Ffa\) and \(\Fb\) belong to \(\RFqab\) having the asserted \tRabMa{s}.
 Using the representation above, we obtain furthermore
\[\begin{split}
 \rk{\obg-z}\rk{z-\ug}F(z)
 &=\int_\ab\frac{\rk{\obg-z}\rk{t-\ug}}{t-z}\rabmFa{\dif t}-\rk{\obg-z}\rabmFA{\ab}\\
 &=\int_\ab\rk*{\frac{\obg-t}{t-z}+1}\rk{t-\ug}\rabmFa{\dif t}-\rk{\obg-z}\rabmFA{\ab}\\
 &=\int_\ab\frac{\rk{\obg-t}\rk{t-\ug}}{t-z}\rabmFa{\dif t}+\int_\ab t\rabmFa{\dif t}-\rk{\ug+\obg-z}\rabmFA{\ab}.
\end{split}\]
 Hence, \(\Fcv{z}=\int_\ab\frac{1}{t-z}\rabmFcv{\dif t}\).
 By virtue of \rthm{ab.T1614}, thus \(\Fc\) belongs to \(\RFqab\) having the asserted \tRabMa{}.
\eproof

 The combination of \rprop{ab.L1401} with \rremss{ab.R0857}{A.R.XRIX} yields:

\breml{ab.R0922}
 Let \(F\in\RFqab\).
 For all \(z\in\Cab\), then
\begin{align*}
 \im\ek*{\rk{z-\ug}F(z)}&=\im\Fav{z}=\im\rk{z}\int_\ab\frac{t-\ug}{\abs{t-z}^2}\rabmFa{\dif t},\\
 \im\ek*{\rk{\obg-z}F(z)}&=\im\Fbv{z}=\im\rk{z}\int_\ab\frac{\obg-t}{\abs{t-z}^2}\rabmFa{\dif t},
\shortintertext{and}
 \im\ek*{\rk{\obg-z}\rk{z-\ug}F(z)}&=\im\Fcv{z}+\im\rk{z}\rabmFA{\ab}\\
 &=\im\rk{z}\ek*{\rabmFA{\ab}+\int_\ab\frac{\rk{\obg-t}\rk{t-\ug}}{\abs{t-z}^2}\rabmFa{\dif t}}.
\end{align*}
\erem

 Using \rrem{ab.R0922}, the following result is readily checked:

\breml{ab.L1005}
 If \(F\in\RFqab\), then \(\frac{1}{\im z}\im\ek{\rk{\obg-z}\rk{z-\ug}F(z)}\lgeq\rabmFa{\ab}\lgeq\Oqq\) for all \(z\in\C\setminus\R\).
\erem

\section{An equivalent problem in the class \(\RFqab\)}\label{S1723}

 To describe the solution set of moment problems on the real axis, the transition to holomorphic functions by means of the \tSTion{} considered in detail in \rapp{B.s1.HN} has turned out to be very helpful.
 For the sake of a simpler description of the relation to \tNF{s} treated in the previous section, as usual, we choose in the case \(\Omega=\R\) for the \tST{} \(S\) the domain \(\uhp\) instead of \(\C\setminus\R\):

\bdefnl{H.D.ST}
 Let \(\sigma\in\MggqR\).
 Then we call the matrix-valued function \(\STHu{\sigma}\colon\uhp\to\Cqq\) defined by
\[
 \STHua{\sigma}{z}
 \defeq\int_\R\frac{1}{t-z}\sigma(\dif t)
\]
 the \emph{\tSTHv{\(\sigma\)}}.
\edefn

 From \rthm{T8-7} and \rdefnss{D1441}{H.D.ST} we immediately see the well-known connection of the \tSTHion{} to the class \(\RFuq{0}\):

\bpropl{H.L.ST-R0}
\benui
 \il{H.L.ST-R0.a} If \(F\in\RFuq{0}\), then there exists a unique \(\sigma\in\MggqR\) fulfilling \(F=\STHu{\sigma}\), namely \(\sigma=\sigmaF\).
 \il{H.L.ST-R0.b} If \(\sigma\in\MggqR\), then the \tSTH{} \(\STHu{\sigma}\) of \(\sigma\) belongs to \(\RFuq{0}\).
\eenui
\eprop

 According to our interest in the matricial Hausdorff moment problem, we consider the integral transformation \eqref{B.G.STO} for the particular case of \tnnH{} measures \(\sigma\) belonging to \(\MggqF\):

\bdefnl{F.D.ST}
 Let \(\sigma\in\MggqF\).
 Then we call the matrix-valued function \(\STFu{\sigma}\colon\Cab\to\Cqq\) defined by
\beql{F.D.ST.1}
 \STFua{\sigma}{z}
 \defeq\int_\ab\frac{1}{t-z}\sigma(\dif t)
\eeq
 the \emph{\tSTFv{\(\sigma\)}}.
\edefn

 The \tSTF{} of a \tnnH{} measure from \(\MggqF\) admits a power series representation at \(z_0=\infc\) involving the corresponding moments:

\bpropnl{\zitaa{CR01}{\cSatz{1.2.16}{34}}}{CR01.S1.2.16}
 Let \(\sigma\in\MggqF\).
 Then the moments \(\mpm{\sigma}{j}\defeq\int_\ab x^j\sigma(\dif x)\) exist for all \(j\in\NO\).
 For each \(z\in\C\) with \(\abs{z}>\max\set{\abs{\ug},\abs{\obg}}\), furthermore \(z\in\Cab\) and
\[
  \STFua{\sigma}{z}
  =-\sum_{j=0}^\infi z^{-(j+1)}\mpm{\sigma}{j}.
\]
\eprop

 The following reformulation of \rthm{ab.T1614} describes the relation between \tSTF{} \(\STFu{\sigma}\) and \tRabMa{} \(\rabmF\):

\bpropl{F.P.STbij}
 The mapping \(\Lambda_\ab\colon\MggqF\to\RFqab\) given by \(\sigma\mapsto\STFu{\sigma}\), where \(\STFu{\sigma}\) is given by \eqref{F.D.ST.1}, is well defined and bijective.
 Its inverse \(\Lambda_\ab^\inv\colon\RFqab\to\MggqF\) is given by \(F\mapsto\rabmF\), where \(\rabmF\) denotes the \tRabMav{F}.
\eprop

 By virtue of \rprop{F.P.STbij}, the \rmprobm{\ab}{\kappa}{=} admits a reformulation as an equivalent problem for functions belonging to the class \(\RFqab\):

\begin{Problem}[\rabproblem{\kappa}]
 Given a sequence \(\seqska \) of complex \tqqa{matrices}, parametrize the set \(\RFqabskag\) of all \(F\in\RFqab\) with \tRabMa{} \(\rabmF\) belonging to \(\MggqFksg\).
\end{Problem}

 In particular, \rprobrab{\kappa} has a solution if and only if the \rmprobm{\ab}{\kappa}{=} has a solution.
 From \rthm{I.P.ab} we can therefore conclude:

\bpropl{F.P.FPsolv}
 Let \(\seqska\) be a sequence of complex \tqqa{matrices}.
 Then the set \(\RFqabskag\) is non-empty if and only if the sequence \(\seqska\) belongs to \(\Fggqka\).
\eprop

 In view of \rprop{F.P.STbij}, the solution set \(\RFqabskag\) of \rprobrabs{\kappa} can also be described in the following way:

\breml{F.R.SF-ST}
 \(\RFqabskag=\setaca{\STFu{\sigma}}{\sigma\in\MggqFksg}\).
\erem

 In combination with \rpropss{I.P.ab8}{CR01.S1.2.16}, we can infer from \rrem{F.R.SF-ST} in particular:

\bpropl{F.P.FP8}
 Let \(\seqsinf\in\Fggqinf\).
 Then the set \(\RFqabsg{\infi}\) consists of exactly one element \(S\).
 For all \(z\in\C\) with \(\abs{z}>\max\set{\abs{\ug},\abs{\obg}}\), furthermore \(z\in\Cab\) and
\beql{F.P.FP8.B}
 S(z)
 =-\sum_{j=0}^\infi z^{-(j+1)}\su{j}.
\eeq
\eprop

 In addition, we have the following result:

\bpropl{F.P.R-F}
 Let \(\seqsinf\in\Fggqinf\), let \(F\colon\Cab\to\Cqq\) be holomorphic, and let \(\rho\in\R\) with \(\rho\geq\max\set{\abs{\ug},\abs{\obg}}\).
 Suppose that \(F(z)=-\sum_{j=0}^\infi z^{-(j+1)}\su{j}\) holds true for all \(z\in\C\) with \(\abs{z}>\rho\).
 Then \(F\in\RFqabsg{\infi}\).
\eprop
\bproof
 Due to \rprop{F.P.FP8}, the set \(\RFqabsg{\infi}\) consists of exactly one element \(S\) and \eqref{F.P.FP8.B} holds true for all \(z\in\C\) with \(\abs{z}>\max\set{\abs{\ug},\abs{\obg}}\).
 In particular, \(F(z)=S(z)\) for all \(z\in\C\) with \(\abs{z}>\rho\) follows.
 Consequently, the application of the identity theorem for holomorphic functions yields \(F=S\).
 Therefore, \(F\) belongs to \(\RFqabsg{\infi}\).
\eproof

 In view of \rprop{ab.L1401}, the sequences \(\seqsa{\kappa-1}\), \(\seqsb{\kappa-1}\), and \(\seqsab{\kappa-2}\) introduced in \rnota{F.N.sa} consist of the first power moments of the \tRabMa{s} \(\rabmFfa\), \(\rabmFb\), and \(\rabmFc\) of the matrix-valued functions \(\Ffa\), \(\Fb\), and \(\Fc\) built, according to \rnota{ab.N1537}, from a given function \(F\in\RFqabskag\):
 
\breml{ab.L0921}
 Let \(\seqska\in\Fggqka\) and let \(F\in\RFqabskag\).
 If \(\kappa\geq1\), then \(\Ffa\in\RFqabg{\seqsa{\kappa-1}}\) and \(\Fb\in\RFqabg{\seqsb{\kappa-1}}\).
 If \(\kappa\geq2\), then \(\Fc\in\RFqabg{\seqsab{\kappa-2}}\).
\erem

\section{A Schur--Nevanlinna type algorithm in the class \(\RFqab\)}\label{S1323}

 On the background of \rprop{I.P.ab8Fgg} and \rthm{F.T.Fggcia}, we parametrized in~\zitaa{MR3979701}{\csec{8}} the set \(\MggqF\) of \tnnH{} \tqqa{measures} on \(\ab\) and generalized several results from the scalar theory of canonical moments (\tcf{}~\zita{MR1468473}) to the matrix case.
 To that end, we associated to such a measure the sequences built via \rdefnss{D1861}{D0754} from its \tfpmf{}:

\bdefnl{F.D.meacia}
 Let \(\sigma\in\MggqF\) with \tfpmf{} \(\seqmpm{\sigma}\).
 Denote by \(\seqmcm{\sigma}\) the \tfcfa{\(\seqmpm{\sigma}\)} and by \(\seqmdm{\sigma}\) the \tfdfa{\(\seqmpm{\sigma}\)}.
 Then we call \(\seqmcm{\sigma}\) the \emph{\tfmcmfa{\(\sigma\)}} and we say that \(\seqmdm{\sigma}\) is the \emph{\tfmdmfa{\(\sigma\)}}.
\edefn

\bthmnl{\zitaa{MR3979701}{\cthm{8.2}{2160}}}{F.T.Mabcia}
 The mapping \(\Pi_\ab\colon\MggqF\to\es{q}{\infi}{\bam}\) given by \(\sigma\mapsto\seqmcm{\sigma}\) is well defined and bijective.
\ethm

 On the basis of \rprop{F.P.STbij}, a parametrization of the class \(\RFqab\) immediately follows:

\bdefnl{F.D.funcia}
 Let \(F\in\RFqab\) with \tRabMa{} \(\rabmF\).
 Denote by \(\seqRabSp{F}\) the \tfmcmfa{\(\rabmF\)} and by \(\seqRabil{F}\) the \tfmdmfa{\(\rabmF\)}.
 Then we call \(\seqRabSp{F}\) the \emph{\tfRabSp{\(F\)}} and we say that \(\seqRabil{F}\) is the \emph{\tfRabil{\(F\)}}.
\edefn

\bthml{F.T.Rabcia}
 The mapping \(\Delta_\ab\colon\RFqab\to\es{q}{\infi}{\bam}\) given by \(F\mapsto\seqRabSp{F}\) is well defined and bijective.
\ethm
\bproof
 Use \rprop{F.P.STbij} and \rthm{F.T.Mabcia}.
\eproof

 By means of this one-to-one correspondence, results obtained in~\zitas{MR3979701,MR4051874} on \tmcm{s} associated with \tnnH{} measures from \(\MggqF\) carry over to matrix-valued functions belonging to the class \(\RFqab\) and their \tRabSp{s}.

 By virtue of \rpropss{I.P.ab8}{ab.P1030}, the \tFTion{} (see~ \rdefnss{ab.N0940}{ab.N1020}) for \tFnnd{} sequences of matrices gave rise to a transformation considered in~\zitaa{MR4051874}{\cdefn{10.1}{202}} for \tnnH{} measures from \(\MggqF\):

\bdefnl{F.D.SN-M}
 Let \(\sigma\in\MggqF\) with \tfpmf{} \(\seqsinf\) and let \(k\in\NO\).
 Denote by \(\seq{\su{j}^\FTa{k}}{j}{0}{\infi}\) the \tnFTv{k}{\(\seqsinf\)} and by \(\sigma^\FTa{k}\) the uniquely determined element in \(\MggqFag{s^\FTa{k}}{\infi}\).
 Then we call \(\sigma^\FTa{k}\) the \emph{\tnFmTv{k}{\(\sigma\)}}.
\edefn

\breml{R0713} 
 Let \(\sigma\in\MggqF\).
 According to \rdefnss{F.D.SN-M}{ab.N1020}, then \(\sigma^\FTa{0}=\sigma\) and \(\sigma^\FTa{k}\) is exactly the first \tFmTv{\(\sigma^\FTa{k-1}\)} for each \(k\in\N\).
\erem

\breml{R0718} 
 Let \(\seqska\in\Fggqka\), let \(\sigma\in\MggqFksg\), and let \(k\in\mn{0}{\kappa}\).
 Then, in view of \rdefn{F.D.SN-M} and \rrem{F.R.kFTtr}, it is readily checked that \(\sigma^\FTa{k}\) belongs to \(\MggqFag{s^\FTa{k}}{\kappa-k}\).
\erem

 In view of \rprop{F.P.STbij}, we can define a corresponding transformation for functions belonging to the class \(\RFqab\):

\bdefnl{F.D.SN-F}
 Let \(F\in\RFqab\) with \tRabMa{} \(\sigma\) and let \(k\in\NO\).
 Denote by \(\sigma^\FTa{k}\) the \tnFmTv{k}{\(\sigma\)} and by \(F^\FTa{k}\) the \tSTFv{\(\sigma^\FTa{k}\)}.
 Then we call \(F^\FTa{k}\) the \emph{\tnSfTv{k}{\(F\)}}.
\edefn

\breml{R1615}
 In the situation of \rdefn{F.D.SN-F} we see from \rprop{F.P.STbij} that \(F^\FTa{k}\) belongs to \(\RFqab\) and that \(\sigma^\FTa{k}\) is the \tRabMav{\(F^\FTa{k}\)}.
\erem

\breml{R0730} 
 Let \(F\in\RFqab\).
 Then, regarding \rdefn{F.D.SN-F}, \rrem{R0713}, and \rprop{F.P.STbij}, it is readily checked that \(F^\FTa{0}=F\) and that \(F^\FTa{k}\) is exactly the first \tSfTv{\(F^\FTa{k-1}\)} for each \(k\in\N\).
\erem

\breml{R0734} 
 Let \(\seqska\in\Fggqka\), let \(F\in\RFqabsg{\kappa}\), and let \(k\in\mn{0}{\kappa}\).
 Because of \rremss{R0718}{R1615}, then \(F^\FTa{k}\in\RFqabag{s^\FTa{k}}{\kappa-k}\).
\erem

 One of the results in~\zita{MR4051874} states that the \tFmT{ation} of a \tnnH{} measure from \(\MggqF\) is essentially equivalent to left shifting its \tfmcmf{}:
 
\bpropnl{\zitaa{MR4051874}{\cprop{10.4}{203}}}{F.R.mdealg}
 Let \(k\in\NO\) and let \(\sigma\in\MggqF\) with \tnFmT{k} \(\mu\).
 Then \(\mcm{\mu}{0}=\ba^{k-1}\mdm{\sigma}{k}\) and \(\mcm{\mu}{j}=\mcm{\sigma}{k+j}\) for all \(j\in\N\).
 Furthermore, \(\mdm{\mu}{j}=\ba^k\mdm{\sigma}{k+j}\) for all \(j\in\NO\).
 In particular, \(\mu\rk{\ab}=\ba^{k-1}\mdm{\sigma}{k}\).
\eprop
 
 The following analogous result for matrix-valued functions from the class \(\RFqab\) justifies the notions \emph{\tRabSp{s}} and \emph{\tSfT{}} chosen in \rdefnss{F.D.funcia}{F.D.SN-F}, \tresp{}
 
\bpropl{F.R.fdealg}
 Let \(k\in\NO\) and let \(F\in\RFqab\) with \tnSfT{k} \(G\).
 Then \(\RabSp{G}{0}=\ba^{k-1}\Rabil{F}{k}\) and \(\RabSp{G}{j}=\RabSp{F}{k+j}\) for all \(j\in\N\).
 Furthermore, \(\Rabil{G}{j}=\ba^k\Rabil{F}{k+j}\) for all \(j\in\NO\).
\eprop
\bproof
 In view of \rdefnss{F.D.funcia}{F.D.SN-F} and \rrem{R1615}, this is an immediate consequence of \rprop{F.R.mdealg}.
\eproof
 
 Let \(\Omega\in\BsAR\setminus\set{\emptyset}\).
 A \tnnH{} measure \(\sigma\in\MggqO\) is said to be \emph{molecular} if there exists a finite subset \(B\) of \(\Omega\) satisfying \(\sigma\rk{\Omega\setminus B}=\Oqq\).
 Obviously, this is equivalent to the existence of an \(m\in\N\) and sequences \(\seq{\xi_\ell}{\ell}{1}{m}\) and \(\seq{A_\ell}{\ell}{1}{m}\) from \(\Omega\) and \(\Cggq\), \tresp{}, such that \(\sigma=\sum_{\ell=1}^m\Kronu{\xi_\ell}A_\ell\), where \(\Kronu{\xi_\ell}\) is the Dirac measure on \(\rk{\ab,\BsAF}\) with unit mass at \(\xi_\ell\).

 It was shown in~\zitaa{MR4051874}{\cprop{10.5}{203}} that \(\sigma\in\MggqF\) is molecular if and only if for some \(k\in\NO\) its \tnFmT{k} \(\sigma^\FTa{k}\) coincides with the \tqqa{zero} measure on \(\rk{\ab,\BsAF}\).
 This leads to a characterization of rational matrix-valued functions from \(\RFqab\) in terms of their \tnSfT{k}s:
 
\bpropl{P1620}
 Let \(F\in\RFqab\).
 Then the following statements are equivalent:
\baeqi{0}
 \il{P1620.i} There exist complex \tqqa{matrix} polynomials \(P\) and \(Q\) such that \(\det Q\) does not vanish identically and that \(F\) coincides with the restriction of \(PQ^\inv\) onto \(\Cab\).
 \il{P1620.ii} There exists an integer \(k\in\NO\) such that \(F^\FTa{k}\) coincides with the constant function on \(\Cab\) with value \(\Oqq\).
\eaeqi
\eprop
\bproof 
 First observe that the \tRabMa{} \(\sigma\defeq\rabmF\) given via \rdefn{F.D.RabMass} belongs to \(\MggqF\).
 According to \rlem{F.L.RabSp}, the restriction \(f\) of \(F\) onto \(\uhp\) belongs to \(\RFOq\).
 Furthermore, the \tSpMa{} \(\mu\defeq\sigmau{f}\) of \(f\) given via \rdefn{D1441} belongs to \(\MggqR\) and, in view of \rlem{F.L.RabSp}, fulfills \(\mu\rk{\cab}=\Oqq\) and \(\mu\rk{B}=\sigma\rk{B}\) for all \(B\in\BsAF\).
 In particular, \(\mu\) is molecular if and only if \(\sigma\) is molecular.
 From~\zitaa{MR4051874}{\cprop{10.5}{203}} we see that \(\sigma\) is molecular if and only if, for some \(k\in\NO\), the \tnFmT{k} \(\sigma^\FTa{k}\) of \(\sigma\) given via \rdefn{F.D.SN-M} coincides with the \tqqa{zero measure} on \(\rk{\ab,\BsAF}\).
 For an arbitrary \(k\in\NO\), by virtue of \rdefnss{F.D.SN-F}{F.D.ST} and \rprop{F.P.STbij}, we infer  that \(\sigma^\FTa{k}\) is the \tqqa{zero measure} on \(\rk{\ab,\BsAF}\) if and only if \(F^\FTa{k}\) is the constant function on \(\Cab\) with value \(\Oqq\).
 Consequently, we have shown that~\rstat{P1620.ii} is equivalent to the following statement:
\baeqi{2}
 \il{P1620.iii} \(\mu\) is molecular.
\eaeqi

\bimp{P1620.i}{P1620.iii}
 Suppose there exist complex \tqqa{matrix} polynomials \(P\) and \(Q\) such that \(\det Q\) does not vanish identically and \(F\) is the restriction of \(PQ^\inv\) onto \(\Cab\).
 Then \(f\) coincides with the restriction of \(PQ^\inv\) onto \(\uhp\).
 Hence, the application of~\zitaa{MR2570113}{\clem{B.4}{823}} yields~\ref{P1620.iii}.
\eimp

\bimp{P1620.iii}{P1620.i}
 Suppose that \(\mu\) is molecular.
 From \rprop{F.P.STbij} we see that \(F\) is exactly the \tSTF{} \(\STFu{\sigma}\) of \(\sigma\) given via \rdefn{F.D.ST}.
 If \(\mu\) is the \tqqa{zero measure} on \(\rk{\R,\BsAR}\), then \(\sigma\) coincides with the \tqqa{zero measure} on \(\rk{\ab,\BsAF}\), hence \(\STFu{\sigma}\) is, by \eqref{F.D.ST.1}, the constant function on \(\Cab\) with value \(\Oqq\), and, regarding \(F=\STFu{\sigma}\), thus~\rstat{P1620.i} obviously holds true.
 Now assume that \(\mu\) is not the \tqqa{zero measure} on \(\rk{\R,\BsAR}\), \tie{}\ \(\mu\rk{\R}\neq\Oqq\).
 \rprop{CR01.S1.2.16} shows that the moments \(\su{j}\defeq\int_\ab x^j\sigma(\dif x)\) exist for all \(j\in\NO\) and that we have, for each \(z\in\C\) with \(\abs{z}>\max\set{\abs{\ug},\abs{\obg}}\), furthermore \(z\in\Cab\) and
\(
  \STFua{\sigma}{z}
  =-\sum_{j=0}^\infi z^{-(j+1)}\su{j}
\).
 For all \(j\in\NO\), obviously \(\su{j}=\int_\R x^j\mu(\dif x)\).
 In view of~\rstat{P1620.iii} and \(\mu\rk{\R}\neq\Oqq\), then~\zitaa{MR2570113}{\crem{4.6}{792}} shows that \(\seqsinf\) is, using the terminology of~\zita{MR2570113}, a completely degenerate Hankel non-negative definite sequence of order \(n\) for some \(n\in\N\).
 Thus, from~\zitaa{MR2570113}{\cprop{9.2}{820} and \crem{3.5}{781}} we obtain the existence of a constant \(\rho\in[0,\infp)\) and specific complex \tqqa{matrix} polynomials \(a_n\) and \(b_n\) such that \(\det b_n\) does not vanish identically and \(\sum_{j=0}^\infi z^{-(j+1)}\su{j}=\rk{a_nb_n^\inv}\rk{z}\) holds true for all \(z\in\C\) with \(\abs{z}>\rho\).
Setting \(P\defeq-a_n\) and \(Q\defeq b_n\), then \(\det Q\) does not vanish identically and
\[
 F(z)
 =\STFua{\sigma}{z}
 =-\sum_{j=0}^\infi z^{-(j+1)}\su{j}
 =-\rk{a_nb_n^\inv}\rk{z}
 =\rk{PQ^\inv}\rk{z}
\]
 is valid for all \(z\in\C\) with \(\abs{z}>\max\set{\abs{\ug},\abs{\obg},\rho}\).
 Since \(F\) is holomorphic in \(\Cab\), then poles of the matrix-valued rational function \(PQ^\inv\) can only occur in \(\ab\) and \(F\) coincides with the restriction of \(PQ^\inv\) onto \(\Cab\).
 Consequently,~\ref{P1620.i} is valid.
\eimp
\eproof

\section{The class \(\PRFabq\)}\label{F.s2.PRF}
 In the scalar case \(q=1\), the set of all solutions of problem~\rabproblem{2n+1} can be parametrized with functions of the class \(\RFab{1}\) augmented by the constant function with value \(\infc\) defined on \(\Cab\) (\tcf{}~\zitaa{MR0458081}{\cthm{7.2}{147}}).
 The corresponding approach for the matricial situation \(q\geq1\) consists of extending the class \(\RFqab\) of holomorphic matrix-valued functions according to \rapp{A.S.cp} to some class of regular \tcp{q}{q}s of meromorphic matrix-valued functions.
 Such a class was already considered in~\zitaa{MR2222521}{\csec{5}}.
 As a first step, we extend the class \(\NFq\), using the terminology from \rapp{A.S.cp} and the end of \rapp{B.S.hol}, without explaining these notations here.
 We only recall that a \tcp{p}{q} \(\copa{P}{Q}\) is said to be regular if it satisfies \(\rank\tmatp{P}{Q}=q\).
 Furthermore, we observe that the set \(\pol{F}\) of poles of any meromorphic matrix-valued function \(F\) is discrete.

\bnotal{F.N.PRF}
 Denote by \(\PRFq\) the set of all ordered pairs \(\copa{P}{Q}\) consisting of \(\Cqq\)\nobreakdash-valued functions \(P\) and \(Q\) which are meromorphic in \(\uhp\), such that a discrete subset \(\mathcal{D}\) of \(\uhp\) exists, satisfying the following three conditions:
 \bAeqi{0}
  \il{F.N.PRF.I} \(\pol{P}\cup\pol{Q}\subseteq\mathcal{D}\).
  \il{F.N.PRF.II} \(\rank\smatp{P(z)}{Q(z)}=q\) for all \(z\in\uhp\setminus\mathcal{D}\).
  \il{F.N.PRF.III} \(\im\rk{\ek{Q(z)}^\ad \ek{P(z)}}\in\Cggq\) for all \(z\in\uhp\setminus\mathcal{D}\).
 \eAeqi
\enota

 Using a continuity argument, the following result is readily checked:

\breml{F.R.PRFim>}
 If \(\copa{P}{Q}\in\PRFq\), then \(\im\rk{\ek{Q(z)}^\ad \ek{P(z)}}\in\Cggq\) for all \(z\in\uhp\setminus\ek{\pol{P}\cup\pol{Q}}\).
\erem

 Now we supplement \rnota{F.N.PRF} in the following way:

\bnotal{F.N.PRFex}
 For each \(\copa{P}{Q}\in\PRFq\), denote by \(\PRFex{P}{Q}\) the set of all \(z\in\uhp\setminus\ek{\pol{P}\cup\pol{Q}}\) satisfying \(\rank\smatp{P(z)}{Q(z)}\neq q\).
\enota

 Regarding \rdefn{A.D.cp}, for each \(\copa{P}{Q}\in\PRFq\), we see that \(\PRFex{P}{Q}\) is exactly the set of all points \(z\in\uhp\) at which \(P\) and \(Q\) are both defined and the \tcp{q}{q} \(\copa{P(z)}{Q(z)}\) is not regular.
 In general, the linear subspace \(\gracp{P(z)}{Q(z)}\) depends on \(z\), whereas its dimension as well as the linear subspaces \(\domcp{P(z)}{Q(z)}\), \(\rancp{P(z)}{Q(z)}\), \(\nulcp{P(z)}{Q(z)}\), \(\mulcp{P(z)}{Q(z)}\), and the difference \(\rankcp{P(z)}{Q(z)}\) are essentially independent of \(z\):

\bpropl{ab.P1525}
 Let \(\copa{P}{Q}\in\PRFq\).
 Then \(\mathcal{P}\defeq\pol{P}\cup\pol{Q}\) is a discrete subset of \(\uhp\) and \(\mathcal{E}\defeq\PRFex{P}{Q}\) is a discrete subset of \(\dom\defeq\uhp\setminus\mathcal{P}\) admitting the representation \(\mathcal{E}=\setaca{z\in\dom}{\det\ek{Q(z)-\iu P(z)}=0}\).
 The set \(\exset\defeq\mathcal{P}\cup\mathcal{E}\) is the smallest discrete subset of \(\uhp\) satisfying the conditions~\ref{F.N.PRF.I}--\ref{F.N.PRF.III} in \rnota{F.N.PRF}.
 For all \(z\in\uhp\setminus\exset\), the \tcp{q}{q} \(\copa{P(z)}{Q(z)}\) is regular.
 For all \(z,w\in\uhp\setminus\exset\), furthermore
 \begin{align*}
  \Domcp{P(z)}{Q(z)}&=\Domcp{P(w)}{Q(w)},&\Rancp{P(z)}{Q(z)}&=\Rancp{P(w)}{Q(w)},\\
  \Nulcp{P(z)}{Q(z)}&=\Nulcp{P(w)}{Q(w)},&\Mulcp{P(z)}{Q(z)}&=\Mulcp{P(w)}{Q(w)},
 \end{align*}
 and \(\rankcp{P(z)}{Q(z)}=\rankcp{P(w)}{Q(w)}\) hold true.
\eprop
\bproof
 Observe that the matrix-valued functions \(P\) and \(Q\) are both meromorphic in \(\uhp\).
 Hence, the sets \(\pol{P}\) and \(\pol{Q}\) of poles as well as their union \(\mathcal{P}\) are discrete subsets of \(\uhp\).
 Consider an arbitrary discrete subset \(\mathcal{D}\) of \(\uhp\), satisfying the conditions~\ref{F.N.PRF.I}--\ref{F.N.PRF.III} in \rnota{F.N.PRF}.
 Such a subset exists by virtue of \rnota{F.N.PRF}.
 In view of \rnota{F.N.PRF}\ref{F.N.PRF.II}, and \rnota{F.N.PRFex}, then the set \(\mathcal{E}\) is a subset of \(\mathcal{D}\) and hence discrete.
 In particular, \(\mathcal{E}\) is a discrete subset of \(\dom\).
 Because of \rnota{F.N.PRFex} and \rrem{F.R.PRFim>}, the conditions~\ref{F.N.PRF.I}--\ref{F.N.PRF.III} in \rnota{F.N.PRF} are fulfilled where the set \(\mathcal{D}\) is substituted by \(\exset\) .
 Due to \rnota{F.N.PRF}\ref{F.N.PRF.I}, we have \(\mathcal{P}\subseteq\mathcal{D}\).
 Taking additionally into account \(\mathcal{E}\subseteq\mathcal{D}\), we see that the set \(\exset\) is a subset of \(\mathcal{D}\) and thus a discrete subset of \(\uhp\).
 Therefore, the set \(\exset\) is the smallest discrete subset \(\mathcal{D}\) of \(\uhp\) satisfying the conditions~\ref{F.N.PRF.I}--\ref{F.N.PRF.III} in \rnota{F.N.PRF}.
 Obviously, the matrix-valued functions \(F\defeq Q+\iu P\) and \(G\defeq Q-\iu P\) are both holomorphic in \(\dom\).
 From \rlem{A.L.detB} we infer that, for all \(z\in\dom\) with \(\det G(z)\neq0\), the \tcp{q}{q} \(\copa{P(z)}{Q(z)}\) is regular, implying \(z\notin\mathcal{E}\).
 In view of \rnota{F.N.PRFex} and \rdefn{A.D.cp}, we can conclude that, for all \(z\in\dom\setminus\mathcal{E}\), the \tcp{q}{q} \(\copa{P(z)}{Q(z)}\) is regular and fulfills \(\im\rk{\ek{Q(z)}^\ad \ek{P(z)}}\in\Cggq\), by virtue of \rrem{F.R.PRFim>}.
 Because of \rlem{ab.L1508}, we conversely have then \(\det G(z)\neq0\) for all \(z\in\dom\) with \(z\notin\mathcal{E}\).
 Consequently, \(\mathcal{E}=\setaca{z\in\dom}{\det G(z)=0}\).
 Hence, \(S\colon\dom\setminus\mathcal{E}\to\Cpq\) defined by \(S(z)\defeq\ek{F(z)}\ek{G(z)}^\inv\) is a holomorphic matrix-valued function.
 Furthermore, due to \rlem{ab.L1508}, we have \(\normS{S(z)}\leq1\) for all \(z\in\dom\setminus\mathcal{E}\).
 In particular, \(S\) and \(-S\) both belong to the class \(\SchF{q}{q}{\dom\setminus\mathcal{E}}\) of Schur functions (in \(\dom\setminus\mathcal{E}\)) introduced in \rnota{A.N.SF}.
 In view of \rlem{ab.P2028}, thus \(\ran{\Iq\pm S(z)}=\ran{\Iq\pm S(w)}\) and \(\nul{\Iq\pm S(z)}=\nul{\Iq\pm S(w)}\) hold true for all \(z,w\in\dom\setminus\mathcal{E}\).
 Regarding additionally \(\uhp\setminus\exset=\dom\setminus\mathcal{E}\), the application of \rlem{A.L.detB} completes the proof.
\eproof

 After transition to an appropriate equivalence relation, we can identify the class \(\NFq\) of matrix-valued functions with the set of equivalence classes of pairs \(\copa{P}{Q}\in\PRFq\) for which \(\det Q\) does not identically vanish in \(\uhp\).
 The analogous considerations are worked out in detail below for the following extension \(\PRFabq\) of the class \(\RFqab\) in the above mentioned sense.

\bnotanl{\tcf{}~\zitaa{MR2222521}{\cdefn{5.2}{144}}}{F.N.PRFab}
 Denote by \(\PRFabq\) the set of all ordered pairs \(\copa{P}{Q}\) consisting of \(\Cqq\)\nobreakdash-valued functions \(P\) and \(Q\) which are meromorphic in \(\Cab\), for which a discrete subset \(\mathcal{D}\) of \(\Cab\) exists, satisfying the following conditions:
\bAeqi{0}
 \il{F.N.PRFab.I} \(\pol{P}\cup\pol{Q}\subseteq\mathcal{D}\).
 \il{F.N.PRFab.II} \(\rank\smatp{P(z)}{Q(z)}=q\) for all \(z\in\C\setminus\rk{\ab\cup\mathcal{D}}\).
 \il{F.N.PRFab.III} \(\frac{1}{\im z}\im\rk{\rk{z-\ug}\ek{Q(z)}^\ad \ek{P(z)}}\in\Cggq\) for all \(z\in\C\setminus\rk{\R\cup\mathcal{D}}\).
 \il{F.N.PRFab.IV} \(\frac{1}{\im z}\im\rk{\rk{\obg-z}\ek{Q(z)}^\ad \ek{P(z)}}\in\Cggq\) for all \(z\in\C\setminus\rk{\R\cup\mathcal{D}}\).
\eAeqi
\enota

 Again using a continuity argument, the following result is readily checked:

\breml{F.R.PRFabim>}
 If \(\copa{P}{Q}\in\PRFabq\), then \(\frac{1}{\im z}\im\rk{\rk{z-\ug}\ek{Q(z)}^\ad \ek{P(z)}}\in\Cggq\) and \(\frac{1}{\im z}\im\rk{\rk{\obg-z}\ek{Q(z)}^\ad \ek{P(z)}}\in\Cggq\) for all \(z\in\C\setminus\ek{\R\cup\pol{P}\cup\pol{Q}}\).
\erem

 As done for \rnota{F.N.PRFex} above, we analogously supplement \rnota{F.N.PRFab} in the following way:
 
\bnotal{F.N.PRFabex}
 For each \(\copa{P}{Q}\in\PRFabq\) denote by \(\PRFabex{P}{Q}\) the set of all \(z\in\C\setminus\rk{\ab\cup\pol{P}\cup\pol{Q}}\) satisfying \(\rank\tmatp{P(z)}{Q(z)}\neq q\).
\enota

 Regarding \rdefn{A.D.cp}, we see that, for each \(\copa{P}{Q}\in\PRFabq\), the set \(\PRFabex{P}{Q}\) is exactly the set of all points \(z\in\Cab\) at which \(P\) and \(Q\) are both defined and for which the \tcp{q}{q} \(\copa{P(z)}{Q(z)}\) is not regular.
 The pairs belonging to \(\PRFabq\) fulfill conditions analogous to those in \rprop{F.P.RabH} for matrix-valued functions belonging to \(\RFqab\):
 
\bleml{ab.L1252}
 Let \(\copa{P}{Q}\in\PRFabq\) and let \(\mathcal{P}\defeq\pol{P}\cup\pol{Q}\).
 Then \(\frac{1}{\im z}\im\rk{\ek{Q(z)}^\ad\ek{P(z)}}\in\Cggq\) for all \(z\in\C\setminus\rk{\R\cup\mathcal{P}}\) and \(\re\rk{\ek{Q(w)}^\ad\ek{P(w)}}\in\Cggq\) for all \(w\in\ek{\Hla}\setminus\mathcal{P}\) and \(-\re\rk{\ek{Q(w)}^\ad\ek{P(w)}}\in\Cggq\) for all \(w\in\ek{\Hrb}\setminus\mathcal{P}\).
 Furthermore, \(\ek{Q(x)}^\ad\ek{P(x)}\in\Cggq\) for all \(x\in\crhl\setminus\mathcal{P}\) and \(-\ek{Q(x)}^\ad\ek{P(x)}\in\Cggq\) for all \(x\in\clhl\setminus\mathcal{P}\).
\elem
\bproof
 Consider an arbitrary \(z\in\C\setminus\rk{\R\cup\mathcal{P}}\).
 We have
\[
 \im\rk*{\rk{\obg-z}\ek*{Q(z)}^\ad \ek*{P(z)}}+\im\rk*{\rk{z-\ug}\ek*{Q(z)}^\ad \ek*{P(z)}}
 =\rk{\obg-\ug}\im\rk*{\ek*{Q(z)}^\ad\ek*{P(z)}}.
\]
 Regarding \(\ug<\obg\), we obtain, by virtue of \rremss{F.R.PRFabim>}{A.R.kK}, consequently \(\frac{1}{\im z}\im\rk{\ek{Q(z)}^\ad\ek{P(z)}}\in\Cggq\).
 Let \(A\defeq\ek{Q(z)}^\ad\ek{P(z)}\).
 Using \rremsss{A.R.XRIX}{A.R.kK}{F.R.PRFabim>}, we infer in the case \(\re z<\ug\) then
\[
 \re A
 =\frac{1}{\im z}\im\rk{zA}-\frac{\re z}{\im z}\im A
 \lgeq\frac{1}{\im z}\im\rk{zA}-\frac{\ug}{\im z}\im A
 =\frac{1}{\im z}\im\ek*{\rk{z-\ug}A}
 \lgeq\Oqq,
\]
 \tie{}, \(\re\rk{\ek{Q(z)}^\ad\ek{P(z)}}\in\Cggq\).
 In the case \(\re z>\obg\), we can conclude analogously
\begin{multline*}
 -\re A
 =\frac{\re z}{\im z}\im\rk{A}-\frac{1}{\im z}\im\rk{zA}\\
 \lgeq\frac{\obg}{\im z}\im\rk{A}-\frac{1}{\im z}\im\rk{zA}
 =\frac{1}{\im z}\im\ek*{\rk{\obg-z}A}
 \lgeq\Oqq,
\end{multline*}
 \tie{}, \(-\re\rk{\ek{Q(z)}^\ad\ek{P(z)}}\in\Cggq\).
 Observe that the matrix-valued functions \(P\) and \(Q\) are both meromorphic in \(\Cab\).
 Since the set \(\mathcal{P}\) is the union of the poles of \(P\) and \(Q\), it is a discrete subset of \(\Cab\).
 Furthermore, \(P\) and \(Q\) are both holomorphic in \(\C\setminus\rk{\ab\cup\mathcal{P}}\).
 Consequently, a continuity argument shows that we have \(\re\rk{\ek{Q(w)}^\ad\ek{P(w)}}\in\Cggq\) for all \(w\in\ek{\Hla}\setminus\mathcal{P}\) and \(-\re\rk{\ek{Q(w)}^\ad\ek{P(w)}}\in\Cggq\) for all \(w\in\ek{\Hrb}\setminus\mathcal{P}\).
 Regarding the continuity of the function \(S\colon\C\setminus\rk{\ab\cup\mathcal{P}}\to\Cqq\) defined by \(S(z)\defeq\ek{Q(z)}^\ad\ek{P(z)}\), we can conclude as in the proof of \rlem{B.L.STOJ} that \(\im\rk{\ek{Q(x)}^\ad\ek{P(x)}}=\Oqq\) holds true for all \(x\in\R\setminus\rk{\ab\cup\mathcal{P}}\).
 Therefore, we get \(\ek{Q(x)}^\ad\ek{P(x)}=\re\rk{\ek{Q(x)}^\ad\ek{P(x)}}\) for all \(x\in\R\setminus\rk{\ab\cup\mathcal{P}}\).
 Taking into account the already shown inequalities, we can infer then \(\ek{Q(x)}^\ad\ek{P(x)}\in\Cggq\) for all \(x\in\crhl\setminus\mathcal{P}\) and \(-\ek{Q(x)}^\ad\ek{P(x)}\in\Cggq\) for all \(x\in\clhl\setminus\mathcal{P}\).
\eproof

 By virtue of \rlem{ab.L1252}, we can think of \(\PRFabq\) as a subclass of \(\PRFq\) via restricting to the open upper half-plane \(\uhp\).
 Analogous as done for the class \(\PRFq\) in \rprop{ab.P1525} above, we are now going to prove that for each pair \(\copa{P}{Q}\in\PRFabq\) certain linear subspaces associated with the \tcp{q}{q} \(\copa{P(z)}{Q(z)}\) are essentially independent of \(z\).
 This is in accordance with \rprop{ab.P1648L1409}.
 In the proof we will use \rlem{ab.L1252} to reduce the situation to several open half-planes, in order to apply \rprop{ab.P1525}.

\bpropl{ab.P1528}
 Let \(\copa{P}{Q}\in\PRFabq\).
 Let \(\Pi_1\defeq\uhp\), \(\Pi_2\defeq\Hla,\), \(\Pi_3\defeq\lhp\), and \(\Pi_4\defeq\Hrb\).
 Then \(\mathcal{P}\defeq\pol{P}\cup\pol{Q}\) is a discrete subset of \(\Cab\) and \(\mathcal{E}\defeq\PRFabex{P}{Q}\) is a discrete subset of \(\dom\defeq\C\setminus\rk{\ab\cup\mathcal{P}}\) with \(\Pi_k\cap\mathcal{E}=\setaca{z\in\Pi_k\setminus\mathcal{P}}{\det\ek{Q(z)-\iu^k P(z)}=0}\) for each \(k\in\set{1,2,3,4}\).
 The set \(\exset\defeq\mathcal{P}\cup\mathcal{E}\) is the smallest discrete subset \(\mathcal{D}\) of \(\Cab\) satisfying the conditions~\ref{F.N.PRFab.I}--\ref{F.N.PRFab.IV} in \rnota{F.N.PRFab}.
 For all \(z\in\C\setminus\rk{\ab\cup\exset}\), the \tcp{q}{q} \(\copa{P(z)}{Q(z)}\) is regular.
 For every choice of \(z\) and \(w\) in \(\C\setminus\rk{\ab\cup\exset}\), furthermore
 \begin{align}
  \Domcp{P(z)}{Q(z)}&=\Domcp{P(w)}{Q(w)},&\Rancp{P(z)}{Q(z)}&=\Rancp{P(w)}{Q(w)},\label{ab.P1528.A}\\
  \Nulcp{P(z)}{Q(z)}&=\Nulcp{P(w)}{Q(w)},&\Mulcp{P(z)}{Q(z)}&=\Mulcp{P(w)}{Q(w)},\label{ab.P1528.B}
 \end{align}
 and \(\rankcp{P(z)}{Q(z)}=\rankcp{P(w)}{Q(w)}\) hold true.
\eprop
\bproof
 Since \(P\) and \(Q\) are matrix-valued functions meromorphic in \(\Cab\), the union \(\mathcal{P}\) of their poles is a discrete subset of \(\Cab\).
 Consider an arbitrary discrete subset \(\mathcal{D}\) of \(\Cab\) satisfying the conditions~\ref{F.N.PRFab.I}--\ref{F.N.PRFab.IV} in \rnota{F.N.PRFab}.
 Such a subset exists by virtue of \rnota{F.N.PRFab}.
 In view of \rnota{F.N.PRFab}\ref{F.N.PRFab.II} and \rnota{F.N.PRFabex}, then the set \(\mathcal{E}\) is a subset of \(\mathcal{D}\) and hence discrete.
 In particular, \(\mathcal{E}\) is a discrete subset of \(\dom\).
 Because of \rnota{F.N.PRFabex} and \rrem{F.R.PRFabim>}, the conditions~\ref{F.N.PRFab.I}--\ref{F.N.PRFab.IV} in \rnota{F.N.PRFab} are fulfilled with the set \(\exset\) instead of \(\mathcal{D}\).
 Due to \rnota{F.N.PRFab}\ref{F.N.PRFab.I}, we have \(\mathcal{P}\subseteq\mathcal{D}\).
 Taking additionally into account \(\mathcal{E}\subseteq\mathcal{D}\), we see that the set \(\exset\) is a subset of \(\mathcal{D}\) and thus a discrete subset of \(\Cab\).
 Therefore, the set \(\exset\) is the smallest discrete subset \(\mathcal{D}\) of \(\Cab\) satisfying the conditions~\ref{F.N.PRFab.I}--\ref{F.N.PRFab.IV} in \rnota{F.N.PRFab}.
 For all \(z\in\C\setminus\rk{\ab\cup\exset}\), we have furthermore \(z\in\mathcal{G}\) and \(z\notin\mathcal{E}\), implying \(\rank\tmatp{P(z)}{Q(z)}=q\) according to \rnota{F.N.PRFabex}, which, in view of \rdefn{A.D.cp}, shows that the \tcp{q}{q} \(\copa{P(z)}{Q(z)}\) is regular.
 Let \(\phi_1(\omega)\defeq\omega\), \(\phi_2(\omega)\defeq\iu\omega+\ug\), \(\phi_3(\omega)\defeq-\omega\), and \(\phi_4(\omega)\defeq-\iu\omega+\obg\).
 It is readily checked that,
 for each \(k\in\set{1,2,3,4}\), the mapping \(\phi_k\colon\uhp\to\Pi_k\) is bijective and that the union of their images \(\Pi_1,\Pi_2,\Pi_3,\Pi_4\) is exactly the whole domain \(\Cab\).
 
 Consider now an arbitrary \(k\in\set{1,2,3,4}\).
 Since the inverse \(\psi_k\defeq\phi_k^\inv\) of \(\phi_k\) is an affine bijection from \(\Pi_k\) onto \(\uhp\) and since the sets \(\mathcal{P}\) and \(\mathcal{E}\) are discrete, we can infer that \(\mathcal{D}_k\defeq\psi_k(\Pi_k\cap\exset)\) fulfills 
\(
 \mathcal{D}_k
 =\psi_k\rk*{\Pi_k\cap\rk{\mathcal{P}\cup\mathcal{E}}}
 =\psi_k(\Pi_k\cap\mathcal{P})\cup\psi_k(\Pi_k\cap\mathcal{E})
\)
 and is a discrete subset of \(\uhp\). 
 Regarding \rrem{A.R.XRIX} and \rlem{ab.L1252}, it is then readily checked that via
\begin{align}\label{ab.P1528.0}
 P_k(\omega)&\defeq\iu^{k-1}P\rk*{\phi_k(\omega)}&
&\text{and}&
 Q_k(\omega)&\defeq Q\rk*{\phi_k(\omega)}
\end{align}
 matrix-valued functions \(P_k\colon\psi_k\rk{\Pi_k\setminus\pol{P}}\to\Cqq\) and \(Q_k\colon\psi_k\rk{\Pi_k\setminus\pol{Q}}\to\Cqq\) are given, such that the pair \(\copa{P_k}{Q_k}\) consists of \(\Cqq\)\nobreakdash-valued functions, which are meromorphic in \(\uhp\), for which \(\mathcal{P}_k\defeq\pol{P_k}\cup\pol{Q_k}\) fulfills
 \(
  \mathcal{P}_k
  =\psi_k\rk{\Pi_k\cap\pol{P}}\cup\psi_k\rk{\Pi_k\cap\pol{Q}}
  =\psi_k\rk{\Pi_k\cap\mathcal{P}}
  \subseteq\mathcal{D}_k,
 \)
 and for which \(\rank\tmatp{P_k(\omega)}{Q_k(\omega)}=q\) and \(\im\rk{\ek{Q_k(\omega)}^\ad \ek{P_k(\omega)}}\in\Cggq\) hold true for all \(\omega\) in
\(
 \psi_k\rk{\Pi_k\setminus\exset}
 =\psi_k\rk{\Pi_k\setminus\rk{\Pi_k\cap\exset}}
 =\uhp\setminus\mathcal{D}_k
\).
 Consequently, \(\copa{P_k}{Q_k}\in\PRFq\).
 In view of \rprop{ab.P1525}, then \(\mathcal{P}_k\) we see that is a discrete subset of \(\uhp\) and \(\mathcal{E}_k\defeq\PRFex{P_k}{Q_k}\) is a discrete subset of \(\dom_k\defeq\uhp\setminus\mathcal{P}_k\), admitting the representation
 \[
  \mathcal{E}_k
  =\setaca*{\omega\in\dom_k}{\det\ek*{Q_k(\omega)-\iu P_k(\omega)}=0}
  =\setaca*{\omega\in\dom_k}{\det\ek*{Q\rk*{\phi_k(\omega)}-\iu^k P\rk*{\phi_k(\omega)}}=0}.
 \]
 Furthermore, \(\exset_k\defeq\mathcal{P}_k\cup\mathcal{E}_k\) is a discrete subset of \(\uhp\) and the linear subspaces \(\domcp{P_k(\omega)}{Q_k(\omega)}\), \(\rancp{P_k(\omega)}{Q_k(\omega)}\), \(\nulcp{P_k(\omega)}{Q_k(\omega)}\), \(\mulcp{P_k(\omega)}{Q_k(\omega)}\) and the difference of dimensions \(\rankcp{P_k(\omega)}{Q_k(\omega)}\) are independent of \(\omega\in\uhp\setminus\exset_k\).
 Since \(\phi_k\) is an affine bijection from \(\uhp\) onto \(\Pi_k\), we can conclude that \(\mathcal{Q}_k\defeq\Pi_k\cap\mathcal{P}\) fulfills
\(
 \mathcal{Q}_k
 =\phi_k(\mathcal{P}_k)
\) and is a discrete subset of \(\Pi_k\) and that \(\mathcal{F}_k\defeq\phi_k(\mathcal{E}_k)\) is a discrete subset of \(\mathcal{H}_k\defeq\Pi_k\setminus\mathcal{Q}_k\).
 We have
\(
 \mathcal{H}_k
 =\phi_k(\uhp\setminus\mathcal{P}_k)
 =\phi_k(\dom_k)
\) and
\beql{ab.P1528.1}
 \mathcal{F}_k
 =\setaca*{\zeta\in\mathcal{H}_k}{\det\ek*{Q\rk{\zeta}-\iu^k P\rk{\zeta}}=0}.
\eeq
 Moreover, \(\mathcal{B}_k\defeq\mathcal{Q}_k\cup\mathcal{F}_k\) fulfills
 \(
  \mathcal{B}_k
  =\phi_k\rk{\mathcal{P}_k}\cup\phi_k\rk{\mathcal{E}_k}
  =\phi_k\rk{\mathcal{P}_k\cup\mathcal{E}_k}
  =\phi_k\rk{\exset_k}
 \)
 and is a discrete subset of \(\Pi_k\).
 Thus, \(\Pi_k\setminus\mathcal{B}_k=\phi_k\rk{\uhp\setminus\exset_k}\).
 Furthermore, \(\domcp{P(\zeta)}{Q(\zeta)}\), \(\rancp{P(\zeta)}{Q(\zeta)}\), \(\nulcp{P(\zeta)}{Q(\zeta)}\), \(\mulcp{P(\zeta)}{Q(\zeta)}\), and \(\rankcp{P(\zeta)}{Q(\zeta)}\) are, in view of \eqref{ab.P1528.0}, independent of \(\zeta\in\Pi_k\setminus\mathcal{B}_k\), \tie{}, independent of \(\zeta\in\phi_k\rk{\uhp\setminus\exset_k}\).
 
 We are now going to verify \(\Pi_k\cap\mathcal{E}=\mathcal{F}_k\).
 First consider an arbitrary \(\zeta\in\Pi_k\cap\mathcal{E}\).
 Because of
\beql{ab.P1528.2}
 \mathcal{H}_k
 =\Pi_k\setminus\mathcal{Q}_k
 =\Pi_k\setminus\rk{\Pi_k\cap\mathcal{P}}
 =\Pi_k\setminus\mathcal{P}
\eeq
 and \rnota{F.N.PRFabex}, we have \(\zeta\in\mathcal{H}_k\) and \(\rank\tmatp{P(\zeta)}{Q(\zeta)}\neq q\).
 To the contrary, assume \(\zeta\notin\mathcal{F}_k\).
 In view of \eqref{ab.P1528.1}, then \(\det\ek{Q\rk{\zeta}-\iu^k P\rk{\zeta}}\neq0\).
 By virtue of \rlem{A.L.detB}, hence the \tcp{q}{q} \(\copa{P(\zeta)}{Q(\zeta)}\) is regular, \tie{}, \(\rank\tmatp{P(\zeta)}{Q(\zeta)}=q\), according to \rdefn{A.D.cp}.
 Since this is a contradiction, we necessarily have \(\zeta\in\mathcal{F}_k\).
 Conversely, consider an arbitrary \(\zeta\in\mathcal{F}_k\).
 Because of \eqref{ab.P1528.1} and \eqref{ab.P1528.2}, then
 \(
  \zeta
  \in\mathcal{H}_k
  =\Pi_k\setminus\mathcal{P}
  \subseteq\C\setminus\rk{\ab\cup\mathcal{P}}
 \).
 To the contrary, assume \(\zeta\notin\mathcal{E}\).
 In view of \rnota{F.N.PRFabex}, then \(\rank\tmatp{P(\zeta)}{Q(\zeta)}=q\).
 According to \rdefn{A.D.cp}, hence the \tcp{q}{q} \(\copa{P(\zeta)}{Q(\zeta)}\) is regular.
 Regarding \eqref{ab.P1528.0} and \rrem{A.R.PQre}, for \(\omega\defeq\psi_k\rk{\zeta}\), we have
\(
 \omega
 \in\psi_k\rk{\mathcal{H}_k}
 =\dom_k
 =\uhp\setminus\mathcal{P}_k
\) and
 \[
  \det\rk*{\ek*{P_k(\omega)}^\ad\ek*{P_k(\omega)}+\ek*{Q_k(\omega)}^\ad\ek*{Q_k(\omega)}}
  =\det\rk*{\ek*{P(\zeta)}^\ad\ek*{P(\zeta)}+\ek*{Q(\zeta)}^\ad\ek*{Q(\zeta)}}
  \neq0.
 \]
 Consequently, due to \rrem{A.R.PQre}, the \tcp{q}{q} \(\copa{P_k(\omega)}{Q_k(\omega)}\) is regular, \tie{}, \(\rank\tmatp{P_k(\omega)}{Q_k(\omega)}=q\), according to \rdefn{A.D.cp}.
 By virtue of \rnota{F.N.PRFex}, we have then \(\omega\notin\mathcal{E}_k\), implying \(\zeta\notin\mathcal{F}_k\).
 Since this is a contradiction, we see that \(\zeta\) necessarily belongs to \(\mathcal{E}\) and therefore to \(\Pi_k\cap\mathcal{E}\). 
 In view of \eqref{ab.P1528.1} and \eqref{ab.P1528.2}, we obtain the equations 
 \(
  \Pi_k\cap\mathcal{E}
  =\mathcal{F}_k
  =\setaca{\zeta\in\Pi_k\setminus\mathcal{P}}{\det\ek{Q\rk{\zeta}-\iu^k P\rk{\zeta}}=0}
 \).
 As already shown, for each \(k\in\set{1,2,3,4}\) the set \(\mathcal{B}_k\) is a discrete subset of \(\Pi_k\) and the entities \(\domcp{P(\zeta)}{Q(\zeta)}\), \(\rancp{P(\zeta)}{Q(\zeta)}\), \(\nulcp{P(\zeta)}{Q(\zeta)}\), \(\mulcp{P(\zeta)}{Q(\zeta)}\), and \(\rankcp{P(\zeta)}{Q(\zeta)}\) are independent of \(\zeta\in\Pi_k\setminus\mathcal{B}_k\).
 In particular, \(\mathcal{B}\defeq\mathcal{B}_1\cup\dotsb\cup\mathcal{B}_4\) is a discrete subset of \(\Pi_1\cup\dotsb\cup\Pi_4=\Cab\).
 Therefore, the sets \(\rk{\Pi_1\cap\Pi_2}\setminus\mathcal{B}\), \(\rk{\Pi_2\cap\Pi_3}\setminus\mathcal{B}\), and \(\rk{\Pi_3\cap\Pi_4}\setminus\mathcal{B}\) are non-empty.
 Consequently, we can infer that the linear subspaces \(\domcp{P(\zeta)}{Q(\zeta)}\), \(\rancp{P(\zeta)}{Q(\zeta)}\), \(\nulcp{P(\zeta)}{Q(\zeta)}\), \(\mulcp{P(\zeta)}{Q(\zeta)}\), and the difference of dimensions \(\rankcp{P(\zeta)}{Q(\zeta)}\) are independent of \(\zeta\) in \(\C\setminus\rk{\ab\cup\mathcal{B}}\).
 Because of
 \begin{align*}
  \mathcal{P}
  &=\bigcup_{k=1}^4\rk{\Pi_k\cap\mathcal{P}}
  =\bigcup_{k=1}^4\mathcal{Q}_k&
 &\text{and}&
  \mathcal{E}
  &=\bigcup_{k=1}^4\rk{\Pi_k\cap\mathcal{E}}
  =\bigcup_{k=1}^4\mathcal{F}_k,
 \end{align*}
 we have furthermore
\(
 \exset
 =\mathcal{P}\cup\mathcal{E}
 =\rk{\mathcal{Q}_1\cup\mathcal{F}_1}\cup\dotsb\cup\rk{\mathcal{Q}_4\cup\mathcal{F}_4}
 =\mathcal{B}_1\cup\dotsb\cup\mathcal{B}_4
 =\mathcal{B}
\). 
 Thus, \eqref{ab.P1528.A}, \eqref{ab.P1528.B}, and \(\rankcp{P(z)}{Q(z)}=\rankcp{P(w)}{Q(w)}\) follow for all \(z,w\in\C\setminus\rk{\ab\cup\exset}\).
\eproof

 As is easily seen, the class \(\PRFabq\) is closed under right multiplication by meromorphic matrix-valued functions \(R\) with not identically vanishing determinant:

\breml{F.R.PQR}
 Let \(\copa{P}{Q}\in\PRFabq\) and let \(R\) be a \(\Cqq\)\nobreakdash-valued function meromorphic in \(\Cab\) such that \(\det R\) does not vanish identically in \(\Cab\).
 Then \(\copa{PR}{QR}\in\PRFabq\).
\erem

 In view of the remarks on meromorphic matrix-valued functions given at the end of \rapp{B.S.hol}, it is readily checked that an equivalence relation on the set \(\PRFabq\) is given.
 Regarding \rrem{ab.R1457}, this equivalence relation is in accordance with that considered in \rapp{A.S.cp} for arbitrary \tcp{p}{q}s.
 
\bdefnl{ab.N0843}
 Two pairs \(\copa{P}{Q},\copa{S}{T}\in\PRFabq\) are said to be \emph{equivalent} if there exists a \(\Cqq\)\nobreakdash-valued function \(R\) meromorphic in \(\Cab\) such that \(\det R\) does not vanish identically in \(\Cab\) which fulfills \(S=PR\) and \(T=QR\).
 In this case, we write \(\copa{P}{Q}\rpaeq\copa{S}{T}\).
 Furthermore, denote by \(\rpcl{P}{Q}\) the  equivalence class of a pair \(\copa{P}{Q}\in\PRFabq\) and by \(\rsetcl{\mathcal{Q}}\defeq\setaca{\rpcl{S}{T}}{\copa{S}{T}\in\mathcal{Q}}\) the set of equivalence classes of pairs belonging to a subset \(\mathcal{Q}\) of \(\PRFabq\).
\edefn

 Using \rrem{ab.R1652}, the following remark can be easily concluded from \rprop{ab.L1642}:
 
\bremnl{\tcf{}~\zitaa{MR2222521}{\crem{5.4}{145}}}{ab.R1659}
 Let \(F\in\RFqab\) and let the functions \(P,Q\colon\Cab\to\Cqq\) be defined by \(P(z)\defeq F(z)\) and \(Q(z)\defeq\Iq\).
 Then the pair \(\copa{P}{Q}\) belongs to \(\PRFabq\) and \(\det Q(z)\neq0\) holds true for all \(z\in\Cab\).
\erem

 Conversely, we have:

\blemnl{\tcf{}~\zitaa{MR2222521}{\cprop{5.7}{145}}}{ab.P1801}
 Let \(\copa{P}{Q}\in\PRFabq\) be such that \(\det Q\) does not identically vanish in \(\Cab\).
 Then \(F\defeq PQ^\inv\) belongs to \(\RFqab\).
 Furthermore, the pair \(\copa{S}{T}\) consisting of the functions \(S,T\colon\Cab\to\Cqq\) defined by \(S(z)\defeq F(z)\) and \(T(z)\defeq\Iq\) belongs to \(\PRFabq\) and fulfills \(\copa{P}{Q}\rpaeq\copa{S}{T}\) and \(\det T(z)\neq0\) for all \(z\in\Cab\).
\elem
\bproof
 Due to~\zitaa{MR2222521}{\cprop{5.7}{145}}, we have \(F\in\RFqab\).
 In view of \rrem{ab.R1659}, we get then \(\copa{S}{T}\in\PRFabq\) and \(\det T(z)\neq0\) for all \(z\in\Cab\).
 Furthermore, \(R\defeq Q^\inv\) is a \(\Cqq\)\nobreakdash-valued function, which is meromorphic in \(\Cab\), satisfying \(S=PR\) and \(T=QR\).
 Since \(\det R=\rk{\det Q}^\inv\) does not identically vanish in \(\Cab\), thus \(\copa{P}{Q}\rpaeq\copa{S}{T}\) follows.
\eproof

 We end this section with an example of a simple family of pairs belonging to \(\PRFabq\).

\breml{ab.R1338}
 Let \(z\in\C\), let \(x\defeq z-\ug\), and let \(y\defeq\obg-z\).
 Involving \(\ba\defeq\obg-\ug\), it is readily checked that \(\iu\rk{y \ko{x }-x \ko{y }}=2\ba \im\rk{z }\) and \(\abs{y }^2x +\abs{x}^2y=\ba y x\).
\erem

 Given two complex matrices \(A\) and \(B\), we will use the notation
\beql{zdiag}
 \zdiag{A}{B}
 \defeq
 \bMat
  A&\NM\\
  \NM&B
 \eMat.
\eeq
 
\bexal{E1410}
 Let \(X,Y\in\Cqq\) satisfy \(\rank\tmatp{X}{Y}=q\) and \(Y^\ad X\in\Cggq\).
 Let \( P , Q \colon\Cab\to\Cqq\) be defined by \( P (z)\defeq X\) and \( Q (z)\defeq Y\) and let \(g,h\colon\Cab\to\C\) be given by \(g(z)\defeq z-\ug\) and \(h(z)\defeq\obg-z\), \tresp{}
 Denote by \(\IqCab\) and \(\OqqCab\) the constant \tqqa{matrix}-valued functions defined on \(\Cab\) with values \(\Iq\) and \(\Oqq\), \tresp{}
 Then:
\benui
 \il{E1410.a} The pairs \(\copa{ P }{h Q }\) and \(\copa{- P }{g Q }\) belong to \(\PRFabq\).
 \il{E1410.b} If \(X=\Oqq\), then \(\copa{ P }{h Q }\) and \(\copa{- P }{g Q }\) are equivalent to \(\copa{\OqqCab}{\IqCab}\).
 \il{E1410.c} If \(Y=\Oqq\), then \(\copa{ P }{h Q }\) and \(\copa{- P }{g Q }\) are equivalent to \(\copa{\IqCab}{\OqqCab}\).
\eenui
 We verify that statements~\eqref{E1410.a}--\eqref{E1410.c} are true:
 
 \eqref{E1410.a} The functions \( P \), \(h Q \), \(- P \), and \(g Q \) are holomorphic in \(\Cab\).
 Consider an arbitrary \(z\in\Cab\).
 Let \(x\defeq z-\ug\) and let \(y\defeq\obg-z\).
 Observe that the matrices \(\zdiag{\Iq}{\rk{y\Iq}}\) and \(\zdiag{\rk{-\Iq}}{\rk{x\Iq}}\) are invertible.
 Consequently, we get
\(
 \rank\tmatp{ P (z)}{h(z) Q (z)}
 =\rank\tmatp{X}{yY}
 =\rank\rk{\ek{\zdiag{\Iq}{\rk{y\Iq}}}\tmatp{X}{Y}}
 =\rank\tmatp{X}{Y}
 =q
\)
 and
\(
 \rank\tmatp{- P (z)}{g(z) Q (z)}
 =\rank\tmatp{-X}{xY}
 =\rank\rk{\ek{\zdiag{\rk{-\Iq}}{\rk{x\Iq}}}\tmatp{X}{Y}}
 =\rank\tmatp{X}{Y}
 =q
\).
 From the first equation in \rrem{ab.R1338} we can conclude \(\im\rk{x\ko y}=\ba\im z\) and \(\im\rk{y\ko x}=-\ba\im z\).
 Taking additionally into account \(Y^\ad X\in\Cggq\), we thus obtain
\begin{align*}
 \im\rk*{\rk{z-\ug}\ek*{h(z) Q (z)}^\ad\ek*{ P (z)}}
 &=\im\rk{x\ko y Y^\ad X}
 =\im\rk{x\ko y}Y^\ad X
 =\ba\im\rk{z}Y^\ad X,\\
 \im\rk*{\rk{\obg-z}\ek*{h(z) Q (z)}^\ad\ek*{ P (z)}}
 &=\im\rk*{\abs{y}^2Y^\ad X}
 =\Oqq,\\
 \im\rk*{\rk{z-\ug}\ek*{g(z) Q (z)}^\ad\ek*{- P (z)}}
 &=\im\rk*{-\abs{x}^2Y^\ad X}
 =\Oqq,
\intertext{and}
 \im\rk*{\rk{\obg-z}\ek*{g(z) Q (z)}^\ad\ek*{- P (z)}}
 &=\im\rk{-y\ko x Y^\ad X}
 =-\im\rk{y\ko x}Y^\ad X
 =\ba\im\rk{z}Y^\ad X.
\end{align*}
 Now assume in addition \(z\notin\R\).
 Because of \(\ba>0\) and \(Y^\ad X\in\Cggq\), we can infer
\begin{align*} 
 \frac{1}{\im z}\im\rk*{\rk{z-\ug}\ek*{h(z) Q (z)}^\ad\ek*{ P (z)}}
 &\in\Cggq,&
 \frac{1}{\im z}\im\rk*{\rk{\obg-z}\ek*{h(z) Q (z)}^\ad\ek*{ P (z)}}
 &\in\Cggq,
\end{align*}
 and 
\begin{align*} 
 \frac{1}{\im z}\im\rk*{\rk{z-\ug}\ek*{g(z) Q (z)}^\ad\ek*{- P (z)}}
 &\in\Cggq,&
 \frac{1}{\im z}\im\rk*{\rk{\obg-z}\ek*{g(z) Q (z)}^\ad\ek*{- P (z)}}
 &\in\Cggq,
\end{align*}
 by virtue of \rrem{A.R.kK}.
 Hence, \(\copa{ P }{h Q }\) and \(\copa{- P }{g Q }\) belong to \(\PRFabq\).

 \eqref{E1410.b} Assume \(X=\Oqq\).
 Then \(\copa{ P }{h Q }=\copa{\OqqCab}{h Q }\) and \(\copa{- P }{g Q }=\copa{\OqqCab}{g Q }\).
 Since \(\rank\tmatp{\Oqq}{Y}=\rank\tmatp{X}{Y}=q\), we have \(\det Y\neq0\).
 Hence, \(\det\rk{h Q }\) and \(\det\rk{g Q }\) both do not vanish identically in \(\Cab\).
 Because of
\begin{align*}
 \matp{ P \rk{h Q }^\inv}{\rk{h Q }\rk{h Q }^\inv}&=\matp{\OqqCab\rk{h Q }^\inv}{\rk{h Q }\rk{h Q }^\inv}=\matp{\OqqCab}{\IqCab}&
&\text{and}&
 \matp{\rk{- P }\rk{g Q }^\inv}{\rk{g Q }\rk{g Q }^\inv}&=\matp{\OqqCab\rk{g Q }^\inv}{\rk{g Q }\rk{g Q }^\inv}=\matp{\OqqCab}{\IqCab},
\end{align*}
 the pairs \(\copa{ P }{h Q }\) and \(\copa{- P }{g Q }\) are then both equivalent to \(\copa{\OqqCab}{\IqCab}\).
 
 \eqref{E1410.c} Assume \(Y=\Oqq\).
 Then \(\copa{ P }{h Q }=\copa{ P }{\OqqCab}\) and \(\copa{- P }{g Q }=\copa{- P }{\OqqCab}\).
 Since \(\rank\tmatp{X}{\Oqq}=\rank\tmatp{X}{Y}=q\), we have \(\det X\neq0\).
 Hence, \(\det P \) and \(\det\rk{- P }\) both do not vanish identically in \(\Cab\).
 Because of
\begin{align*}
 \matp{ P  P ^\inv}{\rk{h Q } P ^\inv}&=\matp{ P  P ^\inv}{\OqqCab P ^\inv}=\matp{\IqCab}{\OqqCab}&
&\text{and}&
 \matp{\rk{- P }\rk{- P }^\inv}{\rk{g Q }\rk{- P }^\inv}&=\matp{\rk{- P }\rk{- P }^\inv}{\OqqCab\rk{- P }^\inv}=\matp{\IqCab}{\OqqCab},
\end{align*}
 the pairs \(\copa{ P }{h Q }\) and \(\copa{- P }{g Q }\) are then both equivalent to \(\copa{\IqCab}{\OqqCab}\).
\eexa

\section{The class of parameters}\label{F.S.PM}

 The pairs belonging to the subclass of \(\PRFabq\) introduced below generate the equivalence classes, which will be used in \rsec{F.S.glt} as parameters in the description of the set of all solutions to \rprobrabs{m}.

\bnotal{ab.N1534}
 For each \(M\in\Cqp\), let \(\PRFabqa{M}\) be the set of all pairs \(\copa{F}{G}\in\PRFabq\) for which there exists a \(z_0\in\C\setminus\rk{\ab\cup\pol{F}\cup\pol{G}\cup\PRFabex{F}{G}}\) such that \(\rancp{F(z_0)}{G(z_0)}\subseteq\ran{M}\).
\enota

\breml{R1234}
 If \(M\in\Cqp\) fulfills \(\rank M=q\), then \(\PRFabqa{M}=\PRFabq\).
\erem 
 
 The class \(\PRFabqa{M}\) can be characterized by an additional equation involving the transformation matrix \(\OPu{\ran{M}}\) corresponding to the orthogonal projection onto the column space \(\ran{M}\):
 
\bleml{F.L.A<P}
 Let \(M\in\Cqp\) and let \(\copa{F}{G}\in\PRFabq\).
 Then \(\copa{F}{G}\in\PRFabqa{M}\) if and only if \(\OPu{\ran{M}}F=F\).
 In this case, \(\ran{F(z)}\subseteq\ran{M}\) for all \(z\in\C\setminus\rk{\ab\cup\pol{F}}\).
\elem
\bproof
 First observe that \(F\) is a matrix-valued function meromorphic in \(\Cab\) and that \(\dom\defeq\C\setminus\rk{\ab\cup\pol{F}}\) is exactly the set of points at which \(F\) is holomorphic.
 Therefore, \(\OPu{\ran{M}}F\) is a matrix-valued function meromorphic in \(\Cab\) and \(\dom\) is exactly the set of points at which \(\OPu{\ran{M}}F\) is holomorphic.
 In view of \rprop{ab.P1528}, the set \(\exset\defeq\pol{F}\cup\pol{G}\cup\PRFabex{F}{G}\) is a discrete subset of \(\Cab\) and \(\rancp{F(z)}{G(z)}=\rancp{F(w)}{G(w)}\) holds true for all \(z,w\in\C\setminus\rk{\ab\cup\exset}\).
 
 Assume \(\copa{F}{G}\in\PRFabqa{M}\).
 Then there exists some \(z_0\in\C\setminus\rk{\ab\cup\exset}\) with \(\rancp{F(z_0)}{G(z_0)}\subseteq\ran{M}\).
 As a subset of \(\exset\) the set \(\mathcal{D}\defeq\pol{G}\cup\PRFabex{F}{G}\) is discrete and the function \(F\) which is holomorphic in \(\dom\) fulfills, for all \(w\in\dom\setminus\mathcal{D}=\C\setminus\rk{\ab\cup\exset}\), furthermore
\(
 \Ran{F(w)}
 =\Rancp{F(z_0)}{G(z_0)}
 \subseteq\ran{M}
\),
 implying \(\OPu{\ran{M}}F(w)=F(w)\).
 Since \(\mathcal{D}\) is discrete, the set \(\dom\setminus\mathcal{D}\) has an accumulation point in \(\dom\).
 Therefore, the identity theorem for holomorphic functions yields \(\OPu{\ran{M}}F(z)=F(z)\) for all \(z\in\dom\), implying \(\OPu{\ran{M}}F=F\) and \(\ran{F(z)}\subseteq\ran{M}\) for all \(z\in\dom\).
 
 Conversely, assume \(\OPu{\ran{M}}F=F\).
 Since the set \(\exset\) is discrete, there exists some \(z_0\in\C\setminus\rk{\ab\cup\exset}\).
 In particular, \(z_0\in\dom\) and, consequently, \(\OPu{\ran{M}}F(z_0)=F(z_0)\).
 Thus \(\ran{F(z_0)}\subseteq\ran{M}\)%
 , implying \(\copa{F}{G}\in\PRFabqa{M}\).
\eproof

 Using \rlem{F.L.A<P}, it is readily checked that the equivalence relation on the set \(\PRFabq\) introduced in \rdefn{ab.N0843} is compatible with the here considered subclass \(\PRFabqa{M}\) in the following sense:

\breml{F.R.PMaeq}
 Let \(M\in\Cqp\) and let \(\copa{F}{G}\in\PRFabqa{M}\).
 Then \(\copa{\tilde F}{\tilde G}\in\PRFabqa{M}\) for all \(\copa{\tilde F}{\tilde G}\in\rpcl{F}{G}\).
\erem

 We can obtain a description of the set of equivalence classes \(\rsetcl{\PRFabqa{M}}\), depending on the rank \(r\) of the matrix \(M\) in terms of equivalence classes of pairs belonging to \(\PRFab{r}\):

\bleml{F.L.PRFabr}
 Let \(M\in\Cqp\) and let \(r\defeq\rank M\):
\benui
 \il{F.L.PRFabr.a} If \(r=0\), then \(\rsetcl{\PRFabqa{M}}=\set{\rpcl{F_0}{G_0}}\) where \(F_0,G_0\colon\Cab\to\Cqq\) are defined by \(F_0(z)\defeq\Oqq\) and \(G_0(z)\defeq\Iq\).
 \il{F.L.PRFabr.b} Assume \(r\geq1\).
 Let \(u_1,u_2,\dotsc,u_r\) be an arbitrary orthonormal basis of \(\ran{M}\), let \(U\defeq\mat{u_1,u_2,\dotsc,u_r}\), and let \(\Gamma_U\colon\rsetcl{\PRFab{r}}\to\rsetcl{\PRFabqa{M}}\) be defined by \(\Gamma_U(\rpcl{f}{g})\defeq\rpcl{UfU^\ad}{UgU^\ad+\OPu{\ek{\ran{M}}^\orth}}\).
 Then \(\Gamma_U\) is well defined and bijective.
\eenui
\elem
\bproof
 First assume \(r=0\), \tie{}, \(M=\Oqp\).
 Hence, \(\OPu{\ran{M}}=\Oqq\).
 Consider now an arbitrary pair \(\copa{F_1}{G_1}\in\PRFabqa{M}\).
 By virtue of \rlem{F.L.A<P}, then \(F_1=\OPu{\ran{M}}F_1=F_0\) follows.
 Observe that \(\copa{F_1}{G_1}\) belongs to \(\PRFabq\).
 In view of \rnota{F.N.PRFab}, thus \(\rank\tmatp{F_1(z_0)}{G_1(z_0)}=q\) for some \(z_0\in\C\setminus\rk{\ab\cup\pol{F_1}\cup\pol{G_1}}\).
 Because of \(F_1(z_0)=F_0(z_0)=\Oqq\), then necessarily \(\det G_1(z_0)\neq0\) holds true.
 In particular, \(\det G_1\) does not vanish identically in \(\Cab\).
 Consequently, the application of \rlem{ab.P1801} to the pair \(\copa{F_1}{G_1}\) and regarding \(F_1G_1^\inv=F_0G_1^\inv=F_0\), we obtain \(\copa{F_1}{G_1}\rpaeq\copa{F_0}{G_0}\).
 Therefore, \(\rsetcl{\PRFabqa{M}}\subseteq\set{\rpcl{F_0}{G_0}}\) is verified.
 Using \rrem{ab.R1659} and \rlem{F.L.A<P}, we can easily infer \(\copa{F_0}{G_0}\in\PRFabqa{M}\).
 Hence, \(\rsetcl{\PRFabqa{M}}=\set{\rpcl{F_0}{G_0}}\).
 
 Now assume \(r\geq1\).
 We have \(U^\ad U=\Iu{r}\) and \(\ran{U}=\ran{M}\).
 Let \(N\defeq\OPu{\ek{\ran{M}}^\orth}\).
 Using \rremss{R.P}{A.R.P=UU*}, we immediately obtain the representations
\beql{F.L.PRFabr.1}
 N
 =\OPu{\ek{\ran{U}}^\orth}
 =\Iq-\OPu{\ran{U}}
 =\Iq-UU^\ad.
\eeq
 In particular, we see
\begin{align}\label{F.L.PRFabr.2}
 NU
 &=U-UU^\ad U
 =\Ouu{q}{r}&
&\text{and}&
 U^\ad N
 &=U^\ad-U^\ad UU^\ad
 =\Ouu{r}{q}.
\end{align}
 
 We first consider an arbitrary pair \(\copa{f}{g}\in\PRFab{r}\) and show that \(\Gamma_U(\rpcl{f}{g})\) belongs to \(\rsetcl{\PRFabqa{M}}\):
 According to \rnota{F.N.PRFab}, the \(\Coo{r}{r}\)\nobreakdash-valued functions \(f\) and \(g\) are meromorphic in \(\Cab\) and there exists a discrete subset \(\mathcal{D}\) of \(\Cab\) with \(\pol{f}\cup\pol{g}\subseteq\mathcal{D}\) such that \(\rank\tmatp{f(z)}{g(z)}=r\) for all \(z\in\C\setminus\rk{\ab\cup\mathcal{D}}\) and furthermore \(\frac{1}{\im z}\im\rk{\rk{z-\ug}\ek{g(z)}^\ad \ek{f(z)}}\in\Cggq\) and \(\frac{1}{\im z}\im\rk{\rk{\obg-z}\ek{g(z)}^\ad \ek{f(z)}}\in\Cggq\) for all \(z\in\C\setminus\rk{\R\cup\mathcal{D}}\) hold true.
 Obviously,
\begin{align}\label{F.L.PRFabr.4}
 F&\defeq UfU^\ad&
&\text{and}&
 G&\defeq UgU^\ad+\OPu{\ek{\ran{M}}^\orth}
\end{align}
 are \(\Cqq\)\nobreakdash-valued functions meromorphic in \(\Cab\) with \(\pol{F}\subseteq\pol{f}\) and \(\pol{G}\subseteq\pol{g}\).
 Thus, \(\pol{F}\cup\pol{G}\subseteq\mathcal{D}\) follows.
 Consider an arbitrary \(z\in\C\setminus\rk{\ab\cup\mathcal{D}}\).
 In view of \rdefn{A.D.cp}, the \tcp{r}{r} \(\copa{f(z)}{g(z)}\) is regular.
 Furthermore, we have \(F(z)=U\ek{f(z)}U^\ad\) and \(G(z)\defeq U\ek{g(z)}U^\ad+\OPu{\ek{\ran{M}}^\orth}\).
 By virtue of \rlem{A.L.cpblo}, then \(\copa{F(z)}{G(z)}\) is a regular \tcp{q}{q} fulfilling \(\rancp{F(z)}{G(z)}\subseteq\ran{U}\) and \(\ek{G(z)}^\ad\ek{F(z)}=U\rk{\ek{g(z)}^\ad\ek{f(z)}}U^\ad\).
 According to \rdefn{A.D.cp}, in particular \(\rank\tmatp{F(z)}{G(z)}=q\).
 In the case \(z\notin\R\), using \rremss{A.R.XRIX}{A.R.XAX}, we can infer \(\frac{1}{\im z}\im\rk{\rk{z-\ug}\ek{G(z)}^\ad \ek{F(z)}}=U\ek{\frac{1}{\im z} \im\rk{\rk{z-\ug}\ek{g(z)}^\ad \ek{f(z)}}}U^\ad\in\Cggq\) and similarly \(\frac{1}{\im z}\im\rk{\rk{\obg-z}\ek{G(z)}^\ad \ek{F(z)}}\in\Cggq\).
 So \(\mathcal{D}\) is a discrete subset of \(\Cab\) such that the conditions~\ref{F.N.PRFab.I}--\ref{F.N.PRFab.IV} in \rnota{F.N.PRFab} are fulfilled for \(\copa{P}{Q}=\copa{F}{G}\).
 Consequently, \(\copa{F}{G}\) belongs to \(\PRFabq\).
 Because of \(\rank\tmatp{F(z)}{G(z)}=q\), we see from \rnota{F.N.PRFabex} that \(z\notin\PRFabex{F}{G}\).
 Summarizing, we have \(z\in\C\setminus\rk{\ab\cup\pol{F}\cup\pol{G}\cup\PRFabex{F}{G}}\) and \(\rancp{F(z)}{G(z)}\subseteq\ran{U}=\ran{M}\).
 Therefore, we obtain \(\copa{F}{G}\in\PRFabqa{M}\) and, regarding \eqref{F.L.PRFabr.4}, thus \(\Gamma_U(\rpcl{f}{g})\in\rsetcl{\PRFabqa{M}}\) follows.
 
 Next we are going to show that \(\Gamma_U(\rpcl{f}{g})\) is independent of the choice of the particular  representative \(\copa{f}{g}\) of the equivalence class \(\rpcl{f}{g}\in\rsetcl{\PRFab{r}}\):
 To that end, consider two arbitrary pairs \(\copa{f_1}{g_1}\) and \(\copa{f_2}{g_2}\) from \(\PRFab{r}\) satisfying \(\copa{f_1}{g_1}\rpaeq\copa{f_2}{g_2}\).
 For each \(j\in\set{1,2}\), let
\begin{align}\label{F.L.PRFabr.5}
 F_j&\defeq Uf_jU^\ad&
&\text{and}&
 G_j&\defeq Ug_jU^\ad+N.
\end{align}
 According to \rdefn{ab.N0843}, there is a \(\Coo{r}{r}\)\nobreakdash-valued function \(\rho\) meromorphic in \(\Cab\) such that \(\det\rho\) does not vanish identically in \(\Cab\), for which \(f_2=f_1\rho\) and \(g_2=g_1\rho\).
 Let \(R\defeq U\rho U^\ad+N\).
 Then \(R\) is a \(\Cqq\)\nobreakdash-valued function meromorphic in \(\Cab\).
 Regarding \(U^\ad U=\Iu{r}\), \(N^2=N\), and \eqref{F.L.PRFabr.2}, we get
\(
 F_1R
 =Uf_1U^\ad U\rho U^\ad+Uf_1U^\ad N
 =Uf_1\rho U^\ad
 =Uf_2 U^\ad
 =F_2
\)
 and
\[
 G_1R
 =Ug_1U^\ad U\rho U^\ad+Ug_1U^\ad N+NU\rho U^\ad+N^2
 =Ug_1\rho U^\ad+N^2
 =Ug_2U^\ad+N
 =G_2.
\]
 Furthermore, there exists some \(z_0\in\Cab\) such that \(\rho\) is holomorphic at \(z_0\) with \(\det\rho(z_0)\neq0\).
 In addition, we are going to check now that \(\det R(z_0)\neq0\).
 For this reason, consider an arbitrary vector \(v\in\nul{R(z_0)}\).
 Let \(w\defeq U^\ad v\).
 Then we have
\beql{F.L.PRFabr.3}
 U\rho(z_0)w+Nv
 =\ek*{U\rho(z_0)U^\ad+N}v
 =R(z_0)v
 =\Ouu{q}{1},
\eeq
 implying \(U^\ad U\rho(z_0)w+U^\ad Nv=\Ouu{r}{1}\).
 In view of \(U^\ad U=\Iu{r}\) and \eqref{F.L.PRFabr.2}, thus \(\rho(z_0)w=\Ouu{r}{1}\) follows.
 Because of \(\det \rho(z_0)\neq0\), then necessarily \(w=\Ouu{r}{1}\) holds true.
 Substituting this into \eqref{F.L.PRFabr.3}, we get \(Nv=\Ouu{q}{1}\).
 Regarding \eqref{F.L.PRFabr.1}, hence
\(
 v
 =UU^\ad v
 =Uw
 =\Ouu{q}{1}
\).
 Therefore, the linear subspace \(\nul{R(z_0)}\) is trivial, implying \(\det R(z_0)\neq0\).
 In particular, \(\det R\) does not vanish identically in \(\Cab\).
 According to \rdefn{ab.N0843}, then \(\copa{F_1}{G_1}\rpaeq\copa{F_2}{G_2}\).
 Consequently, \(\Gamma_U(\rpcl{f_1}{g_1})=\Gamma_U(\rpcl{f_2}{g_2})\).
 Thus, the mapping \(\Gamma_U\) is well defined. 
 We are now going to show that the mapping \(\Gamma_U\) is injective.
 To that end, consider two arbitrary pairs \(\copa{f_1}{g_1}\) and \(\copa{f_2}{g_2}\) from \(\PRFab{r}\) satisfying \(\Gamma_U(\rpcl{f_1}{g_1})=\Gamma_U(\rpcl{f_2}{g_2})\).
 Let \(\copa{F_1}{G_1}\) and \(\copa{F_2}{G_2}\) be given via \eqref{F.L.PRFabr.5}.
 Then \(\copa{F_1}{G_1}\rpaeq\copa{F_2}{G_2}\).
 Hence, according to \rdefn{ab.N0843}, there exists a \(\Cqq\)\nobreakdash-valued function \(R\) meromorphic in \(\Cab\) such that \(\det R\) does not vanish identically in \(\Cab\), fulfilling \(F_2=F_1R\) and \(G_2=G_1R\).
 Let \(\rho\defeq U^\ad RU\).
 Then \(\rho\) is a \(\Coo{r}{r}\)\nobreakdash-valued function meromorphic in \(\Cab\).
 Regarding \(U^\ad U=\Iu{r}\), \eqref{F.L.PRFabr.5}, and \eqref{F.L.PRFabr.2}, we have
\(
 f_1\rho
 =U^\ad Uf_1U^\ad RU
 =U^\ad F_1RU
 =U^\ad F_2U
 =U^\ad Uf_2U^\ad U
 =f_2
\)
 and
\[
 g_1\rho
 =U^\ad Ug_1U^\ad RU
 =U^\ad\rk{G_1-N}RU
 =U^\ad G_1RU
 =U^\ad G_2U
 =U^\ad Ug_2U^\ad U
 =g_2.
\]
 Furthermore, there exists some \(z_0\in\Cab\) such that \(R\) is holomorphic at \(z_0\) with \(\det R(z_0)\neq0\).
 In addition, we now prove that \(\det \rho(z_0)\neq0\).
 To do this, we consider an arbitrary vector \(w\in\nul{\rho(z_0)}\).
 Let \(v\defeq Uw\).
 Then
\(
 U^\ad R(z_0)v
 =U^\ad R(z_0)Uw
 =\rho(z_0)w
 =\Ouu{r}{1}
\).
 Because of \eqref{F.L.PRFabr.5} and \eqref{F.L.PRFabr.2}, we have \(NG_j=N^2=N\) for each \(j\in\set{1,2}\).
 In view of \(G_2=G_1R\) and \eqref{F.L.PRFabr.1}, we hence get \(N=NR=R-UU^\ad R\).
 From \eqref{F.L.PRFabr.2} we infer \(Nv=NUw=\Ouu{q}{1}\).
 Taking additionally into account \(U^\ad R(z_0)v=\Ouu{r}{1}\), we can then conclude \(R(z_0)v=\Ouu{q}{1}\).
 Since \(\det R(z_0)\neq0\) holds true, necessarily \(v=\Ouu{q}{1}\) follows.
 Regarding \(U^\ad U=\Iu{r}\), we thus obtain \(w=U^\ad v=\Ouu{r}{1}\).
 Therefore, the linear subspace \(\nul{\rho(z_0)}\) is trivial, implying \(\det \rho(z_0)\neq0\).
 In particular, \(\det\rho\) does not vanish identically in \(\Cab\).
 According to \rdefn{ab.N0843}, consequently \(\copa{f_1}{g_2}\rpaeq\copa{f_2}{g_2}\), \tie{}, \(\rpcl{f_1}{g_1}=\rpcl{f_2}{g_2}\).
 
 We finish the proof by showing that the mapping \(\Gamma_U\) is surjective.
 To that end, consider an arbitrary pair \(\copa{F_1}{G_1}\) from \(\PRFabqa{M}\).
 Then \(\copa{F_1}{G_1}\in\PRFabq\).
 According to \rnota{F.N.PRFab}, the \(\Cqq\)\nobreakdash-valued functions \(F_1\) and \(G_1\) are meromorphic in \(\Cab\).
 Consequently, \(B\defeq G_1-\iu F_1\) is a \(\Cqq\)\nobreakdash-valued function meromorphic in \(\Cab\).
 Due to \rprop{ab.P1528}, the set \(\exset\defeq\pol{F_1}\cup\pol{G_1}\cup\PRFabex{F_1}{G_1}\) is a discrete subset of \(\Cab\).
 Consider an arbitrary \(w\in\uhp\setminus\exset\).
 In view of \rprop{ab.P1528}, the \tcp{q}{q} \(\copa{F_1(w)}{G_1(w)}\) is regular.
 Using \rlem{ab.L1252} and \rrem{A.R.kK}, we get \(\im\rk{\ek{G(w)}^\ad\ek{F(w)}}\in\Cggq\).
 We see from \rlem{F.L.A<P} moreover \(\rancp{F_1(w)}{G_1(w)}\subseteq\ran{M}\).
 Therefore, \rprop{A.L.cpred} applies to the \tcp{q}{q} \(\copa{F_1(w)}{G_1(w)}\) and we get \(\det B(w)\neq0\).
 In particular, \(\det B\) does not vanish identically in \(\Cab\).
 Hence, \(R\defeq B^\inv\) is a \(\Cqq\)\nobreakdash-valued function meromorphic in \(\Cab\).
 Consequently,
\begin{align}\label{F.L.PRFabr.6}
 f&\defeq U^\ad F_1RU&
&\text{and}&
 g&\defeq U^\ad G_1RU
\end{align}
 are \(\Coo{r}{r}\)\nobreakdash-valued functions meromorphic in \(\Cab\).
 In addition, \(f\) and \(g\) are both holomorphic at \(w\) with \(f(w)=U^\ad\ek{F_1(w)}\ek{B(w)}^\inv U\) and \(g(w)=U^\ad\ek{G_1(w)}\ek{B(w)}^\inv U\).
 Thus,
\begin{align}\label{F.L.PRFabr.7}
 F_2&\defeq UfU^\ad&
&\text{and}&
 G_2&\defeq UgU^\ad+N
\end{align}
 are \(\Cqq\)\nobreakdash-valued functions, which are meromorphic in \(\Cab\) and holomorphic at \(w\) with \(F_2(w)=U\ek{f(w)}U^\ad\) and \(G_2(w)=U\ek{g(w)}U^\ad+\OPu{\ek{\ran{M}}^\orth}\).
 Due to \rprop{A.L.cpred}, the \tcp{q}{q} \(\copa{F_2(w)}{G_2(w)}\) is regular and satisfies \(F_2(w)=\ek{F_1(w)}\ek{B(w)}^\inv\) and \(G_2(w)=\ek{G_1(w)}\ek{B(w)}^\inv\).
 Let \(\mathcal{H}\rk{F_1}\), \(\mathcal{H}\rk{G_1}\), \(\mathcal{H}\rk{R}\), \(\mathcal{H}\rk{F_2}\), and \(\mathcal{H}\rk{G_2}\) be the sets of complex numbers at which \(F_1\), \(G_1\), \(R\), \(F_2\), and \(G_2\) are holomorphic, \tresp{}
 Taking into account the arbitrary choice of \(w\in\uhp\setminus\exset\), we can infer that \(\uhp\setminus\exset\) is a subset of each of the sets \(\mathcal{H}\rk{F_1}\), \(\mathcal{H}\rk{G_1}\), \(\mathcal{H}\rk{R}\), \(\mathcal{H}\rk{F_2}\), and \(\mathcal{H}\rk{G_2}\).
 Since \(\exset\) is discrete, the set \(\uhp\setminus\exset\) has in particular an accumulation point in \(\mathcal{H}\rk{F_1}\cap\mathcal{H}\rk{G_1}\cap\mathcal{H}\rk{R}\cap\mathcal{H}\rk{F_2}\cap\mathcal{H}\rk{G_2}\).
 Using the identity theorem for holomorphic functions, we thus can conclude
\begin{align}\label{F.L.PRFabr.8}
 F_2&=F_1R&
&\text{and}&
 G_2&=G_1R.
\end{align}
 Observe that \(\det R=\rk{\det B}^\inv\) does not vanish identically in \(\Cab\).
 Consequently, \rrem{F.R.PQR} yields \(\copa{F_2}{G_2}\in\PRFabq\) and \(\copa{F_1}{G_1}\rpaeq\copa{F_2}{G_2}\).
 We are now going to show that \(\copa{f}{g}\) belongs to \(\PRFab{r}\):
 Since \(\det B\) does not vanish identically in \(\Cab\), we obtain from the identity theorem for holomorphic functions that \(\mathcal{N}\defeq\setaca{\zeta\in\C\setminus\rk{\ab\cup\pol{\det B}}}{\det B(\zeta)=0}\) is a discrete subset of \(\Cab\).
 As already mentioned, \(\exset\) is a discrete subset of \(\Cab\).
 Therefore, \(\mathcal{D}\defeq\exset\cup\pol{R}\cup\pol{B}\cup\mathcal{N}\) is a discrete subset of \(\Cab\) as well.
 In view of \eqref{F.L.PRFabr.6}, we have \(\pol{f}\subseteq\pol{F_1}\cup\pol{R}\) and \(\pol{g}\subseteq\pol{G_1}\cup\pol{R}\).
 Hence, \(\pol{f}\cup\pol{g}\subseteq\mathcal{D}\) follows.
 Consider an arbitrary \(z\in\C\setminus\rk{\ab\cup\mathcal{D}}\).
 Let \(P\defeq U\ek{f(z)}U^\ad\) and \(Q\defeq U\ek{g(z)}U^\ad+\OPu{\ek{\ran{U}}^\orth}\).
 By virtue of \eqref{F.L.PRFabr.1}, \eqref{F.L.PRFabr.7}, and \eqref{F.L.PRFabr.8}, we immediately see that \(P=F_2(z)=\ek{F_1(z)}\ek{B(z)}^\inv\) and \(Q=G_2(z)=\ek{G_1(z)}\ek{B(z)}^\inv\).
 Observe that, due to \rprop{ab.P1528}, the \tcp{q}{q} \(\copa{F_1(z)}{G_1(z)}\) is regular.
 Because of \rrem{A.R.PQV}, thus the \tcp{q}{q} \(\copa{P}{Q}\) is regular and \(Q^\ad P=\ek{B(z)}^\invad\rk{\ek{G_1(z)}^\ad\ek{F_1(z)}}\ek{B(z)}^\inv\) holds true.
 Regarding \(U^\ad U=\Iu{r}\), we can now apply \rlem{A.L.cpblo} to the \tcp{r}{r} \(\copa{f(z)}{g(z)}\) to see that \(\copa{f(z)}{g(z)}\) is regular and that \(Q^\ad P=U\rk{\ek{g(z)}^\ad\ek{f(z)}}U^\ad\) is fulfilled.
 In particular, we obtain \(\rank\tmatp{f(z)}{g(z)}=r\), according to \rdefn{A.D.cp}, and furthermore
\begin{multline*}
 \ek*{g(z)}^\ad\ek*{f(z)}
 =U^\ad U\ek*{g(z)}^\ad\ek*{f(z)}U^\ad U
 =U^\ad Q^\ad P U\\
 =U^\ad \ek*{B(z)}^\invad\rk*{\ek*{G_1(z)}^\ad\ek*{F_1(z)}}\ek*{B(z)}^\inv U
 =\rk*{\ek*{B(z)}^\inv U}^\ad\rk*{\ek*{G_1(z)}^\ad\ek*{F_1(z)}}\rk*{\ek*{B(z)}^\inv U}.
\end{multline*}
 In the case \(z\notin\R\), due to \rprop{ab.P1528}, we have \(\frac{1}{\im z}\im\rk{\rk{z-\ug}\ek{G_1(z)}^\ad \ek{F_1(z)}}\in\Cggq\) and \(\frac{1}{\im z}\im\rk{\rk{\obg-z}\ek{G_1(z)}^\ad \ek{F_1(z)}}\in\Cggq\), implying, by virtue of \rremss{A.R.XRIX}{A.R.XAX}, then \(\frac{1}{\im z}\im\rk{\rk{z-\ug}\ek{g(z)}^\ad \ek{f(z)}}=\rk{\ek{B(z)}^\inv U}^\ad\ek{\frac{1}{\im z} \im\rk{\rk{z-\ug}\ek{G_1(z)}^\ad \ek{F_1(z)}}}\rk{\ek{B(z)}^\inv U}\in\Cggo{r}\) and, similarly, \(\frac{1}{\im z}\im\rk{\rk{\obg-z}\ek{g(z)}^\ad \ek{f(z)}}\in\Cggo{r}\).
 According to \rnota{F.N.PRFab}, hence \(\copa{f}{g}\) belongs to \(\PRFab{r}\).
 Applying \(\Gamma_U\) to the equivalence class of \(\copa{f}{g}\), we get with \eqref{F.L.PRFabr.7} and \(\copa{F_1}{G_1}\rpaeq\copa{F_2}{G_2}\) then \(\Gamma_U(\rpcl{f}{g})=\rpcl{F_2}{G_2}=\rpcl{F_1}{G_1}\).
\eproof

\bexal{E1147} 
 Let \(M\in\Cqp\) and let \(F,G\colon\Cab\to\Cqq\) be defined by \(F(z)\defeq\Oqq\) and \(G(z)\defeq\Iq\).
 Then \(\copa{F}{G}\in\PRFabqa{M}\).
\eexa

\bexal{FR}%
 Let \(M\in\CHq\) and let \(f,g\colon\Cab\to\C\) be both holomorphic and not identically vanishing.
 Let \(\cU\) be a linear subspace of \(\ran{M}\) and let \(P\defeq\OPu{\cU}\).
 Let \(F,G\colon\Cab\to\Cqq\) be defined by \(F(z)\defeq f\rk{z}MPM\) and \(G(z)\defeq g\rk{z}\rk{\Iq-M^\mpi PM}\).
 Then \(\copa{F}{G}\) belongs to \(\PRFabqa{M}\). 
 Indeed, since \(f\) and \(g\) are both holomorphic and not identically vanishing, the matrix-valued functions \( F \) and \(G\) are holomorphic and in particular meromorphic in \(\Cab\) and \(\mathcal{D}\defeq\setaca{z\in\Cab}{f\rk{z}=0\text{ or }g\rk{z}=0}\) is a discrete subset of \(\Cab\).
 Observe that \(M^\ad=M\) implies \(\rk{M^\mpi}^\ad=M^\mpi\) and \(M^\mpi M=MM^\mpi\), by virtue of \rremss{A.R.A++*}{ab.R1052}.
 Consider an arbitrary \(z\in\C\setminus\rk{\ab\cup\mathcal{D}}\).
 Let \(v\in\nul{\tmatp{F(z)}{G(z)}}\), \tie{}, \(f\rk{z}MPMv=\Ouu{q}{1}\) and \(g\rk{z}\rk{\Iq-M^\mpi PM}v=\Ouu{q}{1}\).
 Since \(f\rk{z}\neq0\) and \(g\rk{z}\neq0\), hence \(MPMv=\Ouu{q}{1}\) and \(v=M^\mpi PMv\).
 Thus, taking additionally into account \eqref{mpi}, we can conclude \(v=M^\mpi PMv=M^\mpi MM^\mpi PMv=M^\mpi M^\mpi M PMv=\Ouu{q}{1}\).
 This shows \(\rank\tmatp{F(z)}{G(z)}=q\).
 Observe that \(\cU\subseteq\ran{M}\) implies \(MM^\mpi P=P\), by virtue of \rremss{R.P}{R.AA+B.A}.
 Taking additionally into account \(P^\ad=P=P^2\), we can furthermore conclude
\[\begin{split}
 \rk{\Iq-M^\mpi PM}^\ad\rk{MPM}
 &=\rk{\Iq-MPM^\mpi}MPM
 =MPM-MPM^\mpi MPM\\
 &=MPM-MPMM^\mpi PM
 =MPM-MP^2M
 =\Oqq.
\end{split}\]
 Consequently, \(\ek{G(z)}^\ad\ek{F(z)}=\ko{g\rk{z}}f\rk{z}\rk{\Iq-M^\mpi PM}^\ad\rk{MPM}=\Oqq\).
 If \(z\notin\R\), thus the matrices \(\frac{1}{\im z}\im\rk{\rk{z-\ug}\ek{G(z)}^\ad\ek{F(z)}}\) and \(\frac{1}{\im z}\im\rk{\rk{\obg-z}\ek{G(z)}^\ad\ek{F(z)}}\) are both \tnnH{}.
 Hence, \(\copa{F}{G}\) belongs to \(\PRFabq\).
 Since obviously \(\OPu{\ran{M}}F=F \), then \(\copa{F}{G}\in\PRFabqa{M}\) follows by virtue of \rlem{F.L.A<P}.
\eexa

\section{The \hFATion{} and its inverse}\label{S1229}
 Our next considerations are aimed at preparing the foundations for the desired function-theoretic \tSchur{}--\tNevanlinna{} type algorithm.
 This algorithm consists of two different instances, because the first step differs from the remaining ones.
 In this section, we treat the algebraic formalism for the  first step.
 Doing this, we take into account as well the forward as the backward form of the algorithm.
 
 We are now going to introduce a transformation of matrix-valued functions, which is intimately connected with the \tFTion{} for sequences of complex matrices (see \rdefn{ab.N0940}).
 
 In this section, for an arbitrarily given complex matrix \(E\), we write
\begin{align*}
 \OPr{E}&\defeq\OPu{\ran{E}}&
&\text{and}&
 \OPn{E}&\defeq\OPu{\nul{E}}
\end{align*}
 for the transformation matrix corresponding to the orthogonal projection onto \(\ran{E}\) and \(\nul{E}\), \tresp{}
 In view of \rremsss{R.P}{ab.R1842*}{ab.R1052}, we have
\begin{align}
 \ran{\OPr{E}}&=\ran{E},&\nul{\OPr{E}}&=\nul{E^\ad},&\ran{\OPn{E}}&=\nul{E},&\nul{\OPn{E}}&=\ran{E^\ad},\notag\\
 \OPr{E}^2&=\OPr{E},&\OPr{E}^\ad&=\OPr{E},&\OPn{E}^2&=\OPn{E},&\OPn{E}^\ad&=\OPn{E},\label{F.G.PQ2*}
\shortintertext{and}
 \OPr{E}&=EE^\mpi,&\OPr{E^\ad}&=E^\mpi E,&\OPn{E}&=\Iq-E^\mpi E,&\OPn{E^\ad}&=\Ip-E E^\mpi.\label{F.G.PQ}
\end{align}
 In the sequel, we will also use these identities without explicitly mentioning.
 Furthermore, we consider here a complex matrix \(M\), which in the context of the matricial Hausdorff moment problem will later be the \tnnH{} matrix \(\su{0}\) taken from a sequence \(\seqska\) belonging to \(\Fggqka\).

\bdefnl{F.D.FTF}
 Let \(\dom\) be a non-empty subset of \(\C\), let \(F\colon\dom\to\Cpq\) be a matrix-valued function, and let \(M\) be a complex \tpqa{matrix}.
 Then the pair \(\copa{G_1}{G_2}\) built with the functions \(G_1\) and \(G_2\) defined on \(\dom\) by
\begin{align*}
 G_1(z)&\defeq\rk{\obg-z}F(z)-M&
&\text{and}&
 G_2(z)&\defeq\rk{\obg-z}\ek*{\rk{z-\ug}M^\mpi F(z)+\OPr{M^\ad}}+\ba\OPn{M}
\end{align*}
 is called the \emph{\tFatpv{M}{F}}.
\edefn

 In connection with the \tFATion{}, we consider the following quadratic \taaa{(p+q)}{(p+q)}{matrix} polynomial:
 
\bnotal{ab.N1246a}
 Let \(M\in\Cpq\).
 Then let \(\mFTu{M}\colon\C\to\Coo{(p+q)}{(p+q)}\) be defined by
\[
 \mFTua{M}{z}
 \defeq
 \begin{pmat}[{|}]
  -\rk{\obg-z}\OPr{M}&M\cr\-
  -\rk{\obg-z}\rk{z-\ug}M^\mpi&-\rk{\obg-z}\OPr{M^\ad}-\ba\OPn{M}\cr
 \end{pmat}.
\]
\enota

 In what follows, we will use the notation given via \eqref{zdiag} to calculate, in view of \rrem{ab.R1543}, certain forms involving the signature matrix \(\Jimq\) given by \eqref{J}.
 For an arbitrarily given \(z\in\C\), we will furthermore write abbreviatory \(x\defeq z-\ug\) and \(y\defeq\obg-z\).
 Obviously, we have
\begin{align}\label{F.G.y+x}
 y+x&=\obg-\ug=\ba&
&\text{and}&
 \ug y+\obg x
 =\ug\obg-\ug z+\obg z-\obg\ug
 =\rk{\obg-\ug}z
 =\ba z
\end{align}%
 as well as
\beql{F.G.ax+by}
 \ug x+\obg y
 =\ug z-\ug^2+\obg^2-\obg z
 =\rk{\obg+\ug}\rk{\obg-\ug}-\rk{\obg-\ug}z
 =\ba\rk{\obg+\ug-z}.
\eeq
 
\bleml{F.L.WJ}
 Let \(M\in\CHq\).
 Let \(z\in\C\), let \(x\defeq z-\ug\), and let \(y\defeq\obg-z\).
\benui
 \il{F.L.WJ.a} Let \(W_0\defeq\mFTua{M}{z}\).
 Then
\beql{F.L.WJ.A}
 W_0^\ad\Jimq W_0
 =
 \begin{pmat}[{|}]
  -2\abs{y}^2\im\rk{z}M^\mpi&\iu\rk{\ko{yx}+\abs{y}^2}\OPr{M}\cr\-
  -\iu\rk{yx+\abs{y}^2}\OPr{M}&-2\im\rk{z}M\cr
 \end{pmat}.
\eeq
 \il{F.L.WJ.b} Let \(W_1\defeq\ek{\zdiag{\rk{x\Iq}}{\Iq}}W_0\) and let \(W_2\defeq\ek{\zdiag{\rk{y\Iq}}{\Iq}}W_0\).
 Then
\beql{F.L.WJ.Ax}
 W_1^\ad\Jimq W_1
 =\ba
 \bMat
  \Oqq&\iu\ko{yx}\OPr{M}\\
  -\iu yx\OPr{M}&-2\im\rk{z}M
 \eMat
 =
 \begin{pmat}[{|}]
  \Oqq&\iu\rk{\abs{y}^2\ko x+\abs{x}^2\ko y}\OPr{M}\cr\-
  -\iu\rk{\abs{y}^2 x+\abs{x}^2 y}\OPr{M}&-2\ba\im\rk{z}M\cr
 \end{pmat}
\eeq
 and
\beql{F.L.WJ.Ay}
 W_2^\ad\Jimq W_2
 =\ba\abs{y}^2
 \bMat
  -2\im\rk{z}M^\mpi&\iu \OPr{M}\\
  -\iu \OPr{M}&\Oqq
 \eMat.
\eeq
 \il{F.L.WJ.c} Let \(W_3\defeq\ek{\zdiag{\rk{yx\Iq}}{\Iq}}W_0\).
 Then
\beql{F.L.WJ.Ayx}
 W_3^\ad\Jimq W_3
 =\abs{y}^2
 \begin{pmat}[{|}]
  -2\abs{x}^2\im\rk{z}M^\mpi&\iu\rk{\ko{yx}+\abs{x}^2}\OPr{M}\cr\-
  -\iu\rk{yx+\abs{x}^2}\OPr{M}&-2\im\rk{z}M\cr
 \end{pmat}.
\eeq
\eenui
\elem
 In view of \eqref{F.G.PQ2*}, \eqref{F.G.PQ}, and \rremss{A.R.A++*}{ab.R1338}, the proof of \rlem{F.L.WJ} is straightforward.
 We omit the details.

 In addition, we now rewrite the right-hand sides of the equations \eqref{F.L.WJ.A}--\eqref{F.L.WJ.Ayx}, using the signature matrix \(\Jimq\):

\bpropl{ab.L1336}
 Let \(M\in\CHq\).
 Let \(z\in\C\), let \(x\defeq z-\ug\), and let \(y\defeq\obg-z\).
\benui
 \il{ab.L1336.a} Let \(W_0\defeq\mFTua{M}{z}\).
 Then
\begin{multline*}
 W_0^\ad\Jimq W_0
 =\ek*{\zdiag{\rk{yx\OPr{M}}}{\Iq}}^\ad\Jimq\ek*{\zdiag{\rk{yx\OPr{M}}}{\Iq}}\\
 +\abs{y}^2\ek*{\rk{\zdiag{\OPr{M}}{\Iq}}^\ad\Jimq\rk{\zdiag{\OPr{M}}{\Iq}}-2\im\rk{z}\rk{\zdiag{M^\mpi}{\Oqq}}}-2\im\rk{z}\rk{\zdiag{\Oqq}{M}}.
\end{multline*}
 \il{ab.L1336.b} Let \(W_1\defeq\ek{\zdiag{\rk{x\Iq}}{\Iq}}W_0\) and let \(W_2\defeq\ek{\zdiag{\rk{y\Iq}}{\Iq}}W_0\).
 Then
\[\begin{split}
 W_1^\ad\Jimq W_1
 &=\ba\rk*{\ek*{\zdiag{\rk{yx\OPr{M}}}{\Iq}}^\ad\Jimq\ek*{\zdiag{\rk{yx\OPr{M}}}{\Iq}}-2\im\rk{z}\rk{\zdiag{\Oqq}{M}}}\\
 &=\abs{y}^2\ek*{\zdiag{\rk{x\OPr{M}}}{\Iq}}^\ad\Jimq\ek*{\zdiag{\rk{x\OPr{M}}}{\Iq}}\\
 &\qquad+\abs{x}^2\ek*{\zdiag{\rk{y\OPr{M}}}{\Iq}}^\ad\Jimq\ek*{\zdiag{\rk{y\OPr{M}}}{\Iq}}-2\ba\im\rk{z}\rk{\zdiag{\Oqq}{M}}
\end{split}\]
 and
\[
 W_2^\ad\Jimq W_2
 =\ba\abs{y}^2\ek*{\rk{\zdiag{\OPr{M}}{\Iq}}^\ad\Jimq\rk{\zdiag{\OPr{M}}{\Iq}}-2\im\rk{z}\rk{\zdiag{M^\mpi}{\Oqq}}}.
\]
 \il{ab.L1336.c} Let \(W_3\defeq\ek{\zdiag{\rk{yx\Iq}}{\Iq}}W_0\).
 Then
\begin{multline*}
 W_3^\ad\Jimq W_3
 =\abs{y}^2\biggl\{\ek*{\zdiag{\rk{yx\OPr{M}}}{\Iq}}^\ad\Jimq\ek*{\zdiag{\rk{yx\OPr{M}}}{\Iq}}\\
 +\abs{x}^2\ek*{\rk{\zdiag{\OPr{M}}{\Iq}}^\ad\Jimq\rk{\zdiag{\OPr{M}}{\Iq}}-2\im\rk{z}\rk{\zdiag{M^\mpi}{\Oqq}}}-2\im\rk{z}\rk{\zdiag{\Oqq}{M}}\biggr\}.
\end{multline*}
\eenui
\eprop
\bproof
 Taking into account \(\OPr{M}^\ad=\OPr{M}\) and \eqref{J}, the asserted identities immediately follow from \rlem{F.L.WJ}.
\eproof

 It will be clear from \rlemss{F.L.FT-11}{F.L.FT1-1} below that the following transformation for pairs of meromorphic matrix-valued functions is, under certain conditions, essentially the inversion of the \tFATion{}.
 To define this inverse transformation, we use the terminology given at the end of \rapp{B.S.hol}.

\bdefnl{F.D.iFTF}
 Let \(\dom\) be a domain.
 Let \(G_1\) be a \(\Cpq\)\nobreakdash-valued function meromorphic in \(\dom\) and let \(G_2\) be a \(\Cqq\)\nobreakdash-valued function meromorphic in \(\dom\).
 Let \(M\in\Cpq\), let the functions \(g,h\colon\dom\to\C\) be defined by
\begin{align}\label{F.D.iFTF.1}
 g(z)&\defeq z-\ug&
&\text{and}&
 h(z)&\defeq\obg-z,
\end{align}
 \tresp{}, and let
\begin{align}\label{F.D.iFTF.2}
 F_1&\defeq h\OPr{M}G_1+MG_2&
&\text{and}&
 F_2&\defeq-hgM^\mpi G_1+hG_2.
\end{align}
 Suppose that \(\det F_2\) does not identically vanish in \(\dom\).
 Then we call the \(\Cpq\)\nobreakdash-valued function \(F\defeq F_1F_2^\inv\)(, which is meromorphic in \(\dom\),) the \emph{\tiFaTv{M}{\(\copa{G_1}{G_2}\)}}.
\edefn

 To the \tiFaT{M} we can associate the following matrix polynomial:

\bnotal{ab.N1246b}
 Let \(M\in\Cpq\).
 Then let \(\mFTiu{M}\colon\C\to\Coo{(p+q)}{(p+q)}\) be defined by
\[
 \mFTiua{M}{z}
 \defeq
 \begin{pmat}[{|}]
  \rk{\obg-z}\OPr{M}&M\cr\-
  -\rk{\obg-z}\rk{z-\ug}M^\mpi&\rk{\obg-z}\Iq\cr
 \end{pmat}.
\]
\enota

\breml{R1204} 
 Let \(M\in\Cqq\) and let \(z\in\C\).
 Then:
\benui
 \il{R1204.a} If \(M=\Oqq\), then
\(
 \mFTiua{M}{z}
 =\smat{\Oqq&\Oqq\\
  \Oqq&\rk{\obg-z}\Iq}
\).
 \il{R1204.b} If \(M\) is invertible, then
\(
 \mFTiua{M}{z}
 =\smat{\rk{\obg-z}\Iq&M\\
  -\rk{\obg-z}\rk{z-\ug}M^\inv&\rk{\obg-z}\Iq}
\).
\eenui
\erem

 Regarding \(y+x=\ba\) and \(\OPn{M}=\Iq-\OPr{M^\ad}\), it is readily checked that the matrix polynomials \(\mFTiu{M}\) and \(\mFTu{M}\) are connected in the following way by the signature matrix \(\Jabspq\) given in \eqref{J}:
 
\breml{ab.R1025}
 If \(z\in\C\), then \(\ek{\mFTiua{M}{z}}\Jabspq=-\Jabspq\ek{\mFTua{M}{z}+\rk{z-\ug}\rk{\zdiag{\Opp}{\OPn{M}}}}\).
\erem

 Consequently, a result analogous to \rprop{ab.L1336} follows:
 
\bpropl{ab.L1059}
 Let \(M\in\CHq\).
 Let \(z\in\C\), let \(x\defeq z-\ug\), and let \(y\defeq\obg-z\).
\benui
 \il{ab.L1059.a} Let \(V_0\defeq\mFTiua{M}{z}\).
 Then
\begin{multline*}
 V_0^\ad\Jimq V_0
 =\ek*{\zdiag{\rk{yx\OPr{M}}}{\Iq}}^\ad\Jimq\ek*{\zdiag{\rk{yx\OPr{M}}}{\Iq}}\\
 +\abs{y}^2\ek*{\rk{\zdiag{\OPr{M}}{\Iq}}^\ad\Jimq\rk{\zdiag{\OPr{M}}{\Iq}}+2\im\rk{z}\rk{\zdiag{M^\mpi}{\Oqq}}}+2\im\rk{z}\rk{\zdiag{\Oqq}{M}}.
\end{multline*}
 \il{ab.L1059.b} Let \(V_1\defeq\ek{\zdiag{\rk{x\Iq}}{\Iq}}V_0\) and let \(V_2\defeq\ek{\zdiag{\rk{y\Iq}}{\Iq}}V_0\).
 Then
\[\begin{split}
 V_1^\ad\Jimq V_1
 &=\ba\rk*{\ek*{\zdiag{\rk{yx\OPr{M}}}{\Iq}}^\ad\Jimq\ek*{\zdiag{\rk{yx\OPr{M}}}{\Iq}}+2\im\rk{z}\rk{\zdiag{\Oqq}{M}}}\\
 &=\abs{y}^2\ek*{\zdiag{\rk{x\OPr{M}}}{\Iq}}^\ad\Jimq\ek*{\zdiag{\rk{x\OPr{M}}}{\Iq}}\\
 &\qquad+\abs{x}^2\ek*{\zdiag{\rk{y\OPr{M}}}{\Iq}}^\ad\Jimq\ek*{\zdiag{\rk{y\OPr{M}}}{\Iq}}+2\ba\im\rk{z}\rk{\zdiag{\Oqq}{M}}
\end{split}\]
 and
\[
 V_2^\ad\Jimq V_2
 =\ba\abs{y}^2\ek*{\rk{\zdiag{\OPr{M}}{\Iq}}^\ad\Jimq\rk{\zdiag{\OPr{M}}{\Iq}}+2\im\rk{z}\rk{\zdiag{M^\mpi}{\Oqq}}}.
\]
 \il{ab.L1059.c} Let \(V_3\defeq\ek{\zdiag{\rk{yx\Iq}}{\Iq}}V_0\).
 Then
\begin{multline*}
 V_3^\ad\Jimq V_3
 =\abs{y}^2\biggl\{\ek*{\zdiag{\rk{yx\OPr{M}}}{\Iq}}^\ad\Jimq\ek*{\zdiag{\rk{yx\OPr{M}}}{\Iq}}\\
 +\abs{x}^2\ek*{\rk{\zdiag{\OPr{M}}{\Iq}}^\ad\Jimq\rk{\zdiag{\OPr{M}}{\Iq}}+2\im\rk{z}\rk{\zdiag{M^\mpi}{\Oqq}}}+2\im\rk{z}\rk{\zdiag{\Oqq}{M}}\biggr\}.
\end{multline*}
\eenui
\eprop
\bproof
 Consider an arbitrary \(\ell\in\set{0,1,2,3}\).
 Using the notation given in \rprop{ab.L1336}, we obtain, by virtue of \rremss{ab.R1031}{ab.R1025}, then \(V_\ell\Jabsqq=-\Jabsqq\ek{W_\ell+x\rk{\zdiag{\Oqq}{\OPn{M}}}}\) and \(\Jimq\Jabsqq=-\Jabsqq\Jimq\).
 Taking additionally into account \(\Jabsqq^2=\Iu{2q}\) and \(\Jabsqq^\ad=\Jabsqq\) and setting \(U_\ell\defeq W_\ell+x\rk{\zdiag{\Oqq}{\OPn{M}}}\), we can conclude hence
\[
 V_\ell^\ad\Jimq V_\ell
 =\rk{-\Jabsqq U_\ell\Jabsqq}^\ad\Jimq\rk{-\Jabsqq U_\ell\Jabsqq}
 =\Jabsqq U_\ell^\ad\rk{\Jabsqq\Jimq\Jabsqq}U_\ell\Jabsqq
 =-\Jabsqq\rk{U_\ell^\ad\Jimq U_\ell}\Jabsqq.
\]
 Because of \(\OPn{M}^\ad=\OPn{M}\) and \(M^\ad=M\), we have \(\OPn{M}M=\rk{M\OPn{M}}^\ad=\Oqq\).
 Consequently, \(\OPn{M}\OPr{M}=\OPn{M}MM^\mpi=\Oqq\) follows.
 In view of \rnota{ab.N1246a}, thus we obtain
\[
 \rk{\zdiag{\OPn{M}}{\Oqq}}W_0
 =\rk{\zdiag{\OPn{M}}{\Oqq}}
 \bMat
  -y\OPr{M}&M\\
  \ub&\ub
 \eMat
 =
 \bMat
  -y\OPn{M}\OPr{M}&\OPn{M}M\\
  \Oqq&\Oqq
 \eMat
 =\Ouu{2q}{2q}.
\]
 In particular, \(\rk{\zdiag{\OPn{M}}{\Oqq}}W_\ell=\Ouu{2q}{2q}\).
 Using \(\Jimq^\ad=\Jimq\) and \rrem{ab.R1031}, we get then
\begin{multline*}
 U_\ell^\ad\Jimq U_\ell
 =W_\ell^\ad\Jimq W_\ell+2\re\ek*{\ko x\rk{\zdiag{\Oqq}{\OPn{M}}}\Jimq W_\ell}+\ko xx\rk{\zdiag{\Oqq}{\OPn{M}}}\Jimq\rk{\zdiag{\Oqq}{\OPn{M}}}\\
 =W_\ell^\ad\Jimq W_\ell+2\re\ek*{\ko x\Jimq\rk{\zdiag{\OPn{M}}{\Oqq}} W_\ell}+\abs{x}^2\Jimq\rk{\zdiag{\OPn{M}}{\Oqq}}\rk{\zdiag{\Oqq}{\OPn{M}}}
 =W_\ell^\ad\Jimq W_\ell,
\end{multline*}
 implying \(V_\ell^\ad\Jimq V_\ell=-\Jabsqq\rk{W_\ell^\ad\Jimq W_\ell}\Jabsqq\).
 \rrem{ab.R1031} yields \(\Jabsqq\rk{\zdiag{R}{S}}\Jabsqq=\zdiag{R}{S}\) and
\[
 \Jabsqq\rk{\zdiag{R}{S}}^\ad\Jimq\rk{\zdiag{R}{S}}\Jabsqq
 =\rk{\zdiag{R}{S}}^\ad\rk{\Jabsqq\Jimq\Jabsqq}\rk{\zdiag{R}{S}}
 =-\rk{\zdiag{R}{S}}^\ad\Jimq\rk{\zdiag{R}{S}}
\]
 for all \(R,S\in\Cqq\).
 The asserted identities can now be deduced from \rprop{ab.L1336}.
\eproof

 We are now going to consider the composition of the transformations introduced in \rdefnss{F.D.FTF}{F.D.iFTF}.

\bleml{ab.L1403}
 Let \(M\in\Cpq\) and let \(z\in\C\).
 Then
 \begin{align*}
  \ek*{\mFTiua{M}{z}}\ek*{\mFTua{M}{z}}
  =-\rk{\obg-z}\ba\rk{\zdiag{\OPr{M}}{\Iq}}
  =\ek*{\mFTua{M}{z}}\ek*{\mFTiua{M}{z}}.
 \end{align*}
\elem
\bproof
 Let \(x\defeq z-\ug\) and let \(y\defeq\obg-z\).
 We have then
\begin{multline*}
 \ek*{\mFTiua{M}{z}}\ek*{\mFTua{M}{z}}
 =
 \begin{pmat}[{|}]
  y\OPr{M}&M\cr\-
  -yxM^\mpi&y\Iq\cr
 \end{pmat}
 \begin{pmat}[{|}]
  -y\OPr{M}&M\cr\-
  -yxM^\mpi&-y\OPr{M^\ad}-\ba\OPn{M}\cr
 \end{pmat}\\
 =
 \begin{pmat}[{|}]
  -y^2\OPr{M}^2-yxMM^\mpi&y\OPr{M}M-yM\OPr{M^\ad}-\ba M\OPn{M}\cr\-
  y^2xM^\mpi\OPr{M}-y^2xM^\mpi&-yxM^\mpi M-y^2\OPr{M^\ad}-y\ba\OPn{M}\cr
 \end{pmat}
\end{multline*}
 and
\begin{multline*}
 \ek*{\mFTua{M}{z}}\ek*{\mFTiua{M}{z}}
 =
 \begin{pmat}[{|}]
  -y\OPr{M}&M\cr\-
  -yxM^\mpi&-y\OPr{M^\ad}-\ba\OPn{M}\cr
 \end{pmat}
 \begin{pmat}[{|}]
  y\OPr{M}&M\cr\-
  -yxM^\mpi&y\Iq\cr
 \end{pmat}\\
 =
 \begin{pmat}[{|}]
  -y^2\OPr{M}^2-yxMM^\mpi&-y\OPr{M}M+yM\cr\-
  -y^2xM^\mpi\OPr{M}+y^2x\OPr{M^\ad}M^\mpi+yx\ba\OPn{M}M^\mpi&-yxM^\mpi M-y^2\OPr{M^\ad}-y\ba\OPn{M}\cr
 \end{pmat}.
\end{multline*}
 Consequently, in view of \eqref{F.G.PQ2*}, \eqref{F.G.PQ}, and \(y+x=\ba\), the assertion follows.
\eproof

 For a given \tnnH{} matrix \(M\), the condition in \rdefn{F.D.iFTF} is satisfied for pairs belonging to the subclass \(\PRFabqa{M}\) of \(\PRFabq\), introduced in \rsec{F.S.PM}.
 Hence, for suchlike pairs the corresponding \tiFaT{M} exists and can be written as a linear fractional transformation, as considered in \rapp{A.s1.lft}:
 
\bpropl{F.P.Filft}
 Let \(M\in\Cggq\) and let \(\copa{G_1}{G_2}\in\PRFabqa{M}\).
 In view of the functions \(g,h\colon\Cab\to\C\) defined by \eqref{F.D.iFTF.1}, let \(F_1\) and \(F_2\) be given via \eqref{F.D.iFTF.2} as matrix-valued functions meromorphic in \(\Cab\).
 Then \(\det F_2\) does not identically vanish in \(\Cab\).
 Furthermore, \(\det F_2(z)\neq0\) and \(F(z)=\ek{F_1(z)}\ek{F_2(z)}^\inv\) for all \(z\in\C\setminus\rk{\ab\cup\pol{G_1}\cup\pol{G_2}\cup\PRFabex{G_1}{G_2}}\), where \(F\) denotes the \tiFaTv{M}{\(\copa{G_1}{G_2}\)}.
\eprop
\bproof
 According to \rnota{ab.N1534}, the pair \(\copa{G_1}{G_2}\) belongs to \(\PRFabq\).
 Hence, \(G_1\) and \(G_2\) are \(\Cqq\)\nobreakdash-valued functions, which are meromorphic in \(\Cab\).
 Furthermore, by virtue of \rprop{ab.P1528}, the set \(\exset\defeq\pol{G_1}\cup\pol{G_2}\cup\PRFabex{G_1}{G_2}\) is a discrete subset of \(\Cab\).
 Consequently, \(\C\setminus\rk{\ab\cup\exset}\neq\emptyset\).
 Consider an arbitrary \(z\in\C\setminus\rk{\ab\cup\exset}\).
 Then \(G_1\) and \(G_2\) are both holomorphic in \(z\).
 Thus \(F_1\) and \(F_2\) are both holomorphic in \(z\) as well.
 Consider an arbitrary \(v\in\nul{F_2(z)}\).
 Regarding \rrem{A.R.A+>}, we are going to show in a first step that
\beql{F.P.Filft.1}
 \normE*{R\ek*{G_1(z)}v}
 =0
\eeq
 holds true, where \(R\defeq\sqrt{M^\mpi}\). 
 Because of \(z\neq\obg\), we have, according to \eqref{F.D.iFTF.1} and \eqref{F.D.iFTF.2}, the equation
\beql{F.P.Filft.3}
 \rk{z-\ug}M^\mpi\ek*{G_1(z)}v
 =\ek*{G_2(z)}v.
\eeq
 In view of \rrem{A.R.A++*}, hence
\beql{F.P.Filft.2}
 v^\ad\ek*{G_2(z)}^\ad\ek{G_1(z)}v
 =\rk{\ko z-\ug}v^\ad\ek*{G_1(z)}^\ad M^\mpi\ek{G_1(z)}v
 =\rk{\ko z-\ug}\normE*{R\ek{G_1(z)}v}^2.
\eeq
 In the case \(z\in\C\setminus\R\), we see from \rlem{ab.L1252} and \rrem{A.R.XRIX} then that
\[\begin{split}
 0
 \leq v^\ad\rk*{\frac{1}{\im z}\im\rk*{\ek*{G_2(z)}^\ad\ek{G_1(z)}}}v
 &=\frac{1}{\im z}\im\rk*{v^\ad\ek*{G_2(z)}^\ad\ek{G_1(z)}v}\\
 &=\frac{1}{\im z}\im\rk*{\rk{\ko z-\ug}\normE*{R\ek{G_1(z)}v}^2}
 =-\normE*{R\ek{G_1(z)}v}^2
 \leq0,
\end{split}\]
 implying \eqref{F.P.Filft.1}.
 If \(z\in\crhl\), then \(\ko z=z<\ug\) and we obtain, by virtue of \rlem{ab.L1252} and \eqref{F.P.Filft.2}, thus
 \[
  0
  \leq v^\ad\ek*{G_2(z)}^\ad\ek{G_1(z)}v
  =\rk{\ko z-\ug}\normE*{R\ek{G_1(z)}v}^2
  \leq0,
 \]
 implying again \eqref{F.P.Filft.1}.
 In the case \(z\in\clhl\), we have \(\ko z=z>\obg>\ug\) and, because of \rlem{ab.L1252} and \eqref{F.P.Filft.2}, similarly
 \[
  0
  \leq v^\ad\rk*{-\ek*{G_2(z)}^\ad\ek{G_1(z)}}v
  =-v^\ad\ek*{G_2(z)}^\ad\ek{G_1(z)}v
  =\rk{\ug-\ko z}\normE*{R\ek{G_1(z)}v}^2
  \leq0,
 \]
 \tie{}, \eqref{F.P.Filft.1}. 
 Thus, \eqref{F.P.Filft.1} is verified.
 Consequently, we can infer
\[
 \OPr{M}\ek{G_1(z)}v
 =MM^\mpi\ek{G_1(z)}v
 =MR^2\ek*{G_1(z)}v
 =\Ouu{q}{1}.
\]
 In view of \rlem{F.L.A<P}, we have furthermore \(\OPr{M}G_1=G_1\).
 Hence, \(\ek{G_1(z)}v=\Ouu{q}{1}\) follows.
 Because of \eqref{F.P.Filft.3}, this implies \(\ek{G_2(z)}v=\Ouu{q}{1}\).
 Observe that, due to \rprop{ab.P1528}, the \tcp{q}{q} \(\copa{G_1(z)}{G_2(z)}\) is regular.
 According to \rrem{A.R.PQre}, thus necessarily \(v=\Ouu{q}{1}\) holds true.
 Therefore, the linear subspace \(\nul{F_2(z)}\) is trivial, implying \(\det F_2(z)\neq0\).
 In particular, \(\det F_2\) does not identically vanish in \(\Cab\) and \(F(z)=\ek{F_1(z)}\ek{F_2(z)}^\inv\).
\eproof

 For any \tnnH{} matrix \(M\), the \tiFATion{} induces, according to \rdefn{ab.N0843} and \rrem{F.R.PMaeq}, a well-defined transformation for equivalence classes from \(\rsetcl{\PRFabqa{M}}\):

\bcorl{F.C.iFwd}
 Let \(M\in\Cggq\) and let the pairs \(\copa{G_1}{G_2},\copa{\tilde G_1}{\tilde G_2}\in\PRFabqa{M}\) be equivalent.
 Then the \tiFaT{M} \(F\) of \(\copa{G_1}{G_2}\) coincides with the \tiFaT{M} \(\tilde F\) of \(\copa{\tilde G_1}{\tilde G_2}\).
\ecor
\bproof
 Using the functions \(g,h\colon\Cab\to\C\) given via \eqref{F.D.iFTF.1}, we define by \eqref{F.D.iFTF.2} two \(\Cqq\)\nobreakdash-valued functions \(F_1\) and \(F_2\) meromorphic in \(\Cab\).
 Because of \rprop{F.P.Filft}, then \(\det F_2\) does not vanish identically in \(\Cab\).
 By virtue of \rdefn{F.D.iFTF}, we thus have \(F=F_1F_2^\inv\).
 Furthermore, due to \rdefn{ab.N0843}, there exists a \(\Cqq\)\nobreakdash-valued function \(R\) meromorphic in \(\Cab\) such that \(\det R\) does not vanish identically in \(\Cab\) satisfying \(\tilde G_1=G_1R\) and \(\tilde G_2=G_2R\).
 Using again the functions \(g\) and \(h\), we define according to \eqref{F.D.iFTF.2} by \(\tilde F_1\defeq h\OPr{M}\tilde G_1+M\tilde G_2\) and \(\tilde F_2\defeq-hgM^\mpi\tilde G_1+h\tilde G_2\) two \(\Cqq\)\nobreakdash-valued functions \(\tilde F_1\) and \(\tilde F_2\) meromorphic in \(\Cab\).
 Then \(\tilde F_1=F_1R\) and \(\tilde F_2=F_2R\).
 The application of \rprop{F.P.Filft} to the pair \(\copa{\tilde G_1}{\tilde G_2}\) yields furthermore that \(\det\tilde F_2\) does not vanish identically in \(\Cab\).
 By virtue of \rdefn{F.D.iFTF}, hence \(\tilde F=\tilde F_1\tilde F_2^\inv\).
 Consequently, \(\tilde F=\rk{F_1R}\rk{F_2R}^\inv=F\).
\eproof

 Furthermore, after transition to equivalence classes as mentioned above, the \tFATion{} turns out to be inverse to the \tiFATion{}:

\bleml{F.L.FT-11}
 Let \(M\in\Cggq\) and let \(\copa{G_1}{G_2}\in\PRFabqa{M}\) with \tiFaT{M} \(F\).
 Then \(\copa{G_1}{G_2}\) is equivalent to the \tFatpv{M}{\(F\)}.
\elem
\bproof
 Using the functions \(g,h\colon\Cab\to\C\) given via \eqref{F.D.iFTF.1}, we define by \eqref{F.D.iFTF.2} two \(\Cqq\)\nobreakdash-valued functions \(F_1\) and \(F_2\) meromorphic in \(\Cab\).
 Then \(\OPr{M}F_1=F_1\).
 Denote by \(\copa{H_1}{H_2}\) the \tFatpv{M}{\(F\)} and by \(W\) the restriction of the holomorphic \(\Coo{2q}{2q}\)\nobreakdash-valued function \(\mFTu{M}\) onto \(\Cab\).
 In view of \rdefn{F.D.FTF} and \rnota{ab.N1246a}, then \(\tmatp{H_1}{H_2}=-W\tmatp{F}{\Iq}\).
 Due to \rprop{F.P.Filft}, the function \(\det F_2\) does not vanish identically in \(\Cab\).
 According to \rdefn{F.D.iFTF}, we thus have \(F=F_1F_2^\inv\).
 Denote by \(V\) the restriction of the holomorphic \(\Coo{2q}{2q}\)\nobreakdash-valued function \(\mFTiu{M}\) onto \(\Cab\).
 Regarding \rnota{ab.N1246b}, then \(\tmatp{F_1}{F_2}=V\tmatp{G_1}{G_2}\).
 Taken all together, we get
\[
 \matp{H_1}{H_2}
 =-W\matp{F_1F_2^\inv}{\Iq}
 =-W\matp{F_1}{F_2}F_2^\inv
 =-WV\matp{G_1}{G_2}F_2^\inv
 =-WV\matp{G_1F_2^\inv}{G_2F_2^\inv}.
\]
 From \rlem{F.L.A<P} we see \(\OPr{M}G_1=G_1\).
 Using \rlem{ab.L1403}, we thus can infer \(H_1=h\ba\OPr{M}G_1F_2^\inv=h\ba G_1F_2^\inv\) and \(H_2=h\ba G_2F_2^\inv\).
 Observe that \(R\defeq\ba h F_2^\inv\) is a \(\Cqq\)\nobreakdash-valued function, which is meromorphic in \(\Cab\) satisfying \(H_1=G_1R\) and \(H_2=G_2R\).
 Furthermore, because of \(\ba\neq0\), the function \(\det R\) does not vanish identically in \(\Cab\).
 According to \rdefn{ab.N0843}, consequently \(\copa{G_1}{G_2}\rpaeq\copa{H_1}{H_2}\).
\eproof

 Conversely, we have:

\bleml{F.L.FT1-1}
 Let \(M\in\Cggq\) and let \(F\colon\Cab\to\Cqq\) be a matrix-valued function with \tFatp{M} \(\copa{G_1}{G_2}\) such that \(\OPr{M}F=F\) and \(\copa{G_1}{G_2}\in\PRFabqa{M}\) hold true.
 Then the \tiFaTv{M}{\(\copa{G_1}{G_2}\)} coincides with \(F\).
\elem
\bproof
 Since \(\copa{G_1}{G_2}\) belongs to \(\PRFabqa{M}\), we see that \(G_1\) and \(G_2\) are \(\Cqq\)\nobreakdash-valued functions, which are meromorphic in \(\Cab\).
 Using the functions \(g,h\colon\Cab\to\C\) given via \eqref{F.D.iFTF.1}, we can thus define by \eqref{F.D.iFTF.2} two \(\Cqq\)\nobreakdash-valued functions \(F_1\) and \(F_2\), which then are meromorphic in \(\Cab\) as well.
 Denote by \(V\) the restriction of the holomorphic \(\Coo{2q}{2q}\)\nobreakdash-valued function \(\mFTiu{M}\) onto \(\Cab\).
 From \rnota{ab.N1246b} we see \(\tmatp{F_1}{F_2}=V\tmatp{G_1}{G_2}\).
 Denote by \(W\) the restriction of the holomorphic \(\Coo{2q}{2q}\)\nobreakdash-valued function \(\mFTu{M}\) onto \(\Cab\).
 Regarding \rdefn{F.D.FTF}, \rnota{ab.N1246a}, and \(\OPr{M}F=F\), we have furthermore \(W\tmatp{F}{\Iq}=-\tmatp{G_1}{G_2}\).
 Taken all together, we obtain \(\tmatp{F_1}{F_2}=-VW\tmatp{F}{\Iq}\).
 In view of \rlem{ab.L1403}, thus \(F_1=h\ba\OPr{M}F=h\ba F\) and \(F_2=h\ba\Iq\) follow.
 Taking into account \(\copa{G_1}{G_2}\in\PRFabqa{M}\), we see from \rprop{F.P.Filft} that \(\det F_2\) does not vanish identically in \(\Cab\).
 Denote by \(H\) the \tiFaTv{M}{\(\copa{G_1}{G_2}\)}.
 According to \rdefn{F.D.iFTF}, then \(H=F_1F_2^\inv\).
 Consequently, \(H=F\).
\eproof

\section{The \hFABTion{} and its inverse}\label{S1230}
 In this section, we continue the preceding considerations concerning the construction of the function-theoretic version of the \tSchur{}--\tNevanlinna{} type algorithm.
 We prepare the algebraic formalism for the remaining steps after the first one.

 In what follows, we consider two complex \tqqa{matrices} \(A\) and \(M\), which, in the context of the matricial Hausdorff moment problem, will later be the \tnnH{} matrices \(\sau{0}\) and \(\su{0}\) for a given sequence \(\seqska\in\Fggqka\).
 In this reading, the matrices
\begin{align}\label{F.G.BN}
 B&\defeq\ba M-A&
&\text{and}&
 N&\defeq A+\ug M
\end{align}
 correspond to \(\sub{0}\) and \(\su{1}\), \tresp{}, and we have%
\begin{align}\label{F.G.AB}
 A&=-\ug M+N&
&\text{and}&
 B&=\obg M-N,
\end{align}
 according to \rnota{F.N.sa}.
 Consider an arbitrarily given \(z\in\C\) and let \(x\defeq z-\ug\) and \(y\defeq \obg-z\).
 Taking additionally into account \eqref{F.G.y+x} and \eqref{F.G.ax+by}, we then infer
\begin{align}
 yA-xB
 &=\rk{y+x}N-\rk{\ug y+\obg x}M
 =\ba\rk{N-zM}\label{F.G.yA-xB}
\shortintertext{and}
 xA-yB
 &=\rk{x+y}N-\rk{\ug x+\obg y}M 
 =\ba\ek*{N-\rk{\obg+\ug-z}M}.\label{F.G.xA-yB}
\end{align}
 
\bdefnl{F.D.FTFAB}
 Let \(\dom\) be a non-empty subset of \(\C\), let \(F\colon\dom\to\Cpq\) be a matrix-valued function, and let \(A\) and \(M\) be two complex \tpqa{matrices}.
 Then \(G\colon\dom\to\Cpq\) defined by
\[
 G(z)
 \defeq AM^\mpi\ek*{\rk{\obg-z}F(z)-M}\rk*{\rk{\obg-z}\ek*{\rk{z-\ug}F(z)+M}}^\mpi A
\]
 is called the \emph{\tFaaTv{A}{M}{\(F\)}}.
\edefn

 In connection with the \tFABTion{}, we consider the following complex \taaa{(p+q)}{(p+q)}{matrix} polynomial:

\bnotal{ab.N1546a}
 Let \(A\) and \(M\) be two complex \tpqa{matrices}.
 Then let \(\mFTuu{A}{M}\colon\C\to\Coo{(p+q)}{(p+q)}\) be defined by
\[
 \mFTuua{A}{M}{z}
 \defeq
 \begin{pmat}[{|}]
  -\rk{\obg-z}AM^\mpi&A\cr\-
  -\rk{\obg-z}\rk{z-\ug}A^\mpi&-\rk{\obg-z}A^\mpi M-\OPn{A}\cr
 \end{pmat}.
\]
\enota

 Under certain conditions, we can write the \tFaaT{A}{M} as a linear fractional transformation with the generating matrix-valued function \(\mFTuu{A}{M}\).

\bleml{F.L.FAMlft}
 Let \(A,M\in\Cpq\) be such that \(\nul{M}\subseteq\nul{A}\).
 Let \(F\colon\Cab\to\Cpq\) be a matrix-valued function with \tFaaT{A}{M} \(G\) and let \(G_1,G_2\colon\Cab\to\Cqq\) be defined by
\begin{align}
 G_1(w)&\defeq-(\obg-w)AM^\mpi F(w)+A\label{F.L.FAMlft.G1}
\shortintertext{and}
 G_2(w)&\defeq-\rk{\obg-w}\rk{w-\ug}A^\mpi F(w)-\rk{\obg-w}A^\mpi M-\OPn{A}.\label{F.L.FAMlft.G2}
\end{align}
 Let \(z\in\Cab\) be such that \(\ran{F(z)}\subseteq\ran{M}\) and \(\nul{M}\subseteq\nul{F(z)}\) as well as \(\ran{{\rk{z-\ug}}F(z)+M}=\ran{A}\) and \(\nul{\rk{z-\ug}F(z)+M}=\nul{A}\) are fulfilled.
 Then \(\det G_2(z)\neq0\) and \(G(z)=\ek{G_1(z)}\ek{G_2(z)}^\inv\).
\elem
\bproof
 Let \(X\defeq\rk{z-\ug}F(z)+M\) and let \(Y\defeq\rk{\obg-z}F(z)-M\).
 From \rrem{R.AA+B.B} we get \(AM^\mpi M=A\).
 \rrem{ab.R1052} shows that \(\Iq-A^\mpi A=\OPn{A}\).
 Setting \(y\defeq\obg-z\) and \(Z\defeq yX\), we obtain then
\beql{F.L.FAMlft.1}
 -AM^\mpi Y
 =-yAM^\mpi F(z)+AM^\mpi M
 =-yAM^\mpi F(z)+A
 =G_1(z)
\eeq
 and
\beql{F.L.FAMlft.2}
 -A^\mpi Z-\rk{\Iq-A^\mpi A}
 =-yA^\mpi X-\OPn{A}
 =G_2(z).
\eeq
 Using \rrem{A.R.AM+B}, we get \(XM^\mpi Y=YM^\mpi X\).
 In view of \rrem{ab.R1052}, furthermore \(X^\mpi X=A^\mpi A\) holds true.
 Thus, we can infer
\[
 AM^\mpi Y
 =AA^\mpi AM^\mpi Y
 =AX^\mpi XM^\mpi Y
 =AX^\mpi YM^\mpi X
\]
 and, therefore,
\[
 AM^\mpi YA^\mpi A
 =AX^\mpi YM^\mpi XA^\mpi A
 =AX^\mpi YM^\mpi XX^\mpi X
 =AX^\mpi YM^\mpi X
 =AM^\mpi Y.
\]
 Regarding \(y\neq0\), we have \(\ran{Z}=\ran{A}\) and \(\nul{Z}=\nul{A}\).
 By virtue of \rrem{ab.R1052}, hence \(ZZ^\mpi=AA^\mpi\) and \(Z^\mpi Z=A^\mpi A\).
 Consequently,
\[\begin{split}
 \ek*{-Z^\mpi A-\rk{\Iq-A^\mpi A}}&\ek*{-A^\mpi Z-\rk{\Iq-A^\mpi A}}\\
 &=Z^\mpi AA^\mpi Z+Z^\mpi A(\Iq-A^\mpi A)+(\Iq-A^\mpi A)A^\mpi Z+(\Iq-A^\mpi A)^2\\
 &=Z^\mpi ZZ^\mpi Z+\Iq-A^\mpi A
 =Z^\mpi Z+\Iq-Z^\mpi Z
 =\Iq.
\end{split}\]
 Hence, \(\det\ek{-A^\mpi Z-\rk{\Iq-A^\mpi A}}\neq0\) and \(\ek{-A^\mpi Z-\rk{\Iq-A^\mpi A}}^\inv=-Z^\mpi A-\rk{\Iq-A^\mpi A}\).
 Thus,
\[\begin{split}
 -AM^\mpi Y\ek*{-A^\mpi Z-\rk{\Iq-A^\mpi A}}^\inv
 &=AM^\mpi Y\ek*{Z^\mpi A+\rk{\Iq-A^\mpi A}}\\
 &=AM^\mpi YZ^\mpi A+AM^\mpi Y\rk{\Iq-A^\mpi A}
 =AM^\mpi YZ^\mpi A
 =G(z).
\end{split}\]
 In view of \eqref{F.L.FAMlft.1} and \eqref{F.L.FAMlft.2}, the proof is complete.
\eproof

 From \rlemss{F.L.FAM-11}{F.L.FAM1-1} we will see that the following transformation for matrix-valued functions is in a generic situation indeed the inversion of the \tFABTion{}.
 Against this background we introduce the following notation:

\bdefnl{F.D.iFTFAB}
 Let \(\dom\) be a non-empty subset of \(\C\), let \(G\colon\dom\to\Cpq\) be a matrix-valued function, and let \(A\) and \(M\) be two complex \tpqa{matrices}.
 Let \(B\defeq\ba M-A\) and let \(F\colon\dom\to\Cpq\) be defined by
\begin{multline*}
 F(z)
 \defeq-\ek*{\rk{\obg-z}MA^\mpi G(z)+A+M\OPn{A}M^\mpi B}\\
 \times\rk*{\rk{\obg-z}\ek*{\rk{z-\ug}A^\mpi G(z)-M^\mpi A}+\rk{z-\ug}\OPn{A}M^\mpi B}^\mpi.
\end{multline*}
 Then we call the matrix-valued function \(F\) the \emph{\tiFaaTv{A}{M}{\(G\)}}.
\edefn

\bleml{F.L.iFAMhol}
 Let \(A\in\CHq\) and let \(M\in\Cggq\) with \(\ran{A}\subseteq\ran{M}\).
 Let \(G\in\RFqab\) be such that \(\ran{G(z_0)}\subseteq\ran{A}\) holds true for some \(z_0\in\Cab\).
 Let \(B\defeq\ba M-A\) and let \(E_1,E_2\colon\Cab\to\Cqq\) be defined by
\begin{align}
 E_1(w)&\defeq\rk{\obg-w}MA^\mpi G(w)+A+M\OPn{A}M^\mpi B\label{F.L.iFAMhol.E1}
\shortintertext{and}
 E_2(w)&\defeq-\rk{\obg-w}\rk{w-\ug}A^\mpi G(w)+\rk{\obg-w}M^\mpi A-\rk{w-\ug}\OPn{A}M^\mpi B.\label{F.L.iFAMhol.E2}
\end{align}
 For all \(z\in\Cab\), then \(\ran{E_1(z)}\subseteq\ran{M}\), \(\nul{M}\subseteq\nul{E_1(z)}\), \(\ran{E_2(z)}=\ran{M}\), and \(\nul{E_2(z)}=\nul{M}\).
 Furthermore, the \tiFaaTv{A}{M}{\(G\)} is holomorphic in \(\Cab\).
\elem
\bproof
 Consider an arbitrary \(z\in\Cab\).
 Regarding \rrem{A.R.rs+}, we have \(\ran{E_1(z)}\subseteq\ran{MA^\mpi G(z)}+\ran{A}+\ran{M\OPn{A}M^\mpi B}\subseteq\ran{M}\).
 From \rprop{ab.P1648L1409} we infer \(\ran{G(\ko z)}=\ran{G(z_0)}\subseteq\ran{A}\subseteq\ran{M}\).
 Consequently, \(\ek{\ran{M}}^\orth\subseteq\ek{\ran{A}}^\orth\subseteq\ek{\Ran{G(\ko z)}}^\orth\).
 In view of \rremss{ab.R1842*}{ab.R1652}, then
\beql{F.L.iFAMhol.0}
 \nul{M}
 \subseteq\nul{A}
 \subseteq\Nul{\ek*{G(\ko z)}^\ad}
 =\Nul{G(z)}
\eeq
 follows.
 Regarding \(\nul{M}\subseteq\nul{A}\) and \rrem{A.R.rs+}, we can conclude \(\nul{M}\subseteq\nul{B}\) and  furthermore \(\nul{M}\subseteq\nul{MA^\mpi G(z)}\cap\nul{A}\cap\nul{M\OPn{A}M^\mpi B}\subseteq\nul{E_1(z)}\).
 By virtue of \rrem{A.R.A++*}, we have
\(
 \ran{A^\mpi}
 =\ran{A}
 \subseteq\ran{M}
 =\ran{M^\mpi}
\).
 In view of \rrem{A.R.A++*}, then \(M^\mpi M A^\mpi=A^\mpi\) holds true, according to \rrem{R.AA+B.A}.
 Since \rrem{ab.R1052} shows that \(\OPn{A}=\Iq-A^\mpi A\), we infer in particular \(\OPn{A}M^\mpi=M^\mpi\rk{\Iq-MA^\mpi A M^\mpi}\).
 Therefore, we obtain \(\ran{\OPn{A}M^\mpi}\subseteq\ran{M^\mpi}\) and, regarding additionally \rrem{A.R.rs+}, furthermore \(\ran{E_2(z)}\subseteq\ran{A^\mpi G(z)}+\ran{M^\mpi A}+\ran{\OPn{A}M^\mpi B}\subseteq\ran{M^\mpi}=\ran{M}\).
 In view of \rrem{A.R.A+>}, let \(R\defeq\sqrt{M^\mpi}A\).
 Consider an arbitrary \(v\in\nul{E_2(z)}\).
 We are now going to check that
 \beql{F.L.iFAMhol.1}
  \normE{Rv}
  =0
 \eeq
 holds true.
 Obviously,
\beql{F.L.iFAMhol.4}
 \rk{\obg-z}M^\mpi Av
 =\rk{\obg-z}\rk{z-\ug}A^\mpi\ek*{G(z)}v+\rk{z-\ug}\OPn{A}M^\mpi Bv.
\eeq
 In view of \rprop{ab.P1648L1409}, we have \(\ran{G(z)}=\ran{G(z_0)}\subseteq\ran{A}\).
 According to \rrem{R.AA+B.A}, thus \(AA^\mpi G(z)=G(z)\).
 Regarding \(z\neq\obg\), we can multiply equation~\eqref{F.L.iFAMhol.4} from the left by \(\rk{\obg-z}^\inv A\) to obtain then \(AM^\mpi Av=\rk{z-\ug}\ek*{G(z)}v\).
 Left multiplication of this identity by \(\rk{\ko z-\ug}v^\ad\) yields
\beql{F.L.iFAMhol.2}
 \abs{z-\ug}^2v^\ad\ek*{G(z)}v
 =\rk{\ko z-\ug}v^\ad AM^\mpi Av
 =\rk{\ko z-\ug}\normEs{Rv}.
\eeq
 In the case \(z\in\C\setminus\R\), we can infer from \rprop{F.P.RabH} and \rrem{A.R.XRIX} that
\[\begin{split}
 0
 \leq\abs{z-\ug}^2v^\ad\ek*{\frac{1}{\im z}\im G(z)}v
 &=\frac{1}{\im z}\im\rk*{\abs{z-\ug}^2v^\ad\ek*{G(z)}v}\\
 &=\frac{1}{\im z}\im\ek*{\rk{\ko z-\ug}\normEs{Rv}}
 =-\normEs{Rv}
 \leq0,
\end{split}\]
 implying \eqref{F.L.iFAMhol.1}.
 If \(z\in\crhl\), then \(\ko z=z<\ug\) and we obtain, by virtue of \rnota{ab.N1603} and \eqref{F.L.iFAMhol.2}, thus
\(
 0
 \leq\abs{z-\ug}^2v^\ad\ek{G(z)}v
 =\rk{\ko z-\ug}\normEs{Rv}
 \leq0
\),
 implying again \eqref{F.L.iFAMhol.1}.
 In the case \(z\in\clhl\), we have \(\ko z=z>\obg>\ug\) and, because of \rnota{ab.N1603} and \eqref{F.L.iFAMhol.2}, similarly
 \[
  0
  \leq\abs{z-\ug}^2v^\ad\ek*{-G(z)}v
  =-\abs{z-\ug}^2v^\ad\ek*{G(z)}v
  =\rk{\ug-\ko z}\normEs{Rv}
  \leq0,
 \]
 \tie{}\ \eqref{F.L.iFAMhol.1}. 
 Hence, \eqref{F.L.iFAMhol.1} is verified.
 Consequently, using \rrem{R.AA+B.A}, we can infer
\(
 Av
 =MM^\mpi Av
 =M\sqrt{M^\mpi}Rv
 =\Ouu{q}{1}
\).
 Regarding \eqref{F.L.iFAMhol.0}, we thus obtain \(\ek{G(z)}v=\Ouu{q}{1}\).
 Because of \eqref{F.L.iFAMhol.0} and \rrem{R.AA+B.B}, we have \(AM^\mpi M=A\).
 In view of \(\OPn{A}=\Iq-A^\mpi A\), in particular \(\OPn{A}M^\mpi M=M^\mpi M-A^\mpi A\) holds true.
 Taking into account \(Av=\Ouu{q}{1}\) and \(\ek{G(z)}v=\Ouu{q}{1}\), we see from \eqref{F.L.iFAMhol.4} then
\[
 \Ouu{q}{1}
 =\rk{z-\ug}\OPn{A}M^\mpi Bv
 =\rk{z-\ug}\OPn{A}M^\mpi\rk{\ba M}v
 =\rk{z-\ug}\ba M^\mpi Mv.
\]
 Left multiplication of the latter by \(M\) yields \(\rk{z-\ug}\ba Mv=\Ouu{q}{1}\).
 Since \(z\neq\ug\) and \(\ba>0\) hold true, necessarily \(Mv=\Ouu{q}{1}\) follows.
 Hence, \(\nul{E_2(z)}\subseteq\nul{M}\).
 Because of \eqref{F.L.iFAMhol.0} and \(\nul{M}\subseteq\nul{B}\), we obtain, by virtue of \rrem{A.R.rs+}, on the other hand \(\nul{M}\subseteq\nul{A^\mpi G(z)}\cap\nul{M^\mpi A}\cap\nul{\OPn{A}M^\mpi B}\subseteq\nul{E_2(z)}\).
 Consequently, \(\nul{E_2(z)}=\nul{M}\) is verified.
 Using \rrem{A.R.RNT}, we can, in view of \(\ran{E_2(z)}\subseteq\ran{M}\), then easily conclude \(\ran{E_2(z)}=\ran{M}\).
 Observe that the matrix-valued function \(G\) is holomorphic.
 Thus, \(E_1\) and \(E_2\) are holomorphic in \(\Cab\) as well.
 Let \(D_2\colon\Cab\to\Cqq\) be defined by \(D_2(w)\defeq\ek{E_2(w)}^\mpi\).
 As already shown, the linear subspaces \(\ran{E_2(w)}\) and \(\nul{E_2(w)}\) do not depend on the point \(w\in\Cab\).
 Due to \rprop{ab.L0850}, thus the matrix-valued function \(D_2\) is holomorphic.
 Denote by \(F\) the \tiFaaTv{A}{M}{\(G\)}.
 In view of \rremss{A.R.l*A}{A.R.A-1}, we have \(E_1D_2=F\).
 Using \rrem{B.R.cp}, we can conclude then that the matrix-valued function \(F\) is holomorphic in \(\Cab\).
\eproof

 The following complex \taaa{(p+q)}{(p+q)}{matrix} polynomial is intimately connected to the \tiFaaT{A}{M}:

\bnotal{ab.N1546b}
 Let \(A\) and \(M\) be two complex \tpqa{matrices} and let \(B\defeq\ba M-A\).
 Then let \(\mFTiuu{A}{M}\colon\C\to\Coo{(p+q)}{(p+q)}\) be defined by
\[
 \mFTiuua{A}{M}{z}
 \defeq
 \begin{pmat}[{|}]
  \rk{\obg-z}MA^\mpi&A+M\OPn{A}M^\mpi B\cr\-
  -\rk{\obg-z}\rk{z-\ug}A^\mpi&\rk{\obg-z}\rk{\ba\OPn{M}+M^\mpi A}-\rk{z-\ug}\OPn{A}M^\mpi B\cr
 \end{pmat}.
\]
\enota

\breml{R1210} 
 Let \(A,M\in\Cqq\) and let \(z\in\C\).
 Let \(B\defeq\ba M-A\).
 Then:
\benui
 \il{R1210.a} If \(A=\Oqq\), then \(B=\ba M\) and
\(
 \mFTiuua{A}{M}{z}
 =\ba\smat{\Oqq&M\\
  \Oqq&\rk{\obg-z}\Iq-\ba M^\mpi M}
\).
 \il{R1210.b} If \(B=\Oqq\), then \(A=\ba M\) and
\(
 \mFTiuua{A}{M}{z}
 =\smat{\ba^\inv\rk{\obg-z}MM^\mpi&\ba M\\
  -\ba^\inv\rk{\obg-z}\rk{z-\ug}M^\mpi&\ba\rk{\obg-z}\Iq}
\).
 \il{R1210.c} If all the matrices \(M,A,B\) are invertible, then
\(
 \mFTiuua{A}{M}{z}
 =\smat{\rk{\obg-z}MA^\inv&A\\
  -\rk{\obg-z}\rk{z-\ug}A^\inv&\rk{\obg-z}M^\inv A}
\).
\eenui
\erem

 Under certain conditions, we can write the \tiFaaT{A}{M} as a linear fractional transformation with the generating matrix-valued function \(\mFTiuu{A}{M}\).

\bleml{F.L.iFAMlft}
 Let \(A\in\CHq\) and let \(M\in\Cggq\) with \(\ran{A}\subseteq\ran{M}\).
 Let \(G\in\RFqab\) with \tiFaaT{A}{M} \(F\) be such that \(\ran{G(z_0)}\subseteq\ran{A}\) holds true for some \(z_0\in\Cab\).
 Let \(B\defeq\ba M-A\) and let \(F_1,F_2\colon\Cab\to\Cqq\) be defined by
\begin{align}
 F_1(w)&\defeq\rk{\obg-w}MA^\mpi G(w)+A+M\OPn{A}M^\mpi B\label{F.L.iFAMlft.F1}
\shortintertext{and}
 F_2(w)&\defeq-\rk{\obg-w}\rk{w-\ug}A^\mpi G(w)+\rk{\obg-w}\rk{\ba\OPn{M}+M^\mpi A}-\rk{w-\ug}\OPn{A}M^\mpi B.\label{F.L.iFAMlft.F2}
\end{align}
 For all \(z\in\Cab\), then \(\det F_2(z)\neq0\) and \(F(z)=\ek{F_1(z)}\ek{F_2(z)}^\inv\).
\elem
\bproof
 Let \(E_1,E_2\colon\Cab\to\Cqq\) be defined by \eqref{F.L.iFAMhol.E1} and \eqref{F.L.iFAMhol.E2}.
 Consider an arbitrary \(z\in\Cab\).
 We have \(F_1(z)=E_1(z)\) and \(F_2(z)=E_2(z)+\rk{\obg-z}\ba\OPn{M}\).
 From \rlem{F.L.iFAMhol} we get \(\ran{E_2(z)}=\ran{M}\) and \(\nul{E_2(z)}=\nul{M}\).
 \rrem{ab.R1842*} yields then \(\ran{E_2(z)}=\ran{M^\ad}=\ran{\ek{E_2(z)}^\ad}\).
 In view of \(z\neq\obg\) and \(\ba>0\), we thus can apply \rlem{141.B463} with \(\eta\defeq\rk{\obg-z}\ba\) to see that the matrix \(E_2(z)+\eta\OPn{M}=F_2(z)\) is invertible and that \(\ek{E_2(z)}^\mpi=\ek{F_2(z)}^\inv-\eta^\inv\OPn{M}\) holds true.
 By virtue of \rlem{F.L.iFAMhol}, we have \(\nul{M}\subseteq\nul{E_1(z)}\).
 Consequently, we obtain
\[
 \ek*{F_1(z)}\ek*{F_2(z)}^\inv
 =\ek*{E_1(z)}\rk*{\ek*{E_2(z)}^\mpi+\eta^\inv\OPn{M}}
 =\ek*{E_1(z)}\ek*{E_2(z)}^\mpi
 =F(z).\qedhere
\]
\eproof

\bleml{F.L.iFAMrn}
 Let \(A\in\CHq\) and let \(M\in\Cggq\) with \(\ran{A}\subseteq\ran{M}\).
 Let \(G\in\RFqab\) with \tiFaaT{A}{M} \(F\).
 Suppose that \(\ran{G(z_0)}\subseteq\ran{A}\) holds true for some \(z_0\in\Cab\).
 For all \(z\in\Cab\), then \(\ran{F(z)}\subseteq\ran{M}\) and \(\nul{M}\subseteq\nul{F(z)}\) and furthermore \(\ran{\rk{z-\ug}F(z)+M}=\ran{A}\) and \(\nul{\rk{z-\ug}F(z)+M}=\nul{A}\).
\elem
\bproof
 Consider an arbitrary \(z\in\Cab\).
 Let \(E_1,E_2\colon\Cab\to\Cqq\) be defined by \eqref{F.L.iFAMhol.E1} and \eqref{F.L.iFAMhol.E2}.
 According to \rdefn{F.D.iFTFAB}, then \(F(z)=\ek{E_1(z)}\ek{E_2(z)}^\mpi\).
 \rlem{F.L.iFAMhol} yields furthermore \(\ran{E_1(z)}\subseteq\ran{M}\) and \(\ran{E_2(z)}=\ran{M}\).
 By virtue of \rremss{A.R.A++*}{ab.R1842*}, we can infer from the last identity \(\nul{\ek{E_2(z)}^\mpi}=\nul{M^\ad}=\nul{M}\).
 Consequently, we obtain \(\ran{F(z)}\subseteq\ran{M}\) and \(\nul{M}\subseteq\nul{F(z)}\).
 Taking additionally into account \rrem{A.R.A++*}, we can conclude \(M\ek{E_2(z)}\ek{E_2(z)}^\mpi=M\).
 \rrem{R.AA+B.A} yields \(MM^\mpi A=A\).
 Let \(x\defeq z-\ug\), let \(y\defeq\obg-z\), and let \(X\defeq xF(z)+M\).
 Taken all together, we get
\[
 X
 =\ek*{xE_1(z)+ME_2(z)}\ek*{E_2(z)}^\mpi
 =\rk{xA+yMM^\mpi A}\ek*{E_2(z)}^\mpi
 =\ba A\ek*{E_2(z)}^\mpi.
\]
 Analogous to the corresponding considerations in the proof of \rlem{F.L.iFAMlft}, we can show that the matrix \(R\defeq E_2(z)+\OPn{M}\) is invertible and that \(\ek{E_2(z)}^\mpi=R^\inv-\OPn{M}\) holds true.
 As in the proof of \rlem{F.L.iFAMhol}, we can obtain \eqref{F.L.iFAMhol.0}.
 Thus, \(X=\ba A\rk{R^\inv-\OPn{M}}=\ba AR^\inv\) follows.
 Regarding \(\ba>0\), we see from \rrem{A.R.rnLAR} hence \(\ran{X}=\ran{A}\) and \(\nul{X}=R\nul{A}\).
 Let \(B\defeq\ba M-A\).
 In view of \eqref{F.L.iFAMhol.E2} and \eqref{F.L.iFAMhol.0}, each \(v\in\nul{A}\) satisfies \(Rv=-x\OPn{A}M^\mpi Bv+\OPn{M}v\) and thus \(ARv=\Ouu{q}{1}\).
 Consequently, \(\nul{X}\subseteq\nul{A}\) is verified.
 Taking additionally into account \(\ran{X}=\ran{A}\), we infer by virtue of \rrem{A.R.RNT} then easily \(\nul{X}=\nul{A}\).
\eproof

 Now we are going to study the composition of the two transformations introduced in \rdefnss{F.D.FTFAB}{F.D.iFTFAB}.
 Doing this, we will take into account that, in view of \rlemss{F.L.FAMlft}{F.L.iFAMlft}, these transformations can be written under certain conditions as linear fractional transformations of matrices with generating matrix-valued functions \(\mFTuu{A}{M}\) and \(\mFTiuu{A}{M}\), \tresp{}

\bleml{ab.L0852r}
 Let \(A,M\in\Cpq\) with \(\ran{A}\subseteq\ran{M}\) and \(\nul{M}\subseteq\nul{A}\), let \(B\defeq\ba M-A\), and let \(N\defeq A+\ug M\).
 Let \(z\in\C\), let \(x\defeq z-\ug\), and let \(y\defeq\obg-z\).
 Then
\[
 \ek*{\mFTuua{A}{M}{z}}\ek*{\mFTiuua{A}{M}{z}}
 =-\ba\rk*{\zdiag{\ek*{y\OPr{A}}}{\ek*{y\OPr{A^\ad}-\obg\OPn{M}+\OPn{A}\rk{z\Iq-M^\mpi N}}}}
\]
 and
\beql{ab.L0852r.B}
 y\OPr{A^\ad}-\obg\OPn{M}+\OPn{A}\rk{z\Iq-M^\mpi N}
 =\rk{M^\mpi B+\OPr{A^\ad}M^\mpi A}-y\OPn{A}-x\OPr{A^\ad}.
\eeq
\elem
\bproof
 For the \tbr{} \(\ek{\mFTuua{A}{M}{z}}\ek{\mFTiuua{A}{M}{z}}=\tmat{X_{11}&X_{12}\\X_{21}&X_{22}}\) with \tppa{block} \(X_{11}\), we have
\begin{align*}
 X_{11}&=-yAM^\mpi\rk{yMA^\mpi}+A\rk{-yxA^\mpi},\\
 X_{12}&=-yAM^\mpi\rk{A+M\OPn{A}M^\mpi B}+A\ek*{y\rk{\ba\OPn{M}+M^\mpi A}-x\OPn{A}M^\mpi B},\\
 X_{21}&=-yxA^\mpi\rk{yMA^\mpi}-\rk{yA^\mpi M-\OPn{A}}\rk{-yxA^\mpi},
\shortintertext{and}
 X_{22}&=-yxA^\mpi\rk{A+M\OPn{A}M^\mpi B }-\rk{yA^\mpi M-\OPn{A}}\ek*{y\rk{\ba\OPn{M}+M^\mpi A}-x\OPn{A}M^\mpi B}.
\end{align*}
 The application of \rremss{R.AA+B.A}{R.AA+B.B} yields \(MM^\mpi A=A\) and \(AM^\mpi M=A\).
 Regarding \eqref{F.G.PQ2*}, \eqref{F.G.PQ}, and \(y+x=\ba\), we obtain then
\begin{align*}
 X_{11}
 &=-y^2AM^\mpi MA^\mpi-yxAA^\mpi
 =-y\rk{yAA^\mpi+xAA^\mpi}
 =- y\ba\OPr{A},\\
 X_{12}
 &=-yAM^\mpi A-yAM^\mpi M\OPn{A}M^\mpi B +yAM^\mpi A
 =-yA\OPn{A}M^\mpi B 
 =\Opq,\\
 X_{21}
 &=-y^2xA^\mpi MA^\mpi+y^2 xA^\mpi MA^\mpi+yx\OPn{A}A^\mpi
 =yx\rk{\Iq-A^\mpi A}A^\mpi
 =\Oqp,
\end{align*}
 and
\beql{ab.L0852r.1}\begin{split}
 X_{22}
 &=-yxA^\mpi A-yxA^\mpi M\OPn{A}M^\mpi B\\
 &\qquad-y^2 A^\mpi MM^\mpi A+yxA^\mpi M\OPn{A}M^\mpi B+y\ba\OPn{A}\OPn{M}+y\OPn{A}M^\mpi A-x\OPn{A}^2M^\mpi B\\
 &=-yxA^\mpi A-y^2 A^\mpi A+y\ba\OPn{M}+y\OPn{A}M^\mpi A-x\OPn{A}M^\mpi B\\
 &=-y\ba\OPr{A^\ad}+\rk{\obg-z}\ba\OPn{M}+y\OPn{A}M^\mpi A-x\OPn{A}M^\mpi B.
\end{split}\eeq
 In view of \eqref{F.G.yA-xB}, we have
\[
 yM^\mpi A-xM^\mpi B
 =M^\mpi\rk{yA-xB}
 =\ba M^\mpi\rk{N-zM}.
\]
 By virtue of \(\OPn{M}=\Iq-M^\mpi M\), we thus get
\beql{ab.L0852r.2}
  yM^\mpi A-xM^\mpi B-z\ba\OPn{M}
  =\ba\rk{M^\mpi N-zM^\mpi M-z\Iq+zM^\mpi M}
  =\ba\rk{M^\mpi N-z\Iq}.
\eeq
 From \eqref{ab.L0852r.1} we can infer then
\[\begin{split}
 X_{22}+y\ba\OPr{A^\ad}-\obg\ba\OPn{M}
 &=-z\ba\OPn{M}+y\OPn{A}M^\mpi A-x\OPn{A}M^\mpi B\\
 &=y\OPn{A}M^\mpi A-x\OPn{A}M^\mpi B-z\ba\OPn{A}\OPn{M}
 =\ba\OPn{A}\rk{M^\mpi N-z\Iq},
\end{split}\]
 \tie{}, \(X_{22}=-\ba\ek{y\OPr{A^\ad}-\obg\OPn{M}+\OPn{A}\rk{z\Iq-M^\mpi N}}\).
 Because of \eqref{F.G.PQ}, we have
\begin{gather*}
 M^\mpi B
 =M^\mpi\rk{\ba M-A}
 =\ba M^\mpi M-M^\mpi A
 =\ba\OPr{M^\ad}-M^\mpi A,\\
 \OPr{A^\ad}-\OPn{M}
 =A^\mpi A-\rk{\Iq-M^\mpi M}
 =M^\mpi M-\rk{\Iq-A^\mpi A}
 =\OPr{M^\ad}-\OPn{A},
\end{gather*}
 and, taking additionally into account \(AM^\mpi M=A\), furthermore
\[
 \OPn{A}\OPr{M^\ad}
 =\rk{\Iq-A^\mpi A}M^\mpi M
 =M^\mpi M-A^\mpi A
 =\OPr{M^\ad}-\OPr{A^\ad}.
\]
 From \eqref{ab.L0852r.1} we thus can conclude
\[\begin{split}
 X_{22}
 &=-y\ek*{\ba\rk{\OPr{A^\ad}-\OPn{M}}-\OPn{A}M^\mpi A}-x\OPn{A}\rk{\ba\OPr{M^\ad}-M^\mpi A}\\
 &=-y\ek*{\ba\rk{\OPr{M^\ad}-\OPn{A}}-\OPn{A}M^\mpi A}-x\ek*{\ba\rk{\OPr{M^\ad}-\OPr{A^\ad}}-\OPn{A}M^\mpi A}\\
 &=-\rk{y+x}\rk{\ba\OPr{M^\ad}-\OPn{A}M^\mpi A}+y\ba\OPn{A}+x\ba\OPr{A^\ad}\\
 &=-\ba\rk{\ba\OPr{M^\ad}-M^\mpi A+A^\mpi AM^\mpi A}+\ba\rk{y\OPn{A}+x\OPr{A^\ad}}\\
 &=-\ba\ek*{\rk{M^\mpi B+\OPr{A^\ad}M^\mpi A}-y\OPn{A}-x\OPr{A^\ad}}.
\end{split}\] 
 Comparing the two representations of \(X_{22}\), we can infer then \eqref{ab.L0852r.B}.
\eproof

 The next result is concerned with the matrix polynomial from \eqref{ab.L0852r.B}:

\bleml{F.L.yAzP}
 Let \(A,M\in\Cpq\) with \(\nul{M}\subseteq\nul{A}\) and let \(z\in\Cab\).
 Let \(N\defeq A+\ug M\), let \(y\defeq\obg-z\), and let \(H\defeq y\OPr{A^\ad}-\obg\OPn{M}+\OPn{A}\rk{z\Iq-M^\mpi N}\).
 Then \(\det H\neq0\) and \(A H^\inv=y^\inv A\).
\elem
\bproof
 Consider an arbitrary \(v\in\nul{H}\).
 We have then
\beql{F.L.yAzP.4}
 \obg\OPn{M}v-yA^\mpi Av
 =\rk{\obg\OPn{M}-y\OPr{A^\ad}}v
 =\OPn{A}\rk{z\Iq-M^\mpi N}v
 =\rk{\Iq-A^\mpi A}\rk{z\Iq-M^\mpi N}v.
\eeq
 Regarding \(z\neq\obg\), left multiplication of the latter identity by \(-y^\inv A\) yields \(Av=\Ouu{p}{1}\).
 Consequently, \(Nv=\ug Mv\).
 Taking into account \rrem{R.AA+B.B}, thus \(AM^\mpi Nv=\ug AM^\mpi Mv=\ug Av=\Ouu{p}{1}\).
 From \eqref{F.L.yAzP.4} we then infer
\beql{F.L.yAzP.5}
 \obg\OPn{M}v
 =\rk{\Iq-A^\mpi A}\rk{z\Iq-M^\mpi N}v
 =\rk{z\Iq-M^\mpi N}v
 =zv-\ug M^\mpi Mv.
\eeq
 Left multiplying this by \(M\), we get \(\Ouu{p}{1}=\rk{z-\ug}Mv\).
 Since \(z\neq\ug\), then necessarily \(Mv=\Ouu{p}{1}\).
 Substituting this into \eqref{F.L.yAzP.5} and regarding \(\OPn{M}=\Iq-M^\mpi M\), we obtain \(\obg v=zv\).
 Because of \(z\neq\obg\), hence \(v=\Ouu{p}{1}\) follows.
 Consequently, the linear subspace \(\nul{H}\) is trivial, implying \(\det H\neq0\).
 By virtue of \(AH=yA\OPr{A^\ad}=yAA^\mpi A=yA\) and \(z\neq\obg\), we thus have \(AH^\inv=y^\inv A\).
\eproof

 In generic situations, the \tFABTion{} turns out to be inverse to the \tiFABTion{}:

\bleml{F.L.FAM-11}
 Let \(A\in\CHq\) and let \(M\in\Cggq\) with \(\ran{A}\subseteq\ran{M}\).
 Let \(G\in\RFqab\) with \tiFaaT{A}{M} \(F\).
 Suppose that \(\ran{G(z_0)}\subseteq\ran{A}\) holds true for some \(z_0\in\Cab\).
 Then \(G\) is exactly the \tFaaTv{A}{M}{\(F\)}.
\elem
\bproof
 Consider an arbitrary \(z\in\Cab\).
 Let \(F_1,F_2\colon\Cab\to\Cqq\) be defined by \eqref{F.L.iFAMlft.F1} and \eqref{F.L.iFAMlft.F2}.
 Then \(\tmat{F_1(z)\\ F_2(z)}=\ek{\mFTiuua{A}{M}{z}}\tmatp{G(z)}{\Iq}\).
 Due to \rlem{F.L.iFAMlft}, furthermore \(\det F_2(z)\neq0\) and \(F(z)=\ek{F_1(z)}\ek{F_2(z)}^\inv\).
 Denote by \(H\) the \tFaaTv{A}{M}{\(F\)} and let \(H_1,H_2\colon\Cab\to\Cqq\) be defined by
\begin{align*}
 H_1(w)&\defeq-(\obg-w)AM^\mpi F(w)+A
\shortintertext{and}
 H_2(w)&\defeq-\rk{\obg-w}\rk{w-\ug}A^\mpi F(w)-\rk{\obg-w}A^\mpi M-\OPn{A}.
\end{align*}
 Using \rrem{ab.R1842*}, we infer \(\nul{M}\subseteq\nul{A}\).
 In view of \rlem{F.L.iFAMrn}, we thus can apply \rlem{F.L.FAMlft} to \(F\) to obtain \(\det H_2(z)\neq0\) and \(H(z)=\ek{H_1(z)}\ek{H_2(z)}^\inv\).
 Since by construction \(\tmat{H_1(z)\\ H_2(z)}=\ek{\mFTuua{A}{M}{z}}\tmat{F(z)\\ \Iq}\) holds true, we have
\[
 \bMat
  H_1(z)\\
  H_2(z)
 \eMat
 =\ek*{\mFTuua{A}{M}{z}}
 \bMat
  F_1(z)\\
  F_2(z)
 \eMat\ek*{F_2(z)}^\inv
 =\ek*{\mFTuua{A}{M}{z}}\ek{\mFTiuua{A}{M}{z}}
 \bMat
  G(z)\\
  \Iq
 \eMat\ek*{F_2(z)}^\inv.
\]
 Let \(x\defeq z-\ug\), let \(y\defeq\obg-z\), and let \(N\defeq A+\ug M\).
 Taking into account \rlem{ab.L0852r}, then
 \begin{align*}
  H_1(z)&=-\ba y\OPr{A}\ek*{G(z)}\ek*{F_2(z)}^\inv,&
  H_2(z)&=-\ba\ek*{y\OPr{A^\ad}-\obg\OPn{M}+\OPn{A}\rk{z\Iq-M^\mpi N}}\ek*{F_2(z)}^\inv
 \end{align*}
 follow by comparing both sides of the latter identity.
 According to \rprop{ab.P1648L1409}, we get \(\ran{G(z)}=\ran{G(z_0)}\subseteq\ran{A}\).
 Thus, \(\OPr{A}G(z)=G(z)\).
 Using again \rprop{ab.P1648L1409}, we can conclude \(\ran{G(\ko z)}=\ran{G(z_0)}\subseteq\ran{A}\), implying \(\ek{\ran{A}}^\orth\subseteq\ek{\Ran{G(\ko z)}}^\orth\).
 In view of \rremss{ab.R1842*}{ab.R1652}, we then obtain \(\nul{A}\subseteq\nul{\ek{G(\ko z)}^\ad}=\nul{G(z)}\).
 Due to \rrem{R.AA+B.B}, therefore \(\ek{G(z)}A^\mpi A=G(z)\) holds true.
 Regarding \(\det H_2(z)\neq0\) and \rlem{F.L.yAzP}, we have furthermore
\[
 -\ba A\ek*{F_2(z)}^\inv\ek*{H_2(z)}^\inv
 =A\rk*{-\ba^\inv\ek*{H_2(z)}\ek*{F_2(z)}}^\inv
 =y^\inv A.
\]
 Consequently,
\[
 H(z)
 =\ek*{H_1(z)}\ek*{H_2(z)}^\inv
 =-\ba y\ek*{G(z)}A^\mpi A\ek*{F_2(z)}^\inv\ek*{H_2(z)}^\inv
 =\ek*{G(z)}A^\mpi A
 =G(z).\qedhere
\]
\eproof

 Conversely, we have:

\bleml{F.L.FAM1-1}
 Let \(A\in\CHq\) and let \(M\in\Cggq\) with \(\ran{A}\subseteq\ran{M}\).
 Let \(F\colon\Cab\to\Cqq\) with \tFaaT{A}{M} \(G\) and denote by \(H\) the \tiFaaTv{A}{M}{\(G\)}.
 Let \(z\in\Cab\) be such that \(\ran{F(z)}\subseteq\ran{M}\) and \(\nul{M}\subseteq\nul{F(z)}\) as well as \(\ran{\rk{z-\ug}F(z)+M}=\ran{A}\) and \(\nul{\rk{z-\ug}F(z)+M}=\nul{A}\) are fulfilled.
 Suppose that \(G\) belongs to \(\RFqab\) and that \(\ran{G(z_0)}\subseteq\ran{A}\) holds true for some \(z_0\in\Cab\).
 Then \(H(z)=F(z)\).
\elem
\bproof
 Because of \rrem{ab.R1842*}, we have \(\nul{M}\subseteq\nul{A}\), implying \(AM^\mpi M=A\), by virtue of \rrem{R.AA+B.B}.
 Let \(G_1,G_2\colon\Cab\to\Cqq\) be defined by \eqref{F.L.FAMlft.G1} and \eqref{F.L.FAMlft.G2}.
 The application of \rlem{F.L.FAMlft} yields then \(\det G_2(z)\neq0\) and \(G(z)=\ek{G_1(z)}\ek{G_2(z)}^\inv\).
 Let \(x\defeq z-\ug\) and let \(y\defeq\obg-z\).
 Setting \(X\defeq xF(z)+M\) and \(Y\defeq yF(z)-M\), we get
\[
 G_1(z)
 =-yAM^\mpi F(z)+A
 =-yAM^\mpi F(z)+AM^\mpi M
 =-AM^\mpi Y
\]
 and
\[
 G_2(z)
 =-yxA^\mpi F(z)-yA^\mpi M-\OPn{A}
 =-yA^\mpi X-\OPn{A}.
\]
 Taking into account the assumptions, we have \(X\OPn{A}=\Oqq\) and, in view of \rremss{R.AA+B.A}{R.AA+B.B}, furthermore \(MM^\mpi A=A\), \(AA^\mpi X=X\), and \(XA^\mpi A=X\).
 Let \(B\defeq\ba M-A\) and let \(E_1,E_2\colon\Cab\to\Cqq\) be defined by \eqref{F.L.iFAMhol.E1} and \eqref{F.L.iFAMhol.E2}.
 Since \(A\OPn{A}=\Oqq\) holds obviously true, we can infer then
\[\begin{split}
 &F(z)E_2(z)-E_1(z)\\
 &=F(z)\ek*{-yxA^\mpi G(z)+yM^\mpi A-x\OPn{A}M^\mpi B}-\ek*{yMA^\mpi G(z)+A+M\OPn{A}M^\mpi B}\\
 &=-yx\ek*{F(z)}A^\mpi G(z)+y\ek*{F(z)}M^\mpi A-x\ek*{F(z)}\OPn{A}M^\mpi B-yMA^\mpi G(z)-MM^\mpi A-M\OPn{A}M^\mpi B\\
 &=-yXA^\mpi G(z)+YM^\mpi A-X\OPn{A}M^\mpi B
 =-yXA^\mpi\ek*{G_1(z)}\ek*{G_2(z)}^\inv+YM^\mpi A\ek*{G_2(z)}\ek*{G_2(z)}^\inv\\
 &=\ek*{-yXA^\mpi\rk{-AM^\mpi Y}+YM^\mpi A\rk{-yA^\mpi X-\OPn{A}}}\ek*{G_2(z)}^\inv
 =y\rk{XM^\mpi Y-YM^\mpi X}\ek*{G_2(z)}^\inv.
\end{split}\]
 Regarding \rrem{A.R.AM+B}, we see \(XM^\mpi Y=YM^\mpi X\).
 Consequently, \(F(z)E_2(z)=E_1(z)\) follows.
 Observe that \(H(z)=\ek{E_1(z)}\ek{E_2(z)}^\mpi\).
 Due to \rlem{F.L.iFAMhol}, we have \(\ran{E_2(z)}=\ran{M}\).
 Using \rremss{A.R.A++*}{ab.R1842*}, we thus can conclude \(\nul{\ek{E_2(z)}^\mpi}=\nul{M^\ad}=\nul{M}\).
 In view of \(\nul{M}\subseteq\nul{F(z)}\), hence \(\nul{\ek{E_2(z)}^\mpi}\subseteq\nul{F(z)}\).
 Because of \rremss{R.AA+B.B}{A.R.A++*}, then \(F(z)\ek{E_2(z)}\ek{E_2(z)}^\mpi=F(z)\) holds true.
 Consequently, we obtain
\[
 H(z)
 =\ek*{E_1(z)}\ek*{E_2(z)}^\mpi
 =F(z)\ek*{E_2(z)}\ek*{E_2(z)}^\mpi
 =F(z).\qedhere
\]
\eproof

 In the particular completely degenerate situation \(B=\Opq\) we have \(A=\ba M\), according to \eqref{F.G.BN}.
 Because of \eqref{F.G.PQ}, then the matrix polynomials \(\mFTiuu{A}{M}\) and \(\mFTuu{A}{M}\) from \rnotass{ab.N1546b}{ab.N1546a} essentially coincide with \(\mFTiu{M}\) and \(\mFTu{M}\) introduced in \rnotass{ab.N1246b}{ab.N1246a}, \tresp{:}
 
\breml{ab.R1557}
 If \(M\in\Cpq\), then the equations \(\mFTiuu{\ba M}{M}=\mFTiu{M}\ek{\zdiag{\rk{\ba^\inv\Ip}}{\rk{\ba\Iq}}}\) and \(\mFTuu{\ba M}{M}=\ek{\zdiag{\rk{\ba^\inv\Ip}}{\rk{\ba\Iq}}}\mFTu{M}\) hold true.
\erem
 
  In addition to the matrices \(A\) and \(M\) and the matrices \(B\) and \(N\) built from them via \eqref{F.G.BN}, we now consider the matrix
\(
 D
 \defeq AM^\mpi B
\),
 which, in view of \rrem{ab.L0907}, corresponds to \(\dia{1}\).
 Because of \eqref{F.G.AB}, we have
\[
 D
 =\rk{-\ug M+N}M^\mpi\rk{\obg M-N}
 =-\ug\obg M+\ug\OPr{M}N+\obg N\OPr{M^\ad}-NM^\mpi N
 \]
 in analogy to the second equation in \eqref {F.G.d01}.
 Taking into account \eqref{F.G.BN} and \eqref{ps}, we get furthermore
\beql{F.G.D=ApsB}
 D
 =A\ek*{\frac{1}{\ba}\rk{A+B}}^\mpi B
 =\ba\ek*{A\rk{A+B}^\mpi B}
 =\ba\rk{A\ps B}.
\eeq
  
\bnotal{ab.N1409} 
 Let \(A,M\in\Cpq\) and let \(B\defeq\ba M-A\) and \(D\defeq AM^\mpi B\).
 Then let \(\cFTiuu{A}{M}\colon\C\to\Coo{(p+q)}{(p+q)}\) be defined by
\[
 \cFTiuua{A}{M}{z}
 \defeq
 \begin{pmat}[{|}]
  M\ek{\rk{\obg-z}\OPr{A^\ad}M^\mpi B+\rk{z-\ug}\OPn{A}M^\mpi A}D^\mpi&B\cr\-
  -\rk{\obg-z}\rk{z-\ug}M^\mpi AD^\mpi&\rk{\obg-z}\rk{\ba\OPn{M}+M^\mpi A}\cr
 \end{pmat}.
\]
\enota

\breml{R1130} 
 Let \(A,M\in\Cqq\) and let \(z\in\C\).
 Let \(B\defeq\ba M-A\) and \(D\defeq AM^\mpi B\).
 Then:
\benui
 \il{R1130.a} If \(D=\Oqq\), then
\(
 \cFTiuua{A}{M}{z}
 =
 \smat{
  \Oqq&B\\
  \Oqq&\rk{\obg-z}\rk{\ba\OPn{M}+M^\mpi A}}
\).
 \il{R1130.b} If \(A=\Oqq\), then \(B=\ba M\), \(D=\Oqq\), and
\(
 \cFTiuua{A}{M}{z}
 =\ba\smat{\Oqq& M\\\Oqq&\rk{\obg-z}\OPn{M}}
\).
 \il{R1130.c} If \(B=\Oqq\), then \(A=\ba M\), \(D=\Oqq\) and
\(
 \cFTiuua{A}{M}{z}
 =\ba\smat{\Oqq&\Oqq\\\Oqq&\rk{\obg-z}\Iq}
\).
 \il{R1130.d} If the matrices \(M,A,B,D\) are invertible and \(AM^\inv B=BM^\inv A\), then
\(
 \cFTiuua{A}{M}{z}
 =
 \smat{
  \rk{\obg-z}MA^\inv&B\\
  -\rk{\obg-z}\rk{z-\ug}B^\inv&\rk{\obg-z}M^\inv A}
\).
\eenui
\erem

\bleml{ab.R1443} 
 Let \(A,M\in\CHq\) with \(\ran{A}\subseteq\ran{M}\) and let \(z\in\C\).
 Let \(B\defeq\ba M-A\), \(D\defeq AM^\mpi B\), and \(N\defeq A+\ug M\) and let \(x\defeq z-\ug\) and \(y\defeq\obg-z\).
 Then
 \begin{multline}\label{ab.R1443.B}
  \ek*{\mFTiuua{A}{M}{z}}\ek*{\mFTiua{D}{z}}\\
  =-y\ba
  \begin{pmat}[{|}]
  \ek{\rk{z-\ug-\obg}M+N-M\OPn{A}M^\mpi N}D^\mpi&-M\cr\-
  x\rk{y\Iq+\OPn{A}M^\mpi N}D^\mpi&z\Iq-M^\mpi N-\obg\OPn{M}\cr
 \end{pmat}\\
  =\ek*{\mFTiua{M}{z}}\ek*{\cFTiuua{A}{M}{z}}.
 \end{multline}
\elem
\bproof
 Because of \rrem{A.R.rs+}, we have \(\ran{B}\subseteq\ran{M}\) and \(\ran{N}\subseteq\ran{M}\).
 Consequently, \rrem{R.AA+B.A} shows \(MM^\mpi A=A\), \(MM^\mpi B=B\), and \(MM^\mpi N=N\).
 In view of \(\OPn{A}=\Iq-A^\mpi A\), then
\beql{ab.R1443.5}
 M\OPn{A}M^\mpi N
 =MM^\mpi N-MA^\mpi AM^\mpi N
 =N-MA^\mpi AM^\mpi N
\eeq
 and, by virtue of
\beql{ab.R1443.1}
 \OPn{A}M^\mpi B
 =M^\mpi B-A^\mpi AM^\mpi B
 =M^\mpi B-A^\mpi D,
\eeq
 furthermore
\beql{ab.R1443.3}
 M\OPn{A}M^\mpi B
 =MM^\mpi B-MA^\mpi D
 =B-MA^\mpi D
\eeq
 follow.
 Taking into account \(A+B=\ba M\) and \(\ba>0\), we can infer \(\ran{A}\subseteq\ran{A+B}\). %
 Since \(A\) and \(M\) are \tH{}, \rrem{A.R.kK} shows that \(B\) is \tH{} as well.
 Using \zitaa{MR325642}{\cthm{2.2(b)}{93}}, we can thus conclude \(\rk{A\ps B}^\ad=A\ps B\).
 In view of \eqref{F.G.D=ApsB}, then \(D^\ad=D\) follows.
 Furthermore, we have \(\ran{D}\subseteq\ran{A}\subseteq\ran{M}\).
 Using \rrem{A.R.A++*}, we obtain then \(\ran{D^\mpi}\subseteq\ran{A^\mpi}\subseteq\ran{M^\mpi}\).
 From \rremss{R.AA+B.A}{A.R.A++*}, we can thus conclude 
 \begin{align}\label{ab.R1443.6}
  M^\mpi MD^\mpi&=D^\mpi,&
  M^\mpi MA^\mpi&=A^\mpi,&
 &\text{and}&
  A^\mpi AD^\mpi&=D^\mpi.
 \end{align}
 In particular,
 \[
  A^\mpi AM^\mpi MD^\mpi
  =A^\mpi AD^\mpi
  =D^\mpi
  =M^\mpi MD^\mpi
 \]
 follows.
 In view of \eqref{F.G.AB}, we infer then
\beql{ab.R1443.4}
 A^\mpi DD^\mpi
 =A^\mpi AM^\mpi BD^\mpi
 =A^\mpi AM^\mpi\rk{\obg M-N}D^\mpi
 =\rk{\obg M^\mpi M-A^\mpi AM^\mpi N}D^\mpi.
\eeq
 Taking into account \(\OPr{D}=DD^\mpi\) and \rnotass{ab.N1546b}{ab.N1246b}, the computation of the \tqqa{matrices} in  the \tbr{} \(\ek{\mFTiuua{A}{M}{z}}\ek{\mFTiua{D}{z}}=\tmat{X_{11}&X_{12}\\X_{21}&X_{22}}\) yields
\begin{align*}
 X_{11}&=y^2MA^\mpi\OPr{D}-yx\rk{A+M\OPn{A}M^\mpi B}D^\mpi=y\rk{yMA^\mpi D-xA-xM\OPn{A}M^\mpi B}D^\mpi,\\
 X_{12}&=yMA^\mpi D+y\rk{A+M\OPn{A}M^\mpi B}=y\rk{MA^\mpi D+A+M\OPn{A}M^\mpi B},\\
 X_{21}&=-y^2xA^\mpi\OPr{D}-yx\ek*{y\rk{\ba\OPn{M}+M^\mpi A}-x\OPn{A}M^\mpi B}D^\mpi\\
 &=-yx\rk{yA^\mpi D+y\ba\OPn{M}+yM^\mpi A-x\OPn{A}M^\mpi B}D^\mpi,
\intertext{and}
 X_{22}&=-yxA^\mpi D+y\ek*{y\rk{\ba\OPn{M}+M^\mpi A}-x\OPn{A}M^\mpi B}\\
 &=-y\rk{xA^\mpi D-y\ba\OPn{M}-yM^\mpi A+x\OPn{A}M^\mpi B}.
\end{align*}
 By virtue of \eqref{ab.R1443.3}, \eqref{ab.R1443.4}, and \eqref{ab.R1443.5}, we obtain
\[\begin{split}
 X_{11}
 &=y\rk{yMA^\mpi D-xA-xB+xMA^\mpi D}D^\mpi
 =y\rk{\ba MA^\mpi D-x\ba M}D^\mpi\\
 &=-y\ba\rk{xM-MA^\mpi D}D^\mpi
 =-y\ba\ek*{\rk{z-\ug}M-M\rk{\obg M^\mpi M-A^\mpi AM^\mpi N}}D^\mpi\\
 &=-y\ba\ek*{\rk{z-\ug-\obg}M+MA^\mpi AM^\mpi N}D^\mpi
 =-y\ba\ek*{\rk{z-\ug-\obg}M+N-M\OPn{A}M^\mpi N}D^\mpi.
\end{split}\]
 Because of \eqref{ab.R1443.3}, we have
\[
 X_{12}
 =y\rk{MA^\mpi D+A+B-MA^\mpi D}
 =y\ba M
 =-y\ba\rk{-M}.
\]
 Using \eqref{F.G.yA-xB}, we get \eqref{ab.L0852r.2} by the same reasoning as in the proof of \rlem{ab.L0852r}.
 The combination of \eqref{ab.R1443.1}, \eqref{ab.L0852r.2}, \eqref{ab.R1443.4}, and \eqref{F.G.PQ} yields
\[\begin{split}
 X_{21}
 &=-yx\ek*{yA^\mpi D+\rk{\obg-z}\ba\OPn{M}+yM^\mpi A-xM^\mpi B+xA^\mpi D}D^\mpi\\
 &=-yx\rk{\ba A^\mpi D+\obg\ba\OPn{M}+yM^\mpi A-xM^\mpi B-z\ba\OPn{M}}D^\mpi\\
 &=-yx\ek*{\ba A^\mpi D+\obg\ba\OPn{M}+\ba\rk{M^\mpi N-z\Iq}}D^\mpi\\
 &=-yx\ba\ek*{\obg M^\mpi M-A^\mpi AM^\mpi N+\obg\rk{\Iq-M^\mpi M}+M^\mpi N-z\Iq}D^\mpi\\
 &=-yx\ba\rk{y\Iq+\OPn{A}M^\mpi N}D^\mpi.
\end{split}\]
 Taking into account \eqref{ab.R1443.1} and \eqref{ab.L0852r.2}, we get furthermore
\[\begin{split}
 X_{22}
 &=-y\ek*{xA^\mpi D-\rk{\obg-z}\ba\OPn{M}-yM^\mpi A+xM^\mpi B-xA^\mpi D}\\
 &=-y\rk{z\ba\OPn{M}-yM^\mpi A+xM^\mpi B-\obg\ba\OPn{M}}
 =-y\ba\rk{z\Iq-M^\mpi N-\obg\OPn{M}}.
\end{split}\]
 Hence, the first equation in \eqref{ab.R1443.B} is verified.
 
 Because of \eqref{F.G.xA-yB} and \eqref{ab.R1443.6}, we have
\beql{ab.R1443.7}\begin{split}
 MA^\mpi AM^\mpi\rk{yB-xA}D^\mpi
 &=\ba MA^\mpi AM^\mpi\ek*{\rk{\obg+\ug-z}M-N}D^\mpi\\
 &=\ba\ek*{\rk{\obg+\ug-z}M-MA^\mpi AM^\mpi N}D^\mpi.
\end{split}\eeq 
 Taking into account \(\OPr{M}=MM^\mpi\) and \rnotass{ab.N1246b}{ab.N1409}, the computation of the \tqqa{matrices} in the \tbr{} \(\ek{\mFTiua{M}{z}}\ek{\cFTiuua{A}{M}{z}}=\tmat{Y_{11}&Y_{12}\\Y_{21}&Y_{22}}\) yields
\begin{align*}
 Y_{11}
 &=y\OPr{M} M\rk{y\OPr{A^\ad}M^\mpi B+x\OPn{A}M^\mpi A}D^\mpi-yxMM^\mpi AD^\mpi\\
 &=y\rk{yM\OPr{A^\ad}M^\mpi B+xM\OPn{A}M^\mpi A-xMM^\mpi A}D^\mpi,\\
 Y_{12}
 &=y\OPr{M}B+yM\rk{\ba\OPn{M}+M^\mpi A}
 =y\OPr{M}\rk{B+A}
 =y\OPr{M}\rk{\ba M}
 =y\ba M
 =-y\ba\rk{-M},\\
 Y_{21}
 &=-yxM^\mpi M\rk{y\OPr{A^\ad}M^\mpi B+x\OPn{A}M^\mpi A}D^\mpi-y^2xM^\mpi AD^\mpi\\
 &=-yxM^\mpi\rk{yM\OPr{A^\ad}M^\mpi B+xM\OPn{A}M^\mpi A+yA}D^\mpi,
\intertext{and}
 Y_{22}
 &=-yxM^\mpi B+y^2\rk{\ba\OPn{M}+M^\mpi A}
 =-y\ek*{xM^\mpi B-\rk{\obg-z}\ba\OPn{M}-yM^\mpi A}.
\end{align*}
 In view of \(\OPr{A^\ad}=A^\mpi A\) and \(\OPn{A}=\Iq-A^\mpi A\), we can infer from \eqref{ab.R1443.7} and \eqref{ab.R1443.5} then
\[\begin{split}
 Y_{11}
 &=y\rk{yMA^\mpi AM^\mpi B+xMM^\mpi A-xMA^\mpi AM^\mpi A-xMM^\mpi A}D^\mpi\\
 &=yMA^\mpi AM^\mpi\rk{yB-xA}D^\mpi
 =-y\ba\ek*{\rk{z-\ug-\obg}M+MA^\mpi AM^\mpi N}D^\mpi\\
 &=-y\ba\ek*{\rk{z-\ug-\obg}M+N-M\OPn{A}M^\mpi N}D^\mpi.
\end{split}\]
 Because of \(MM^\mpi A=A\) and the identities \eqref{ab.R1443.7} and \eqref{ab.R1443.6}, we have furthermore
\[\begin{split}
 Y_{21}
 &=-yxM^\mpi\rk{yMA^\mpi AM^\mpi B+xMM^\mpi A-xMA^\mpi AM^\mpi A+yA}D^\mpi\\
 &=-yxM^\mpi\ek*{MA^\mpi AM^\mpi\rk{yB-xA}D^\mpi+xAD^\mpi+yAD^\mpi}\\
 &=-yxM^\mpi\rk*{\ba\ek*{\rk{\obg+\ug-z}M-MA^\mpi AM^\mpi N}D^\mpi+\ba AD^\mpi}\\
 &=-yx\ba\ek*{\rk{y+\ug}M^\mpi MD^\mpi-M^\mpi MA^\mpi AM^\mpi ND^\mpi+M^\mpi AD^\mpi}\\
 &=-yx\ba\rk{yD^\mpi+\ug M^\mpi MD^\mpi-A^\mpi AM^\mpi ND^\mpi+M^\mpi AD^\mpi}\\
 &=-yx\ba\ek*{y\Iq+M^\mpi\rk{\ug M+A}-A^\mpi AM^\mpi N}D^\mpi
 =-yx\ba\rk{y\Iq+\OPn{A}M^\mpi N}D^\mpi.
\end{split}\]
 From \eqref{ab.L0852r.2} moreover \(Y_{22}=-y\ba\rk{z\Iq-M^\mpi N-\obg\OPn{M}}\) follows.
 By virtue of \(Y_{12}=-y\ba\rk{-M}\), thus the second equation in \eqref{ab.R1443.B} is verified.
\eproof

\section{On the elementary steps of the forward algorithm}\label{S1531}
 This section is aimed to work out the elementary step of the forward algorithm by applying the transformations studied in the previous section.
  
\bleml{ab.P1505}
 Let  \(\seqs{0}\in\Fggqu{0}\) and let \(F\in\RFqabsg{0}\).
 Then the \tFatpv{\su{0}}{\(F\)} belongs to \(\PRFabqa{\su{0}}\).
\elem
\bproof
 Denote by \(\copa{G_1}{G_2}\) the \tFatpv{\su{0}}{\(F\)}.
 Obviously, \(\mathcal{D}\defeq\emptyset\) is a discrete subset of \(\Cab\).
 Observe that \(F\) is holomorphic.
 In view of \rdefn{F.D.FTF}, then \(G_1\) and \(G_2\) are holomorphic as well.
 In particular, \(G_1\) and \(G_2\) are \(\Cqq\)\nobreakdash-valued functions, which are meromorphic in \(\Cab\) with \(\pol{G_1}\cup\pol{G_2}\subseteq\mathcal{D}\).
 Consequently, condition~\ref{F.N.PRFab.I} in \rnota{F.N.PRFab} is fulfilled with the set \(\mathcal{D}\) for the pair \(\copa{P}{Q}=\copa{G_1}{G_2}\).
 Consider an arbitrary \(z\in\Cab\).
 By assumption, the \tRabMa{} \(\rabmF\) of \(F\) belongs to \(\MggqFsg{0}\), \tie{}, \(\rabmF(\ab)=\su{0}\).
 Taking additionally into account \rprop{ab.P1648L1409}, hence \(\ran{F(z)}=\ran{\su{0}}\) follows.
 Thus, \(\OPu{\ran{\su{0}}}F(z)=F(z)\).
 Regarding \rnota{ab.N1246a}, we can infer then \(\smat{G_1(z)\\ G_2(z)}=-\ek{\mFTua{\su{0}}{z}}\smat{F(z)\\ \Iq}\).
 Using \rlem{ab.L1403}, we obtain
\[
 \ek*{\mFTiua{\su{0}}{z}}
 \bMat
  G_1(z)\\
  G_2(z)
 \eMat
 =-\ek*{\mFTiua{\su{0}}{z}}\ek*{\mFTua{\su{0}}{z}}
 \bMat
  F(z)\\
  \Iq
 \eMat
 =\rk{\obg-z}\ba
 \bMat
  \OPu{\ran{\su{0}}}F(z)\\
  \Iq
 \eMat
 =\rk{\obg-z}\ba
 \bMat
  F(z)\\
  \Iq
 \eMat.
\]
 In view of  \(z\neq\obg\) and \(\ba>0\), we can conclude \(q\geq\rank\smat{G_1(z)\\ G_2(z)}\geq\rank\smat{F(z)\\ \Iq}=q\), implying \(\rank\smatp{G_1(z)}{G_2(z)}=q\).
 Let \(x\defeq z-\ug\) and let \(y\defeq\obg-z\).
 Furthermore, let \(W_0\defeq\mFTua{\su{0}}{z}\), let \(W_1\defeq\ek{\zdiag{\rk{x\Iq}}{\Iq}}W_0\), and let \(W_2\defeq\ek{\zdiag{\rk{y\Iq}}{\Iq}}W_0\).
 As already mentioned above, we have \(-W_0\smat{F(z)\\ \Iq}=\smat{G_1(z)\\ G_2(z)}\).
 Consequently,
\begin{align*}
 -W_1
 \bMat
  F(z)\\
  \Iq
 \eMat
 &=
 \bMat
  xG_1(z)\\
  G_2(z)
 \eMat&
&\text{and}&
 -W_2
 \bMat
  F(z)\\
  \Iq
 \eMat
 &=
 \bMat
  yG_1(z)\\
  G_2(z)
 \eMat.
\end{align*}
 Taking additionally into account \(\OPu{\ran{\su{0}}}F(z)=F(z)\), the application of \rprop{ab.L1336} yields
\begin{align*}
 \bMat
  xG_1(z)\\
  G_2(z)
 \eMat^\ad\Jimq
 \bMat
  xG_1(z)\\
  G_2(z)
 \eMat
 &=\ba\rk*{
 \bMat
  yxF(z)\\
  \Iq
 \eMat^\ad\Jimq
 \bMat
  yxF(z)\\
  \Iq
 \eMat-2\im\rk{z}\su{0}}
\shortintertext{and}
 \bMat
  yG_1(z)\\
  G_2(z)
 \eMat^\ad\Jimq
 \bMat
  yG_1(z)\\
  G_2(z)
 \eMat
 &=\ba\abs{y}^2\rk*{
 \bMat
  F(z)\\
  \Iq
 \eMat^\ad\Jimq
 \bMat
  F(z)\\
  \Iq
 \eMat-2\im\rk{z}\ek*{F(z)}^\ad\su{0}^\mpi\ek*{F(z)}}.
\end{align*}
 Because of \rrem{ab.R1543}, we have
\begin{align*}
 \bMat
  xG_1(z)\\
  G_2(z)
 \eMat^\ad\Jimq
 \bMat
  xG_1(z)\\
  G_2(z)
 \eMat
 &=2\im\rk*{\ek*{G_2(z)}^\ad\ek*{xG_1(z)}}
 =2\im\rk*{\rk{z-\ug}\ek*{G_2(z)}^\ad\ek*{G_1(z)}},\\
 \bMat
  yG_1(z)\\
  G_2(z)
 \eMat^\ad\Jimq
 \bMat
  yG_1(z)\\
  G_2(z)
 \eMat
 &=2\im\rk*{\ek*{G_2(z)}^\ad\ek*{yG_1(z)}}
 =2\im\rk*{\rk{\obg-z}\ek*{G_2(z)}^\ad\ek*{G_1(z)}},
\end{align*}
 and, furthermore,
\begin{align*}
 \bMat
  yxF(z)\\
  \Iq
 \eMat^\ad\Jimq
 \bMat
  yxF(z)\\
  \Iq
 \eMat
 &=2\im\ek*{yxF(z)}&
&\text{and}&
 \bMat
  F(z)\\
  \Iq
 \eMat^\ad\Jimq
 \bMat
  F(z)\\
  \Iq
 \eMat
 &=2\im\ek*{F(z)}.
\end{align*}
 Now assume in addition \(z\notin\R\).
 Taken all together, we get then
\beql{ab.P1505.1}
 \frac{1}{\im z}\im\rk*{\rk{z-\ug}\ek*{G_2(z)}^\ad\ek*{G_1(z)}}
 =\ba\rk*{\frac{1}{\im z}\im\ek*{yxF(z)}-\su{0}}
\eeq
 and
\beql{ab.P1505.2}
 \frac{1}{\im z}\im\rk*{\rk{\obg-z}\ek*{G_2(z)}^\ad\ek*{G_1(z)}}
 =\ba\abs{y}^2\rk*{\frac{1}{\im z}\im\ek*{F(z)}-\ek*{F(z)}^\ad\su{0}^\mpi\ek*{F(z)}}.
\eeq
 \rrem{ab.L1005} provides us \(\frac{1}{\im z}\im\ek{yxF(z)}\lgeq\su{0}\).
 Taking additionally into account \(\ba>0\) and \rrem{A.R.kK}, we conclude from \eqref{ab.P1505.1} then \(\frac{1}{\im z}\im\rk{\rk{z-\ug}\ek{G_2(z)}^\ad\ek{G_1(z)}}\in\Cggq\).
 Moreover, in view of \rthm{ab.T1614}, we  apply \rlem{ab.L1400} to \(F\) and obtain \(\ek{F(\ko z)}\su{0}^\mpi\ek{F(\ko z)}^\ad\lleq\rk{\im\ko z}^\inv\im\ek{F(\ko z)}\).
 Because of \rrem{ab.R1652}, we have \(\rk{\im\ko z}^\inv\im\ek{F(\ko z)}-\ek{F(\ko z)}\su{0}^\mpi\ek{F(\ko z)}^\ad=\frac{1}{\im z}\im\ek{F(z)}-\ek{F(z)}^\ad\su{0}^\mpi\ek{F(z)}\).
 Taking additionally into account \(\ba>0\) and \rrem{A.R.kK}, we infer from \eqref{ab.P1505.2} then similarly \(\frac{1}{\im z}\im\rk{\rk{\obg-z}\ek{G_2(z)}^\ad\ek{G_1(z)}}\in\Cggq\).
 In view of the choice of \(z\) in \(\Cab=\C\setminus\rk{\ab\cup\mathcal{D}}\), we have thus shown that the conditions~\ref{F.N.PRFab.II}--\ref{F.N.PRFab.IV} in \rnota{F.N.PRFab} are fulfilled with the set \(\mathcal{D}\) for the pair \(\copa{P}{Q}=\copa{G_1}{G_2}\).
 Consequently, \(\copa{G_1}{G_2}\in\PRFabq\).
 Furthermore, \(\OPu{\ran{\su{0}}}F=F\) is verified.
 According to \rdefn{F.D.FTF}, then \(\OPu{\ran{\su{0}}}G_1=G_1\).
 By virtue of \rlem{F.L.A<P}, hence \(\copa{G_1}{G_2}\in\PRFabqa{\su{0}}\) follows.
\eproof

\bleml{ab.P1235}
 Assume \(\kappa\geq1\).
 Let \(\seqska\in\Fggqka\) with \tFT{} \(\seqt{\kappa-1}\) and let \(F\in\RFqabskag\).
 Further, let \(\sau{0}\) be given via \rnota{F.N.sa}.
 Then the \tFaaTv{\sau{0}}{\su{0}}{\(F\)} belongs to \(\RFqabg{\seqt{\kappa-1}}\).
\elem
\bproof
 We first consider the case \(\kappa=\infi\).
 Let \(\rho\defeq\max\set{\abs{\ug},\abs{\obg}}\) and let \(\diskc{\rho}\defeq\setaca{z\in\C}{\abs{z}>\rho}\).
 Obviously, \(\diskc{\rho}\subseteq\Cab\) and \(0\notin\diskc{\rho}\).
 According to \rprop{F.L.sabF}, the sequences \(\seqsainf\) and \(\seqsbinf\) introduced in \rnota{F.N.sa} both belong to \(\Fggqinf\).
 Because of \rrem{ab.L0921}, the matrix-valued functions \(\Ffa\) and \(\Fb\) given in \rnota{ab.N1537} fulfill \(\Ffa\in\RFqabg{\seqsainf}\) and \(\Fb\in\RFqabg{\seqsbinf}\).
 Consequently, \(\Fb\) and \(\Ffa\) are holomorphic in \(\Cab\) and \rprop{F.P.FP8} yields, for all \(z\in\diskc{\rho}\), the series expansions
\begin{align*}
 \Fbv{z}&=-\sum_{j=0}^\infi z^{-(j+1)}\sub{j}&
&\text{and}&
 \Fav{z}&=-\sum_{j=0}^\infi z^{-(j+1)}\sau{j}.
\end{align*}
 In accordance with \rdefn{ab.N1137b}, denote by \(\seqapb{\infi}\) the \tbmodv{\seqsainf}.
 Then \(\apb{0}=-\sau{0}\) and, in view of \rnota{F.N.sa}, furthermore \(\apb{j}=\obg\sau{j-1}-\sau{j}\) for all \(j\in\N\).
 The matrix-valued functions \(Y,Z\colon\Cab\to\Cqq\) defined by \(Y(z)\defeq-z\Fbv{z}\) and \(Z(z)\defeq-\rk{\obg-z}\Fav{z}\), \tresp{}, are both holomorphic in \(\Cab\) with series expansions
\begin{align*}
 Y(z)&=\sum_{n=0}^\infi z^{-n}\sub{n}&
&\text{and}&
 Z(z)&=\sum_{n=0}^\infi z^{-n}\apb{n}
\end{align*}
 for all \(z\in\diskc{\rho}\).
 Let \(R\colon\Cab\to\Cqq\) be defined by \(R(z)\defeq\ek{Z(z)}^\mpi\).
 Observe that the \tRabMa{} \(\rabmFfa\) of \(\Ffa\) fulfills \(\rabmFav{\ab}=\sau{0}\).
 Using \rprop{ab.P1648L1409}, we obtain, for all \(z\in\Cab\), thus \(\ran{Z(z)}=\ran{\Fav{z}}=\ran{\sau{0}}=\ran{\apb{0}}\) and, analogously, \(\nul{Z(z)}=\nul{\apb{0}}\).
 Therefore, we see from \rprop{ab.L0850} that the matrix-valued function \(R\) is holomorphic in \(\Cab\).
 Denote by \(\seq{r_j}{j}{0}{\infi}\) the \trFa{\(\seqapb{\infi}\)}.
 The application of \rlem{ab.L0836} yields the series expansion
\[
 R(z)
 =\sum_{n=0}^\infi z^{-n}r_n
\]
 for all \(z\in\diskc{\rho}\).
 Denote by \(\seq{x_j}{j}{0}{\infi}\) the \tCPa{\(\seqsb{\infi}\)}{\(\seq{r_j}{j}{0}{\infi}\)}.
 From \rlem{ab.L1324} we see then that the function \(X\defeq YR\) is holomorphic in \(\Cab\) with series expansion
\[
 X(z)
 =\sum_{n=0}^\infi z^{-n}x_n
\]
 for all \(z\in\diskc{\rho}\).
 Observe that the \tRabMa{} \(\rabmF\) of \(F\) fulfills \(\rabmFa{\ab}=\su{0}\).
 Denote by \(G\) the \tFaaTv{\sau{0}}{\su{0}}{\(F\)}.
 Regarding \rnota{ab.N1537} and \rdefn{F.D.FTFAB}, we obtain
\[
 G(z)
 =\sau{0}\su{0}^\mpi\ek*{\Fbv{z}}\ek*{\rk{\obg-z}\Fav{z}}^\mpi\sau{0}
 =\sau{0}\su{0}^\mpi\ek*{\Fbv{z}}\ek*{-Z(z)}^\mpi\sau{0}
 =-\sau{0}\su{0}^\mpi\ek*{\Fbv{z}}\ek{R(z)}\sau{0}
\]
 for all \(z\in\Cab\).
 In particular, by virtue of \rremss{B.R.cp}{B.R.Tay*}, thus \(G\) is holomorphic in \(\Cab\).
 Let \(H\colon\Cab\to\Cqq\) be defined by \(H(z)\defeq zG(z)\).
 We have
\[
 H(z)
 =-z\sau{0}\su{0}^\mpi\ek*{\Fbv{z}}\ek{R(z)}\sau{0}
 =\sau{0}\su{0}^\mpi\ek*{Y(z)}\ek*{R(z)}\sau{0}
 =\sau{0}\su{0}^\mpi\ek*{X(z)}\sau{0}
\]
 for all \(z\in\Cab\).
 From \rrem{B.R.Lau*} we see then that the matrix-valued function \(H\) is holomorphic in \(\Cab\) with series expansion
\[
 H(z)
 =\sum_{n=0}^\infi z^{-n}\rk{\sau{0}\su{0}^\mpi x_n\sau{0}}
\]
 for all \(z\in\diskc{\rho}\).
 In view of \rdefn{ab.N0940}, consequently
\[
 G(z)
 =\frac{1}{z}H(z)
 =-\sum_{j=0}^\infi z^{-(j+1)}t_j
\]
 for all \(z\in\diskc{\rho}\) follows.
 Due to \rprop{ab.P1030}, the sequence \(\seqt{\infi}\) belongs to \(\Fggqinf\).
 Since \(G\) is holomorphic in \(\Cab\), the application of \rprop{F.P.R-F} thus yields \(G\in\RFqabg{\seqt{\infi}}\), completing the proof in the case \(\kappa=\infi\).

 Now we consider the case \(\kappa<\infi\).
 Then \(m\defeq\kappa\) belongs to \(\N\).
 Regarding \rrem{F.R.RabM8}, denote by \(\hat s_j\defeq\int_\ab \xi^j\rabmFa{\dif\xi}\) for all \(j\in\NO\) the power moments of the \tRabMa{} \(\rabmF\) of \(F\).
 Then \(\rabmF\in\MggqFag{\hat s}{\infi}\), \tie{}, \(F\in\RFqabg{\seq{\hat s_j}{j}{0}{\infi}}\).
 In particular, by virtue of \rprop{F.P.FPsolv}, we therefore have \(\seq{\hat s_j}{j}{0}{\infi}\in\Fggqinf\).
 Denote by \(\seq{\hat t_j}{j}{0}{\infi}\) the \tFTv{\(\seq{\hat s_j}{j}{0}{\infi}\)}.
 Let \(\hat a_0\defeq-\ug\hat s_0+\hat s_1\) and denote by \(G\) the \tFaaTv{\hat a_0}{\hat s_0}{\(F\)}.
 Since the assertion is already proved for \(\kappa=\infi\), we see that \(G\) belongs to \(\RFqabg{\seq{\hat t_j}{j}{0}{\infi}}\).
 Observe that by assumption \(m\geq1\) and \(\su{j}=\hat s_j\) for all \(j\in\mn{0}{m}\) hold true.
 Hence, we have \(\hat a_0=\sau{0}\) and, because of \rrem{F.R.FTtr}, furthermore \(\hat t_j=t_j\) for all \(j\in\mn{0}{m-1}\).
 Consequently, \(G\) is exactly the \tFaaTv{\sau{0}}{\su{0}}{\(F\)} and belongs to \(\RFqabg{\seqt{m-1}}\).
\eproof

\section{On the elementary steps of the backward algorithm}\label{S1141}
 This section can be considered as the analogue of the preceding one for the backward algorithm.
 More precisely, we will work out the elementary step of the backward algorithm by applying the transformation studied in \rsec{S1230}.
 
\bleml{ab.P1448}
 Let \(\seqs{0}\in\Fggqu{0}\) and let \(\copa{G_1}{G_2}\in\PRFabqa{\su{0}}\).
 Then the \tiFaTv{\su{0}}{\(\copa{G_1}{G_2}\)} belongs to \(\RFqabsg{0}\).
\elem
\bproof
 According to \rnota{ab.N1534}, we have \(\copa{G_1}{G_2}\in\PRFabq\).
 In particular, \(G_1\) and \(G_2\) are \(\Cqq\)\nobreakdash-valued functions, which are meromorphic in \(\Cab\).
 Using the functions \(g,h\colon\Cab\to\C\) given via \eqref{F.D.iFTF.1}, we define by \eqref{F.D.iFTF.2} two \(\Cqq\)\nobreakdash-valued functions \(F_1\) and \(F_2\) meromorphic in \(\Cab\).
 From \rrem{F.R.Fgg-s} we get \(\su{0}\in\Cggq\).
 By virtue of \rprop{F.P.Filft}, thus \(\det F_2\) does not vanish identically in \(\Cab\) and the \tiFaT{\su{0}} \(F\) of \(\copa{G_1}{G_2}\) admits the representation \(F=F_1F_2^\inv\), according to \rdefn{F.D.iFTF}.
 
 In a first step, we are now going to show that the pair \(\copa{F_1}{F_2}\) belongs to \(\PRFabq\).
 Due to \rprop{ab.P1528}, the set \(\exset\defeq\pol{G_1}\cup\pol{G_2}\cup\PRFabex{G_1}{G_2}\) is a discrete subset \(\mathcal{D}\) of \(\Cab\), satisfying the conditions~\ref{F.N.PRFab.I}--\ref{F.N.PRFab.IV} in \rnota{F.N.PRFab} for the pair \(\copa{P}{Q}=\copa{G_1}{G_2}\).
 In view of \eqref{F.D.iFTF.2}, we have \(\pol{F_1}\subseteq\pol{G_1}\cup\pol{G_2}\) and \(\pol{F_2}\subseteq\pol{G_1}\cup\pol{G_2}\).
 Consequently, \(\pol{F_1}\cup\pol{F_2}\subseteq\exset\), \tie{}, condition~\ref{F.N.PRFab.I} in \rnota{F.N.PRFab} is fulfilled with the set \(\mathcal{D}=\exset\) for the pair \(\copa{P}{Q}=\copa{F_1}{F_2}\).
 Consider an arbitrary \(z\in\C\setminus\rk{\ab\cup\exset}\).
 Regarding \(\su{0}\in\Cggq\), \rprop{F.P.Filft} yields \(\det F_2(z)\neq0\) and
\beql{ab.P1448.0}
 F(z)
 =\ek*{F_1(z)}\ek*{F_2(z)}^\inv.
\eeq
 In particular, \(\rank\smatp{F_1(z)}{F_2(z)}=q\).
 Let \(x\defeq z-\ug\) and let \(y\defeq\obg-z\).
 Furthermore, let \(V_0\defeq\mFTiua{\su{0}}{z}\), let \(V_1\defeq\ek{\zdiag{\rk{x\Iq}}{\Iq}}V_0\), and let \(V_2\defeq\ek{\zdiag{\rk{y\Iq}}{\Iq}}V_0\).
 According to \rnota{ab.N1246b}, we have \(V_0\smatp{G_1(z)}{G_2(z)}=\smatp{F_1(z)}{F_2(z)}\).
 Hence,
\begin{align*}
 V_1
 \bMat
  G_1(z)\\
  G_2(z)
 \eMat
 &=
 \bMat
  xF_1(z)\\
  F_2(z)
 \eMat&
&\text{and}&
 V_2
 \bMat
  G_1(z)\\
  G_2(z)
 \eMat
 &=
 \bMat
  yF_1(z)\\
  F_2(z)
 \eMat.
\end{align*}
 From \rlem{F.L.A<P} we see \(\OPu{\ran{\su{0}}}G_1=G_1\).
 Using \rprop{ab.L1059}, we thus can infer
\begin{multline*}
 \bMat
  xF_1(z)\\
  F_2(z)
 \eMat^\ad\Jimq
 \bMat
  xF_1(z)\\
  F_2(z)
 \eMat
 =\abs{y}^2
 \bMat
  xG_1(z)\\
  G_2(z)
 \eMat^\ad\Jimq
 \bMat
  xG_1(z)\\
  G_2(z)
 \eMat+\abs{x}^2
 \bMat
  yG_1(z)\\
  G_2(z)
 \eMat^\ad\Jimq
 \bMat
  yG_1(z)\\
  G_2(z)
 \eMat\\
 +2\ba\im\rk{z}\ek*{G_2(z)}^\ad\su{0}\ek*{G_2(z)}
\end{multline*}
 and
\[
 \bMat
  yF_1(z)\\
  F_2(z)
 \eMat^\ad\Jimq
 \bMat
  yF_1(z)\\
  F_2(z)
 \eMat
 =\ba\abs{y}^2\rk*{
 \bMat
  G_1(z)\\
  G_2(z)
 \eMat^\ad\Jimq
 \bMat
  G_1(z)\\
  G_2(z)
 \eMat+2\im\rk{z}\ek*{G_1(z)}^\ad\su{0}^\mpi\ek*{G_1(z)}}.
\]
 Because of \rrem{ab.R1543}, we see that
\begin{align*}
 \bMat
  \xi F_1(z)\\
  F_2(z)
 \eMat^\ad\Jimq
 \bMat
  \xi F_1(z)\\
  F_2(z)
 \eMat
 &=2\im\rk*{\ek*{F_2(z)}^\ad\ek*{\xi F_1(z)}}
 =2\im\rk*{\xi\ek*{F_2(z)}^\ad\ek*{F_1(z)}}
\shortintertext{and}
 \bMat
  \xi G_1(z)\\
  G_2(z)
 \eMat^\ad\Jimq
 \bMat
  \xi G_1(z)\\
  G_2(z)
 \eMat
 &=2\im\rk*{\ek*{G_2(z)}^\ad\ek*{\xi G_1(z)}}
 =2\im\rk*{\xi\ek*{G_2(z)}^\ad\ek*{G_1(z)}}
\end{align*}
 hold true for each \(\xi\in\set{x,y}\) and that
\[
 \bMat
  G_1(z)\\
  G_2(z)
 \eMat^\ad\Jimq
 \bMat
  G_1(z)\\
  G_2(z)
 \eMat
 =2\im\rk*{\ek*{G_2(z)}^\ad\ek*{G_1(z)}}
\]
 is valid.
 Now assume in addition \(z\notin\R\).
 Taken all together, we obtain then
\beql{ab.P1448.1}\begin{split}
 \frac{1}{\im z}\im\rk*{x\ek*{F_2(z)}^\ad\ek*{F_1(z)}}
 &=\abs{y}^2\ek*{\frac{1}{\im z}\im\rk*{x\ek*{G_2(z)}^\ad\ek*{G_1(z)}}}\\
 &\qquad+\abs{x}^2\ek*{\frac{1}{\im z}\im\rk*{y\ek*{G_2(z)}^\ad\ek*{G_1(z)}}}+\ba\ek*{G_2(z)}^\ad\su{0}\ek*{G_2(z)}
\end{split}\eeq
 and
\beql{ab.P1448.2}
 \frac{1}{\im z}\im\rk*{y\ek*{F_2(z)}^\ad\ek*{F_1(z)}}
 =\ba\abs{y}^2\ek*{\frac{1}{\im z}\im\rk*{\ek*{G_2(z)}^\ad\ek*{G_1(z)}}+\ek*{G_1(z)}^\ad\su{0}^\mpi\ek*{G_1(z)}}.
\eeq
 Since the conditions~\ref{F.N.PRFab.III} and \ref{F.N.PRFab.IV} in \rnota{F.N.PRFab} are satisfied with the set \(\mathcal{D}=\exset\) for the pair \(\copa{P}{Q}=\copa{G_1}{G_2}\), we have
\begin{align}\label{ab.P1448.9}
 \frac{1}{\im z}\im\rk*{x\ek*{G_2(z)}^\ad\ek*{G_1(z)}}&\in\Cggq&
&\text{and}&
 \frac{1}{\im z}\im\rk*{y\ek*{G_2(z)}^\ad\ek*{G_1(z)}}&\in\Cggq.
\end{align}
 Taking into account \eqref{ab.P1448.9}, \(\ba>0\), and \(\su{0}\in\Cggq\), we use \rremss{A.R.kK}{A.R.XAX} to infer from \eqref{ab.P1448.1} that \(\frac{1}{\im z}\im\rk{x\ek{F_2(z)}^\ad\ek{F_1(z)}}\in\Cggq\).
 By virtue of \rlem{ab.L1252}, the matrix \(\frac{1}{\im z}\im\rk{\ek{G_2(z)}^\ad\ek{G_1(z)}}\) is \tnnH{}.
 Because of \rrem{A.R.A+>}, we have \(\su{0}^\mpi\in\Cggq\).
 Regarding additionally \(\ba>0\), from \rremss{A.R.kK}{A.R.XAX} and \eqref{ab.P1448.2} we conclude similarly \(\frac{1}{\im z}\im\rk*{y\ek*{F_2(z)}^\ad\ek*{F_1(z)}}\in\Cggq\).
 In view of the choice of \(z\in\C\setminus\rk{\ab\cup\exset}\), we thus have verified that conditions~\ref{F.N.PRFab.II}--\ref{F.N.PRFab.IV} in \rnota{F.N.PRFab} are fulfilled with the set \(\mathcal{D}=\exset\) for the pair \(\copa{P}{Q}=\copa{F_1}{F_2}\).
 Consequently, \(\copa{F_1}{F_2}\in\PRFabq\).
 Therefore, the application of \rlem{ab.P1801} to the pair \(\copa{P}{Q}=\copa{F_1}{F_2}\) yields
\beql{ab.P1448.5}
 F
 =F_1F_2^\inv
 \in\RFqab.
\eeq
 In a second step, we are now going to show that \(\rabmFa{\ab}=\su{0}\).
 First we verify
\beql{ab.P1448.3}
 \rabmFA{\ab}
 \lgeq\su{0}.
\eeq
 Let \(W_0\defeq\mFTua{\su{0}}{z}\).
 Because of \rlem{ab.L1403} and \(\OPu{\ran{\su{0}}}G_1=G_1\), we have
\beql{ab.P1448.7}
 W_0
 \bMat
  F_1(z)\\
  F_2(z)
 \eMat
 =W_0V_0
 \bMat
  G_1(z)\\
  G_2(z)
 \eMat
 =-y\ba
 \bMat
  \OPu{\ran{\su{0}}}G_1(z)\\
  G_2(z)
 \eMat
 =-y\ba
 \bMat
  G_1(z)\\
  G_2(z)
 \eMat.
\eeq
 Setting \(W_1\defeq\ek{\zdiag{\rk{x\Iq}}{\Iq}}W_0\), hence \(
  W_1
  \smat{
   F_1(z)\\
   F_2(z)
  }
  =-y\ba
  \smat{
   xG_1(z)\\
   G_2(z)
  }\) follows.
 In view of \eqref{F.D.iFTF.2}, we have \(\OPu{\ran{\su{0}}}F_1=F_1\).
 Thus, using \rprop{ab.L1336}, we can conclude then
\[
 \abs{y}^2\ba^2
 \bMat
  x G_1(z)\\
  G_2(z)
 \eMat^\ad\Jimq
 \bMat
  x G_1(z)\\
  G_2(z)
 \eMat
 =\ba\rk*{
 \bMat
  yxF_1(z)\\
  F_2(z)
 \eMat^\ad\Jimq
 \bMat
  yxF_1(z)\\
  F_2(z)
 \eMat-2\im\rk{z}\ek*{F_2(z)}^\ad\su{0}\ek*{F_2(z)}}.
\]
 \rrem{ab.R1543} yields \(
   \smat{
    xG_1(z)\\
    G_2(z)
   }^\ad\Jimq
   \smat{
    xG_1(z)\\
    G_2(z)
   }
   =2\im\rk{x\ek{G_2(z)}^\ad\ek{G_1(z)}}\).
 Consequently,
\[
 \abs{y}^2\ba\ek*{\frac{1}{\im z}\im\rk*{x\ek*{G_2(z)}^\ad\ek*{G_1(z)}}}
 =\frac{1}{2\im z}
 \bMat
  yxF_1(z)\\
  F_2(z)
 \eMat^\ad\Jimq
 \bMat
  yxF_1(z)\\
  F_2(z)
 \eMat-\ek*{F_2(z)}^\ad\su{0}\ek*{F_2(z)}.
\]
 Taking into account \(\ba>0\) and \eqref{ab.P1448.9}, we see from \rrem{A.R.kK} that the matrix on the left-hand side of the last equation is \tnnH{}.
 Since \(\Jimq\) and \(\su{0}\) are \tH{} matrices, we can then use \rremss{A.R.kK}{A.R.XAX} to conclude
\[
 \ek*{F_2(z)}^\ad\su{0}\ek*{F_2(z)}
 \lleq\frac{1}{2\im z}
 \bMat
  yxF_1(z)\\
  F_2(z)
 \eMat^\ad\Jimq
 \bMat
  yxF_1(z)\\
  F_2(z)
 \eMat.
\]
 In view of \eqref{ab.P1448.0} and \rremss{A.R.XAX}{ab.R1543}, hence
\[
 \su{0}
 \lleq\frac{1}{2\im z}
 \bMat
  yxF(z)\\
  \Iq
 \eMat^\ad\Jimq
 \bMat
  yxF(z)\\
  \Iq
 \eMat
 =\frac{1}{\im z}\im\ek*{yxF(z)}
 =\frac{1}{\im z}\im\ek*{\rk{\obg-z}\rk{z-\ug}F(z)}
\]
 follows.
 Since \(\exset\) is a discrete set, we in particular infer
 \beql{ab.P1448.4}
  \lim_{\eta\to\infp}\frac{1}{\eta}\im\ek*{\rk{\obg-\iu\eta}\rk{\iu\eta-\ug}F(\iu\eta)}
  \lgeq\su{0}.
 \eeq
 For all \(\eta>0\), we have \(\eta^{-2}\rk{\obg-\iu\eta}\rk{\iu\eta-\ug}=\rk{\obg\eta^\inv-\iu}\rk{\iu-\ug\eta^\inv}\) and therefore \(\lim_{\eta\to\infp}\eta^{-2}{\rk{\obg-\iu\eta}}{\rk{\iu\eta-\ug}}=1\).
 Regarding \eqref{ab.P1448.5} and \rthm{ab.T1614}, we can apply \rlem{Mae05.P2-1-8} to \(F\) and obtain \(\lim_{\eta\to\infp}\iu\eta F(\iu \eta)=-\rabmFa{\ab}\).
 Consequently,
\[
 -\rabmFa{\ab}
 =\ek*{\lim_{\eta\to\infp}\frac{1}{\eta^2}\rk{\obg-\iu\eta}\rk{\iu\eta-\ug}}\ek*{\lim_{\eta\to\infp}\iu\eta F(\iu \eta)}
 =\lim_{\eta\to\infp}\frac{\iu}{\eta}\ek*{\rk{\obg-\iu\eta}\rk{\iu\eta-\ug}F(\iu\eta)}.
\]
 Because of \(\rabmFa{\ab}\in\CHq\), we obtain from \rrem{A.R.XRIX} then
\[
 \rabmFa{\ab}
 =-\re\rk*{\lim_{\eta\to\infp}\frac{\iu}{\eta}\ek*{\rk{\obg-\iu\eta}\rk{\iu\eta-\ug}F(\iu\eta)}}
 =\lim_{\eta\to\infp}\frac{1}{\eta}\im\ek*{\rk{\obg-\iu\eta}\rk{\iu\eta-\ug}F(\iu\eta)}.
\]
 In combination with \eqref{ab.P1448.4}, this implies \eqref{ab.P1448.3}.
 
 Conversely, we now verify \(\rabmFa{\ab}\lleq\su{0}\).
 Because of \eqref{ab.P1448.5} and \rprop{ab.P1648L1409}, we have \(\ran{F(z)}=\ran{\rabmFa{\ab}}\).
 Using \rrem{A.R.rnLAR} and \eqref{ab.P1448.0}, we infer \(\ran{F(z)}=\ran{F_1(z)}\).
 In view of \(\OPu{\ran{\su{0}}}F_1=F_1\), we have \(\ran{F_1(z)}\subseteq\ran{\su{0}}\).
 Consequently, \(\ran{\rabmFa{\ab}}\subseteq\ran{\su{0}}\).
 Let \(W_2\defeq\ek{\zdiag{\rk{y\Iq}}{\Iq}}W_0\).
 In view of \eqref{ab.P1448.7}, then \(
  W_2
  \tmat{
   F_1(z)\\
   F_2(z)
  }
  =-y\ba
  \tmat{
   y G_1(z)\\
   G_2(z)
  }\).
 Taking additionally into account \(\OPu{\ran{\su{0}}}F_1=F_1\), from \rprop{ab.L1336} we get
\[
 \abs{y}^2\ba^2
 \bMat
  yG_1(z)\\
  G_2(z)
 \eMat^\ad\Jimq
 \bMat
  yG_1(z)\\
  G_2(z)
 \eMat
 =\ba\abs{y}^2\rk*{
 \bMat
  F_1(z)\\
  F_2(z)
 \eMat^\ad\Jimq
 \bMat
  F_1(z)\\
  F_2(z)
 \eMat-2\im\rk{z}\ek*{F_1(z)}^\ad\su{0}^\mpi\ek*{F_1(z)}}.
\]
 \rrem{ab.R1543} yields\(
   \smat{
    yG_1(z)\\
    G_2(z)
   }^\ad\Jimq
   \smat{
    yG_1(z)\\
    G_2(z)
   }
   =2\im\rk{y\ek{G_2(z)}^\ad\ek{G_1(z)}}\).
 Consequently,
\[
 \ba\ek*{\frac{1}{\im z}\im\rk*{y\ek*{G_2(z)}^\ad\ek*{G_1(z)}}}
 =\frac{1}{2\im z}
 \bMat
  F_1(z)\\
  F_2(z)
 \eMat^\ad\Jimq
 \bMat
  F_1(z)\\
  F_2(z)
 \eMat-\ek*{F_1(z)}^\ad\su{0}^\mpi\ek*{F_1(z)}.
\]
 Regarding \(\ba>0\) and \eqref{ab.P1448.9}, we see from \rrem{A.R.kK} that the matrix on the left-hand side of the last equation is \tnnH{}.
 Since \(\su{0}\) is \tH{}, \rrem{A.R.A++*} shows that \(\su{0}^\mpi\) is \tH{} as well.
 Taking additionally into account that \(\Jimq\) is \tH{}, we can use \rremss{A.R.kK}{A.R.XAX} to conclude
\[
 \ek*{F_1(z)}^\ad\su{0}^\mpi\ek*{F_1(z)}
 \lleq\frac{1}{2\im z}
 \bMat
  F_1(z)\\
  F_2(z)
 \eMat^\ad\Jimq
 \bMat
  F_1(z)\\
  F_2(z)
 \eMat.
\]
 Because of \eqref{ab.P1448.0}, the application of \rremss{A.R.XAX}{ab.R1543} thus yields
\[
 \ek*{F(z)}^\ad\su{0}^\mpi\ek*{F(z)}
 \lleq\frac{1}{2\im z}
 \bMat
  F(z)\\
  \Iq
 \eMat^\ad\Jimq
 \bMat
  F(z)\\
  \Iq
 \eMat
 =\frac{1}{\im z}\im F(z).
\]
 Regarding that the set \(\exset\) is discrete and that, according to \eqref{ab.P1448.5}, the function \(F\) is holomorphic, we can apply a continuity argument to show that \(\ek{F(\zeta)}^\ad\su{0}^\mpi\ek{F(\zeta)}\lleq\rk{\im \zeta}^\inv\im F(\zeta)\) holds true for all \(\zeta\in\uhp\).
 In view of \eqref{ab.P1448.5} and \rthm{ab.T1614}, we can apply \rlem{B.L.s<B} to \(F\).
 Taking additionally into account \(\su{0}\in\Cggq\) and \(\ran{\rabmFa{\ab}}\subseteq\ran{\su{0}}\), then \(\rabmFa{\ab}\lleq\su{0}\) follows by \rlem{B.L.s<B}.
 In combination with \eqref{ab.P1448.3}, this implies \(\rabmFa{\ab}=\su{0}\).
 Because of \eqref{ab.P1448.5}, thus \(F\) belongs to \(\RFqabsg{0}\).
\eproof

\bleml{ab.P1555}
 Assume \(\kappa\geq1\).
 Let \(\seqska\in\Fggqka\) with \tFT{} \(\seqt{\kappa-1}\) and let \(G\in\RFqabg{\seqt{\kappa-1}}\).
 Then the \tiFaaTv{\sau{0}}{\su{0}}{\(G\)} belongs to \(\RFqabskag\).
\elem
\bproof
 Because of \rprop{ab.P1030}, we have \(\seqt{\kappa-1}\in\Fggqu{\kappa-1}\).
 First we consider the case \(\kappa=\infi\).
 By virtue of \rprop{F.P.FP8}, the set \(\RFqabsg{\infi}\) consists of exactly one element, say \(F\).
 According to \rlem{ab.P1235}, the \tFaaTv{\sau{0}}{\su{0}}{\(F\)} belongs to \(\RFqabg{\seqt{\infi}}\).
 Since, due to \rprop{F.P.FP8}, the set \(\RFqabg{\seqt{\infi}}\) consists of exactly one element, we conclude that \(G\) coincides with the \tFaaTv{\sau{0}}{\su{0}}{\(F\)}.
 Using \rrem{F.R.Fgg-s}, we easily infer \(\sau{0}\in\CHq\) and \(\su{0}\in\Cggq\).
 Due to \rrem{F.R.Fgg-r}, we have \(\ran{\sau{0}}\subseteq\ran{\su{0}}\).
 Consider now an arbitrary \(z\in\Cab\).
 Observe that the \tRabMa{} \(\rabmF\) of \(F\) fulfills \(\rabmFa{\ab}=\su{0}\).
 Taking into account \rprop{ab.P1648L1409}, we obtain thus \(\ran{F(z)}=\ran{\su{0}}\) and \(\nul{F(z)}=\nul{\su{0}}\).
 Because of \rrem{ab.L0921}, the function \(\Ffa\) given in \rnota{ab.N1537} belongs to \(\RFqabg{\seqsainf}\).
 Hence, we analogously get \(\ran{\Fav{z}}=\ran{\sau{0}}\) and \(\nul{\Fav{z}}=\nul{\sau{0}}\).
 In a similar way, we can conclude \(\ran{G(z)}=\ran{t_0}\).
 Regarding \rdefn{ab.N0940}, we see furthermore \(\ran{t_0}\subseteq\ran{\sau{0}}\).
 Consequently, \(\ran{G(z)}\subseteq\ran{\sau{0}}\).
 Thus, we can apply \rlem{F.L.FAM1-1} to the function \(F\) and its \tFaaT{\sau{0}}{\su{0}} \(G\) and obtain with the \tiFaaT{\sau{0}}{\su{0}} \(H\) of \(G\) then \(H(z)=F(z)\).
 Hence, \(H=F\), implying \(H\in\RFqabsg{\infi}\).
 
 Now we consider the case \(\kappa<\infi\).
 Then \(m\defeq\kappa\) belongs to \(\N\).
 Regarding \rrem{F.R.RabM8}, denote by \(\hat t_j\defeq\int_\ab \xi^j\rabmGa{\dif\xi}\) for all \(j\in\NO\) the power moments of the \tRabMa{} \(\rabmG\) of \(G\).
 Then \(\rabmG\in\MggqFag{\hat t}{\infi}\), \tie{}, \(G\in\RFqabg{\seq{\hat t_j}{j}{0}{\infi}}\).
 By virtue of \rprop{F.P.FPsolv}, we have in particular \(\seq{\hat t_j}{j}{0}{\infi}\in\Fggqinf\).
 We are now going to construct the \tfpf{} of a sequence from \(\Fggqinf\) with \tFT{} \(\seq{\hat t_j}{j}{0}{\infi}\):
 Due to \rthm{F.T.FggFP}, the \tfpf{} \(\seq{\hat{\mathfrak{g}}_j}{j}{0}{\infi}\) of \(\seq{\hat t_j}{j}{0}{\infi}\) belongs to the class \(\csqinfd\) introduced in \rnota{F.N.cs}.
 Since \(G\) belongs to \(\RFqabg{\seqt{m-1}}\), we have \(\hat t_j=t_j\) for all \(j\in\mn{0}{m-1}\).
 Denote by \(\fgpseq{2(m-1)}\) the \tfpfa{\(\seqt{m-1}\)}.
 Because of \rrem{F.R.fpftr}, then \(\hat{\mathfrak{g}}_j=\fgpu{j}\) for all \(j\in\mn{0}{2(m-1)}\).
 Denote by \(\fpseq{2m}\) the \tfpfa{\(\seqs{m}\)}.
 Regarding \(\ba>0\), let the sequence \(\seq{\hat{\mathfrak{f}}_j}{j}{0}{\infi}\) be given by
\[
 \hat{\mathfrak{f}}_j
 \defeq
 \begin{cases}
  \fpu{j}\tincase{j\leq2m}\\
  \ba^\inv\hat{\mathfrak{g}}_{j-2}\tincase{j\geq2m+1}
 \end{cases}.
\]
 \rthm{F.T.FggFP} yields \(\fpseq{2m}\in\csqud{m}\).
 Taking additionally into account \(\ba>0\) and \(m\geq1\), we can conclude then that the sequence \(\seq{\hat{\mathfrak{f}}_j}{j}{0}{\infi}\) is a sequence of \tnnH{} matrices fulfilling the relations
\(
 \ba\hat{\mathfrak{f}}_0
 =\ba\fpu{0}
 =\fpu{1}+\fpu{2}
 =\hat{\mathfrak{f}}_1+\hat{\mathfrak{f}}_2
\)
 and
\[
 \ba\rk{\hat{\mathfrak{f}}_{2k-1}\ps\hat{\mathfrak{f}}_{2k}}
 =\ba\rk{\fpu{2k-1}\ps\fpu{2k}}
 =\fpu{2k+1}+\fpu{2k+2}
 =\hat{\mathfrak{f}}_{2k+1}+\hat{\mathfrak{f}}_{2k+2}
\]
 for all \(k\in\N\) with \(k\leq m-1\), and, regarding \rrem{ab.R1132}, furthermore
\[
 \ba\rk{\hat{\mathfrak{f}}_{2k-1}\ps\hat{\mathfrak{f}}_{2k}}
 =\ba\ek*{\rk{\ba^\inv\hat{\mathfrak{g}}_{2k-3}}\ps\rk{\ba^\inv\hat{\mathfrak{g}}_{2k-2}}}
 =\hat{\mathfrak{g}}_{2k-3}\ps\hat{\mathfrak{g}}_{2k-2}
 =\ba^\inv\rk{\hat{\mathfrak{g}}_{2k-1}+\hat{\mathfrak{g}}_{2k}}
 =\hat{\mathfrak{f}}_{2k+1}+\hat{\mathfrak{f}}_{2k+2}
\]
 for all \(k\in\N\) with \(k\geq m+1\).
 In view of \rcor{F.C.Salg1}, in addition \(\fgpu{0}=\ba\rk{\fpu{1}\ps\fpu{2}}\) and \(\fgpu{j}=\ba\fpu{j+2}\) hold true for all \(j\in\mn{1}{2(m-1)}\).
 In the case \(m=1\), we therefore have
\[
 \ba\rk{\fpu{2m-1}\ps\fpu{2m}}
 =\ba\rk{\fpu{1}\ps\fpu{2}}
 =\fgpu{0}
 =\hat{\mathfrak{g}}_{0}
 =\ba^\inv\rk{\hat{\mathfrak{g}}_{1}+\hat{\mathfrak{g}}_{2}}
 =\ba^\inv\rk{\hat{\mathfrak{g}}_{2m-1}+\hat{\mathfrak{g}}_{2m}},
\]
 whereas, because of \rrem{ab.R1132}, in the case \(m\geq2\) then
\[
 \ba\rk{\fpu{2m-1}\ps\fpu{2m}}
 =\rk{\ba\fpu{2m-1}}\ps\rk{\ba\fpu{2m}}
 =\fgpu{2m-3}\ps\fgpu{2m-2}
 =\hat{\mathfrak{g}}_{2m-3}\ps\hat{\mathfrak{g}}_{2m-2}
 =\ba^\inv\rk{\hat{\mathfrak{g}}_{2m-1}+\hat{\mathfrak{g}}_{2m}}
\]
 follows.
 Consequently, \(\ba\rk{\hat{\mathfrak{f}}_{2m-1}\ps\hat{\mathfrak{f}}_{2m}}= \hat{\mathfrak{f}}_{2m+1}+\hat{\mathfrak{f}}_{2m+2}\).
 Hence, the sequence \(\seq{\hat{\mathfrak{f}}_j}{j}{0}{\infi}\) belongs to \(\csqinfd\).
 According to \rthm{F.T.FggFP}, then there exists a sequence \(\seq{\hat s_j}{j}{0}{\infi}\) from \(\Fggqinf\) with \tfpf{} \(\seq{\hat{\mathfrak{f}}_j}{j}{0}{\infi}\).
 Denote by \(\seq{\tilde t_j}{j}{0}{\infi}\) the \tFTv{\(\seq{\hat s_j}{j}{0}{\infi}\)} and by \(\seq{\tilde{\mathfrak{g}}_j}{j}{0}{\infi}\) the \tfpfa{\(\seq{\tilde t_j}{j}{0}{\infi}\)}.
 Using \rcor{F.C.Salg1}, we infer
\begin{align*}
 \tilde{\mathfrak{g}}_{0}
 &=\ba\rk{\hat{\mathfrak{f}}_{1}\ps\hat{\mathfrak{f}}_{2}}
 =\ba\rk{\fpu{1}\ps\fpu{2}}
 =\fgpu{0}
 =\hat{\mathfrak{g}}_{0}&
&\text{and}&
 \tilde{\mathfrak{g}}_{j}
 &=\ba\hat{\mathfrak{f}}_{j+2}
 =\ba\fpu{j+2}
 =\fgpu{j}
 =\hat{\mathfrak{g}}_{j}
\end{align*}
 for all \(j\in\N\) with \(j\leq 2m-2\) and, furthermore
\(
 \tilde{\mathfrak{g}}_{j}
 =\ba\hat{\mathfrak{f}}_{j+2}
 =\hat{\mathfrak{g}}_{j}
\)
 for all \(j\in\N\) with \(j\geq 2m-1\).
 Consequently, the \tfpf{} \(\seq{\hat{\mathfrak{g}}_j}{j}{0}{\infi}\) of \(\seq{\hat t_j}{j}{0}{\infi}\) coincides with the \tfpf{} \(\seq{\tilde{\mathfrak{g}}_j}{j}{0}{\infi}\) of \(\seq{\tilde t_j}{j}{0}{\infi}\).
 Since \rprop{ab.P1030} yields \(\seq{\tilde t_j}{j}{0}{\infi}\in\Fggqinf\), we can conclude from \rthm{F.T.FggFP} that the sequences \(\seq{\hat t_j}{j}{0}{\infi}\) and \(\seq{\tilde t_j}{j}{0}{\infi}\) coincide.
 In particular, \(\seq{\hat t_j}{j}{0}{\infi}\) is the \tFTv{\(\seq{\hat s_j}{j}{0}{\infi}\)}.
 The application of the already for \(\kappa=\infi\) proved assertion of the present lemma yields with \(\hat a_0\defeq-\ug\hat s_0+\hat s_1\) for the \tiFaaT{\hat a_0}{\hat s_0} \(H\) of \(G\) thus \(H\in\RFqabg{\seq{\hat s_j}{j}{0}{\infi}}\).
 In view of \rprop{ab.R0933} and \rrem{F.R.fpftr}, the sequence \(\seq{\hat s_j}{j}{0}{m}\) belongs to \(\Fggqu{m}\) and its \tfpf{} is exactly \(\seq{\hat{\mathfrak{f}}_j}{j}{0}{2m}\).
 Consequently, the \tfpf{s} of \(\seq{\hat s_j}{j}{0}{m}\) and \(\seqs{m}\) coincide as well.
 By virtue of \rthm{F.T.FggFP}, then the sequences \(\seq{\hat s_j}{j}{0}{m}\) and \(\seqs{m}\) coincide. 
 Hence, \(H\) is exactly the \tiFaaTv{\sau{0}}{\su{0}}{\(G\)} and belongs to \(\RFqabsg{m}\).
\eproof

\section{Parametrization of the set of all solutions}\label{S1245}
 We are now going to iterate the \tFATion{} introduced in \rdefn{F.D.FTF} with the \tFABTion{} introduced in \rdefn{F.D.FTFAB}.
 To that end, we use the \tnFT{k} \(\seq{\su{j}^\FTa{k}}{j}{0}{\kappa-k}\) of a sequence \(\seqska\) constructed in \rdefn{ab.N1020} and, in addition, the sequence \(\seq{\sau{j}^\FTa{k}}{j}{0}{\kappa-1-k}\) given by \(\sau{j}^\FTa{k}\defeq-\ug\su{j}^\FTa{k}+\su{j+1}^\FTa{k}\), \tie{}, the sequence built from \(\seq{\su{j}^\FTa{k}}{j}{0}{\kappa-k}\) according to \rnota{F.N.sa}:

\bdefnl{F.D.stepFT}
 Let \(\dom\) be a non-empty subset of \(\C\), let \(F\colon\dom\to\Cpq\) be a matrix-valued function, and let \(\seqs{\kappa}\) be a sequence of complex \tpqa{matrices}.
 Let \(\FTUa{0}{F}{\seqs{\kappa}}\defeq F\).
 Recursively, for all \(k\in\mn{1}{\kappa}\), denote by \(\FTUa{k}{F}{\seqs{\kappa}}\) the \tFaaTv{\sau{0}^\FTa{k-1}}{\su{0}^\FTa{k-1}}{\(\FTUa{k-1}{F}{\seqs{\kappa}}\)}.
 In view of \rdefn{F.D.FTF}, for all \(k\in\mn{0}{\kappa}\), denote by \(\FTPUa{k}{F}{\seqs{\kappa}}\) the \tFatpv{\su{0}^\FTa{k}}{\(\FTUa{k}{F}{\seqs{\kappa}}\)}.
 Then, for all \(m\in\mn{0}{\kappa}\), we call \(\FTPUa{m}{F}{\seq{\su{j}}{j}{0}{\kappa}}\) the \emph{\tnFTPv{m}{\(F\)}{\(\seq{\su{j}}{j}{0}{\kappa}\)}} and we call \(\FTUa{m}{F}{\seq{\su{j}}{j}{0}{\kappa}}\) the \emph{\tnFTFv{m}{\(F\)}{\(\seq{\su{j}}{j}{0}{\kappa}\)}}.
\edefn

\breml{F.R.1FT-FTA}
 The pair \(\FTPUa{0}{F}{\seqs{\kappa}}\) is exactly the \tFatpv{\su{0}}{\(F\)}.
 If \(\kappa\geq1\), then \(\FTUa{1}{F}{\seqs{\kappa}}\) is exactly the \tFaaTv{\sau{0}}{\su{0}}{\(F\)}.
\erem

 Regarding \rprop{F.P.Filft}, we will use the following mappings:

\bnotal{F.N.iFTabb}
 For each matrix \(M\in\Cggq\), let the mapping \(\FTiu{M}\) be defined on the class \(\PRFabqa{M}\) by \(\FTiua{M}{\copa{G_1}{G_2}}\defeq F\), where \(F\) is the \tiFaTv{M}{\(\copa{G_1}{G_2}\)}.
 Furthermore, given two matrices \(A,M\in\Cpq\), let the mapping \(\FTiuu{A}{M}\) be defined on the set of all matrix-valued functions \(G\colon\Cab\to\Cpq\) by \(\FTiuua{A}{M}{G}\defeq F\), where \(F\) is the \tiFaaTv{A}{M}{\(G\)}.
\enota

\bpropl{F.P.P-S0}
 Let \(\seqs{0}\in\Fggqu{0}\).
 Then \(\psi\colon\rsetcl{\PRFabqa{\su{0}}}\to\RFqabsg{0}\) defined by \(\psi(\rpcl{G_1}{G_2})\defeq\FTiua{\su{0}}{\copa{G_1}{G_2}}\) is a bijection with inverse \(\psi^\inv\) given by \(\psi^\inv(F)=\rsetcl{\FTPUa{0}{F}{\seqs{0}}}\).
\eprop
\bproof
 Consider arbitrary \(\copa{G_1}{G_2}\in\PRFabqa{\su{0}}\) and \(F\in\RFqabsg{0}\).
 Due to \rrem{F.R.Fgg-s}, we have \(\su{0}\in\Cggq\).
 According to \rcor{F.C.iFwd}, then \(\psi(\rpcl{G_1}{G_2})\) is independent of the concrete representative of the equivalence class \(\rpcl{G_1}{G_2}\).
 Since, because of \rlem{ab.P1448}, furthermore \(\FTiua{\su{0}}{\copa{G_1}{G_2}}\) belongs to \(\RFqabsg{0}\), the mapping \(\psi\) is  well defined.
 Regarding \rrem{F.R.1FT-FTA}, we obtain from \rlem{ab.P1505} moreover \(\FTPUa{0}{F}{\seqs{0}}\in\PRFabqa{\su{0}}\).
 Consequently, the mapping \(\chi\colon\RFqabsg{0}\to\rsetcl{\PRFabqa{\su{0}}}\) defined by \(\chi(S)\defeq\rsetcl{\FTPUa{0}{S}{\seqs{0}}}\) is well defined as well.
 Using \rlem{F.L.FT-11} and \rrem{F.R.1FT-FTA}, we conclude \((\chi\circ\psi)(\rpcl{G_1}{G_2})=\rpcl{G_1}{G_2}\).
 Because of \rrem{F.R.SF-ST} and \rlem{B.L.rn}, we get \(\ran{F(z)}=\ran{\su{0}}\) for all \(z\in\Cab\).
 Therefore, \(\OPu{\ran{\su{0}}}F=F\).
 Taking additionally into account \rrem{F.R.1FT-FTA} and \(\FTPUa{0}{F}{\seqs{0}}\in\PRFabqa{\su{0}}\), \rlem{F.L.FT1-1} then yields \((\psi\circ\chi)(F)=F\).
 Consequently, \(\psi\) is a bijection with inverse \(\chi\).
\eproof

\bpropl{F.P.S0-S1}
 Let \(m\in\N\) and let \(\seqs{m} \in\Fggqu{m}\) with \tFT{} \(\seqt{m-1}\).
 Then \(\psi\colon\RFqabg{\seqt{m-1}}\to\RFqabsg{m}\) defined by \(\psi(G)\defeq\FTiuua{\sau{0}}{\su{0}}{G}\) is a bijection with inverse \(\psi^\inv\) given by \(\psi^\inv(F)=\FTUa{1}{F}{\seqs{m}}\).
\eprop
\bproof
 In view of \rlem{ab.P1555}, the mapping \(\psi\) is well defined.
 Regarding \rrem{F.R.1FT-FTA}, we see from \rlem{ab.P1235} that the mapping \(\chi\colon\RFqabsg{m}\to\RFqabg{\seqt{m-1}}\) given by \(\chi(F)\defeq\FTUa{1}{F}{\seqs{m}}\) is also well defined.
 Using \rrem{F.R.Fgg-s}, we easily infer \(\sau{0}\in\CHq\) and \(\su{0}\in\Cggq\).
 Because of \rrem{F.R.Fgg-r}, we have \(\ran{\sau{0}}\subseteq\ran{\su{0}}\).
 Consider now an arbitrary \(G\in\RFqabg{\seqt{m-1}}\).
 Then \(G\in\RFqab\).
 Taking into account \rrem{F.R.SF-ST}, we conclude from \rlem{B.L.rn} furthermore \(\ran{G(z)}=\ran{t_0}\) for all \(z\in\Cab\).
 Regarding \rdefn{ab.N0940}, we thus obtain \(\ran{G(z)}\subseteq\ran{\sau{0}}\) for all \(z\in\Cab\).
 By construction, \(\psi(G)\) is the \tiFaaTv{\sau{0}}{\su{0}}{\(G\)}.
 Denote by \(H\) the \tFaaTv{\sau{0}}{\su{0}}{\(\psi(G)\)}.
 In view of \rrem{F.R.1FT-FTA}, then \(H=\chi(\psi(G))\).
 Using \rlem{F.L.FAM-11}, hence \(H=G\) follows.
 Consequently, \((\chi\circ\psi)(G)=G\).
 Now we consider an arbitrary \(F\in\RFqabsg{m}\).
 By virtue of \rlem{ab.P1235}, the \tFaaT{\sau{0}}{\su{0}} \(\chi(F)\) of \(F\) then belongs to \(\RFqabg{\seqt{m-1}}\).
 As above, we thus have \(\chi(F)\in\RFqab\) and \(\ran{[\chi(F)](z)}\subseteq\ran{\sau{0}}\) for all \(z\in\Cab\).
 Observe that the \tRabMa{} \(\rabmF\) of \(F\) satisfies \(\rabmFa{\ab}=\su{0}\).
 Using \rprop{ab.P1648L1409}, we get \(\ran{F(z)}=\ran{\su{0}}\) and \(\nul{F(z)}=\nul{\su{0}}\) for all \(z\in\Cab\).
 Because of \rrem{ab.L0921}, the function \(\Ffa\) given in \rnota{ab.N1537} belongs to \(\RFqabg{\seqsa{m-1}}\).
 Therefore, we obtain analogously \(\ran{\Fav{z}}=\ran{\sau{0}}\) and \(\nul{\Fav{z}}=\nul{\sau{0}}\) for all \(z\in\Cab\).
 Hence, we can apply \rlem{F.L.FAM1-1} to the function \(F\) and its \tFaaT{\sau{0}}{\su{0}} \(\chi(F)\) and obtain with the \tiFaaT{\sau{0}}{\su{0}} \(H=\psi(\chi(F))\) of \(\chi(F)\) then \(H(z)=F(z)\) for all \(z\in\Cab\).
 Thus, \((\psi\circ\chi)(F)=F\).
 Consequently, \(\psi\) is bijective with inverse \(\chi\).
\eproof

 The combination of \rpropss{F.P.P-S0}{F.P.S0-S1} now yields a first parametrization of the solution set of the matricial \tHausdorff{} \rmprobm{\ab}{m}{=}, where, however, the set of parameters still depends on the given data.
 
\bthml{F.T.Phi}
 Let \(m\in\NO\) and let \(\seqs{m}\in\Fggqu{m}\).
 Let \(\psi_m\colon\rsetcl{\PRFabqa{\su{0}^\FTa{m}}}\to\RFqabg{\seq{\su{j}^\FTa{m}}{j}{0}{0}}\) be defined by \(\psi_m(\rpcl{G_1}{G_2})\defeq\FTiua{\su{0}^\FTa{m}}{\copa{G_1}{G_2}}\).
 In the case \(m\geq1\) let, for all \(k\in\mn{0}{m-1}\), furthermore \(\psi_k\colon\RFqabg{\seq{\su{j}^\FTa{k+1}}{j}{0}{m-k-1}}\to\RFqabg{\seq{\su{j}^\FTa{k}}{j}{0}{m-k}}\) be given by \(\psi_k(G)\defeq\FTiuua{\sau{0}^\FTa{k}}{\su{0}^\FTa{k}}{G}\).
 Then \(\Psi_m\colon\rsetcl{\PRFabqa{\su{0}^\FTa{m}}}\to\RFqabsg{m}\) defined by \(\Psi_m(\rpcl{G_1}{G_2})\defeq(\psi_0\circ\psi_1\circ\dotsm\circ\psi_m)(\rpcl{G_1}{G_2})\) is a bijection with inverse \(\Psi_m^\inv\) given by \(\Psi_m^\inv(F)=\rsetcl{\FTPUa{m}{F}{\seqs{m}}}\).
\ethm
\bproof
 Because of \rprop{ab.P1030}, we have \(\seq{\su{j}^\FTa{k}}{j}{0}{m-k}\in\Fggqu{m-k}\) for all \(k\in\mn{0}{m}\).
 According to \rprop{F.P.P-S0}, then \(\psi_m\) is a bijection with inverse \(\psi_m^\inv\) given by \(\psi_m^\inv(F)=\rsetcl{\FTPUa{0}{F}{\seq{\su{j}^\FTa{m}}{j}{0}{0}}}\).
 Regarding \rdefn{ab.N1020}, we infer in the case \(m\geq1\) for all \(k\in\mn{0}{m-1}\) from \rprop{F.P.S0-S1} that \(\psi_k\) is a bijection with inverse \(\psi_k^\inv\) given by \(\psi_k^\inv(F)=\FTUa{1}{F}{\seq{\su{j}^\FTa{k}}{j}{0}{m-k}}\).
 In view of \rrem{F.R.1FT-FTA}, we see, for all \(F\in\RFqabg{\seq{\su{j}^\FTa{m}}{j}{0}{0}}\), that \(\psi_m^\inv(F)\) is exactly the equivalence class of the \tFatnpv{\su{0}^\FTa{m}}{\(F\)}.
 In the case \(m\geq1\), for all \(k\in\mn{0}{m-1}\) and all \(F\in\RFqabg{\seq{\su{j}^\FTa{k}}{j}{0}{m-k}}\), furthermore \(\psi_k^\inv(F)\) coincides with the equivalence class of the \tFaaTv{\sau{0}^\FTa{k}}{\su{0}^\FTa{k}}{\(F\)}.
 Regarding \rdefn{F.D.stepFT}, we obtain in the case \(m\geq1\), for all \(F\in\RFqabg{\seqs{m}}\), the equation \((\psi_{m-1}^\inv\circ\dotsm\circ\psi_0^\inv)(F)=\FTUa{m}{F}{\seqs{m}}\) and hence \((\psi_m^\inv\circ\psi_{m-1}^\inv\circ\dotsm\circ\psi_0^\inv)(F)=\rsetcl{\FTPUa{m}{F}{\seqs{m}}}\), implying that \(\Psi_{m}\) is a bijection with inverse \(\Psi_{m}^\inv\) given by \(\Psi_{m}^\inv(F)=\rsetcl{\FTPUa{m}{F}{\seqs{m}}}\).
\eproof

 We end this section by mentioning a relation between the \tnFT{k} given in \rdefn{F.D.stepFT} of a matrix-valued function with respect to a sequence of matrices and the \tnSfT{k} introduced in \rdefn{F.D.SN-F}.
 We start with the case \(k=1\):
 
\bleml{L0951}
 Suppose \(\kappa\geq1\).
 Let \(\seqska\in\Fggqka\) and let \(F\in\RFqabsg{\kappa}\) with \tFaaT{\sau{0}}{\su{0}} \(G\) and first \tSfT{} \(F^\FTa{1}\).
 Then \(G=F^\FTa{1}\).
\elem
\bproof 
 By assumption we have \(F\in\RFqab\) with \tRabMa{} \(\sigma\defeq\rabmF\) belonging to \(\MggqFksg\).
 Observe that \(\sigma\in\MggqinfF\) according to \rrem{F.R.RabM8}.
 Setting \(\su{j}\defeq\int_\ab x^j\sigma\rk{\dif x}\) for all \(j\in\minf{\kappa+1}\), we have then \(\sigma\in\MggqFsg{\infi}\) and hence \(F\in\RFqabsg{\infi}\).
 In particular, \rthm{I.P.ab} shows \(\seqska\in\Fggqinf\).
 Regarding \(\kappa\geq1\), the matrix \(\sau{0}=-\ug\su{0}+\su{1}\) is not affected by the above extension of the sequence \(\seqska\).
 Thus, we can apply \rlem{ab.P1235} to obtain \(G\in\RFqabg{\seqt{\infi}}\), where \(\seqt{\infi}\) is the \tFTv{\(\seqsinf\)}.
 Consequently, we have \(G\in\RFqab\) with \tRabMa{} \(\nu\defeq\rabmG\) belonging to \(\MggqFag{t}{\infi}\).
 By virtue of \rprop{F.P.STbij}, in particular \(G=\STFu{\nu}\).
 According to \rdefn{F.D.SN-M}, the first \tFmT{} \(\mu\defeq\sigma^\FTa{1}\) of \(\sigma\) belongs to \(\MggqFag{t}{\infi}\) as well.
 Due to \rprop{I.P.ab8} and \rthm{I.P.ab}, the set \(\MggqFag{t}{\infi}\) consists of at most one element.
 Hence, \(\nu=\mu\) follows.
 Since \(F^\FTa{1}=\STFu{\mu}\) by \rdefn{F.D.SN-F}, we infer then \(G=\STFu{\nu}=\STFu{\mu}=F^\FTa{1}\).
\eproof

\bpropl{P0941}
 Let \(\seqska\in\Fggqka\), let \(k\in\mn{0}{\kappa}\), and let \(F\in\RFqabsg{\kappa}\) with \tnFT{k} \(\FTUa{k}{F}{\seq{\su{j}}{j}{0}{\kappa}}\) of \(F\) with respect to \(\seqska\) and \tnSfT{k} \(F^\FTa{k}\) of \(F\).
 Then \(\FTUa{k}{F}{\seq{\su{j}}{j}{0}{\kappa}}=F^\FTa{k}\).
\eprop
\bproof 
 We use mathematical induction.
 According to \rdefn{F.D.stepFT}, we have \(\FTUa{0}{F}{\seqs{\kappa}}=F\), whereas \rrem{R0730} shows \(F^\FTa{0}=F\).
 Hence, the assertion holds true for \(k=0\).
 Now assume that \(\kappa\geq1\) and that \(\FTUa{k-1}{F}{\seq{\su{j}}{j}{0}{\kappa}}=F^\FTa{k-1}\) is valid for some \(k\in\mn{1}{\kappa}\).
 Setting \(T\defeq F^\FTa{k-1}\) and \(\tu{j}\defeq\su{j}^\FTa{k-1}\) for all \(j\in\mn{0}{\kappa-\rk{k-1}}\), then \(\FTUa{k}{F}{\seqs{\kappa}}\) is, by \rdefn{F.D.stepFT}, exactly the \tFaaTv{-\ug\tu{0}+\tu{1}}{\tu{0}}{\(T\)}.
 Furthermore, \rrem{R0734} yields \(T\in\RFqabag{t}{\kappa-\rk{k-1}}\).
 From \rprop{ab.P1030} we infer \(\seqa{t}{\kappa-\rk{k-1}}\in\Fggqu{\kappa-\rk{k-1}}\).
 Taking additionally into account \(\kappa-\rk{k-1}\geq1\), we can apply \rlem{L0951} to the sequence \(\seqa{t}{\kappa-\rk{k-1}}\) and the function \(T\) to see that the \tFaaTv{-\ug\tu{0}+\tu{1}}{\tu{0}}{\(T\)} coincides with the first \tSfTv{\(T\)}, \tie{}, \(\FTUa{k}{F}{\seqs{\kappa}}=T^\FTa{1}\).
 Since \rrem{R0730} provides \(F^\FTa{k}=T^\FTa{1}\), the proof is complete.
\eproof

 In view of \rlem{L0951}, a more explicit description of the Schur--Nevanlinna type algorithm for the class \(\RFqab\) considered in \rsec{S1323} can be given by means of matricial linear fractional transformations.
 For the sake of simplicity, we illustrate this for the scalar case \(q=1\), where this amounts to a scalar linear fractional transformation or a continued fraction expansion of functions belonging to \(\RFab{1}\).
 These considerations are along the lines of the classical results by Schur~\zitas{MR1580958,MR1580948} for the class \(\SchF{1}{1}{\mathbb{D}}\) of holomorphic functions mapping the open unit disc \(\mathbb{D}\defeq\setaca{z\in\C}{\abs{z}<1}\) into the closed unit disc \(\overline{\mathbb{D}}\) (\tcf{}~\rnota{A.N.SF}) and by Nevanlinna~\zita{zbMATH02604576} for the class \(\RFuu{0}{1}\) introduced in \rnota{N1409}:
 
 Let \(f\in\RFab{1}\) with \tRabMa{} \(\rabmfu{f}\).
 Then \(f\in\RFuabg{1}{\seqsinf}\) with the sequence \(\seqsinf\) of power moments \(\su{j}\defeq\int_\ab x^j\rabmfu{f}\rk{\dif x}\) associated with \(\rabmfu{f}\) belonging to \(\Fgguuuu{1}{\infi}{\ug}{\obg}\), by virtue of \rprop{I.P.ab8Fgg}.
 We see from \rlem{L0951} that the first \tSfT{} \(f^\FTa{1}\) of \(f\) is exactly the \tFaaTv{\sau{0}}{\su{0}}{\(f\)}.
 Consequently, \rlem{ab.P1235} yields \(f^\FTa{1}\in\RFuabg{1}{\seqa{s^\FTa{1}}{\infi}}\) with the \tFT{} \(\seqa{s^\FTa{1}}{\infi}\) of \(\seqsinf\) belonging to \(\Fgguuuu{1}{\infi}{\ug}{\obg}\), by virtue of \rprop{ab.P1030}.
 Hence, we can conclude from \rprop{F.P.FP8} that the expansions
\begin{align*}
 f(z)
 &=-\frac{\su{0}}{z}-\frac{\su{1}}{z^2}-\frac{\su{2}}{z^3}-\dotsb&
&\text{and}&
 f^\FTa{1}(z)
 &=-\frac{\su{0}^\FTa{1}}{z}-\frac{\su{1}^\FTa{1}}{z^2}-\frac{\su{2}^\FTa{1}}{z^3}-\dotsb
\end{align*}
 are valid for all \(z\in\C\) with \(\abs{z}>\max\set{\abs{\ug},\abs{\obg}}\). 
 Recall that \(\sau{0}=-\ug\su{0}+\su{1}\) and \(\sub{0}=\obg\su{0}-\su{1}\), according to \rnota{F.N.sa}.
 \rlem{F.R.Fgg-s} yields \(\su{0}\geq0\) and \(\sau{0}\geq0\).
 In what follows, we assume \(\su{0}>0\).
 Because of \eqref{F.G.d01} and \rcor{F.C.diaFT}, we have then
\begin{align*}
 \dia{0}&=\ba\su{0}&
&\text{and}&
 \su{0}^\FTa{1}
 =\dia{1}
 &=-\ug\obg\su{0}+(\ug+\obg)\su{1}-\frac{\su{1}^2}{\su{0}}
 =\frac{\sau{0}\sub{0}}{\su{0}}.
\end{align*}
 In view of \(\ba=\obg-\ug>0\), thus \(\dia{0}>0\).
 Regarding \rrem{A.R.A-1} and \(\sau{0}+\sub{0}=\ba\su{0}\), we obtain, by virtue of \rdefn{D0754} and \eqref{F.G.f012}, hence
\begin{align*}
 \cia{0}&=\fpu{0}=\su{0}&
&\text{and}&
 \cia{1}
 &=\frac{\fpu{2}}{\dia{0}}
 =\frac{\sub{0}}{\ba\su{0}}
 =\frac{\ba\su{0}-\sau{0}}{\ba\su{0}}
 =1-\frac{\sau{0}}{\ba\su{0}}.
\end{align*}
 Taken all together, we can infer by direct calculation \(\su{0}=\cia{0}\), \(\sau{0}=\ba\cia{0}\rk{1-\cia{1}}\), \(\sub{0}=\ba\cia{0}\cia{1}\), and \(\su{0}^\FTa{1}=\ba^2\cia{0}\cia{1}\rk{1-\cia{1}}\).
 In view of \(\su{0}=\rabmfua{f}{\ab}\), we can conclude from \rrem{ab.L0921} that the function \(\Fka{f}\colon\Cab\to\C\) given, according to \rnota{ab.N1537}, by \(\Fkav{f}{z}\defeq\rk{z-\ug}f\rk{z}+\su{0}\) belongs to \(\RFab{1}\) with \tRabMa{} \(\rabmFfa\) fulfilling \(\rabmFav{\ab}=\sau{0}\).
 If \(\sau{0}=0\), then \rprop{ab.P1648L1409} yields \(\Fkav{f}{z}=0\), \tie{}, \(f\rk{z}=\rk{\ug-z}^\inv\su{0}\) for all \(z\in\Cab\), implying \(\rabmfu{f}=\su{0}\Kronu{\ug}\), by virtue of \rprop{F.P.STbij}, where \(\Kronu{\ug}\) is the Dirac measure on \(\rk{\ab,\BsAF}\) with unit mass at \(\ug\).
 Now assume \(\sau{0}>0\).
 \rprop{ab.P1648L1409} yields, for all \(z\in\Cab\), then \(\Fkav{f}{z}\neq0\).
 In view of \rdefn{F.D.FTFAB} and \rrem{A.R.A-1}, we thus can infer the representation 
\beql{g=f/f}
 f^\FTa{1}\rk{z}
 =\frac{\sau{0}^2/\su{0}}{\obg-z}\cdot\frac{\rk{\obg-z}f\rk{z}-\su{0}}{\rk{z-\ug}f\rk{z}+\su{0}}
 =\frac{\ba^2\cia{0}\rk{1-\cia{1}}^2}{\obg-z}\cdot\frac{\rk{\obg-z}f\rk{z}-\cia{0}}{\rk{z-\ug}f\rk{z}+\cia{0}}
\eeq
 for all \(z\in\Cab\). 
 As seen above, we have \(f^\FTa{1}\in\RFqab\).
 Taking additionally into account \(\Fkav{f}{z}\neq0\) for all \(z\in\Cab\) and the assumptions \(\su{0}>0\) and \(\sau{0}>0\), the conditions of \rlem{F.L.FAM1-1} are fulfilled.
 So its application shows that \(f\) coincides with the \tiFaaTv{\sau{0}}{\su{0}}{\(f^\FTa{1}\)}.
 Observe that \(\OPn{\su{0}}=0\) and \(\OPn{\sau{0}}=0\) by \eqref{F.G.PQ} and the assumptions  \(\su{0}>0\) and \(\sau{0}>0\).
 From \rlem{F.L.iFAMlft} we can conclude for all \(z\in\Cab\) then \(\sau{0}^2/\su{0}-\rk{z-\ug}f^\FTa{1}\rk{z}\neq0\) and the representation
\[
 f\rk{z}
 =\frac{\su{0}}{\obg-z}\cdot\frac{\sau{0}^2/\su{0}+\rk{\obg-z}f^\FTa{1}\rk{z}}{\sau{0}^2/\su{0}-\rk{z-\ug}f^\FTa{1}\rk{z}}
 =\frac{\cia{0}}{\obg-z}\cdot\frac{\ba^2\cia{0}\rk{1-\cia{1}}^2+\rk{\obg-z}f^\FTa{1}\rk{z}}{\ba^2\cia{0}\rk{1-\cia{1}}^2-\rk{z-\ug}f^\FTa{1}\rk{z}},
\]
 which also follows by direct calculation from \eqref{g=f/f}.
 Regarding \(\rk{\obg-z}^\inv+\rk{z-\ug}^\inv=\ba\rk{\obg-z}^\inv\rk{z-\ug}^\inv\), we can rewrite this, for all \(z\in\Cab\), as
\[\begin{split}
 f\rk{z}
 &=\frac{\su{0}}{z-\ug}\cdot\frac{\frac{\sau{0}^2/\su{0}}{\obg-z}+f^\FTa{1}\rk{z}}{\frac{\sau{0}^2/\su{0}}{z-\ug}-f^\FTa{1}\rk{z}}
 =\frac{-\su{0}}{\ug-z}\cdot\frac{\frac{\ba\sau{0}^2/\su{0}}{\rk{\obg-z}\rk{z-\ug}}+f^\FTa{1}\rk{z}-\frac{\sau{0}^2/\su{0}}{z-\ug}}{\frac{\sau{0}^2/\su{0}}{z-\ug}-f^\FTa{1}\rk{z}}\\
 &=\frac{\su{0}\ek*{\frac{\sau{0}^2/\su{0}}{z-\ug}-f^\FTa{1}\rk{z}}+\rk{\ug-z}\frac{\ba\sau{0}^2}{\rk{\obg-z}\rk{z-\ug}^2}}{\rk{\ug-z}\ek*{\frac{\sau{0}^2/\su{0}}{z-\ug}-f^\FTa{1}\rk{z}}}
 =\frac{\su{0}}{\ug-z}+\frac{\frac{\ba\sau{0}^2}{\rk{\obg-z}\rk{z-\ug}^2}}{\frac{\sau{0}^2/\su{0}}{z-\ug}-f^\FTa{1}\rk{z}},
\end{split}\]
 giving rise to a continued fraction expansion of functions \(f\in\RFab{1}\).

\section{Description via linear fractional transformation}\label{F.S.glt}
 We are now going to write the parametrization \(\Psi_m\) from \rthm{F.T.Phi} as a matricial linear fractional transformation, as considered in \rapp{A.s1.lft}.
 The generating matrix-valued function of this transformation is a composition of certain instances of the matrix polynomials introduced in \rnotass{ab.N1246b}{ab.N1546b}:

\bnotal{F.N.Vsm}
 Let \(\seqs{\kappa}\) be a sequence of complex \tpqa{matrices} and let \(m\in\mn{0}{\kappa}\).
 Then let \(\rmif{m}\defeq V_0V_1\dotsm V_m\), where \(V_k\defeq\mFTiuu{\sau{0}^\FTa{k}}{\su{0}^\FTa{k}}\) for all \(k\in\mn{0}{m-1}\) and \(V_m\defeq\mFTiu{\su{0}^\FTa{m}}\).
\enota
 
 Regarding \rnotass{ab.N1246b}{ab.N1546b}, we see that \(\rmif{m}\) is a complex \taaa{(p+q)}{(p+q)}{matrix} polynomial with \(\deg\rmif{m}\leq2(m+1)\).
 As a main result of the present paper, we now obtain a description of the set of all solutions to \rprobrab{m} via a matricial linear fractional transformation generated by that matrix polynomial:
 
\bthml{F.P.PVR}
 Let \(m\in\NO\) and let \(\seqs{m}\in\Fggqu{m}\).
 Denote by \(\smat{\rmifrnw{m}&\rmifrne{m}\\\rmifrsw{m}&\rmifrse{m}}\) the \tqqa{\tbr{}} of the restriction of \(\rmif{m}\) onto \(\Cab\):
 \benui
  \il{F.P.PVR.a} Let \(\Gamma\in\rsetcl{\PRFabqa{\su{0}^\FTa{m}}}\) and let \(\copa{G_1}{G_2}\in\Gamma\).
  Then \(\det\rk{\rmifrsw{m}G_1+\rmifrse{m}G_2}\) does not vanish identically in \(\Cab\) and the function
  \beql{F.P.PVR.B1}
   F
   =\rk{\rmifrnw{m}G_1+\rmifrne{m}G_2}\rk{\rmifrsw{m}G_1+\rmifrse{m}G_2}^\inv
  \eeq
  belongs to \(\RFqabsg{m}\).
  \il{F.P.PVR.b} For each \(F\in\RFqabsg{m}\), there exists a unique equivalence class \(\Gamma\in\rsetcl{\PRFabqa{\su{0}^\FTa{m}}}\) such that \eqref{F.P.PVR.B1} holds true for all \(\copa{G_1}{G_2}\in\Gamma\), namely the equivalence class \(\rsetcl{\FTPUa{m}{F}{\seqs{m}}}\) of the \tnFTP{m} \(\FTPUa{m}{F}{\seqs{m}}\) of \(F\) with respect to \(\seqs{m}\).
 \eenui
\ethm
\bproof
 Let the mappings \(\psi_0,\psi_1,\dotsc,\psi_m\) be defined as in \rthm{F.T.Phi}.
 Let \(V_m\defeq\mFTiu{\su{0}^\FTa{m}}\) and, in the case \(m\geq1\), let \(V_k\defeq\mFTiuu{\sau{0}^\FTa{k}}{\su{0}^\FTa{k}}\) for all \(k\in\mn{0}{m-1}\).
 For all \(k\in\mn{0}{m}\), let \(\smat{\tilde a_k&\tilde b_k\\ \tilde c_k& \tilde d_k}\) be the \tqqa{\tbr{}} of the restriction of \(V_k\) onto \(\Cab\).
 Because of \rprop{ab.P1030}, we have \(\seq{\su{j}^\FTa{k}}{j}{0}{m-k}\in\Fggqu{m-k}\) for all \(k\in\mn{0}{m}\).
 Using \rremss{F.R.Fgg-s}{F.R.Fgg-r}, we thus easily obtain \(\sau{0}^\FTa{k}\in\CHq\) and \(\ran{\sau{0}^\FTa{k}}\subseteq\ran{\su{0}^\FTa{k}}\) for all \(k\in\mn{0}{m-1}\) and, furthermore, \(\su{0}^\FTa{k}\in\Cggq\) for all \(k\in\mn{0}{m}\).
 Consider now an arbitrary pair \(\copa{G_1}{G_2}\in\PRFabqa{\su{0}^\FTa{m}}\).
 According to \rnota{ab.N1534}, then \(\copa{G_1}{G_2}\) belongs to \(\PRFabq\).
 In particular, \(G_1\) and \(G_2\) are \(\Cqq\)\nobreakdash-valued functions meromorphic in \(\Cab\).
 Because of \rprop{ab.P1528}, the set \(\exset\defeq\pol{G_1}\cup\pol{G_2}\cup\PRFabex{G_1}{G_2}\) is a discrete subset of \(\Cab\).
 Hence, \(\C\setminus\rk{\ab\cup\exset}\neq\emptyset\).
 Let \(H_m\defeq\psi_m(\rpcl{G_1}{G_2})\).
 In view of \rdefn{F.D.iFTF} and \rnotass{ab.N1246b}{F.N.iFTabb}, we conclude from \rprop{F.P.Filft} then \(\det\ek{\tilde c_m(z)G_1(z)+\tilde d_m(z)G_2(z)}\neq0\) and
\[
 \ek*{\tilde a_m(z)G_1(z)+\tilde b_m(z)G_2(z)}\ek*{\tilde c_m(z)G_1(z)+\tilde d_m(z)G_2(z)}^\inv
 =\ek*{\psi_m\rk*{\rpcl{G_1}{G_2}}}\rk{z}
 =H_m(z)
\]
 for all \(z\in\C\setminus\rk{\ab\cup\exset}\).
 In the case \(m\geq1\), let \(H_k\defeq\psi_k(H_{k+1})\) for all \(k\in\mn{0}{m-1}\).
 By virtue of \rthm{F.T.Phi}, we have \(H_m\in\RFqabg{\seq{\su{j}^\FTa{m}}{j}{0}{0}}\) and, in the case \(m\geq1\), moreover \(H_k\in\RFqabg{\seq{\su{j}^\FTa{k}}{j}{0}{m-k}}\) for all \(k\in\mn{0}{m-1}\).
 
 In the case \(m\geq1\), we now consider an arbitrary \(\ell\in\mn{0}{m-1}\).
 Then \(H_{\ell+1}\in\RFqab\) and the \tRabMa{} \(\rabmfu{\ell+1}\) of \(H_{\ell+1}\) fulfills \(\rabmfua{\ell+1}{\ab}=\su{0}^\FTa{\ell+1}\).
 Using \rprop{ab.P1648L1409}, we obtain in particular \(\ran{H_{\ell+1}(z)}=\ran{\su{0}^\FTa{\ell+1}}\) for all \(z\in\Cab\).
 Regarding \rdefnss{ab.N1020}{ab.N0940}, we infer furthermore \(\ran{\su{0}^\FTa{\ell+1}}\subseteq\ran{\sau{0}^\FTa{\ell}}\).
 Consequently, \(\ran{H_{\ell+1}(z)}\subseteq\ran{\sau{0}^\FTa{\ell}}\) follows for all \(z\in\Cab\).
 
 Thus, in the case \(m\geq1\), \rlem{F.L.iFAMlft} yields, in view of \rnotass{ab.N1546b}{F.N.iFTabb}, for all \(k\in\mn{0}{m-1}\), then \(\det\ek{\tilde c_k(z)H_{k+1}(z)+\tilde d_k(z)}\neq0\) and
\[
 \ek*{\tilde a_k(z)H_{k+1}(z)+\tilde b_k(z)}\ek*{\tilde c_k(z)H_{k+1}(z)+\tilde d_k(z)}^\inv
 =\ek*{\psi_k(H_{k+1})}(z)
 =H_k(z)
\]
 for all \(z\in\Cab\).
 By virtue of \(\rmif{m}=V_0V_1\dotsm V_m\), we can conclude from \rprop{M.P.gltP163} hence \(\det\ek{\rmifrswa{m}{z}G_1(z)+\rmifrsea{m}{z}G_2(z)}\neq0\) and
\begin{multline*}
 H_0(z)
 =\ek*{(\psi_0\circ\psi_1\circ\dotsm\circ\psi_m)\rk*{\rpcl{G_1}{G_2}}}(z)\\
 =\ek*{\rmifrnwa{m}{z}G_1(z)+\rmifrnea{m}{z}G_1(z)}\ek*{\rmifrswa{m}{z}G_1(z)+\rmifrsea{m}{z}G_2(z)}^\inv
\end{multline*}
 for all \(z\in\C\setminus\rk{\ab\cup\exset}\).
 In particular, \(\det\rk{\rmifrsw{m}G_1+\rmifrse{m}G_2}\) does not vanish identically in \(\Cab\).
 Consequently, \(\rk{\rmifrnw{m}G_1+\rmifrne{m}G_2}\rk{\rmifrsw{m}G_1+\rmifrse{m}G_2}^\inv\) is a \(\Cqq\)\nobreakdash-valued function meromorphic in \(\Cab\).
 By virtue of \(H_0\in\RFqabg{\seqs{m}}\), the matrix-valued function \(H_0\) is holomorphic in \(\Cab\).
 Since the set \(\exset\) is discrete, we can conclude from the identity theorem for holomorphic functions then
 \[
  (\psi_0\circ\psi_1\circ\dotsm\circ\psi_m)\rk*{\rpcl{G_1}{G_2}}
  =H_0
  =\rk{\rmifrnw{m}G_1+\rmifrne{m}G_2}\rk{\rmifrsw{m}G_1+\rmifrse{m}G_2}^\inv.
 \]
 The application of \rthm{F.T.Phi} completes the proof.
\eproof

\breml{R1335} 
 Let \(m\in\NO\) and let \(\seqs{m}\in\Fggqu{m}\).
 In view of \(\ba>0\) and \rcor{F.C.diaFT}, we have \(\ba^{m-1}\dia{m}=\su{0}^\FTa{m}\), and consequently \(\ran{\dia{m}}=\ran{\su{0}^\FTa{m}}\) and \(\nul{\dia{m}}=\nul{\su{0}^\FTa{m}}\).
 According to \rnota{ab.N1534}, hence \(\PRFabqa{\dia{m}}=\PRFabqa{\su{0}^\FTa{m}}\).
\erem

 Now we study separately three cases depending on the rank \(r\) of the matrix \(\dia{m}=\ba^{-\rk{m-1}}\su{0}^\FTa{m}\) determining the amount of determinacy of the truncated moment problem in question.
 We distinguish between the case \(r=q\), the case \(1\leq r\leq q-1\), and the case \(r=0\).
 
 First we consider the so-called non-degenerate case \(r=q\).
 In view of \rrem{P1338}, this is exactly the case of \(\seqs{m}\in\Fgqu{m}\).
 For this situation \rprobrab{m} was already considered in \zitas{MR2222521,MR2342899}, applying Potapov's method of fundamental matrix inequalities separately for even and odd \tnn{} integers \(m\), \tresp{}
 The generating matrix-valued functions of the linear fractional transformation obtained there are matrix polynomials.
 For the scalar case \(q=1\) we refer to Krein/Nudelman \zitaa{MR0458081}{\cch{IV}, \S~7}, where the generating matrix function of the linear fractional transformation is built from orthogonal polynomials of first kind and second kind.

\bthml{R1245} 
 Let \(m\in\NO\) and let \(\seqs{m}\in\Fggqu{m}\) be such that the matrix \(\dia{m}\) is non-singular.
 Then all statements of \rthm{F.P.PVR} are valid with the class \(\PRFabq\) instead of \(\PRFabqa{\su{0}^\FTa{m}}\).
\ethm
\bproof
 Use \rremss{R1335}{R1234} and \rthm{F.P.PVR}.
\eproof

 Now we turn our attention to the degenerate, but not completely degenerate case \(1\leq r\leq q-1\):

\bthml{F.P.RFabmred} 
 Assume \(q\geq2\).
 Let \(m\in\NO\) and let \(\seqs{m}\in\Fggqu{m}\).
 Let \(r\defeq\rank\dia{m}\) and assume \(1\leq r\leq q-1\).
 Denote by \(\smat{\rmifrnw{m}&\rmifrne{m}\\\rmifrsw{m}&\rmifrse{m}}\) the \tqqa{\tbr{}} of the restriction of \(\rmif{m}\) onto \(\Cab\).
 Let \(u_1,u_2,\dotsc,u_q\) be an orthonormal basis of \(\Cq\) with \(\set{u_1,u_2,\dotsc,u_r}\subseteq\ran{\dia{m}}\) and let \(W\defeq\mat{u_1,u_2,\dotsc,u_q}\).
 Then:
\benui
 \il{F.P.RFabmred.a} Let \(\copa{g_1}{g_2}\in\PRFab{r}\).
 Then \(\det\ek{\rmifrsw{m} W\rk{\zdiag{g_1}{\Ouu{\rk{q-r}}{\rk{q-r}}}}+\rmifrse{m}W\rk{\zdiag{g_2}{\Iu{q-r}}}}\) does not vanish identically in \(\Cab\).
 \il{F.P.RFabmred.b} For each \(\copa{g_1}{g_2}\in\PRFab{r}\) let \(S_{W,\copa{g_1}{g_2}}\colon\Cab\to\Cqq\) be defined by
\begin{multline*}
 S_{W,\copa{g_1}{g_2}}
 \defeq\ek*{\rmifrnw{m} W\rk{\zdiag{g_1}{\Ouu{\rk{q-r}}{\rk{q-r}}}}+\rmifrne{m}W\rk{\zdiag{g_2}{\Iu{q-r}}}}\\
 \times\ek*{\rmifrsw{m} W\rk{\zdiag{g_1}{\Ouu{\rk{q-r}}{\rk{q-r}}}}+\rmifrse{m}W\rk{\zdiag{g_2}{\Iu{q-r}}}}^\inv.
\end{multline*}
  Then \(\Sigma_W\colon\rsetcl{\PRFab{r}}\to\RFqabsg{m}\) defined by \(\Sigma_W\rk{\rpcl{g_1}{g_2}}\defeq S_{W,\copa{g_1}{g_2}}\) is well defined and bijective.
 \eenui
\ethm
\bproof
 Obviously, \(W\) is a unitary \tqqa{matrix} and \(U\defeq\mat{u_1,u_2,\dotsc,u_r}\) is the left \taaa{q}{r}{block} of \(W\).
 In view of \(\dim\ran{\dia{m}}=r\), we furthermore see that \(u_1,u_2,\dotsc,u_r\) is an orthonormal basis of \(\ran{\dia{m}}\).
 Due to \rlem{F.L.PRFabr}, the mapping \(\Gamma_U\colon\rsetcl{\PRFab{r}}\to\rsetcl{\PRFabqa{\dia{m}}}\) defined by \(\Gamma_U\rk{\rpcl{f_1}{f_2}}\defeq\rpcl{Uf_1U^\ad}{Uf_2U^\ad+\OPu{\ek{\ran{\dia{m}}}^\orth}}\) is thus well defined and bijective.
 By virtue of \rrem{R1335}, we have \(\PRFabqa{\dia{m}}=\PRFabqa{\su{0}^\FTa{m}}\).
 Consider now an arbitrary pair \(\copa{g_1}{g_2}\in\PRFab{r}\).
 Then \(g_1\) and \(g_2\) are \(\Coo{r}{r}\)\nobreakdash-valued functions, which are meromorphic in \(\Cab\). Observe that \(\rmifrnw{m}\), \(\rmifrne{m}\), \(\rmifrsw{m}\), and \(\rmifrse{m}\), as restrictions of \tqqa{matrix} polynomials, are \(\Cqq\)\nobreakdash-valued functions, which are holomorphic in \(\Cab\).
 Consequently, we can easily conclude that
\begin{align*}
 X_{U,\copa{g_1}{g_2}}&\defeq\rmifrnw{m}Ug_1U^\ad+\rmifrne{m}\rk{Ug_2U^\ad+\OPu{\ek{\ran{\dia{m}}}^\orth}}
\shortintertext{and}
 Y_{U,\copa{g_1}{g_2}}&\defeq\rmifrsw{m}Ug_1U^\ad+\rmifrse{m}\rk{Ug_2U^\ad+\OPu{\ek{\ran{\dia{m}}}^\orth}}
\end{align*} 
 are \(\Cqq\)\nobreakdash-valued functions which are meromorphic in \(\Cab\) with \(\pol{X_{U,\copa{g_1}{g_2}}}\subseteq\pol{g_1}\cup\pol{g_2}\) and \(\pol{Y_{U,\copa{g_1}{g_2}}}\subseteq\pol{g_1}\cup\pol{g_2}\).
 Therefore, \(\det Y_{U,\copa{g_1}{g_2}}\) is a complex-valued function which is meromorphic in \(\Cab\) with \(\pol{\det Y_{U,\copa{g_1}{g_2}}}\subseteq\pol{g_1}\cup\pol{g_2}\).
 Similarly,
\begin{align*}
 T_{W,\copa{g_1}{g_2}}&\defeq\rmifrnw{m}W\rk{\zdiag{g_1}{\Ouu{\rk{q-r}}{\rk{q-r}}}}+\rmifrne{m}W\rk{\zdiag{g_2}{\Iu{q-r}}}
\shortintertext{and}
 R_{W,\copa{g_1}{g_2}}&\defeq\rmifrsw{m}W\rk{\zdiag{g_1}{\Ouu{\rk{q-r}}{\rk{q-r}}}}+\rmifrse{m}W\rk{\zdiag{g_2}{\Iu{q-r}}}
\end{align*}
 are \(\Cqq\)\nobreakdash-valued functions which are meromorphic in \(\Cab\) with \(\pol{T_{W,\copa{g_1}{g_2}}}\subseteq\pol{g_1}\cup\pol{g_2}\) and \(\pol{R_{W,\copa{g_1}{g_2}}}\subseteq\pol{g_1}\cup\pol{g_2}\), and \(\det R_{W,\copa{g_1}{g_2}}\) is a complex-valued function which is meromorphic in \(\Cab\) with \(\pol{\det R_{W,\copa{g_1}{g_2}}}\subseteq\pol{g_1}\cup\pol{g_2}\).
 In view of \(\copa{g_1}{g_2}\in\PRFab{r}\), we have \(\Gamma_U\rk{\rpcl{g_1}{g_2}}\in\rsetcl{\PRFabqa{\dia{m}}}=\rsetcl{\PRFabqa{\su{0}^\FTa{m}}}\).
 Due to \rthmp{F.P.PVR}{F.P.PVR.a}, hence \(\det Y_{U,\copa{g_1}{g_2}}\) does not vanish identically in \(\Cab\).
 Regarding the identity theorem for holomorphic functions, then \(\mathcal{N}\defeq\setaca{\zeta\in\C\setminus\rk{\ab\cup\pol{\det Y_{U,\copa{g_1}{g_2}}}}}{\det Y_{U,\copa{g_1}{g_2}}(\zeta)=0}\) is a discrete subset of \(\Cab\).
 Consequently, \(\mathcal{D}\defeq\pol{g_1}\cup\pol{g_2}\cup\mathcal{N}\) is a discrete subset of \(\Cab\), which fulfills \(\pol{X_{U,\copa{g_1}{g_2}}}\cup\pol{Y_{U,\copa{g_1}{g_2}}}\cup\pol{\det Y_{U,\copa{g_1}{g_2}}}\subseteq\mathcal{D}\) and \(\pol{T_{W,\copa{g_1}{g_2}}}\cup\pol{R_{W,\copa{g_1}{g_2}}}\cup\pol{\det R_{W,\copa{g_1}{g_2}}}\subseteq\mathcal{D}\).
 In particular, \(\C\setminus\rk{\ab\cup\mathcal{D}}\) is non-empty.
 Consider now an arbitrary \(z\in\C\setminus\rk{\ab\cup\mathcal{D}}\).
 Then \(\det Y_{U,\copa{g_1}{g_2}}(z)\neq0\) and, furthermore,
\begin{align*}
 X_{U,\copa{g_1}{g_2}}\rk{z}&=\rmifrnwa{m}{z}\ek*{Ug_1\rk{z}U^\ad}+\rmifrnea{m}{z}\ek*{Ug_2\rk{z}U^\ad+\OPu{\ek{\ran{\dia{m}}}^\orth}},\\
 Y_{U,\copa{g_1}{g_2}}\rk{z}&=\rmifrswa{m}{z}\ek*{Ug_1\rk{z}U^\ad}+\rmifrsea{m}{z}\ek*{Ug_2\rk{z}U^\ad+\OPu{\ek{\ran{\dia{m}}}^\orth}},\\
 T_{W,\copa{g_1}{g_2}}\rk{z}&=\rmifrnwa{m}{z} W\ek*{\zdiag{g_1\rk{z}}{\Ouu{\rk{q-r}}{\rk{q-r}}}}+\rmifrnea{m}{z}W\ek*{\zdiag{g_2\rk{z}}{\Iu{q-r}}},
\shortintertext{and}
 R_{W,\copa{g_1}{g_2}}\rk{z}&=\rmifrswa{m}{z}W\ek*{\zdiag{g_1\rk{z}}{\Ouu{\rk{q-r}}{\rk{q-r}}}}+\rmifrsea{m}{z}W\ek*{\zdiag{g_2\rk{z}}{\Iu{q-r}}}.
\end{align*} 
 Consequently, using \rlemp{A.L.red}{A.L.red.b}, we obtain \(\det R_{W,\copa{g_1}{g_2}}(z)\neq0\) and
\[
 \ek*{X_{U,\copa{g_1}{g_2}}\rk{z}}\ek*{Y_{U,\copa{g_1}{g_2}}\rk{z}}^\inv
 =\ek*{T_{W,\copa{g_1}{g_2}}\rk{z}}\ek*{R_{W,\copa{g_1}{g_2}}\rk{z}}^\inv.
\]
 In particular, \(\det R_{W,\copa{g_1}{g_2}}\) does not vanish identically in \(\Cab\), showing~\eqref{F.P.RFabmred.a}.
 Therefore, \(S_{W,\copa{g_1}{g_2}}=T_{W,\copa{g_1}{g_2}}R_{W,\copa{g_1}{g_2}}^\inv\) is a \(\Cqq\)\nobreakdash-valued function which is meromorphic in \(\Cab\) and holomorphic at \(z\) with
\beql{F.P.RFabmred.0}\begin{split}
 S_{W,\copa{g_1}{g_2}}\rk{z}
 &=\ek*{X_{U,\copa{g_1}{g_2}}\rk{z}}\ek*{Y_{U,\copa{g_1}{g_2}}\rk{z}}^\inv\\
 &=\rk*{\rmifrnwa{m}{z}\ek*{Ug_1\rk{z}U^\ad}+\rmifrnea{m}{z}\ek*{Ug_2\rk{z}U^\ad+\OPu{\ek{\ran{\dia{m}}}^\orth}}}\\
 &\qquad\times\rk*{\rmifrswa{m}{z}\ek*{Ug_1\rk{z}U^\ad}+\rmifrsea{m}{z}\ek*{Ug_2\rk{z}U^\ad+\OPu{\ek{\ran{\dia{m}}}^\orth}}}^\inv.
\end{split}\eeq
 Since \eqref{F.P.RFabmred.0} holds true for all \(z\in\C\setminus\rk{\ab\cup\mathcal{D}}\) and the set \(\mathcal{D}\) is discrete, we can conclude from the identity theorem for holomorphic functions that
\begin{multline*} 
 S_{W,\copa{g_1}{g_2}}
 =\ek*{\rmifrnw{m} \rk{Ug_1U^\ad}+\rmifrne{m}\rk{Ug_2U^\ad+\OPu{\ek{\ran{\dia{m}}}^\orth}}}\\
 \times\ek*{\rmifrsw{m}\rk{Ug_1U^\ad}+\rmifrse{m}\rk{Ug_2U^\ad+\OPu{\ek{\ran{\dia{m}}}^\orth}}}^\inv.
\end{multline*}
 Regarding \rthmp{F.P.PVR}{F.P.PVR.a}, for each pair \(\copa{P}{Q}\in\PRFabqa{\dia{m}}=\PRFabqa{\su{0}^\FTa{m}}\), we can consider the function \(F_{\copa{P}{Q}}\defeq\rk{\rmifrnw{m} P+\rmifrne{m}Q}\rk{\rmifrsw{m} P+\rmifrse{m}Q}^\inv\).
 In view of \rdefn{ab.N0843}, it is readily checked that, given two arbitrary pairs \(\copa{P_1}{Q_1},\copa{P_2}{Q_2}\in\PRFabqa{\dia{m}}\), the equivalence \(\copa{P_2}{Q_2}\rpaeq\copa{P_2}{Q_2}\) implies \(F_{\copa{P_1}{Q_1}}=F_{\copa{P_2}{Q_2}}\).
 According to \rthm{F.P.PVR}, thus the mapping \(\Pi\colon\rsetcl{\PRFabqa{\dia{m}}}\to\RFqabsg{m}\) defined by \(\Pi\rk{\rpcl{P}{Q}}\defeq F_{\copa{P}{Q}}\) is well defined and bijective.
 By virtue of that we have already shown, we get
\[\begin{split}
 \Sigma_W\rk*{\rpcl{g_1}{g_2}}
 =S_{W,\copa{g_1}{g_2}}
 &=F_{\copa{Ug_1U^\ad}{Ug_2U^\ad+\OPu{\ek{\ran{\dia{m}}}^\orth}}}\\
 &=\Pi\rk*{\rpcl{Ug_1U^\ad}{Ug_2U^\ad+\OPu{\ek{\ran{\dia{m}}}^\orth}}}
 =\Pi\rk*{\Gamma_U\rk*{\rpcl{g_1}{g_2}}}.
\end{split}\]
 Since this holds true for all pairs \(\copa{g_1}{g_2}\in\PRFab{r}\), we thus verified that \(\Sigma_W=\Pi\circ\Gamma_U\) is well defined and bijective, \tie{},~\eqref{F.P.RFabmred.b} holds true.
\eproof

 Now we treat the completely degenerate case \(r=0\).
 We will see in particular that in this situation the solution is unique.

\bthml{F.P.RFabunq}
 Let \(m\in\NO\) and let \(\seqs{m}\in\Fggqu{m}\) with \tfdf{} \(\seqdia{m}\).
 Then \(\RFqabsg{m}\) consists of exactly one element if and only if \(\dia{m}=\Oqq\).
 In this case, \(\RFqabsg{m}=\set{\rmifrne{m}\rmifrse{m}^\inv}\), where \(\smat{\rmifrnw{m}&\rmifrne{m}\\\rmifrsw{m}&\rmifrse{m}}\) denotes the \tqqa{\tbr{}} of the restriction of \(\rmif{m}\) onto \(\Cab\).
\ethm
\bproof
 Assume \(\dia{m}\neq\Oqq\).
 For each \(\ell\in\set{0,1}\) let \(\seq{\su{\ell,j}}{j}{0}{m+1}\) be given by
\[
 \su{\ell,j}
 \defeq
 \begin{cases}
  \su{j}\tincase{0\leq j\leq m}\\
  \umg{m}+\ell\dia{m}\tincase{j=m+1}
 \end{cases}.
\]
 Due to \rcor{ab.R1011}, then \(\set{\seq{\su{0,j}}{j}{0}{m+1},\seq{\su{1,j}}{j}{0}{m+1}}\subseteq\Fggqu{m+1}\).
 By virtue of \rprop{F.P.FPsolv}, thus \(\RFqabg{\seq{\su{\ell,j}}{j}{0}{m+1}}\) is non-empty for each \(\ell\in\set{0,1}\).
 Consequently, we can choose, for each \(\ell\in\set{0,1}\), a function \(F_\ell\in\RFqabg{\seq{\su{\ell,j}}{j}{0}{m+1}}\).
 Observe that, for each \(\ell\in\set{0,1}\), the \tRabMa{} \(\sigmau{\ell}\) of \(F_\ell\) belongs to \(\Mggouaag{q}{m+1}{\ab}{\seq{\su{\ell,j}}{j}{0}{m+1}}\), implying \(\int_\ab x^{m+1}\sigmaua{\ell}{\dif x}=\su{\ell,m+1}\).
 Because of \(\su{1,m+1}-\su{0,m+1}=\dia{m}\neq\Oqq\), we have \(\sigmau{0}\neq\sigmau{1}\) and hence \(F_0\neq F_1\).
 Since the functions \(F_0\) and \(F_1\) both belong to \(\RFqabsg{m}\), therefore this set consists of at least two elements.
 Thus, we can reversely conclude that if the set \(\RFqabsg{m}\) consists of exactly one element, then necessarily \(\dia{m}=\Oqq\) follows.

 Assume \(\dia{m}=\Oqq\).
 By virtue of \rlem{F.L.PRFabr}, then \(\rsetcl{\PRFabqa{\dia{m}}}\) consists of exactly one element, namely the equivalence class \(\rpcl{G_1}{G_2}\) of the pair \(\copa{G_1}{G_2}\) built from the functions \(G_1,G_2\colon\Cab\to\Cqq\) defined by \(G_1(z)\defeq\Oqq\) and \(G_2(z)\defeq\Iq\).
 Due to \rrem{R1335}, we have \(\PRFabqa{\dia{m}}=\PRFabqa{\su{0}^\FTa{m}}\).
 Consequently,
\beql{F.P.RFabunq.1}
 \rsetcl*{\PRFabqa{\su{0}^\FTa{m}}}
 =\rsetcl*{\PRFabqa{\dia{m}}}
 =\set*{\rpcl{G_1}{G_2}}
\eeq
 follows.
 Because of \rprop{F.P.FPsolv}, the set \(\RFqabsg{m}\) is non-empty.
 Consider an arbitrary \(F\in\RFqabsg{m}\).
 Taking into account \eqref{F.P.RFabunq.1}, we can infer from \rthm{F.P.PVR} then
\[
 F
 =\rk{\rmifrnw{m}G_1+\rmifrne{m}G_2}\rk{\rmifrsw{m}G_1+\rmifrse{m}G_2}^\inv
 =\rmifrne{m}\rmifrse{m}^\inv,
\]
 implying \(\RFqabsg{m}=\set{\rmifrne{m}\rmifrse{m}^\inv}\).
 In particular, \(\RFqabsg{m}\) consists of exactly one element.
\eproof

 In view of \rpropss{F.P.STbij}{ab.R13371422} and \rrem{F.R.ABL}, we obtain from \rthm{F.P.RFabunq} immediately the following two results:

\bcorl{F.P.MFunq2n+1}
 Let \(n\in\NO\) and let \(\seqs{2n+1}\in\Fggqu{2n+1}\).
 Then \(\MggqFsg{2n+1}\) consists of exactly one element if and only if \(\ran{\Lau{n}}\cap\ran{\Lub{n}}=\set{\Ouu{q}{1}}\).
\ecor

\bcorl{F.P.MFunq2n}
 Let \(n\in\N\) and let \(\seqs{2n}\in\Fggqu{2n}\).
 Then \(\MggqFsg{2n}\) consists of exactly one element if and only if \(\ran{\Lu{n}}\cap\ran{\Lab{n-1}}=\set{\Ouu{q}{1}}\).
\ecor

 We end this section with a necessary condition for unique solvability of \rprobm{\ab}{m}{=}:

\bcorl{F.C.>RFabunq}
 Let \(m\in\N\) and let \(\seqs{m}\) be a sequence of complex \tqqa{matrices} such that \(\MggqFsg{m}\) consists of exactly one element.
 Then \(\rank\usc{m}+\rank\osc{m}\leq q\).
\ecor
\bproof
 Due to \rthm{I.P.ab}, we have \(\seqs{m}\in\Fggqu{m}\).
 In view of \rprop{F.P.STbij}, we obtain from \rthm{F.P.RFabunq} thus \(\dia{m}=\Oqq\), \tie{}\ \(\rank\dia{m}=0\).
 Since \(\rank\dia{m-1}\leq q\) holds true, we can infer by virtue of \rcor{ab.R0938} then \(\rank\usc{m}+\rank\osc{m}\leq q\).
\eproof

  We are now going to factorize \(\rmif{m}\) in a way alternative to \rnota {F.N.Vsm}.
 In a first step, we derive by virtue of \rlem {ab.R1443} a connection between \(\rmif{m}\) and \(\rmif{m-1}\), which against the background of \rthm{F.P.PVR} correlates the solution sets:
 \(\RFqabsg{m}\subseteq\RFqabsg{m-1}\)
 To that end, we make use of the sequences \(\seq{\sau{j}^\FTa{k}}{j}{0}{\kappa-1-k}\) and \(\seq{\sub{j}^\FTa{k}}{j}{0}{\kappa-1-k}\) given by \(\sau{j}^\FTa{k}\defeq-\ug\su{j}^\FTa{k}+\su{j+1}^\FTa{k}\) and \(\sub{j}^\FTa{k}\defeq\obg\su{j}^\FTa{k}-\su{j+1}^\FTa{k}\), \tresp{}, \tie{}, the sequences built from the \tnFT{k} \(\seq{\su{j}^\FTa{k}}{j}{0}{\kappa-k}\) according to \rnota{F.N.sa}:

\bleml{ab.P1348} 
 Suppose \(\kappa\geq1\).
 Let \(\seqska\in\Fggqka\).
 Then \(\rmif{m}=\rmif{m-1}\cFTiuu{\sau{0}^\FTa{m-1}}{\su{0}^\FTa{m-1}}\) for all \(m\in\mn{1}{\kappa}\).
\elem
\bproof
 Consider an arbitrary \(m\in\mn{1}{\kappa}\).
 In view of \rnota{F.N.Vsm} it is sufficient to verify the identity \(\mFTiuu{A}{M}\mFTiu{D}=\mFTiu{M}\cFTiuu{A}{M}\) with \(A\defeq\sau{0}^\FTa{m-1}\), \(M\defeq\su{0}^\FTa{m-1}\), and \(D\defeq\su{0}^\FTa{m}\).
 
 Due to \rprop{ab.P1030}, the sequence \(\seq{\su{j}^\FTa{m-1}}{j}{0}{\kappa-m+1}\) belongs to \(\Fggqu{\kappa-m+1}\).
 Using \rremss{F.R.Fgg-s}{F.R.Fgg-r}, we can hence infer \(A,M\in\CHq\) and \(\ran{A}\subseteq\ran{M}\).
 Let \(B\defeq\ba M-A\).
 Because of \rrem{F.R.012}, then \(B=\ba\su{0}^\FTa{m-1}-\sau{0}^\FTa{m-1}=\sub{0}^\FTa{m-1}\).
 According to \rdefnss{D1861}{ab.N0940}, denote by \(\seq{\dia{j}^\FTa{m-1}}{j}{0}{\kappa-m+1}\) the \tfdfa{\(\seq{\su{j}^\FTa{m-1}}{j}{0}{\kappa-m+1}\)} and by \(\seqt{\kappa-m}\) the \tFTv{\(\seq{\su{j}^\FTa{m-1}}{j}{0}{\kappa-m+1}\)}, \tresp{}
 \rrem{F.R.Fgg<D} yields \(\seq{\su{j}^\FTa{m-1}}{j}{0}{\kappa-m+1}\in\Dqqu{\kappa-m+1}\).
 By virtue of \rrem{ab.L0907} and \rlem{ab.P1728a}, we can hence infer \(AM^\mpi B=\dia{1}^\FTa{m-1}=t_0\).
 In view of \rdefn{ab.N1020}, we have \(t_0=\su{0}^\FTa{m}=D\).
 Consequently, \(D= AM^\mpi B\) follows.
 Thus, we can apply \rlem{ab.R1443} to obtain \(\ek{\mFTiuua{A}{M}{z}}\ek{\mFTiua{D}{z}}=\ek{\mFTiua{M}{z}}\ek{\cFTiuua{A}{M}{z}}\) for all \(z\in\C\).
\eproof 
 
  The above mentioned factorization alternative to \rnota{F.N.Vsm} results from \rlem{ab.P1348} now by means of mathematical induction:
\bpropl{F.P.V=VU} 
 Let \(\seqska\in\Fggqka\) and let \(m\in\mn{0}{\kappa}\).
 Then \(\rmif{m}=U_0U_1\dotsm U_m\), where \(U_\ell\defeq\cFTiuu{\sau{0}^\FTa{\ell-1}}{\su{0}^\FTa{\ell-1}}\) for all \(\ell\in\mn{1}{m}\) and \(U_0\defeq\mFTiu{\su{0}}\).
\eprop

\bleml{ab.L1652} 
 Suppose \(\kappa\geq1\).
 Let \(\seqska\in\Fggqka\) with \tfpf{} \(\fpseqka\) and \tfdf{} \(\seqdiaka\).
 Let \(k\in\mn{0}{\kappa-1}\) and let \(z\in\C\).
 Then \(\cFTiuua{\sau{0}^\FTa{k}}{\su{0}^\FTa{k}}{z}=\tmat{U_{11}&U_{12}\\U_{21}&U_{22}}\), where
 \begin{align*}
  U_{11}&=\dia{k}\ek*{\rk{\obg-z}\fpu{2 k +1}^\mpi\fpu{2 k +1}\dia{k}^\mpi\fpu{2 k +2}+\rk{z-\ug}\rk{\Iq-\fpu{2 k +1}^\mpi\fpu{2 k +1}}\dia{k}^\mpi \fpu{2 k +1}}\dia{ k +1}^\mpi,\qquad
  U_{12}=\ba^ k \fpu{2 k +2},\\
  U_{21}&=-\rk{\obg-z}\rk{z-\ug}\ba^{- k +1}\dia{k}^\mpi\fpu{2 k +1}\dia{ k +1}^\mpi,\hspace{60pt}
  U_{22}=\rk{\obg-z}\ba\ek*{\rk{\Iq-\dia{k}^\mpi\dia{k}}+\dia{k}^\mpi\fpu{2 k +1}}.
 \end{align*}
\elem
\bproof
 Let \(A\defeq\sau{0}^\FTa{k}\), let \(M\defeq\su{0}^\FTa{k}\), and let \(B\defeq\ba M-A\).
 Because of \rrem{F.R.012}, then \(B=\ba\su{0}^\FTa{k}-\sau{0}^\FTa{k}=\sub{0}^\FTa{k}\).
 \rcor{F.C.diaFT} shows \(M=\ba^{k-1}\dia{k}\).
 Furthermore, \rprop{F.P.FPFT} yields \(A=\sau{0}^\FTa{k}=\ba^ k \fpu{2 k +1}\) and \(B=\ba^ k \fpu{2 k +2}\).
 Let \(D\defeq AM^\mpi B\) and denote by \(\seq{\dia{j}^\FTa{k}}{j}{0}{\kappa- k }\) the \tfdfa{\(\seq{\su{j}^\FTa{k}}{j}{0}{\kappa- k }\)}.
 Due to \rprop{ab.P1030} we have \(\seq{\su{j}^\FTa{k}}{j}{0}{\kappa- k }\in\Fggqu{\kappa- k }\).
 Because of \rrem{F.R.Fgg<D}, in particular \(\seq{\su{j}^\FTa{k}}{j}{0}{\kappa- k }\in\Dqqu{\kappa- k }\).
 Consequently, from \rrem{ab.L0907} we conclude \(\dia{1}^\FTa{k}=D\).
 \rprop{F.P.diaalg}, then implies \(D=\ba^ k \dia{ k +1}\).
 Using \rremss{A.R.l*A}{A.R.A-1} and taking into account \(\ba>0\), \rnota{ab.N1409}, and \eqref{F.G.PQ}, the assertion follows.
\eproof

\section{On the sets \(\RFqabsg{m+1}\) in the case of \habHd{} extensions of a sequence \(\seqs{m}\in\Fggqu{m}\)}\label{S0948}%
 In this section, we study \tabHd{} extensions of a sequence \(\seqs{m}\in\Fggqu{m}\).
 First we recall the notion of \tabHd{} sequences belonging to \(\Fggqu{m}\) and a characterization of this class of sequences.

\bdefnnl{\zitaa{MR3775449}{\cdefn{10.24}{203}}}{D0927}
 Let \(\ell\in\NO\) and let \(\seqs{\ell}\in\Fggqu{\ell}\) with \tfdf{} \(\seqdia{\ell}\) given in \rdefn{D1861}.
 Then \(\seqs{\ell}\) is called \emph{\tabHd{}} if \(\dia{\ell}=\Oqq\).
 We denote by \(\Fggdqu{\ell}\) the set of all sequences \(\seqs{\ell}\in\Fggqu{\ell}\) which are \tabHd{}.
\edefn

\bpropnl{\tcf{}~\zitaa{MR3979701}{\cprop{6.38}{2156}}}{F.P.abCDe}
 Let \(\ell\in\N\) and let \(\seqs{\ell}\in\Fggqu{\ell}\) with \tfcf{} \(\seqcia{\ell}\) given in \rdefn{D0754}.
 Then \(\seqs{\ell}\) is \tabHd{} if and only if \(\cia{\ell}^2=\cia{\ell}\).
\eprop

 Observe that in the situation of \rprop{F.P.abCDe}, due to \(\cia{\ell}\lgeq\Oqq\), we have \(\cia{\ell}^\ad=\cia{\ell}\) and thus the condition \(\cia{\ell}^2=\cia{\ell}\) is equivalent to \(\cia{\ell}\) being a transformation matrix corresponding to an orthogonal projection, \tie{}, \(\cia{\ell}=\OPu{\ran{\cia{\ell}}}\).
 
 Against the background of \rprop{F.P.abCDe}, we are looking now for a description of the set
\[
 \setaca{\su{m+1}\in\Cqq}{\seqs{m+1}\in\Fggdqu{m+1}}.
\]
 We will show that this set stands in a bijective correspondence to the set of all linear subspaces of \(\ran{\dia{m}}\).

\bnotal{N1336}
 Let \(m\in\NO\), let \(\seqs{m}\in\Fggqu{m}\) with \tfdf{} \(\seqdia{m}\), and let \(\cU\) be a linear subspace of \(\ran{\dia{m}}\).
 Then let \(\su{m,\cU}\defeq\omg{m}-\dia{m}^\varsqrt\OPu{\cU}\dia{m}^\varsqrt\) if \(m\) is even, and \(\su{m,\cU}\defeq\umg{m}+\dia{m}^\varsqrt\OPu{\cU}\dia{m}^\varsqrt\) if \(m\) is odd.
\enota

\bexal{R1403}
 Let \(m\in\NO\) and let \(\seqs{m}\in\Fggqu{m}\) with \tfdf{} \(\seqdia{m}\).
 Let \(\cU_0\defeq\set{\Ouu{q}{1}}\) and let \(\cU_1\defeq\ran{\dia{m}}\).
 Then \(\su{m,\cU_0}=\omg{m}\) and \(\su{m,\cU_1}=\umg{m}\) if \(m\) is even, and \(\su{m,\cU_0}=\umg{m}\) and \(\su{m,\cU_1}=\omg{m}\) if \(m\) is odd.

 Indeed, we have \(\OPu{\cU_0}=\Oqq\) and, in view of \(\ran{\dia{m}^\varsqrt}=\ran{\dia{m}}\), furthermore \(\OPu{\cU_1}=\OPu{\ran{\dia{m}^\varsqrt}}\).
 Consequently, \(\dia{m}^\varsqrt\OPu{\cU_0}\dia{m}^\varsqrt=\Oqq\) and \(\dia{m}^\varsqrt\OPu{\cU_1}\dia{m}^\varsqrt=\dia{m}\).
 Since, according to \rdefn{D1861}, we have \(\dia{m}=\omg{m}-\umg{m}\), the assertions follow by virtue of \rnota{N1336}.
\eexa

\bpropl{L1253}
 Let \(m\in\NO\) and let \(\seqs{m}\in\Fggqu{m}\).
 Then:
\benui
 \il{L1253.a} Let \(\cU\) be a linear subspace of \(\ran{\dia{m}}\) and let \(\su{m+1}\defeq\su{m,\cU}\).
 Then \(\seqs{m+1}\in\Fggdqu{m+1}\).
 Furthermore, \(\fpu{2m+1}=\dia{m}^\varsqrt\rk{\Iq-\OPu{\cU}}\dia{m}^\varsqrt\), \(\fpu{2m+2}=\dia{m}^\varsqrt\OPu{\cU}\dia{m}^\varsqrt\), and \(\cia{m+1}=\OPu{\cU}\).
 \il{L1253.b} Let \(\su{m+1}\in\Cqq\) be such that \(\seqs{m+1}\in\Fggdqu{m+1}\).
 Then there exists a linear subspace \(\cU\) of \(\ran{\dia{m}}\) such that \(\su{m+1}=\su{m,\cU}\), namely \(\cU=\ran{\cia{m+1}}\).
 \il{L1253.c} Let \(\cU\) and \(\cV\) be linear subspaces of \(\ran{\dia{m}}\).
 Then \(\cU=\cV\) if and only if \(\su{m,\cU}=\su{m,\cV}\).
\eenui
\eprop
\bproof
 First observe that \(\dia{m}\) belongs to \(\Cggq\), due to \rprop{ab.C0929}.
 Let \(D\defeq\dia{m}^\varsqrt\).

 \eqref{L1253.a} By virtue of \rdefnss{D0752}{D2972} and \rnota{N1336}, we have
\[
 \fpu{2m+2}
 =\fpu{4n+2}
 =\osc{2n+1}
 =\osc{m+1}
 =\omg{m}-\su{m+1}
 =D\OPu{\cU}D
\]
 in the case \(m=2n\) for some \(n\in\NO\), and
\[
 \fpu{2m+2}
 =\fpu{4n+4}
 =\usc{2n+2}
 =\usc{m+1}
 =\su{m+1}-\umg{m}
 =D\OPu{\cU}D
\]
 in the case \(m=2n+1\) for some \(n\in\NO\).
 Using \rrem{F.R.f2n-1}, we can infer then
\(
 \fpu{2m+1}
 =\dia{m}-\fpu{2m+2}
 =DD-D\OPu{\cU}D
 =D\rk{\Iq-\OPu{\cU}}D
\).
 Because of \(\Oqq\lleq\OPu{\cU}\lleq\Iq\) and \rrem{A.R.XAX}, consequently the matrices \(\fpu{2m+2}\) and \(\fpu{2m+1}\) are both \tnnH{}.
 From \rprop{F.P.FggFP}, we can conclude now \(\seqs{m+1}\in\Fggqu{m+1}\).
 According to \rdefn{D0754}, then \(\cia{m+1}=D^\mpi\fpu{2m+2}D^\mpi=D^\mpi D\OPu{\cU}DD^\mpi\).
 Because of \(\ran{\dia{m}}=\ran{D}\), we have \(\cU\subseteq\ran{D}\).
 \rremss{R.P}{R.AA+B.A} then yield \(DD^\mpi\OPu{\cU}=\OPu{\cU}\).
 Furthermore \(D^\ad=D\) implies \(D^\mpi D=DD^\mpi\), by virtue of \rrem{ab.R1052}.
 Hence, \(D^\mpi D\OPu{\cU}=\OPu{\cU}\) follows.
 Taking account \(\OPu{\cU}^\ad=\OPu{\cU}\) and \eqref{mpi}, we furthermore obtain \(\OPu{\cU}=\rk{DD^\mpi\OPu{\cU}}^\ad=\OPu{\cU}DD^\mpi\).
 Consequently, we conclude \(\cia{m+1}=D^\mpi D\OPu{\cU}DD^\mpi=\OPu{\cU}\).
 In view of \(\OPu{\cU}^2=\OPu{\cU}\), the application of \rprop{F.P.abCDe} provides then \(\seqs{m+1}\in\Fggdqu{m+1}\).
 
 \eqref{L1253.b} According to \rprop{F.P.abCDe}, we have \(\cia{m+1}^2=\cia{m+1}\).
 From \rprop{F.P.ed} and \eqref{matint}, we get \(\cia{m+1}^\ad=\cia{m+1}\) and \(\Oqq\lleq\cia{m+1}\lleq\OPu{\ran{\dia{m}}}\).
 Consequently, \(\cia{m+1}=\OPu{\cU}\) with \(\cU\defeq\ran{\cia{m+1}}\).
 Furthermore, \(\OPu{\cU}\lleq\OPu{\ran{\dia{m}}}\), implying \(\cU\subseteq\ran{\dia{m}}\).
 \rlem{F.L.cia-fp} yields \(\fpu{2m+2}=D\cia{m+1}D\).
 By virtue of \rdefnss{D0752}{D2972}, we have then
\[
 D\OPu{\cU}D
 =\fpu{2m+2}
 =\fpu{4n+2}
 =\osc{2n+1}
 =\osc{m+1}
 =\omg{m}-\su{m+1}
\]
 in the case \(m=2n\) for some \(n\in\NO\), and
\[
 D\OPu{\cU}D
 =\fpu{2m+2}
 =\fpu{4n+4}
 =\usc{2n+2}
 =\usc{m+1}
 =\su{m+1}-\umg{m}
\]
 in the case \(m=2n+1\) for some \(n\in\NO\).
 In view of \rnota{N1336}, hence \(\su{m+1}=\su{m,\cU}\) follows.
 
 \eqref{L1253.c} Obviously, \(\cU=\cV\) implies \(\su{m,\cU}=\su{m,\cV}\), according to \rnota{N1336}.
 Conversely, suppose \(\su{m,\cU}=\su{m,\cV}\).
 From \rnota{N1336}, then \(D\OPu{\cU}D=D\OPu{\cV}D\) follows.
 By the same reasoning as in the proof of \rpart{L1253.a}, we can infer \(D^\mpi D\OPu{\cU}DD^\mpi=\OPu{\cU}\) and \(D^\mpi D\OPu{\cV}DD^\mpi=\OPu{\cV}\).
 Consequently, \(\OPu{\cU}=\OPu{\cV}\), implying \(\cU=\cV\).
\eproof

\bnotal{N1338}
 Let \(m\in\NO\), let \(\seqs{m}\in\Fggqu{m}\) with \tfdf{} \(\seqdia{m}\), and let \(\cU\) be a linear subspace of \(\ran{\dia{m}}\).
 Then let \(\dphi{m}{\cU},\dpsi{m}{\cU}\colon\Cab\to\Cqq\) be defined by
\begin{align*}
 \dphia{m}{\cU}{z}&\defeq\ba^{m-1}\dia{m}^\varsqrt\OPu{\cU}\dia{m}^\varsqrt,&
 \dpsia{m}{\cU}{z}&\defeq\rk{\obg-z}\ek{\Iq-\rk{\dia{m}^\varsqrt}^\mpi\OPu{\cU}\dia{m}^\varsqrt}.
\end{align*}
\enota

\bexal{E1341}
 Let \(m\in\NO\) and let \(\seqs{m}\in\Fggqu{m}\) with \tfdf{} \(\seqdia{m}\).
 Let \(\cU_0\defeq\set{\Ouu{q}{1}}\) and let \(\cU_1\defeq\ran{\dia{m}}\).
 For all \(z\in\Cab\), then \(\dphia{m}{\cU_0}{z}=\Oqq\) and \(\dpsia{m}{\cU_0}{z}\defeq\rk{\obg-z}\Iq\) as well as \(\dphia{m}{\cU_1}{z}\defeq\ba^{m-1}\dia{m}\) and \(\dpsia{m}{\cU_1}{z}\defeq\rk{\obg-z}\OPu{\nul{\dia{m}}}\).
 Indeed, as in the proof of \rexam{R1403}, we obtain \(\dia{m}^\varsqrt\OPu{\cU_0}\dia{m}^\varsqrt=\Oqq\) and \(\dia{m}^\varsqrt\OPu{\cU_1}\dia{m}^\varsqrt=\dia{m}\).
 Taking into account \(\dia{m}^\mpi\dia{m}^\varsqrt=\rk{\dia{m}^\varsqrt}^\mpi\), then \(\rk{\dia{m}^\varsqrt}^\mpi\OPu{\cU_0}\dia{m}^\varsqrt=\Oqq\) and \(\rk{\dia{m}^\varsqrt}^\mpi\OPu{\cU_1}\dia{m}^\varsqrt=\dia{m}^\mpi\dia{m}\) follow.
 Hence, \(\Iq-\rk{\dia{m}^\varsqrt}^\mpi\OPu{\cU_0}\dia{m}^\varsqrt=\Iq\) and, in view of \rrem{ab.R1052}, furthermore \(\Iq-\rk{\dia{m}^\varsqrt}^\mpi\OPu{\cU_1}\dia{m}^\varsqrt=\OPu{\nul{\dia{m}}}\).
 Now, the assertions follow by virtue of \rnota{N1338}.
\eexa

\breml{R1414}
 Let \(m\in\NO\), let \(\seqs{m}\in\Fggqu{m}\) with \tfdf{} \(\seqdia{m}\), and let \(\cU\) be a linear subspace of \(\ran{\dia{m}}\).
 In view of \rnotass{N1338}{ab.N1534} and \(\ran{\dia{m}}=\ran{\dia{m}^\varsqrt}\) as well as \rexam{FR}, then \(\copa{\dphi{m}{\cU}}{\dpsi{m}{\cU}}\in\PRFabqa{\dia{m}}\).
\erem

\bleml{R1419}
 Let \(m\in\NO\), let \(\seqs{m}\in\Fggqu{m}\) with \tfdf{} \(\seqdia{m}\), and let \(\cU\) be a linear subspace of \(\ran{\dia{m}}\).
 Denote by \(\smat{\rmifrnw{m}&\rmifrne{m}\\\rmifrsw{m}&\rmifrse{m}}\) the \tqqa{\tbr{}} of the restriction of \(\rmif{m}\) onto \(\Cab\).
 Then the function \(\det\rk{\rmifrsw{m}\dphi{m}{\cU}+\rmifrse{m}\dpsi{m}{\cU}}\) does not vanish identically.
\elem
\bproof
 In view of \rremss{R1414}{R1335}, this is a consequence of \rthmp{F.P.PVR}{F.P.PVR.a}.
\eproof

\rlem{R1419} shows that the following notation is correct.

\bnotal{N1417}
 Let \(m\in\NO\), let \(\seqs{m}\in\Fggqu{m}\) with \tfdf{} \(\seqdia{m}\), and let \(\cU\) be a linear subspace of \(\ran{\dia{m}}\)
 Denote by \(\smat{\rmifrnw{m}&\rmifrne{m}\\\rmifrsw{m}&\rmifrse{m}}\) the \tqqa{\tbr{}} of the restriction of \(\rmif{m}\) onto \(\Cab\).
 Then let
\[ 
 \dOSt{m}{\cU}
 \defeq\rk{\rmifrnw{m}\dphi{m}{\cU}+\rmifrne{m}\dpsi{m}{\cU}}\rk{\rmifrsw{m}\dphi{m}{\cU}+\rmifrse{m}\dpsi{m}{\cU}}^\inv.
\]
\enota

\bleml{R1426}
 Let \(m\in\NO\), let \(\seqs{m}\in\Fggqu{m}\) with \tfdf{} \(\seqdia{m}\), and let \(\cU\) be a linear subspace of \(\ran{\dia{m}}\).
 Then \(\dOSt{m}{\cU}\in\RFqabsg{m}\).
\elem
\bproof
 In view of \rremss{R1414}{R1335} and \rnota{N1417}, this is a consequence of \rthmp{F.P.PVR}{F.P.PVR.a}.
\eproof

 Given a sequence \(\seqs{m+1}\in\Fggqu{m+1}\), now we look for the \tqqa{matrix} polynomials in the four \tqqa{blocks} of the \taaa{2q}{2q}{matrix} polynomial \(\rmif{m+1}\) defined in \rnota{F.N.Vsm}.
 
\bleml{L1614}
 Let \(m\in\NO \) and let \(\seqs{m+1}\in\Fggqu{m+1}\) with \tfpf{} \(\fpseq{2m+2}\) and \tfdf{} \(\seqdia{m+1}\). 
 Denote by \(\smat{\rmifnw{m}&\rmifne{m}\\\rmifsw{m}&\rmifse{m}}\) and \(\smat{\rmifnw{m+1}&\rmifne{m+1}\\\rmifsw{m+1}&\rmifse{m+1}}\) the \tqqa{\tbr{s}} of \(\rmif{m}\) and \(\rmif{m+1}\), \tresp{}
 For all \(z\in\C\), then
\begin{align*}
 \rmifnwa{m+1}{z}&=\biggl\{\rmifnwa{m}{z}\dia{m}\ek*{\rk{\obg-z}\fpu{2m+1}^\mpi\fpu{2m+1}\dia{m}^\mpi\fpu{2m+2}+\rk{z-\ug}\rk{\Iq-\fpu{2m+1}^\mpi\fpu{2m+1}}\dia{m}^\mpi\fpu{2m+1}}\\
 &\qquad-\rk{\obg-z}\rk{z-\ug}\ba^{-m+1}\rmifnea{m}{z}\dia{m}^\mpi\fpu{2m+1}\biggr\}\dia{m+1}^\mpi,\\
 \rmifnea{m+1}{z}&=\ba\rk*{\ba^{m-1}\rmifnwa{m}{z}\fpu{2m+2}+\rk{\obg-z}\rmifnea{m}{z}\ek*{\rk{\Iq-\dia{m}^\mpi\dia{m}}+\dia{m}^\mpi\fpu{2m+1}}},\\
 \rmifswa{m+1}{z}&=\biggl\{\rmifswa{m}{z}\dia{m}\ek*{\rk{\obg-z}\fpu{2m+1}^\mpi\fpu{2m+1}\dia{m}^\mpi\fpu{2m+2}+\rk{z-\ug}\rk{\Iq-\fpu{2m+1}^\mpi\fpu{2m+1}}\dia{m}^\mpi\fpu{2m+1}}\\
 &\qquad-\rk{\obg-z}\rk{z-\ug}\ba^{-m+1}\rmifsea{m}{z}\dia{m}^\mpi\fpu{2m+1}\biggr\}\dia{m+1}^\mpi,
\intertext{and}
 \rmifsea{m+1}{z}&=\ba\rk*{\ba^{m-1}\rmifswa{m}{z}\fpu{2m+2}+\rk{\obg-z}\rmifsea{m}{z}\ek*{\rk{\Iq-\dia{m}^\mpi\dia{m}}+\dia{m}^\mpi\fpu{2m+1}}}.
\end{align*}
\elem
\bproof
 From \rlem{ab.P1348} we obtain \(\rmif{m+1}= \rmif{m}\cFTiuu{\sau{0}^\FTa{m}}{\su{0}^\FTa{m}}\).
 Using the \tqqa{\tbr{s}} of \(\rmif{m+1}\) and \(\rmif{m}\) as well as \rlem{ab.L1652}, a straightforward calculation completes the proof.
\eproof

 Now we specify \rlem{L1614} for the case of a sequence \(\seqs{m+1}\in\Fggdqu{m+1}\).
 
\bleml{L1504}
 Let \(m\in\NO \) and let \(\seqs{m+1}\in\Fggdqu{m+1}\) with \tfpf{} \(\fpseq{2m+2}\) and \tfdf{} \(\seqdia{m+1}\). 
 Denote by \(\smat{\rmifnw{m}&\rmifne{m}\\\rmifsw{m}&\rmifse{m}}\) and \(\smat{\rmifnw{m+1}&\rmifne{m+1}\\\rmifsw{m+1}&\rmifse{m+1}}\) the \tqqa{\tbr{s}} of \(\rmif{m}\) and \(\rmif{m+1}\), \tresp{}
 For all \(z\in\C\), then \(\rmifnwa{m+1}{z}=\Oqq\) and \(\rmifswa{m+1}{z}=\Oqq\) as well as
\begin{align*}
 \rmifnea{m+1}{z}&=\ba\rk*{\ba^{m-1}\rmifnwa{m}{z}\fpu{2m+2}+\rk{\obg-z}\rmifnea{m}{z}\ek*{\rk{\Iq-\dia{m}^\mpi\dia{m}}+\dia{m}^\mpi\fpu{2m+1}}}
\intertext{and}
 \rmifsea{m+1}{z}&=\ba\rk*{\ba^{m-1}\rmifswa{m}{z}\fpu{2m+2}+\rk{\obg-z}\rmifsea{m}{z}\ek*{\rk{\Iq-\dia{m}^\mpi\dia{m}}+\dia{m}^\mpi\fpu{2m+1}}}.
\end{align*}
\elem
\bproof
 According to \rdefn{D0927}, we have \(\seqs{m+1}\in\Fggqu{m+1}\) and \(\dia{m+1}=\Oqq\).
 The application of \rlem{L1614} completes the proof.
\eproof

 Let \(\seqs{m}\in\Fggqu{m}\) with \tfdf{} \(\seqdia{m}\) and let \(\cU\) be a linear subspace of \(\ran{\dia{m}}\).
 Using \rnota{N1336} to define \(\su{m+1}\defeq\su{m,\cU}\), we determine now the set \(\RFqabsg{m+1}\).
 
\bpropl{L1237}
 Let \(m\in\NO\) and let \(\seqs{m}\in\Fggqu{m}\) with \tfdf{} \(\seqdia{m}\).
 Let \(\cU\) be a linear subspace of \(\ran{\dia{m}}\) and let \(\su{m+1}\defeq\su{m,\cU}\).
 Then \(\RFqabsg{m+1}=\set{\dOSt{m}{\cU}}\), where \(\dOSt{m}{\cU}\) is given via \rnota{N1417}.
\eprop
\bproof
 Denote by \(\fpseq{2m+2}\) the \tfpfa{\(\seqs{m+1}\)}.
 Because of \rpropp{L1253}{L1253.a}, we have \(\seqs{m+1}\in\Fggdqu{m+1}\) and furthermore \(\fpu{2m+1}=\dia{m}^\varsqrt\rk{\Iq-\OPu{\cU}}\dia{m}^\varsqrt\) and \(\fpu{2m+2}=\dia{m}^\varsqrt\OPu{\cU}\dia{m}^\varsqrt\).
 In particular, \(\dia{m}^\mpi\fpu{2m+1}=\dia{m}^\mpi\dia{m}-\dia{m}^\mpi\dia{m}^\varsqrt\OPu{\cU}\dia{m}^\varsqrt\) and, in view of \(\dia{m}^\mpi\dia{m}^\varsqrt=\rk{\dia{m}^\varsqrt}^\mpi\), consequently \(\rk{\Iq-\dia{m}^\mpi\dia{m}}+\dia{m}^\mpi\fpu{2m+1}=\Iq-\rk{\dia{m}^\varsqrt}^\mpi\OPu{\cU}\dia{m}^\varsqrt\). 
 By virtue of \rnota{N1338}, for all \(z\in\Cab\), hence \(\ba^{m-1}\fpu{2m+2}=\dphia{m}{\cU}{z}\) and \(\rk{\obg-z}\ek{\rk{\Iq-\dia{m}^\mpi\dia{m}}+\dia{m}^\mpi\fpu{2m+1}}=\dpsia{m}{\cU}{z}\). 
 Denote by \(\smat{\rmifrnw{m}&\rmifrne{m}\\\rmifrsw{m}&\rmifrse{m}}\) and \(\smat{\rmifrnw{m+1}&\rmifrne{m+1}\\\rmifrsw{m+1}&\rmifrse{m+1}}\) the \tqqa{\tbr{s}} of the restrictions of \(\rmif{m}\) and \(\rmif{m+1}\), \tresp{}, onto \(\Cab\).
 From \rlem{L1504} we can infer then \(\rmifrnea{m+1}{z}=\ba\ek{\rmifrnwa{m}{z}\dphia{m}{\cU}{z}+\rmifrnea{m}{z}\dpsia{m}{\cU}{z}}\) and \(\rmifrsea{m+1}{z}=\ba\ek{\rmifrswa{m}{z}\dphia{m}{\cU}{z}+\rmifrsea{m}{z}\dpsia{m}{\cU}{z}}\) for all \(z\in\Cab\). 
 In view of \(\ba>0\), \rlem{R1419}, and \rnota{N1417}, we see that \(\det\rmifrse{m+1}\) does not identically vanish and that \(\dOSt{m}{\cU}=\rmifrne{m+1}\rmifrse{m+1}^\inv\).
 By virtue of \rdefn{D0927}, we can apply \rthm{F.P.RFabunq} to complete the proof.
\eproof

 Now we recall some facts and notions from \zita{MR3775449}.

\bpropnl{\zitaa{MR3775449}{\cprop{11.4}{208}}}{L1702}
 Let \(m\in\NO\), let \(\seqs{m}\in\Fggqu{m}\), and let \(\su{m+1}\in\set{\umg{m},\omg{m}}\).
 Then \(\seqs{m+1}\in\Fggdqu{m+1}\).
\eprop

\bdefnnl{\tcf{}~\zitaa{MR3775449}{\cdefn{11.5}{209}}}{D0842}
 Let \(m\in\NO\) and let \(\seqs{m}\in\Fggqu{m}\).
 Let the sequence \((\su{j})_{j=m+1}^\infi\) be recursively defined by \(\su{j}\defeq\umg{j-1}\) (\tresp{}\ \(\su{j}\defeq\omg{j-1}\)).
 Then \(\seqsinf\) is called the \emph{lower} (\tresp{}\ \emph{upper}) \tabHdso{\(\seqs{m}\)}.
\edefn

\bpropnl{\tcf{}~\zitaa{MR3775449}{\cprop{}{12.3(b)}{211}}}{P1809}
 Let \(m\in\NO\) and let \(\seqs{m}\in\Fggqu{m}\).
 Denote by \(\seqsl\) and \(\seqsu\) the lower and upper \tabHdso{\(\seqs{m}\)}, \tresp{}
 Then the set \(\Mggouaag{q}{\infi}{\ab}{\seqsl}\) contains exactly one element \(\cdablm{m}\) and the set \(\Mggouaag{q}{\infi}{\ab}{\seqsu}\) contains exactly one element \(\cdabum{m}\).
\eprop

\bdefnnl{\zitaa{MR3775449}{\cdefn{12.4}{211}}}{D1824}
 Let \(m\in\NO\) and let \(\seqs{m}\in\Fggqu{m}\).
 Then the \tnnH{} \tqqa{measure} \(\cdablm{m}\) (\tresp{}\ \(\cdabum{m}\)) is called the \emph{lower} (\tresp{}\ \emph{upper}) \tCDabmo{\(\seqs{m}\)}. 
\edefn

 Now we are interested in the \tSTF{s} of \(\cdablm{m}\) and \(\cdabum{m}\), \tresp{}

\bdefnl{D0923}
 Let \(m\in\NO\) and let \(\seqs{m}\in\Fggqu{m}\).
 Denote by \(\cdablm{m}\) (\tresp{}\ \(\cdabum{m}\)) the lower (\tresp{}\ upper) \tCDabmo{\(\seqs{m}\)}.
 Let \(\lOSt{m}\) be the \tSTFv{\(\cdablm{m}\)} and let \(\uOSt{m}\) be the \tSTFv{\(\cdabum{m}\)}.
 Then we call \(\lOSt{m}\) (\tresp{}\ \(\uOSt{m}\)) the \emph{lower} (\tresp{}\ \emph{upper}) \emph{\tOStv{\(\seqs{m}\)}}.
\edefn

 In \rthm{F.P.PVR}, we obtained a complete description of the set \(\RFqabsg{m}\) of \tSTF{s} of measures belonging to \(\MggqFsg{m}\).
 Now we are interested in the position of \(\lOSt{m}\) and \(\uOSt{m}\) in the set \(\RFqabsg{m}\).
 In particular, we determine the pairs \(\copa{\lphi{m}}{\lpsi{m}}\in\PRFabqa{\su{0}^\FTa{m}}\) and \(\copa{\uphi{m}}{\upsi{m}}\in\PRFabqa{\su{0}^\FTa{m}}\) which correspond to \(\lOSt{m}\) and \(\uOSt{m}\), \tresp{}, according to \rthmp{F.P.PVR}{F.P.PVR.b}.
 It can be expected that these pairs possess certain extremal properties within the set \(\PRFabqa{\su{0}^\FTa{m}}\).
 
 The preceding considerations lead us now quickly to explicit expressions for \(\lOSt{m}\) and \(\uOSt{m}\).
 
\bpropl{R1603}
 Let \(m\in\NO\) and let \(\seqs{m}\in\Fggqu{m}\) with \tfdf{} \(\seqdia{m}\).
 Let \(\cU_0\defeq\set{\Ouu{q}{1}}\) and let \(\cU_1\defeq\ran{\dia{m}}\).
 Then \(\lOSt{m}=\dOSt{m}{\cU_1}\) and \(\uOSt{m}=\dOSt{m}{\cU_0}\) if \(m\) is even, and \(\lOSt{m}=\dOSt{m}{\cU_0}\) and \(\uOSt{m}=\dOSt{m}{\cU_1}\) if \(m\) is odd.
\eprop
\bproof
 Denote by \(\seqsl\) and \(\seqsu\) the lower and upper \tabHdso{\(\seqs{m}\)}, \tresp{}
 From \rprop{P1809} we can conclude \(\cdablm{m}\in\Mggouaag{q}{m+1}{\ab}{\seq{\slu{j}}{j}{0}{m+1}}\) and \(\cdabum{m}\in\Mggouaag{q}{m+1}{\ab}{\seq{\supu{j}}{j}{0}{m+1}}\).
 \rrem{F.R.SF-ST} then shows \(\lOSt{m}\in\RFqabg{\seq{\slu{j}}{j}{0}{m+1}}\) and \(\uOSt{m}\in\RFqabg{\seq{\supu{j}}{j}{0}{m+1}}\).
 
 First consider the case \(m=2n\) with some \(n\in\NO\).
 Because of \rexam{R1403}, we have then \(\su{m,\cU_0}=\omg{m}\) and \(\su{m,\cU_1}=\umg{m}\).
 In view of \rdefn{D0842}, we can thus apply \rprop{L1237} to the sequences \(\seq{\supu{j}}{j}{0}{m+1}\) and \(\seq{\slu{j}}{j}{0}{m+1}\), \tresp{}, to obtain \(\RFqabg{\seq{\supu{j}}{j}{0}{m+1}}=\set{\dOSt{m}{\cU_0}}\) and \(\RFqabg{\seq{\slu{j}}{j}{0}{m+1}}=\set{\dOSt{m}{\cU_1}}\). 
 Thus, \(\uOSt{m}=\dOSt{m}{\cU_0}\) and \(\lOSt{m}=\dOSt{m}{\cU_1}\) follow.
 
 Now consider the case \(m=2n+1\) with some \(n\in\NO\).
 Because of \rexam{R1403}, we have then \(\su{m,\cU_0}=\umg{m}\) and \(\su{m,\cU_1}=\omg{m}\).
 In view of \rdefn{D0842}, we can thus apply \rprop{L1237} to the sequences \(\seq{\slu{j}}{j}{0}{m+1}\) and \(\seq{\supu{j}}{j}{0}{m+1}\), \tresp{}, to obtain \(\RFqabg{\seq{\slu{j}}{j}{0}{m+1}}=\set{\dOSt{m}{\cU_0}}\) and \(\RFqabg{\seq{\supu{j}}{j}{0}{m+1}}=\set{\dOSt{m}{\cU_1}}\). 
 Thus, \(\lOSt{m}=\dOSt{m}{\cU_0}\) and \(\uOSt{m}=\dOSt{m}{\cU_1}\) follow.
\eproof

 Finally, we want to indicate the announced extremal properties of the pairs \(\copa{\lphi{m}}{\lpsi{m}}\) and \(\copa{\uphi{m}}{\upsi{m}}\) from \(\PRFabqa{\dia{m}}\) which correspond to \(\lOSt{m}\) and \(\uOSt{m}\) according to \rthmp{F.P.PVR}{F.P.PVR.b}.

 \breml{R1150}
 If we look back to \rprop{L1237} and \rexam{E1341} and consider the corresponding pairs \(\copa{\lphi{m}}{\lpsi{m}}\) and \(\copa{\uphi{m}}{\upsi{m}}\) belonging to \(\PRFabqa{\dia{m}}\), then it should be mentioned that these pairs consist of \(\Cqq\)\nobreakdash-valued functions in \(\Cab\), which have extremal rank properties.
 Indeed, the function \(\uphi{m}\) satisfies \(\rank\uphi{m}=\rank\dia{m}\), which is the maximal possible rank of a \tqqa{matrix-valued} function \(X\) with \(\ran{X\rk{z}}\subseteq\ran{\dia{m}}\) for all points \(z\in\Cab\) which are points of holomorphy of \(X\), whereas the function \(\lphi{m}\) has rank \(0\) which is clearly the minimal possible rank.
\erem

\section{On the \hSTF{} of the central solution corresponding to a sequence \(\seqs{m}\in\Fggqu{m}\)}\label{S1014}
 At the beginning of this section we state the necessary background information.
 
 Recall that the sequences \(\seq{\umg{j}}{j}{0}{\kappa}\) and \(\seq{\omg{j}}{j}{0}{\kappa}\) were introduced in \rdefn{D0750}.

\bdefnnl{\tcf{}~\zitaa{MR3775449}{\cdefn{10.11}{195}}}{D0756}
 If \(\seqska \) is a sequence of complex \tpqa{matrices}, then we call \(\seq{\mi{j}}{j}{0}{\kappa}\) given by \(\mi{j}\defeq\frac{1}{2}\rk{\umg{j}+\omg{j}}\) the \emph{\tfmfa{\(\seqska \)}}.
\edefn
 
\bdefnnl{\tcf{}~\zitaa{MR3775449}{\cdefn{10.33}{205}}}{F.D.abZ}
 Let \(\seqska\) be a sequence of complex \tpqa{matrices} with \tfmf{} \(\seqmika\).
 Assume \(\kappa\geq1\) and let \(k\in\mn{1}{\kappa}\).
 Then \(\seqska\) is said to be \emph{\tabZo{k}} if \(\su{j}=\mi{j-1}\) for all \(j\in\mn{k}{\kappa}\).
\edefn

\bdefnnl{\zitaa{MR3775449}{\cdefn{11.9}{209}}}{D1527}
 Let \(m\in\NO\) and let \(\seqs{m}\in\Fggqu{m}\).
 Let the sequence \((\su{j})_{j=m+1}^\infi\) be recursively defined by \(\su{j}\defeq\mi{j-1}\), where \(\mi{j-1}\) is given by \rdefn{D0756}.
 Then \(\seqsinf\) is called the \emph{\tabHcso{\(\seqs{m}\)}}.
\edefn

\bpropnl{\zitaa{MR3775449}{\cprop{11.10}{209}}}{R1532}
 Let \(m\in\NO\) and let \(\seqs{m}\in\Fggqu{m}\).
 Then the \tabHcso{\(\seqs{m}\)} is \tFnnd{} and \tabZo{m+1}.
\eprop

\bpropl{P1337}
 Let \(m\in\NO\) and let \(\seqs{m}\in\Fggqu{m}\).
 Denote by \(\seqscent\) the \tabHcso{\(\seqs{m}\)}.
 Then the set \(\MggqFag{\scent}{\infp}\) contains exactly one element \(\centmu{m}\).
\eprop
\bproof
 Combine \rpropss{R1532}{I.P.ab8}.
\eproof

 \rprop{P1337} leads us to the following notion.
 
\bdefnl{D1346}
 Let \(m\in\NO\) and let \(\seqs{m}\in\Fggqu{m}\).
 Then the \tnnH{} \tqqa{measure} \(\centmu{m}\) mentioned in \rprop{P1337} is called the \emph{\tabHcmo{\(\seqs{m}\)}}.
\edefn

\bdefnl{D1555}%
 Let \(m\in\NO\) and let \(\seqs{m}\in\Fggqu{m}\).
 Denote by \(\centmu{m}\) the \tabHcmo{\(\seqs{m}\)}.
 Then the \tSTF{} \(\centS{m}\) of \(\centmu{m}\) is call the \emph{\tabHcfo{\(\seqs{m}\)}}.
\edefn

 Our next goal is now to determine the position of the \tSTFv{\(\centmu{m}\)} within the parametrization of \(\RFqabsg{m}\) obtained in \rthm{F.P.PVR}.
 In order to realize this plan, we continue our investigations in \zitaa{MR4051874}{\csec{10}} where we studied a \tSchur{} type transformation for matrix measures on \(\ab\) which transforms the concrete matrix measure under consideration in accordance with the \tSchur{} type algorithm considered in \rdefn{ab.N1020}, which has to be applied to the corresponding moment sequence.

\bdefnnl{\zitaa{MR4051874}{\cdefn{10.6}{203}}}{F.D.Zmo}
 Let \(\sigma\in\MggqF\) with \tfpmf{} \(\seqmpm{\sigma}\) and let \(k\in\N\).
 Then \(\sigma\) is called \emph{\tZmo{k}} if \(\seqmpm{\sigma}\) is \tabZo{k}.
\edefn

 Against the background of centrality of measures on \(\ab\), we consider now the scalar case, in particular the following object discussed in \zitaa{MR4051874}{\csec{10}}.

\bnotal{N1322}%
 Let \(a,b\in\R\) with \(a<b\) and let \(\centabm{[a,b]}\colon\BsAu{[a,b]}\to[0,\infp)\) be the \emph{arcsine distribution} on \([a,b]\) given by \(\centabma{[a,b]}{B}\defeq\int_Bh\dif\Leb\), where \(\Leb\colon\BsAu{[a,b]}\to[0,\infp)\) is the Lebesgue measure on \([a,b]\) and \(h\colon[a,b]\to[0,\infp)\) is defined by \(h\rk{x}\defeq0\) if \(x\in\set{a,b}\) and by \(h\rk{x}\defeq\ek{\pi\sqrt{\rk{x-a}\rk{b-x}}}^\inv\) if \(x\in(a,b)\).
\enota

 Now we turn our attention to the ordinary and canonical moments of \(\centabm{[a,b]}\).

\bexal{R1524}%
 Let \(a,b\in\R\) with \(a<b\).
 Then \(\centabm{[a,b]}\in\Mggoa{1}{[a,b]}\).
 Denote by \(\seqsinf\) the \tfpmfa{\(\centabm{[a,b]}\)} and by \(\seqciainf\) the \tfocmfa{\(\centabm{[a,b]}\)} via \rdefn{F.D.meacia}.
 Then \(\su{j}=\sum_{k=0}^j\binom{j}{k}\binom{2k}{k}2^{-2k}\rk{b-a}^ka^{j-k}\) for all \(j\in\NO\).
 In particular, \(\centabma{[a,b]}{[a,b]}=1\) and \(\int_{[a,b]}t\centabma{[a,b]}{\dif t}=\rk{a+b}/2\).
 Furthermore, \(\cia{0}=1\) and \(\cia{j}=1/2\) for all \(j\in\N\).
 
 Indeed, the measure \(\mu\defeq\centabm{[0,1]}\) is a probability measure on \([0,1]\) with moments
\(
 \int_{[0,1]}x^k\mu\rk{\dif x}
 =\binom{2k}{k}2^{-2k}
\)
 for all \(k\in\NO\) (see, \teg{}~\zitaa{MR59329}{\cfrml{(25.1)}{75}}) and (classical) canonical moments \(p_k=1/2\) for all \(k\in\N\) (see, \teg{}~\zitaa{MR1468473}{\cexa{1.3.6}{15}}).
 By virtue of \eqref{F.G.epq}, then the \tfmcmfa{\(\mu\)} via \rdefn{F.D.meacia} fulfills \(\mcm{\mu}{0}=1\) and \(\mcm{\mu}{j}=1/2\) for all \(j\in\N\).
 With \(d\defeq b-a\) let \(T\colon[0,1]\to[a,b]\) be defined by \(T\rk{x}=dx+a\).
 Then it is readily checked that \(\centabm{[a,b]}\) is the image measure of \(\mu\) under \(T\).
 Consequently, we can infer \(\centabm{[a,b]}\in\Mggoa{1}{[a,b]}\) and
\begin{multline*}
 \int_{[a,b]}t^j\centabma{[a,b]}{\dif t}
 =\int_{[0,1]}\ek*{T\rk{x}}^j\mu\rk{\dif x}
 =\int_{[0,1]}\ek*{\sum_{k=0}^j\binom{j}{k}d^k x^k a^{j-k}}\mu\rk{\dif x}\\
 =\sum_{k=0}^j\binom{j}{k}d^k\ek*{\int_{[0,1]}x^k\mu\rk{\dif x}}a^{j-k}
 =\sum_{k=0}^j\binom{j}{k}\binom{2k}{k}2^{-2k}d^ka^{j-k}
\end{multline*}
 for all \(j\in\NO\).
 Furthermore, the \tfmcmfa{\(\centabm{[a,b]}\)} coincides, according to~\zitaa{MR3979701}{\cprop{8.12}{40}}, with \(\seqmcm{\mu}\).
\eexa
 
 We reformulate now \rexamp{E1749}{E1749.b} in the language of measures.

\bpropl{L1331}%
 Suppose \(\ba=2\).
 Denote by \(\mu\) the first \tFmTv{\(\centm\)}.
 Then \(\mu=\centm\), \tie{}, the measure \(\centm\) is a fixed point of the \tFmT{ation}.
 In particular, the measure \(\centabm{[-1,1]}\) is a fixed point of the \tFmabT{-1}{1}{ation}.
\eprop
\bproof
 Regarding \rexam{R1524}, denote by \(\seqsinf\) the \tfpmfa{\(\centm\)}.
 According to \rprop{I.P.ab8Fgg}, then \(\seqsinf\in\Fgguuuu{1}{\infi}{\ug}{\obg}\).
 Denote by \(\seqtinf\) the \tFTv{\(\seqsinf\)} and by \(\seqciainf\) the \tfcfa{\(\seqsinf\)} given in \rdefnss{ab.N0940}{D0754}, \tresp{}
 Taking into account \rrem{F.R.FT1} and \rdefn{F.D.SN-M}, we infer that \(\seqtinf\) is the \tfpmfa{\(\mu\)}.
 By virtue of \rdefn{F.D.meacia}, we see that \(\seqciainf\) is the \tfmcmfa{\(\centm\)}.
 From \rexam{R1524} we thus obtain \(\cia{0}=1\) and \(\cia{j}=1/2\) for all \(j\in\N\).
 Using \rexamp{E1749}{E1749.b}, we can conclude then that \(\seqtinf\) coincides with \(\seqsinf\).
 The application of \rprop{I.P.ab8Fgg} hence yields \(\mu=\centm\).
\eproof

 The following result indicates that the notion of \tabZ{ity} of order \(k\) of matrix measures is intimately connected via \tST{} with the scalar probability measure \(\centm\) introduced in \rnota{N1322}.
 More precisely, this property is characterized by the fact that the \tnFmTv{\rk{k-1}}{the} matrix measure under consideration is a \(q\)\nobreakdash-dimensional inflation of \(\centm\), where the corresponding matrix coefficient is a multiple of the \(\rk{k-1}\)\nobreakdash-th matricial interval length.

\bthmnl{\tcf{}~\zitaa{MR4051874}{\cthm{10.9}{204}}}{F.P.algZmo}
 Let \(\sigma\in\MggqF\) and let \(k\in\N\).
 Denote by \(\seqmdm{\sigma}\) the \tfmdmfa{\(\sigma\)} given in \rdefn{F.D.meacia} and by \(\sigma^\FTa{k-1}\) the \tnFmTv{\rk{k-1}}{\(\sigma\)} given via \rdefn{F.D.SN-M}.
 Let \(M\defeq\ba^{k-2}\mdm{\sigma}{k-1}\) and let \(\mu\colon\BsAF\to\Cqq\) be defined by \(\mu\rk{B}\defeq\ek{\centma{B}}M\).
 Then \(\sigma\) is \tZmo{k} if and only if \(\sigma^\FTa{k-1}=\mu\).
\ethm

 A closer look at the proof of \rthm{F.P.algZmo} given in \zita{MR4051874} shows that one of the central points of it is \zitaa{MR4051874}{\cexa{10.8}{204}}, where we took from \zitaa{MR1468473}{\cexa{1.3.6}{15}} the observation that the sequence \(\seq{p_k}{k}{1}{\infi}\) of canonical moments of \(\centabm{[0,1]}\) is the constant sequence with value \(1/2\).
 This result originates in Karlin/Shapley~\zitaa{MR59329}{\csec{25}}.
 For an updated presentation, we refer also to Karlin/Studden~\zitaa{MR0204922}{\cch{4}, \S~4}.
 The essential method used by Karlin and Shapley is a careful study of the geometry of Chebychev polynomials.
 
\bexal{R1431}%
 Let \(a,b\in\R\) with \(a<b\).
 Then \(\centabm{[a,b]}\) is \tZmo{1}.

 Indeed, let \(\sigma\defeq\centabm{[a,b]}\) and let \(M\defeq\rk{b-a}^\inv\mdm{\sigma}{0}\).
 By virtue of \rexam{R1524}, we see \(\sigma\in\Mggoa{1}{[a,b]}\) and \(\sigma\rk{[a,b]}=1\).
 According to \rrem{R0713}, we have \(\sigma^\FTa{0}=\sigma\).
 \rprop{F.R.mdealg} yields \(\sigma^\FTa{0}\rk{[a,b]}=\rk{b-a}^\inv\mdm{\sigma}{0}\).
 Consequently, we can infer \(M=1\) and hence \(\sigma^\FTa{0}\rk{B}=\ek{\centabma{[a,b]}{B}}M\) for all \(B\in\BsAu{[a,b]}\) follows.
 Applying \rthm{F.P.algZmo} shows that \(\centabm{[a,b]}\) is \tZmo{1}.
\eexa

 Now we turn our attention via \tSTF{} to functions belonging  to \(\RFqab\).

\bdefnl{F.D.ZFo}
 Let \(F\in\RFqab\) with \tRabMa{} \(\sigma\) and let \(k\in\N\).
 We call \(F\) \emph{\tZFo{k}} if \(\sigma\) is \tZmo{k}.
\edefn

\bnotal{N1340}%
 Let \(a,b\in\R\) with \(a<b\).
 Then denote by \(\centabf{[a,b]}\) the \txSTv{[a,b]}{\(\centabm{[a,b]}\)}.
\enota

 The following observation is an easy consequence of the construction of the objects under consideration.
 
\breml{R1425}
 Let \(M\in\Cggq\) and let \(\nu\in\Mggoa{1}{\ab}\) with \tSTF{} \(f\).
 Then \(\mu\colon\BsAF\to\Cqq\) defined by \(\mu\rk{B}\defeq\ek{\nu\rk{B}}M\) belongs to \(\MggqF\) and \(G\colon\Cab\to\Cqq\) defined by \(G\rk{z}\defeq f\rk{z}M\) coincides with the \tSTFv{\(\mu\)}.
\erem

 Now we translate \rthm{F.P.algZmo} into the language of functions belonging to \(\RFqab\).

\bthml{F.L.algZFo}
 Let \(F\in\RFqab\) and let \(k\in\N\).
 Denote by \(\seqRabil{F}\) the \tfRabil{\(F\)} given via \rdefn{F.D.funcia} and by \(F^\FTa{k-1}\) the \tnSfTv{\rk{k-1}}{\(F\)} given in \rdefn{F.D.SN-F}.
 Let \(N\defeq\ba^{k-2}\Rabil{F}{k-1}\) and let \(G\colon\Cab\to\Cqq\) be defined by \(G\rk{z}\defeq\centfa{z}N\).
 Then \(F\) is \tZFo{k} if and only if \(F^\FTa{k-1}=G\).
\ethm
\bproof
 Let \(\sigma\) be the \tRabMav{\(F\)}.
 Denote by \(\seqmdm{\sigma}\) the \tfmdmfa{\(\sigma\)} and by \(\sigma^\FTa{k-1}\) the \tnFmTv{\rk{k-1}}{\(\sigma\)}.
 Let \(M\defeq\ba^{k-2}\mdm{\sigma}{k-1}\) and let \(\mu\colon\BsAF\to\Cqq\) be defined by \(\mu\rk{B}\defeq\ek{\centma{B}}M\).
 According to \rdefnss{F.D.funcia}{F.D.meacia}, we  have \(\Rabil{F}{k-1}=\mdm{\sigma}{k-1}\), implying \(N=M\).
 From \rprop{F.R.mdealg} we can infer \(M=\sigma^\FTa{k-1}\rk{\ab}\).
 In particular, \(M\in\Cggq\).
 Taking additionally into account \rexam{R1524} and \rnota{N1340}, the application of \rrem{R1425} now shows that \(\mu\) belongs to \(\MggqF\) and that \(G\) coincides with the \tSTFv{\(\mu\)}.
 Consequently, in view of \rdefn{F.D.SN-F}, we can conclude from \rprop{F.P.STbij} that \(\sigma^\FTa{k-1}=\mu\) if and only if \(F^\FTa{k-1}=G\).
 By virtue of \rdefn{F.D.ZFo}, the application of \rthm{F.P.algZmo} completes the proof.
\eproof

 \rthmss{F.P.algZmo}{F.L.algZFo} contain further results, which indicate the importance of the arcsine distribution introduced in \rnota{N1322}.
 For other topics in which the arcsine distribution plays a significant role, we refer to sum rules for Jacobi matrices, which are compact perturbations of the free Jacobi matrix associated with the arcsine distribution (see Killip/Simon \zita{MR1999923} and Simon \zita{MR2743058}) and free probability and random matrices (see Hiai/Petz \zita{MR1746976})
 
 Now we are going to determine the \tabHcfo{a} sequence \(\seqs{m}\in\Fggqu{m}\) within the description of \(\RFqabsg{m}\) given in \rthm{F.P.PVR}.

\bnotal{N1519}
 Let \(m\in\NO\) and let \(\seqs{m}\in\Fggqu{m}\) with \tfdf{} \(\seqdia{m}\).
 Then let \(\centG{m},\centX{m},\centY{m}\colon\Cab\to\Cqq\) be defined by \(\centGa{m}{z}\defeq\ba^{m-1}\centfa{z}\dia{m}\), \(\centXa{m}{z}\defeq\ba^{m-1}\ek{\rk{\obg-z}\centfa{z}-1}\dia{m}\), and \(\centYa{m}{z}\defeq\rk{\obg-z}\ek{\rk{z-\ug}\centfa{z}+1}\OPu{\ran{\dia{m}}}+\ba\OPu{\nul{\dia{m}}}\).
\enota

\bpropl{R1527}
 Let \(m\in\NO\), let \(\seqs{m}\in\Fggqu{m}\) with \tfdf{} \(\seqdia{m}\).
 Then:
\benui
 \il{R1527.a} \(\centG{m}\in\RFqabg{\seq{\su{j}^\FTa{m}}{j}{0}{0}}\) and \(\copa{\centX{m}}{\centY{m}}\) is the \tFatpv{\su{0}^\FTa{m}}{\(\centG{m}\)}.
 \il{R1527.b} \(\copa{\centX{m}}{\centY{m}}\in\PRFabqa{\su{0}^\FTa{m}}\) and the \tiFaTv{\su{0}^\FTa{m}}{\(\copa{\centX{m}}{\centY{m}}\)} coincides with \(\centG{m}\).
\eenui
\eprop
\bproof
 Setting \(M\defeq\ba^{m-1}\dia{m}\), we have \(\centG{m}=\centf M\), according to \rnota{N1519}.
 From \(\ba>0\) and \rprop{ab.C0929}, we can infer \(M\in\Cggq\).
 \rrem{R1335} furthermore yields \(\su{0}^\FTa{m}=M\).
 Denote by \(\copa{G_1}{G_2}\) the \tFatpv{M}{\(\centG{m}\)}.
 
 \eqref{R1527.a} Taking into account \rexam{R1524} and \rnota{N1340}, the application of \rrem{R1425} shows that \(\mu\colon\BsAF\to\Cqq\) defined by \(\mu\rk{B}\defeq\ek{\centma{B}}M\) belongs to \(\MggqF\) and that \(\centG{m}\) coincides with the \tSTFv{\(\mu\)}.
 In view of \rexam{R1524}, we have \(\mu\rk{\ab}=M\). 
 Consequently, \(\mu\in\MggqFag{s^\FTa{m}}{0}\) follows.
 \rrem{F.R.SF-ST} then shows \(\centG{m}\in\RFqabg{\seq{\su{j}^\FTa{m}}{j}{0}{0}}\).
 Taking into account \(M^\ad=M\), we can infer from \rrem{R1335} moreover \(\OPu{\ran{M^\ad}}=\OPu{\ran{\dia{m}}}\) and \(\OPu{\nul{M}}=\OPu{\nul{\dia{m}}}\).
 Using \rrem{ab.R1052}, in particular \(M^\mpi M=\OPu{\ran{\dia{m}}}\) follows.
 By virtue of \rdefn{F.D.FTF} and \rnota{N1519} we have, for all \(z\in\Cab\), then
\[
 G_1(z)
 =\rk{\obg-z}\centGa{m}{z}-M
 =\rk{\obg-z}\centfa{z}M-M
 =\ek*{\rk{\obg-z}\centfa{z}-1}M
 =\centXa{m}{z}
\]
 and
\[\begin{split}
 G_2(z)
 &=\rk{\obg-z}\ek*{\rk{z-\ug}M^\mpi\centGa{m}{z}+\OPu{\ran{M^\ad}}}+\ba\OPu{\nul{M}}\\
 &=\rk{\obg-z}\ek*{\rk{z-\ug}\centfa{z}M^\mpi M+\OPu{\ran{M^\ad}}}+\ba\OPu{\nul{M}}\\
 &=\rk{\obg-z}\ek*{\rk{z-\ug}\centfa{z}\OPu{\ran{\dia{m}}}+\OPu{\ran{\dia{m}}}}+\ba\OPu{\nul{\dia{m}}}\\
 &=\rk{\obg-z}\ek*{\rk{z-\ug}\centfa{z}+1}\OPu{\ran{\dia{m}}}+\ba\OPu{\nul{\dia{m}}}
 =\centYa{m}{z}.
\end{split}\]

 \eqref{R1527.b} By virtue of \rprop{ab.P1030}, we have \(\seq{\su{j}^\FTa{m}}{j}{0}{0}\in\Fggqu{0}\).
 Taking additionally into account \rpart{R1527.a}, \rlem{ab.P1505} yields \(\copa{\centX{m}}{\centY{m}}\in\PRFabqa{\su{0}^\FTa{m}}\).
 Obviously, \(\OPu{\ran{M}}\centG{m}=\centf\OPu{\ran{M}}M=\centG{m}\).
 Because of \(M\in\Cggq\) and \rpart{R1527.a}, the application of \rlem{F.L.FT1-1} shows that the \tiFaTv{M}{\(\copa{G_1}{G_2}\)} coincides with \(\centG{m}\).
\eproof

\bpropl{L1355}
 Let \(m\in\NO\), let \(\seqs{m}\in\Fggqu{m}\) with \tfdf{} \(\seqdia{m}\), and let \(F\in\RFqabsg{m}\).
 Denote by \(\FTUa{m}{F}{\seqs{m}}\) the \tnFTFv{m}{\(F\)}{\(\seqs{m}\)} and by \(\FTPUa{m}{F}{\seqs{m}}\) the \tnFTPv{m}{\(F\)}{\(\seqs{m}\)} given in \rdefn{F.D.stepFT}.
 Regarding \rnota{N1519} and \rdefn{ab.N0843}, then the following statements are equivalent:
\baeqi{0}
 \il{L1355.i} \(F\) is \tZFo{m+1}.
 \il{L1355.ii} \(\FTUa{m}{F}{\seqs{m}}=\centG{m}\).
 \il{L1355.iii} \(\FTPUa{m}{F}{\seqs{m}}\rpaeq\copa{\centX{m}}{\centY{m}}\).
\eaeqi
\eprop
\bproof
 Denote by \(F^\FTa{m}\) the \tnSfTv{m}{\(F\)}.
 Due to \rprop{P0941}, we have \(\FTUa{m}{F}{\seqs{m}}=F^\FTa{m}\).
 
\baeq{L1355.i}{L1355.ii}
 Denote by \(\sigma\) the \tRabMav{\(F\)} and by \(\seqmpm{\sigma}\) and \(\seqmdm{\sigma}\) the \tfpmf{} and the \tfmdmfa{\(\sigma\)}, \tresp{}
 Then \(\sigma\in\MggqFsg{m}\) and hence \(\mpm{\sigma}{j}=\su{j}\) for all \(j\in\mn{0}{m}\).
 By virtue of \rrem{F.R.diatr}, consequently \(\mdm{\sigma}{m}=\dia{m}\).
 Denote by \(\seqRabil{F}\) the \tfRabil{\(F\)}.
 Taking additionally into account \rdefnss{F.D.funcia}{F.D.meacia}, then \(\Rabil{F}{m}=\mdm{\sigma}{m}=\dia{m}\) follows.
 According to \rnota{N1519}, thus \(\centG{m}=\centf\ba^{m-1}\Rabil{F}{m}\).
 Now, in view of \(\FTUa{m}{F}{\seqs{m}}=F^\FTa{m}\), the application of \rthm{F.L.algZFo} yields the equivalence of~\ref{L1355.i} and~\ref{L1355.ii}.
\eaeq

\bimp{L1355.ii}{L1355.iii}
 From \rpropp{R1527}{R1527.a} we see that \(\copa{\centX{m}}{\centY{m}}\) is the \tFatpv{\su{0}^\FTa{m}}{\(\centG{m}\)}.
 In view of~\ref{L1355.ii} and \rdefn{F.D.stepFT}, then \(\FTPUa{m}{F}{\seqs{m}}=\copa{\centX{m}}{\centY{m}}\) follows.
 In particular,~\ref{L1355.iii} holds true.
\eimp

\bimp{L1355.iii}{L1355.ii}
 Setting \(M\defeq\su{0}^\FTa{m}\), the combination of \rprop{ab.P1030} and \rlem{F.R.Fgg-s} yields \(M\in\Cggq\).
 From \rpropp{R1527}{R1527.b} we see \(\copa{\centX{m}}{\centY{m}}\in\PRFabqa{M}\) and that \(\centG{m}\) is the \tiFaTv{M}{\(\copa{\centX{m}}{\centY{m}}\)}.
 Because of~\ref{L1355.iii} and \rrem{F.R.PMaeq}, in particular \(\FTPUa{m}{F}{\seqs{m}}\in\PRFabqa{M}\) follows.
 Denote by \(G\) the \tiFaTv{M}{\(\FTPUa{m}{F}{\seqs{m}}\)}.
 Taking into account~\ref{L1355.iii}, we can infer from \rcor{F.C.iFwd} that \(G=\centG{m}\).
 Observe that \(\FTPUa{m}{F}{\seqs{m}}\) is the \tFatpv{M}{\(\FTUa{m}{F}{\seqs{m}}\)}, according to \rdefn{F.D.stepFT}.
 \rrem{R0734} yields \(F^\FTa{m}\in\RFqabag{s^\FTa{m}}{0}\).
 Consequently, the \tRabMa{} \(\mu\) of \(F^\FTa{m}\) belongs to \(\MggqFag{s^\FTa{m}}{0}\), \tie{}, \(\mu\rk{\ab}=M\).
 Taking additionally into account \rpropp{ab.P1648L1409}{ab.P1648L1409.a}, hence \(\ran{F^\FTa{m}\rk{z}}=\ran{M}\) for all \(z\in\Cab\) follows.
 In view of \(\FTUa{m}{F}{\seqs{m}}=F^\FTa{m}\), hence \(\OPu{\ran{M}}\FTUa{m}{F}{\seqs{m}}=\FTUa{m}{F}{\seqs{m}}\).
 Thus, we can apply \rlem{F.L.FT1-1} to obtain \(G=\FTUa{m}{F}{\seqs{m}}\).
 Therefore,~\ref{L1355.ii} holds true.
\eimp
\eproof

 Now we are able to determine the \tSTFv{the} \tabHcmo{a} sequence \(\seqs{m}\in\Fggqu{m}\).

\bpropl{C1542}%
 Let \(m\in\NO\) and let \(\seqs{m}\in\Fggqu{m}\).
 Denote by \(\smat{\rmifrnw{m}&\rmifrne{m}\\\rmifrsw{m}&\rmifrse{m}}\) the \tqqa{\tbr{}} of the restriction of \(\rmif{m}\) onto \(\Cab\).
 Then \(\det\rk{\rmifrsw{m}\centX{m}+\rmifrse{m}\centY{m}}\) does not vanish identically in \(\Cab\) and \(\centS{m}=\rk{\rmifrnw{m}\centX{m}+\rmifrne{m}\centY{m}}\rk{\rmifrsw{m}\centX{m}+\rmifrse{m}\centY{m}}^\inv\).
\eprop
\bproof
 \rremp{R1527}{R1527.b} shows  \(\copa{\centX{m}}{\centY{m}}\in\PRFabqa{\su{0}^\FTa{m}}\).
 Consequently, we can apply \rthmp{F.P.PVR}{F.P.PVR.a} to see that \(\det\rk{\rmifrsw{m}\centX{m}+\rmifrse{m}\centY{m}}\) does not vanish identically in \(\Cab\) and that the matrix-valued function \(F\defeq\rk{\rmifrnw{m}\centX{m}+\rmifrne{m}\centY{m}}\rk{\rmifrsw{m}\centX{m}+\rmifrse{m}\centY{m}}^\inv\) belongs to \(\RFqabsg{m}\).
 \rthmp{F.P.PVR}{F.P.PVR.b} then yields \(\copa{\centX{m}}{\centY{m}}\in\rsetcl{\FTPUa{m}{F}{\seqs{m}}}\).
 From \rprop{L1355} we can thus conclude that \(F\) is \tZFo{m+1}.
 Denote by \(\sigma\) the \tRabMav{\(F\)}.
 Then \(\sigma\in\MggqFsg{m}\).
 According to \rdefn{F.D.ZFo}, furthermore \(\sigma\) is \tZmo{m+1}.
 In view of \rdefn{F.D.Zmo}, this means that the \tfpmf{} \(\seqmpm{\sigma}\) associated with \(\sigma\) is \tabZo{m+1}.
 Since \(\mpm{\sigma}{j}=\su{j}\) for all \(j\in\mn{0}{m}\), thus \rdefnss{F.D.abZ}{D1527} show that \(\seqmpm{\sigma}\) coincides with the \tabHcs{} \(\seqscent\) associated with \(\seqs{m}\).
 Consequently, \(\sigma\in\MggqFag{\scent}{\infp}\) follows.
 \rprop{P1337} and \rdefn{D1346} then yield \(\sigma=\centmu{m}\).
 From \rprop{F.P.STbij} we can furthermore conclude that \(F\) is the \tSTFv{\(\sigma\)}.
 In view of \rdefn{D1555}, the proof is complete. 
\eproof

 In our following considerations, we concentrate on the case of a sequence \(\seqs{m}\in\Fgqu{m}\).
 Before doing that we state some elementary preparations in the scalar case, which are of own interest.
 
\bleml{R1705}%
 Let \(a,b\in\R\) with \(a<b\) and let \(z\in\C\setminus[a,b]\).
 Then there exists a unique \(w\in\C\) satisfying \(w^2=\rk{z-a}\rk{z-b}\) and \(\abs{w-z+c}<d\), where \(c\defeq\rk{a+b}/2\) and let \(d\defeq\rk{b-a}/2\).
\elem
\bproof
 First observe that \(\rk{z-a}\rk{z-b}=\rk{z-c}^2-d^2\) and that there exists either a single one or two different solutions \(w\in\C\) satisfying \(w^2=\rk{z-a}\rk{z-b}\).
 We choose a particular solution \(w_0\).
 If \(w_0=0\), then \(z=a\) or \(z=b\), contradicting \(z\notin[a,b]\).
 Thus, we have \(w_0\neq-w_0\) and hence \(w_1\defeq w_0\) and \(w_2\defeq-w_0\) are the only solutions of the equation \(w^2=\rk{z-a}\rk{z-b}\).
 Consequently, \(t_1\defeq w_1-z+c\) and \(t_2\defeq w_2-z+c\) fulfill \(t_1\neq t_2\) and \(\rk{t_{1,2}+z-c}^2=\rk{\pm w_0}^2=\rk{z-c}^2-d^2\).
 Therefore, \(t_1\) and \(t_2\) are the two solutions of the equation \(t^2+2\rk{z-c}t+d^2=0\).
 Hence, \(t_1+t_2=-2\rk{z-c}\) and \(t_1t_2=d^2\).
 In particular, \(\abs{t_1}\cdot\abs{t_2}=d^2\).
 We are now going to show \(\abs{t_1}\neq\abs{t_2}\).
 Assume to the contrary \(\abs{t_1}=\abs{t_2}\).
 In view of \(t_1t_2=d^2\) and \(d>0\), then \(t_2=\ko{t_1}\) and \(\abs{t_1}=d\) follow.
 Using \(t_1+t_2=-2\rk{z-c}\), we can thus infer \(z=c-\re t_1\in\R\) and furthermore \(\abs{z-c}=\abs{\re t_1}\leq\abs{t_1}\).
 Taking additionally into account \(\abs{t_1}=d\), then \(-d\leq z-c\leq d\) follows, contradicting \(z\notin[a,b]\).
 Thus, we have shown \(\abs{t_1}\neq\abs{t_2}\).
 Since \(\abs{t_1}\cdot\abs{t_2}=d^2\), then either \(\abs{t_1}<d\) and \(\abs{t_2}>d\) or \(\abs{t_1}>d\) and \(\abs{t_2}<d\).
 Consequently, exactly one of the two solutions of the equation \(w^2=\rk{z-a}\rk{z-b}\) fulfills \(\abs{w-z+c}<d\).
\eproof

\bleml{R1526}%
\benui
 \il{R1526.a} The function \(\centsnf\) belongs to \(\RFx{1}{\C\setminus\snab}\) with \tRxMa{\snab} \(\centsnm\) and is \tZFo{1}.
 \il{R1526.b} Let \(z\in\Csnab\).
 Then \(\centsnfa{z}\neq0\) and \(w_z\defeq-1/\centsnfa{z}\) is the unique complex number \(w\) satisfying \(w^2=z^2-1\) and \(\abs{w-z}<1\).
\eenui
\elem
\bproof
\eqref{R1526.a} In view of \rexam{R1524} and \rnota{N1340}, we can infer from \rprop{F.P.STbij} that \(\centsnf\) belongs to \(\RFx{1}{\Csnab}\) and that \(\centsnm\) is the \tRxMav{\snab}{\(\centsnf\)}.
 According to \rdefn{F.D.ZFo} and \rexam{R1431}, thus \(\centsnf\) is \tZFo{1}.
 
\eqref{R1526.b} Because of \rpart{R1526.a}, we can apply \rpropp{ab.P1648L1409}{ab.P1648L1409.a} to obtain \(\nul{\centsnfa{z}}=\nul{\centsnma{\snab}}\).
 From \rexam{R1524} we know \(\centsnma{\snab}=1\).
 Hence, \(\centsnfa{z}\neq0\) follows.
 According to \zitaa{MR1468473}{\cpage{125}, especially \cfrml{(4.5.4)}{}}, the remaining assertion of \rpart{R1526.b} holds true.
 (Observe that in \zita{MR1468473}, for probability measures \(\mu\) on \(\snab\), the integral \(S\rk{z,\mu}=\int_{-1}^1\rk{z-x}^\inv\dif\mu\rk{x}=-\STFua{\mu}{z}\) is considered.)
\eproof

\bpropl{R0956}%
\benui
 \il{R0956.a} The function \(\centf\) belongs to \(\RFab{1}\) with \tRabMa{} \(\centm\) and is \tZFo{1}.
 \il{R0956.b} Let \(z\in\Cab\).
 Then \(\centfa{z}\neq0\) and \(w_z\defeq-1/\centfa{z}\) is the unique complex number \(w\) satisfying \(w^2=\rk{z-\ug}\rk{z-\obg}\) and \(\abs{w-z+\rk{\ug+\obg}/2}<\rk{\obg-\ug}/2\).
\eenui
\eprop
\bproof
 By the same reasoning as in the proof of \rlem{R1526}, we can conclude that \(\centf\) belongs to \(\RFx{1}{\Cab}\) and is \tZFo{1}, that \(\centm\) is the \tRabMav{\(\centf\)}, and that \(\centfa{z}\neq0\).
 Let \(c\defeq\rk{\ug+\obg}/2\) and let \(d\defeq\rk{\obg-\ug}/2\).
 Let \(T\colon\snab\to\ab\) be defined by \(T\rk{x}=dx+c\).
 Then it is readily checked that \(\centm\) is the image measure of \(\centsnm\) under \(T\).
 Observe that \(d>0\) and that \(z\notin\ab\) implies \(\zeta\notin\snab\) for \(\zeta\defeq\rk{z-c}/d\).
 By virtue of \rnota{N1340} and \rdefn{F.D.ST}, we thus can infer
\[\begin{split}
 \centfa{z}
 &=\int_\ab\frac{1}{t-z}\centma{\dif t}
 =\int_{\snab}\frac{1}{T\rk{x}-z}\centsnma{\dif x}
 =\int_{\snab}\frac{1}{dx-z+c}\centsnma{\dif x}\\
 &=\frac{1}{d}\int_{\snab}\frac{1}{x-\rk{z-c}/d}\centsnma{\dif x}
 =\frac{1}{d}\centsnfa{\zeta}.
\end{split}\]
 In view of \rlemp{R1526}{R1526.b}, consequently \(w_z=d\omega_\zeta\), where \(\omega_\zeta\) is the unique complex number \(\omega\) satisfying \(\omega^2=\zeta^2-1\) and \(\abs{\omega-\zeta}<1\).
 Hence, \(w_z^2=d^2\rk{\zeta^2-1}=\rk{z-c}^2-d^2=\rk{z-\ug}\rk{z-\obg}\) and \(\abs{w_z-z+\rk{\ug+\obg}/2}=\abs{d\omega_\zeta-z+c}=d\abs{\omega_\zeta-\zeta}<d=\rk{\obg-\ug}/2\).
 By virtue of \rlem{R1705}, the proof is complete.
\eproof

\bleml{L1335}%
 Suppose \(\ba=2\).
 Denote by \(f\) the first \tSfTv{\(\centf\)}.
 Then \(f=\centf\), \tie{}, the function \(\centf\) is a fixed point of the \tSfT{ation}.
 In particular, the function \(\centabf{[-1,1]}\) is a fixed point of the \tSfabT{-1}{1}{ation}.
\elem
\bproof
 In view of \rexam{R1524} and \rnota{N1340}, we can infer from \rprop{F.P.STbij} that \(\centf\) belongs to \(\RFx{1}{\Cab}\) and that \(\centm\) is the \tRabMav{\(\centf\)}.
 Denote by \(\mu\) the first \tFmTv{\(\centm\)}.
 According to \rdefn{F.D.SN-F}, then \(f\) is the \tSTFv{\(\mu\)}.
 Since \rprop{L1331} yields \(\mu=\centm\), hence \(f\) is the \tSTFv{\(\centm\)}.
 Regarding \rnota{N1340}, the proof is complete.
\eproof

\bpropl{R1101}%
 Let \(m\in\NO\) and let \(\seqs{m}\in\Fgqu{m}\) with \tfdf{} \(\seqdia{m}\).
 Let \(G_1,G_2\colon\Cab\to\Cqq\) be defined by \(G_1\rk{z}\defeq\ba^{m-1}\centfa{z}\dia{m}\) and \(G_2\rk{z}\defeq\Iq\).
 For all \(z\in\Cab\), then
\begin{align}\label{R1101.1}
 \centXa{m}{z}&=\ba^{m-1}\ek*{\rk{\obg-z}\centfa{z}-1}\dia{m},&
 \centYa{m}{z}&=\rk{\obg-z}\ek*{\rk{z-\ug}\centfa{z}+1}\Iq.
\end{align}
 Regarding \rdefn{ab.N0843}, furthermore \(\copa{G_1}{G_2}\in\PRFabq\) and \(\copa{\centX{m}}{\centY{m}}\rpaeq\copa{G_1}{G_2}\).
\eprop
\bproof
 \rrem{P1338} yields \(\det\dia{m}\neq0\).
 Hence, \(\ran{\dia{m}}=\Cq\) and \(\nul{\dia{m}}=\set{\Ouu{q}{1}}\), which imply \(\OPu{\ran{\dia{m}}}=\Iq\) and \(\OPu{\nul{\dia{m}}}=\Oqq\).
 For all \(z\in\Cab\), now \eqref{R1101.1} follows immediately from \rnota{N1519}.
 By virtue of \rpropp{R1527}{R1527.b} and \rnota{ab.N1534}, we get \(\copa{\centX{m}}{\centY{m}}\in\PRFabq\).
 Let \(f\colon\Cab\to\C\) be defined by \(f\rk{z}\defeq\rk{z-\ug}\centfa{z}+1\).
 From \rpropp{R0956}{R0956.a} we see that \(\centf\) belongs to \(\RFqab\) and that \(\centm\) is the \tRabMav{\(\centf\)}.
 \rexam{R1524} furthermore shows \(\centma{\ab}=1\) and \(\int_\ab t\centabma{\ab}{\dif t}=\rk{\ug+\obg}/2\).
 Taking into account \rnota{ab.N1537}, we can thus infer from \rprop{ab.L1401} that \(f\) belongs to \(\RFab{1}\) and that the \tRabMa{} \(\sigma\) of \(f\) fulfills \(\sigma\rk{\ab}=\int_\ab\rk{t-\ug}\centma{\dif t}=\rk{\ug+\obg}/2-\ug=\rk{\obg-\ug}/2\neq0\).
 Consider now an arbitrary \(z\in\Cab\).
 Using \rpropp{ab.P1648L1409}{ab.P1648L1409.a} we obtain then \(\nul{f\rk{z}}=\nul{\sigma\rk{\ab}}=\set{0}\) and hence \(f\rk{z}\neq0\).
 Regarding \(\centYa{m}{z}=\rk{\obg-z}f\rk{z}\Iq\), then \(\det\centYa{m}{z}\neq0\) follows.
 From \rpropp{R0956}{R0956.b} we see \(\centfa{z}\neq0\) and that \(w_z\defeq-1/\centfa{z}\) satisfies \(w_z^2=\rk{z-\ug}\rk{z-\obg}\).
 Hence, \(\ek{\rk{\obg-z}+w_z}w_z=\rk{\obg-z}\ek{w_z-\rk{z-\ug}}\) and, in view of \(w_z\neq0\), therefore \(\rk{\obg-z}+w_z=\rk{\obg-z}\ek{1-\rk{z-\ug}/w_z}=\rk{\obg-z}f\rk{z}\).
 Multiplication by \(\centfa{z}\) yields \(\rk{\obg-z}\centfa{z}-1=\rk{\obg-z}f\rk{z}\centfa{z}\).
 Consequently, we get the identity \(\ek{\rk{\obg-z}\centfa{z}-1}/\ek{\rk{\obg-z}f\rk{z}}=\centfa{z}\), from which we can conclude \(\ek{\centXa{m}{z}}\ek{\centYa{m}{z}}^\inv=G_1\rk{z}\). 
 The application of \rlem{ab.P1801} completes the proof.
\eproof

\bpropl{C1647}%
 Let \(m\in\NO\) and let \(\seqs{m}\in\Fgqu{m}\).
 Let \(\smat{\rmifrnw{m}&\rmifrne{m}\\\rmifrsw{m}&\rmifrse{m}}\) be the \tqqa{\tbr{}} of the restriction of \(\rmif{m}\) onto \(\Cab\).
 Then \(\det\rk{\ba^{m-1}\centf\rmifrsw{m}\dia{m}+\rmifrse{m}}\) does not vanish identically in \(\Cab\) and \(\centS{m}=\rk{\ba^{m-1}\centf\rmifrnw{m}\dia{m}+\rmifrne{m}}\rk{\ba^{m-1}\centf\rmifrsw{m}\dia{m}+\rmifrse{m}}^\inv\).
\eprop
\bproof
 From \rprop{C1542} we see that \(\det\rk{\rmifrsw{m}\centX{m}+\rmifrse{m}\centY{m}}\) does not vanish identically in \(\Cab\) and that \(\centS{m}=\rk{\rmifrnw{m}\centX{m}+\rmifrne{m}\centY{m}}\rk{\rmifrsw{m}\centX{m}+\rmifrse{m}\centY{m}}^\inv\).
 \rremp{R1527}{R1527.b} shows  \(\copa{\centX{m}}{\centY{m}}\in\PRFabqa{\su{0}^\FTa{m}}\).
 Let \(G_1,G_1\colon\Cab\to\Cqq\) be defined by \(G_1\rk{z}\defeq\ba^{m-1}\centfa{z}\dia{m}\) and \(G_2\rk{z}\defeq\Iq\).
 According to \rprop{R1101}, then \(\copa{G_1}{G_2}\in\PRFabq\) and \(\copa{\centX{m}}{\centY{m}}\rpaeq\copa{G_1}{G_2}\).
 The application of \rdefn{ab.N0843} completes the proof.
\eproof

\appendix
\section{Some facts from matrix theory}\label{A0828}
 
 This appendix contains a summary of results from matrix theory, which are used in this paper.
 What concerns results on the Moore--Penrose inverse \(A^\mpi\) of a complex matrix \(A\), we refer, \teg{}, to~\zitaa{MR1152328}{\csec{1}}.

\breml{A.R.AB=BA}
 Let \(m,n\in\N\) and let \(A_1,A_2,\dotsc,A_m\) and \(B_1,B_2,\dotsc,B_n\) be complex \tpqa{matrices}.
 If \(M\in\Cqp\) is such that \(A_jMB_k=B_kMA_j\) holds true for all \(j\in\mn{1}{m}\) and all \(k\in\mn{1}{n}\), then \(\rk{\sum_{j=1}^m\eta_jA_j}M\rk{\sum_{k=1}^n\theta_kB_k}=\rk{\sum_{k=1}^n\theta_kB_k}M\rk{\sum_{j=1}^m\eta_jA_j}\) for all complex numbers \(\eta_1,\eta_2,\dotsc,\eta_m\) and \(\theta_1,\theta_2,\dotsc,\theta_n\).
\erem

\breml{A.R.XRIX}
\benui
 \il{A.R.XRIX.a} If \(Z\in\Cqq\) and \(\eta\in\C\), then \(\re\rk{\eta Z}=\re\rk{\eta}\re\rk{Z}-\im\rk{\eta}\im\rk{Z}\) and \(\im\rk{\eta Z}=\re\rk{\eta}\im\rk{Z}+\im\rk{\eta}\re\rk{Z}\).
 \il{A.R.XRIX.b} If \(Z\in\Cqq\) and \(X\in\Cqp\), then \(\re\rk{X^\ad ZX}=X^\ad\re\rk{Z}X\) and \(\im\rk{X^\ad ZX}=X^\ad\im\rk{Z}X\).
\eenui
\erem

\breml{A.R.RNT}%
 If \(A\in\Cpq\), then \(\dim\ran{A}+\dim\nul{A}=q\).
\erem

\breml{ab.R1842x}
 If \(A\in\Cpq\), then \(\ran{A}=\ran{AA^\ad}\) and \(\nul{A}=\nul{A^\ad A}\).
\erem

\breml{ab.L1044}
 Let \(A\in\Coo{p}{r}\), let \(B\in\Coo{p}{s}\), and let \(C\in\Coo{q}{r}\).
 In view of \rrem{ab.R1842x}, then \(\ran{A}+\ran{B}=\ran{\mat{A,B}}=\ran{AA^\ad+BB^\ad}\) and \(\nul{A}\cap\nul{C}=\nul{\tmatp{A}{C}}=\nul{A^\ad A+C^\ad C}\).
\erem

\breml{A.R.rs+}
 Let \(n\in\N\) and let \(A_1,A_2,\dotsc,A_n\in\Cpq\).
 For all \(\eta_1,\eta_2,\dotsc,\eta_n\in\C\), then \(\ran{\sum_{j=1}^n\eta_jA_j}\subseteq\sum_{j=1}^n\ran{A_j}\) and \(\bigcap_{j=1}^n\nul{A_j}\subseteq\nul{\sum_{j=1}^n\eta_jA_j}\).
\erem

\breml{A.R.rgcolreg}
 If \(A\in\Cpq\) has rank \(q\) and \(B\in\Coo{q}{s}\) has  rank \(s\), then \(\rank\rk{AB}=s\).
\erem

\breml{A.R.rnLAR}
 Let \(L\in\Cpp\) and \(R\in\Cqq\) be both invertible.
 Let \(A\in\Cpq\) and let \(X\defeq LAR\).
 Then \(\ran{X}=L\ran{A}\) and \(\nul{X}=R^\inv\nul{A}\).
\erem

 We think that the following result is well-know.
 However, we did not succeed in finding a reference.
 
\bleml{ab.L1217}
 Let \(A\in\Cpq\) and let \(R\in\Coo{q}{r}\).
 Then \(\ran{AR}=\ran{A}\) if and only if \(\nul{A}+\ran{R}=\Cq\).
\elem
\bproof
 Observe that \(U\defeq\nul{A}+\ran{R}\) is a linear subspace of the \(\C\)\nobreakdash-vector space \(\Cq\).
 Let \(\phi\colon U\to\Cp\) be defined by \(\phi(x)\defeq Ax\).
 Then \(\phi\) is linear with \(\ker\phi=U\cap\nul{A}=\nul{A}\) and \(\phi(U)=A\ran{R}=\ran{AR}\).
 Regarding \(\dim\ker\phi+\dim\phi(U)=\dim U\), then
\(
 \dim U
 =\dim\nul{A}+\dim\ran{AR}
 =q-\rank A+\rank\rk{AR}
\)
 follows.
 Hence, \(U=\Cq\) if and only if \(\rank\rk{AR}=\rank A\).
 Because of \(\ran{AR}\subseteq\ran{A}\), the latter is equivalent to \(\ran{AR}=\ran{A}\).
\eproof

\breml{ab.R1842*}
 If \(A\in\Cpq\), then \(\ran{A^\ad}=\ek{\nul{A}}^\orth\) and  \(\nul{A^\ad}=\ek{\ran{A}}^\orth\).
\erem

 We write \(\OPu{\mathcal{U}}\) for the transformation matrix corresponding to the orthogonal projection onto a linear subspace \(\mathcal{U}\) of the unitary space \(\Cp\) with respect to the standard basis.

\breml{R.P}
 Let \(\mathcal{U}\) be a linear subspace of \(\Cp\).
 Then \(\OPu{\mathcal{U}}\) is the unique complex \tppa{matrix} satisfying \(\OPu{\mathcal{U}}^2=\OPu{\mathcal{U}}\), \(\OPu{\mathcal{U}}^\ad=\OPu{\mathcal{U}}\), and \(\ran{\OPu{\mathcal{U}}}=\mathcal{U}\).
 Furthermore, \(\nul{\OPu{\mathcal{U}}}=\mathcal{U}^\orth\) 
 and \(\OPu{\mathcal{U}}+\OPu{\mathcal{U}^\orth}=\Ip\).
\erem

\breml{A.R.P=UU*}
 If \(\mathcal{U}\) is a linear subspace of the unitary space \(\Cp\) with dimension \(d\defeq\dim\mathcal{U}\geq1\) and some orthonormal basis \(u_1,u_2,\dotsc,u_d\), then \(\OPu{\mathcal{U}}=UU^\ad\), where \(U\defeq\mat{u_1,u_2,\dotsc,u_d}\).
\erem

\breml{A.R.A-1}
 If \(A\in\Cqq\) fulfills \(\det A\neq0\), then \(A^\mpi=A^\inv\).
\erem

\breml{A.R.A++*}
 If \(A\in\Cpq\), then \(\rk{A^\mpi}^\mpi=A\), \(\rk{A^\ad}^\mpi=\rk{A^\mpi}^\ad\), \(\ran{A^\mpi}=\ran{A^\ad}\), and \(\nul{A^\mpi}=\nul{A^\ad}\).
\erem

\breml{A.R.l*A}
 Let \(\eta\in\C\) and let \(A\in\Cpq\).
 Then \(\rk{\eta A}^\mpi=\eta^\mpi A^\mpi\).
\erem

\breml{A.R.A+>}
 If \(A\in\Cggq\), then \(A^\mpi\in\Cggq\) and \(\rk{A^\mpi}^\varsqrt=\rk{A^\varsqrt}^\mpi\).
\erem

 Regarding \eqref{ps}, we easily obtain with \rrem{A.R.l*A}:

\breml{ab.R1132}
 Let \(\eta\in\C\) and let \(A,B\in\Cpq\).
 Then \((\eta A)\ps(\eta B)=\eta(A\ps B)\).
\erem

\breml{ab.R1052}
 If \(A\in\Cpq\), then \(\OPu{\ran{A}}=AA^\mpi\), \(\OPu{\nul{A}}=\Iq-A^\mpi A\), \(\OPu{\ran{A^\ad}}=A^\mpi A\), and \(\OPu{\nul{A^\ad}}=\Ip-AA^\mpi\).
\erem

\bleml{141.B463}
 Let \(A\in\Cpq\) and let \(B\in\Cqq \) be such that \(\ran{B}\subseteq\ran{A^\ad}\subseteq\ran{B^\ad}\).
 Let \(\eta\in\C\setminus\set{0}\).
 Then the matrix \(B+\eta\OPu{\nul{A}}\) is invertible and \(B^\mpi=(B+\eta\OPu{\nul{A}})^\inv-\eta^\inv\OPu{\nul{A}}\).
\elem
\begin{proof}
 First observe that \(\ran{B}=\ran{A^\ad}=\ran{B^\ad}\) follows from the assumption, since \(\rank(B^\ad)=\rank B\) holds true.
 In view of \rrem{ab.R1052}, we thus obtain \(BB^\mpi=A^\mpi A=B^\mpi B\) and \(\OPu{\nul{A}}=\Iq-A^\mpi A\).
 Taking additionally into account \rrem{R.P}, we infer then
\[\begin{split}
 (B+\eta\OPu{\nul{A}})(B^\mpi+\eta^\inv\OPu{\nul{A}})
 &=BB^\mpi+\eta^\inv B\OPu{\nul{A}}+\eta\OPu{\nul{A}}B^\mpi+\OPu{\nul{A}}^2\\
 &=A^\mpi A+\eta^\inv B(\Iq-A^\mpi A)+\eta(\Iq-A^\mpi A)B^\mpi+\OPu{\nul{A}}\\
 &=\Iq+\eta^\inv B(\Iq-B^\mpi B)+\eta(\Iq-B^\mpi B)B^\mpi
 =\Iq.\qedhere
\end{split}\]
\end{proof}

\breml{R.AA+B.A}
 Let \(A\in\Cpq\) and let \(B\in\Coo{p}{m}\).
 Then \(\ran{B}\subseteq\ran{A}\) if and only if \(AA^\mpi B=B\).
\erem

\breml{R.AA+B.B}
 Let \(A\in\Cpq\) and let \(C\in\Coo{n}{q}\).
 Then \(\nul{A}\subseteq\nul{C}\) if and only if \(CA^\mpi A=C\).
\erem

 The combination of \rremsss{R.AA+B.A}{R.AA+B.B}{A.R.AB=BA} yields:

\breml{A.R.AM+B}
 Let \(A\in\Cpq\) and \(M\in\Cqp\) be such that \(\ran{A}\subseteq\ran{M}\) and \(\nul{M}\subseteq\nul{A}\).
 For all \(\eta_1,\eta_2,\theta_1,\theta_2\in\C\), then \(\rk{\eta_1A+\eta_2M}M^\mpi\rk{\theta_1A+\theta_2M}=\rk{\theta_1A+\theta_2M}M^\mpi\rk{\eta_1A+\eta_2M}\).
\erem

 Regarding \rrem{R.AA+B.A}, we can easily conclude from \rlem{ab.L1217}:

\breml{ab.L1418}
 Let \(A\in\Cpq\) and let \(B\in\Coo{p}{r}\).
 Then \(\ran{B}=\ran{A}\) if and only if there exists a complex \taaa{q}{r}{matrix} \(R\) fulfilling \(\nul{A}+\ran{R}=\Cq\) and \(B=AR\).
\erem

\breml{A.R.kK}
 The set \(\CHq\) is an \(\R\)\nobreakdash-vector space and \(\Cggq\) is a convex cone in \(\CHq\).
\erem

\breml{A.R.XAX}
 Let \(A\in\CHq\) and let \(X\in\Cqp\).
 Then \(X^\ad AX\in\CHp\).
 If \(A\in\Cggq\), then \(X^\ad AX\in\Cggp\).
\erem

 We now state a well-known characterization of \tnnH{} block matrices in terms of the Schur complement \eqref{E/A}:

\blemnl{\tcf{}~\zitaa{MR1152328}{\clemss{1.1.9}{18}{1.1.10}{19}}}{L.AEP}
 Let \(\tmat{A&B\\C&D}\) be the \tbr{} of a complex \taaa{(p+q)}{(p+q)}{matrix} \(M\) with \tppa{block} \(A\).
 Then \(M\) is \tnnH{} if and only if \(A\) and \(M\schca A\defeq D-CA^\mpi B\) are both \tnnH{} and furthermore \(\ran{B}\subseteq\ran{A}\) and \(C=B^\ad\) are fulfilled.
 In this case, \(\normSs{B}\leq\normS{A}\normS{D}\).
\elem

\blemnl{\tcf{}~\zitaa{MR3979701}{\clem{A.13}{2168}}}{A.R.rA<rB}
 Let \(A,B\in\CHq\) with \(\Oqq\lleq A\lleq B\).
 Then \(\ran{A}\subseteq\ran{B}\) and \(\nul{B}\subseteq\nul{A}\).
 Furthermore, \(\Oqq\lleq\OPu{\ran{A}}B^\mpi\OPu{\ran{A}}\lleq A^\mpi\).
\elem

 We continue with some observations on the classes of \tEP{matrices} and \tAD{} matrices given in \rdefn{D.EPAD}. 
 From \rrem{A.R.rs+} we easily see:
 
\breml{R1421}
 If \(A\in\CEPq\), then \(\ran{\re A}\subseteq\ran{A}\) and \(\ran{\im A}\subseteq\ran{A}\).
\erem

\breml{R1749}
 If \(A\in\CADq\), then \(\eta A\in\CADq\) for all \(\eta\in\C\).
\erem

\breml{R.AD<}
 If \(A\in\CADq\), then \(\nul{A^\ad}=\nul{A}\) and \(\ran{A^\ad}=\ran{A}\).
\erem

\breml{L.AD<EP}
 Taking into account \rrem{R.AD<}, one can easily check that \(\setaca{M\in\Cqq}{\eta M\in\Cggq}\subseteq\CADq\subseteq\CEPq\) for all \(\eta\in\C\).
\erem

\bleml{L1423}
 Let \(A\in\Cqq\) satisfy \(\im A\in\Cggq\).
 Then \(A\in\CEPq\).
\elem
\bproof
 Let \(x\in\nul{A}\).
 Then \(x^\ad\im\rk{A}x=0\).
 Since, by virtue of \rrem{L.AD<EP}, we have \(\im A\in\CADq\), then \(\im\rk{A}x=\NM\) follows.
 Consequently, \(A^\ad x=Ax=\NM\).
 Hence, \(\nul{A}\subseteq\nul{A^\ad}\), implying \(\ran{A}=\ran{A^\ad}\), \tie{}, \(A\in\CEPq\).
\eproof

\bleml{L1305}
 Let \(A\in\Cqq\) satisfy \(\im A\in\Cggq\) and \(\rank\rk{\im A}=\rank A\).
 Then \(A\in\CADq\) and \(\ran{\im A}=\ran{A}\).
\elem
\bproof
 \rlem{L1423} yields \(
 \ran{A^\ad}
 =\ran{A}
\).
 In view of \rrem{R1421} and \(\rank\rk{\im A}=\rank A\), we get
\(
 \ran{\im A}
 =\ran{A}
\).
 Consider an arbitrary \(x\in\Cq\) with \(x^\ad Ax=0\).
 Then \(x^\ad A^\ad x=\ko{x^\ad Ax}=0\).
 Consequently, \(x^\ad\im\rk{A}x=0\).
 Since \rrem{L.AD<EP} yields \(\im A\in\CADq\), we have \(\im\rk{A}x=\Ouu{q}{1}\).
 Because of \(\nul{\im A}=\ek{\ran{\rk{\im A}^\ad}}^\orth=\ek{\ran{\im A}}^\orth=\ek{\ran{A}}^\orth=\ek{\ran{A^\ad}}^\orth=\nul{A}\), we get \(Ax=\Ouu{q}{1}\) and, thus, \(A\in\CADq\).
\eproof

 A complex \tpqa{matrix} \(K\) is said to be \emph{contractive} if \(\normS{K}\leq1\).

\breml{DFK.L1-1-12}%
 Let \(K\in\Cpq\).
 Then the matrix \(K\) is contractive if and only if \(\Iq-K^\ad K\) is \tnnH{}.
\erem

\breml{A.L.KAB0D}
 Let \(\tmat{A&B\\C&D}\) be the \tbr{} of a contractive complex \taaa{(p+q)}{(p+q)}{matrix} \(K\) with \taaa{p}{p}{block} \(A\).
 Suppose that \(A\) or \(D\) is unitary.
 In view of \rrem{DFK.L1-1-12} then one can easily see that \(B=\Ouu{p}{q}\) and \(C=\Ouu{q}{p}\).
\erem

 A complex square matrix \(J\) is called a \emph{signature matrix} if it satisfies \(J^\ad=J\) and \(J^2=\EM\).
 In this paper, we focus on the particular signature matrices
\begin{align}\label{J}
 \Jimq
 &\defeq
 \bMat
  \Oqq&\iu\Iq\\
  -\iu\Iq&\Oqq
 \eMat&
&\text{and}&
 \Jabspq
 &\defeq
 \bMat
  -\Ip&\Opq\\
  \Oqp&\Iq
 \eMat.
\end{align}
 What concerns several aspects of the so-called \(J\)\nobreakdash-Theory for arbitrary signature matrices \(J\), we refer to~\zitaa{MR1152328}{\csecsspp{1.3}{1.4}{26}{44}}.

\breml{ab.R1543}
 If \(P,Q\in\Cqp\), then \(\tmatp{P}{Q}^\ad\Jimq\tmatp{P}{Q}=2\im\rk{Q^\ad P}\) and \(\tmatp{P}{Q}^\ad\Jabsqq\tmatp{P}{Q}=Q^\ad Q-P^\ad P\).
\erem

\breml{ab.R1031}
 Let \(A\in\Cpp\) and let \(D\in\Cqq\).
 Then
\(
 \smat{
  A&\Opq\\
  \Oqp& D
 }\Jabspq
 =
 \smat{
  -A&\Opq\\
  \Oqp& D
 }
 =\Jabspq
 \smat{
  A&\Opq\\
  \Oqp& D
 }
\).
 Furthermore, if \(p=q\), then
\(
 \smat{A&\Oqq\\
  \Oqp& D
 }\Jimq
 =
 \smat{
  \Oqq&\iu A\\
  -\iu D& \Oqq
 }
 =\Jimq
 \smat{
  D&\Oqq\\
  \Oqp& A
 }
\).
 In particular, \(\Jimq\Jabsqq=-\Jabsqq\Jimq\).
\erem

\section{Integration with respect to \tnnH{} measures}\label{B.se.IT}

 Consider a measurable space \((\dbX,\sigA)\) consisting of a non-empty set \(\dbX\) and a \(\sigma\)\nobreakdash-algebra \(\sigA\) on \(\dbX\).
 A mapping \(\mu\colon\sigA\to\Cggq\) is called \emph{\tnnH{} \tqqa{measure} on \((\dbX,\sigA)\)} if it is \(\sigma\)\nobreakdash-additive, \tie{}\ \(\mu(\bigcup_{n=1}^\infi X_n)=\sum_{n=1}^\infi\mu(X_n)\) holds true for all sequences \(\seq{X_n}{n}{1}{\infi}\) of pairwise disjoint sets from \(\sigA\). 
 We are using the integration theory of pairs of matrix-valued functions with respect to \tnnH{} measures, developed independently by Kats~\zita{MR0080280} and Rosenberg~\zitas{MR0163346,MR0378068,MR0436302}(\tcf{}~\zitaa{MR1152328}{\csec{2.2}}):
 
 First consider an arbitrary ordinary measure \(\nu\) on a measurable space \((\dbX,\sigA)\).
 A measurable function \(F\colon\dbX\to\Cpq\) is said to be \emph{integrable with respect to \(\nu\)} if \(F=\mat{f_{jk}}_{\substack{j=1,\dotsc,p\\k=1,\dotsc,q}}\) belongs to \(\ek{\Loa{1}{\nu}}^\xx{p}{q}\), \tie{}, all entries \(f_{jk}\) belong to the class \(\Loa{1}{\nu}\) of functions  \(f\colon\dbX\to\C\), which are integrable with respect to \(\nu\).
 In this case, let \(\int_X F\dif\nu\defeq\mat{\int_Xf_{jk}\dif\nu}_{\substack{j=1,\dotsc,p\\k=1,\dotsc,q}}\) for all \(X\in\sigA\).
 
 Throughout this part of the appendix, let an arbitrary \tnnH{} \tqqa{measure} \(\mu=\matauuo{\mu_{jk}}{j,k}{1}{q}\) be given.
 Then all the entries \(\mu_{jk}\) are complex-valued measures and all the diagonal entries \(\mu_{jj}\) are ordinary measures with values in \([0,\infp)\).
 Consequently, the \emph{trace measure} \(\trm\defeq\Spur\mu=\sum_{j=1}^q\mu_{jj}\)  of \(\mu\) is an ordinary measure with values in \([0,\infp)\).
 For each \(X\in\sigA\), from \(\trm (X)=0\) necessarily \(\mu(X)=\Oqq\) follows.
 Consequently, there exist \(\trm\)\nobreakdash-a.\,e.\ uniquely determined Radon--Nikodym derivatives \(\dif\mu_{jk}/\dif\trm\).
 The \(\trm\)\nobreakdash-a.\,e.\ uniquely determined measurable mapping \(\mu_\trm'\defeq\matauuo{\dif\mu_{jk}/\dif\trm}{j,k}{1}{q}\) is called \emph{trace derivative of \(\mu\)} and satisfies \(\mu(X)=\int_X\mu_\trm'\dif\trm\) for all \(X\in\sigA\) and \(\Oqq\lleq\mu_\trm'(x)\lleq\Iq\) for \(\trm\)\nobreakdash-a.\,a.\ \(x\in\dbX\).
 
  A measurable function \(f\colon\dbX\to\C\) is said to be \emph{integrable with respect to \(\mu\)} if \(\int_\dbX\abs{f}\dif\nu_{jk}<\infp\) holds true for all \(j,k\in\mn{1}{q}\), where \(\nu_{jk}\) denotes the variation of the complex measure \(\mu_{jk}\).
 In this case, let
\(
 \int_Xf\dif\mu
 \defeq\matauuo{\int_Xf\dif\mu_{jk}}{j,k}{1}{q}
\)
 for all \(X\in\sigA\).
 Denote by \(\LnnH{\mu}\) the set of all such functions \(f\), which are integrable with respect to \(\mu\) in this sense.
 
\breml{142.R1141}
 Let \(u\in\Cq\) and let \(\nu\defeq u^\ad\mu u\).
 Then \(\nu\) is a bounded measure on \((\dbX,\sigA)\), which is absolutely continuous with respect to \(\trm\).
 For all \(f\in\LnnH{\mu}\), furthermore \(\int_\dbX\abs{f}\dif\nu<\infp\) and \(\int_\dbX f\dif\nu=u^\ad\rk{\int_\dbX f\dif\mu}u\).
\erem

\breml{B.R.int+}
 The mapping defined on the \(\C\)\nobreakdash-vector space \(\LnnH{\mu}\) by \(f\mapsto\int_\dbX f\dif\mu\) is \(\C\)\nobreakdash-linear.
\erem

\breml{B.R.int*}
 If \(f\in\LnnH{\mu}\), then \(\ko f\in\LnnH{\mu}\) and \(\int_\dbX\ko f\dif\mu=\rk{\int_\dbX f\dif\mu}^\ad\).
\erem

 An ordered pair \((\Phi,\Psi)\) consisting of measurable functions \(\Phi\colon\dbX\to\Cpq\) and \(\Psi\colon\dbX\to\Coo{r}{q}\) is said to be \emph{\tli{\mu}} if the matrix-valued function \(\Phi\mu_\trm'\Psi^\ad\) belongs to \(\ek{\Loa{1}{\trm}}^\xx{p}{r}\).
 In this case, let \(\int_X\Phi\dif\mu\Psi^\ad\defeq\int_X\Phi\mu_\trm'\Psi^\ad\dif\trm\) for all \(X\in\sigA\). 
 In particular, denote by \(\pqLsqa{\mu}\) the set of all measurable functions \(\Phi\colon\dbX\to\Cpq\) for which the pair \((\Phi,\Phi)\) is \tli{\mu}.

\breml{DFK2.2.3}
 If \(\Phi\in\pqLsqa{\mu}\) and \(\Psi\in\aaLsqa{r}{q}{\mu}\), then \((\Phi,\Psi)\) is \tli{\mu}.
\erem

\breml{B.R.L2HR}
 If \(\Phi,\Psi\in\pqLsqa{\mu}\), then \(\int_X\Phi\dif\mu\Psi^\ad=\rk{\int_X\Psi\dif\mu\Phi^\ad}^\ad\) and \(\int_X\Phi\dif\mu\Phi^\ad\in\Cggp\) for all \(X\in\sigA\).
\erem

\bleml{ab.L1459}
 Let \(\Phi\in\pqLsqa{\mu}\) and let \(\Psi\in\rqLsqa{\mu}\).
 For all \(X\in\sigA\), then
\begin{align*}
 \Ran{\int_X\Phi\dif\mu\Psi^\ad}&\subseteq\Ran{\int_X\Phi\dif\mu\Phi^\ad},&
 \Nul{\int_X\Phi\dif\mu\Phi^\ad}&\subseteq\Nul{\int_X\Psi\dif\mu\Phi^\ad},
\end{align*}
 and
\[
 \rk*{\int_X\Psi\dif\mu\Phi^\ad}\rk*{\int_X\Phi\dif\mu\Phi^\ad}^\mpi\rk*{\int_X\Phi\dif\mu\Psi^\ad}
 \lleq\int_X\Psi\dif\mu\Psi^\ad.
\]
\elem
\bproof
 Consider an arbitrary \(X\in\sigA\).
 First observe that \(\Theta\defeq\tmat{\Phi\\ \Psi}\) is a measurable function satisfying
\[
 \Theta\mu_\trm'\Theta^\ad
 =
 \bMat
  \Phi\\
  \Psi
 \eMat
 \mu_\trm'
 \mat{\Phi^\ad,\Psi^\ad}
 =
 \bMat
  \Phi\mu_\trm'\Phi^\ad&\Phi\mu_\trm\Psi^\ad\\
  \Psi\mu_\trm\Phi^\ad&\Psi\mu_\trm'\Psi^\ad
 \eMat.
\]
 By assumption, we have \(\Phi\mu_\trm'\Phi^\ad\in\ek{\Loa{1}{\trm}}^\x{p}\) and \(\Psi\mu_\trm'\Psi^\ad\in\ek{\Loa{1}{\trm}}^\x{r}\).
 Due to \rrem{DFK2.2.3}, the pairs \((\Phi,\Psi)\) and \((\Psi,\Phi)\) are both \tli{\mu}, \tie{}, \(\Phi\mu_\trm\Psi^\ad\in\ek{\Loa{1}{\trm}}^\xx{p}{r}\) and \(\Psi\mu_\trm\Phi^\ad\in\ek{\Loa{1}{\trm}}^\xx{r}{p}\) are true.
 Thus, we infer \(\Theta\mu_\trm'\Theta\in\ek{\Loa{1}{\trm}}^\x{(p+r)}\), \tie{}, \(\Theta\in\aaLsqa{(p+r)}{q}{\mu}\).
 Setting \(A\defeq\int_X\Phi\dif\mu\Phi^\ad\), \(B\defeq\int_X\Phi\dif\mu\Psi^\ad\), \(C\defeq\int_X\Psi\dif\mu\Phi^\ad\), and \(D\defeq\int_X\Psi\dif\mu\Psi^\ad\), we get
\[%
 \bMat
  A&B\\
  C&D
 \eMat
 =\int_X
 \bMat
  \Phi\mu_\trm'\Phi^\ad&\Phi\mu_\trm'\Psi^\ad\\
  \Psi\mu_\trm'\Phi^\ad&\Psi\mu_\trm'\Psi^\ad
 \eMat\dif\trm
 =\int_X\Theta\mu_\trm'\Theta^\ad\dif\trm
 =\int_X\Theta\dif\mu\Theta^\ad.
\]
 Because of \rrem{B.R.L2HR}, the matrix \(\int_X\Theta\dif\mu\Theta^\ad\) is \tnnH{}.
 Using \rlem{L.AEP}, we can conclude then \(\ran{B}\subseteq\ran{A}\), \(C=B^\ad\), and that the matrices \(A\), \(D\), and \(E\schca A=D-CA^\mpi B\) are \tnnH{}.
 By virtue of \rrem{A.R.A++*}, thus \(D\) and \(CA^\mpi B\) are \tH{} matrices, which satisfy \(CA^\mpi B\lleq D\).
 Taking into account \rrem{ab.R1842*}, we can furthermore infer \(\nul{A}\subseteq\nul{C}\).
\eproof

\breml{DFK2.2.5}
 Let \(f\colon\dbX\to\C\) be a measurable function and let \(\Phi,\Psi\colon\dbX\to\Cqq\) be defined by \(\Phi(x)\defeq f(x)\Iq\) and \(\Psi(x)\defeq\Iq\).
 Then \(f\in\LnnH{\mu}\) if and only if \((\Phi,\Psi)\) is \tli{\mu}.
 In this case, \(f\mu_\trm'\in\ek{\Loa{1}{\trm}}^\xx{q}{q}\) and
\(
 \int_Xf\dif\mu
 =\int_X\rk{f\mu_\trm'}\dif\trm
 =\int_X\Phi\dif\mu\Psi^\ad
\)
 for all \(X\in\sigA\), where \(\mu_\trm'\) denotes the trace derivative of \(\mu\).
\erem

 Using \rrem{DFK2.2.5} and \rlem{ab.L1459}, the following result can be easily verified:

\blemnl{\tcf{}~\zitaa{MR2988005}{\clem{B.2(b)}{1787}}}{L0813}
 If \(f\in\LnnH{\mu}\), then \(\ran{\int_{\dbX} f\dif\mu}\subseteq\ran{\mu (\dbX)}\) and \(\nul{\mu (\dbX)}\subseteq\nul{\int_{\dbX} f\dif\mu}\).
\elem

\bleml{L0818}
 Let \(g\colon\dbX\to\R\) satisfy \(g\in\LnnH{\mu}\) and \(\mu(\set{g\leq0})=\Oqq\).
 Then \(\ran{\int_{\dbX} g\dif\mu}=\ran{\mu (\dbX)}\) and \(\nul{\int_{\dbX} g\dif\mu}=\nul{\mu (\dbX)}\).
\elem
\bproof
 Consider an arbitrary \(u\in\nul{\int_{\dbX}g\dif\mu}\).
 In view of \rrem{142.R1141}, then \(\nu\defeq u^\ad\mu u\) is a bounded measure with \(\int_{\dbX}g\dif\nu=u^\ad\rk{\int_{\dbX}g\dif\mu}u=0\) and \(\nu(\set{g\leq0})=0\).
 Thus, \(\nu(\dbX)=0\).
 Therefore, \(u^\ad\mu(\dbX)u=0\).
 Since, by virtue of \(\mu(\dbX)\in\Cggq\) and \rrem{L.AD<EP}, we have \(\mu(\dbX)\in\CADq\), we conclude then \(u\in\nul{\mu(\dbX)}\).
 So we have \(\nul{\int_{\dbX} g\dif\mu}\subseteq\nul{\mu (\dbX)}\).
 Observe that, due to \rlem{L0813}, furthermore \(\nul{\mu (\dbX)}\subseteq\nul{\int_{\dbX} g\dif\mu}\) holds true.
 Hence, \(\nul{\int_{\dbX} g\dif\mu}=\nul{\mu (\dbX)}\).
 Using \rremss{B.R.int*}{ab.R1842*}, we can then easily infer \(\ran{\int_{\dbX} g\dif\mu}=\ran{\mu (\dbX)}\).
\eproof

 We will particularly apply the following result on integrable functions \(f\) satisfying \(\re f>0\) or \(\im f>0\).

\bleml{L0824}
 Let \(f\in\LnnH{\mu}\), let \(\eta,\theta\in\R\), and let \(g\defeq\eta\re f+\theta\im f\).
 Suppose that \(\mu(\set{g\leq0})=\Oqq\).
 Then \(g\in\LnnH{\mu}\) with \(\ran{\rk{\int_{\dbX} f\dif\mu}^\ad}=\ran{\int_{\dbX} f\dif\mu}\),
\begin{align}
 \Ran{\int_{\dbX} g\dif\mu}&=\Ran{\mu (\dbX)},&&&\Nul{\int_{\dbX} g\dif\mu}&=\Nul{\mu (\dbX)},\label{L0824.A}\\
 \Ran{\int_{\dbX} f\dif\mu}&=\Ran{\mu (\dbX)},&&\text{and}&\Nul{\int_{\dbX} f\dif\mu}&=\Nul{\mu (\dbX)}.\label{L0824.B}
\end{align}
\elem
\bproof
 Because of \rremss{B.R.int*}{B.R.int+}, the real-valued function \(g\) belongs to \(\LnnH{\mu}\).
 Hence, \rlem{L0818} yields \eqref{L0824.A}.
 Let \(h\in\set{f,\ko f}\).
 By virtue of \rrem{B.R.int*}, we have \(h\in\LnnH{\mu}\).
 Consider an arbitrary \(u\in\nul{\int_{\dbX}h\dif\mu}\).
 In view of \rrem{142.R1141}, then \(\nu\defeq u^\ad\mu u\) is a bounded measure and \(\int_{\dbX}h\dif\nu=u^\ad\rk{\int_{\dbX}h\dif\mu}u=0\).
 \rremss{B.R.int*}{B.R.int+} provide us \(\int_{\dbX}\re h\dif\nu=0\) and \(\int_{\dbX}\im h\dif\nu=0\).
 Regarding \(\re h=\re f\) as well as \(\im h=\im f\) if \(h=f\) and \(\im h=-\im f\) if \(h=\ko f\), we thus obtain \(\int_{\dbX}\re f\dif\nu=0\) and \(\int_{\dbX}\im f\dif\nu=0\).
 \rrem{B.R.int+} implies \(\int_{\dbX}g\dif\nu=0\).
 The assumption \(\nu(\set{g\leq0})=0\) yields then \(\nu(\dbX)=0\).
 So we have \(u^\ad\mu(\dbX)u=0\).
 Since, by virtue of \(\mu(\dbX)\in\Cggq\) and \rrem{L.AD<EP}, we have \(\mu(\dbX)\in\CADq\), hence \(u\in\nul{\mu(\dbX)}\).
 Observe that, due to \rlem{L0813}, furthermore \(\nul{\mu (\dbX)}\subseteq\nul{\int_{\dbX} h\dif\mu}\) holds true.
 Consequently, \(\nul{\int_{\dbX} h\dif\mu}=\nul{\mu (\dbX)}\).
 Using \rremss{B.R.int*}{ab.R1842*}, we easily conclude \(\ko h\in\LnnH{\mu}\) and \(\ran{\int_{\dbX}\ko h\dif\mu}=\ran{\rk{\int_{\dbX}h\dif\mu}^\ad}=\ran{\mu (\dbX)}\).
 Choosing the appropriate \(h\), we thus obtain \eqref{L0824.B} and \(\ran{\rk{\int_{\dbX}f\dif\mu}^\ad}=\ran{\mu (\dbX)}=\ran{\int_{\dbX}f\dif\mu}\).
\eproof

 We end this section with a matricial version of  Lebesgue's \emph{dominated convergence theorem}:

\bpropnl{\zitaa{MR3644521}{\cprop{A.6}{348}}}{142.P1102}
 Let \(f,f_1,f_2,\dotsc\colon\dbX\to\C\) be measurable functions and let \(g\in\LnnH{\mu}\) be such that \(\lim_{n\to\infi}f_n(x)=f(x)\) for \(\trm\)\nobreakdash-almost all \(x\in\dbX\) and \(\abs{f_n(x)}\leq\abs{g(x)}\) for all \(n\in\N\) and \(\trm\)\nobreakdash-almost all \(x\in\dbX\) hold true.
 Then the functions \(f,f_1,f_2,\dotsc\) belong to \(\LnnH{\mu}\) and \(\lim_{n\to\infi}\int_\dbX f_n\dif\mu=\int_\dbX f\dif\mu\).
\eprop

\section{The \tST{} of \tnnH{} measures}\label{B.s1.HN}
 In this section, we consider a non-empty \emph{closed} subset \(\Omega\) of \(\R\) and a \tnnH{} \tqqa{measure} \(\sigma\) on \((\Omega,\BsAu{\Omega})\).
 So \(\Omega\) is also a closed subset of \(\C\), whereas \(\dom\defeq\C\setminus\Omega\) is an open subset of \(\C\).
 Observe that, for each \(t\in\Omega\), the function \(h_t\colon\dom\to\C\) defined by \(h_t (z)\defeq1/(t-z)\) is holomorphic.
 Consider an arbitrary \(z\in\dom\) and let \(d_z\defeq\inf_{x\in\Omega}\abs{x-z}\).
 Then \(d_z>0\) and \(1/\abs{t-z}\leq1/d_z\) for all \(t\in\Omega\).
 Consequently, the function \(g_z \colon\Omega\to\C\) defined by \(g_z (t)\defeq1/(t-z)\) belongs to \(\LnnH{\sigma}\).
 For each closed disk \(K\subseteq\dom\), we have with \(d_K\defeq\inf_{(x,w)\in\Omega\times K}\abs{x-w}\), furthermore \(d_K>0\) and \(1/\abs{t-z}\leq1/d_K\).
 By means of that, one can check that the matrix-valued function \(\STOu{\sigma}\colon\dom\to\Cqq\) defined by
\beql{B.G.STO}
 \STOua{\sigma}{z}
 \defeq\int_\Omega\frac{1}{t-z}\sigma(\dif t)
\eeq
 is holomorphic (see, \teg{}~\zitaa{MR2257838}{\cSatz{5.8}{147}, \cKap{IV}}).
 The mapping \(\sigma\mapsto\STOu{\sigma}\) is called \emph{\tSTion{}}.
 Accordingly, the function \(\STOu{\sigma}\) itself is called \emph{\tSTv{\(\sigma\)}}.
 For \(\Omega=\R\), the restriction of \(\STOu{\sigma}\) onto \(\uhp\) is exactly the \tSTH{} \(\STHu{\sigma}\) of \(\sigma\) introduced in \rdefn{H.D.ST}.
 Thus, we obtain:

\bleml{B.L.STOR0}
 Denote by \(F\) the restriction of \(\STOu{\sigma}\) onto \(\uhp\).
 Then \(F\in\RFuq{0}\) and the \tSpMa{} \(\sigmaF\) of \(F\) fulfills \(\sigmaFa{\R\setminus\Omega}=\Oqq\) and \(\sigmaFa{B}=\sigma(B)\) for all \(B\in\BsAu{\Omega}\).
\elem
\bproof
 Let \(\chi\colon\BsAR\to\Cqq\) be defined by \(\chi(B)\defeq\sigma(B\cap\Omega)\).
 Then \(\chi\) is a \tnnH{} \tqqa{measure} on \((\R,\BsAR)\) satisfying \(F(z)=\int_\R\rk{t-z}^\inv\chi(\dif t)\) for all \(z\in\uhp\).
 By virtue of \rprop{H.L.ST-R0}, then \(F\in\RFuq{0}\) and \(\chi=\sigmaF\) follow.
\eproof

 For \(\Omega=\ab\) the function \(\STOu{\sigma}\) coincides with the \tSTF{} \(\STFu{\sigma}\) of \(\sigma\) introduced in \rdefn{F.D.ST}.
 By virtue of \eqref{B.G.STO}, we easily see:

\breml{B.R.STORJ}
 For all \(z\in\C\setminus\Omega\), we have
\begin{align*}
 \re \STOua{\sigma}{z}&=\int_\Omega\frac{t-\re z}{\abs{t-z}^2}\sigma(\dif t)&
&\text{and}&
 \im \STOua{\sigma}{z}&=\int_\Omega\frac{\im z}{\abs{t-z}^2}\sigma(\dif t).
\end{align*}
\erem

\bleml{B.L.STOJ}
 For all \(z\in\C\setminus\R\), the matrix \(\frac{1}{\im z}\im \STOua{\sigma}{z}\) is \tnnH{}.
 Furthermore, \(\im \STOua{\sigma}{x}=\Oqq\) for all \(x\in\R\setminus\Omega\).
 Moreover, for all \(w\in\C\) with \(\re w<\inf\Omega\), the matrix \(\re \STOua{\sigma}{w}\) is \tnnH{} and, for all \(w\in\C\) with \(\re w>\sup\Omega\), the matrix \(-\re \STOua{\sigma}{w}\) is \tnnH{} .
\elem
\bproof
 Except for \(\im \STOua{\sigma}{x}=\Oqq\) for all \(x\in\R\setminus\Omega\), the assertions are a consequence of \rrem{B.R.STORJ}.
 Consider now an arbitrary \(x\in\R\setminus\Omega\).
 The matrix-valued function \(\STOu{\sigma}\) is holomorphic.
 Hence, the two sequences \(\seq{\pm\im \STOua{\sigma}{x\pm\iu/n}}{n}{1}{\infi}\) converge to \(\pm\im \STOua{\sigma}{x}\), \tresp{}
 As already mentioned, we have \(\frac{1}{\im z}\im \STOua{\sigma}{z}\in\Cggq\) for all \(z\in\C\setminus\R\).
 In particular, the sequences \(\seq{\pm\im \STOua{\sigma}{x\pm\iu/n}}{n}{1}{\infi}\) both consist of \tnnH{} matrices.
 Consequently, we obtain \(\pm\im \STOua{\sigma}{x}\in\Cggq\) for their limits, implying \(\im \STOua{\sigma}{x}=\Oqq\).
\eproof

 Using \rprop{142.P1102}, one can easily deduce the following representations for \(\sigma(\Omega)\) via limits along the imaginary axis (\tcf{}~\zitaa{MR3644521}{\clem{A.8(c)}{349}}):
 
\bleml{Mae05.P2-1-8}
 The following equations hold true:
\begin{align*}
 \lim_{y\to\infp}y\re \STOua{\sigma}{\iu y}&=\Oqq,&
 \lim_{y\to\infp}y\im \STOua{\sigma}{\iu y}&=\sigma(\Omega),& 
&\text{and}&
 \lim_{y\to\infp}\iu y\STOua{\sigma}{\iu y}&=-\sigma(\Omega).
\end{align*}
\elem

 The following lemma in particular shows that the matrix-valued function \(\STOu{\sigma}\) has constant column and null space on \(\C\setminus [\inf \Omega, \sup \Omega]\).

\blemnl{\tcf{}~\zitaa{MR3644521}{\clem{A.8(b)}{349}}}{B.L.rn}
 For all \(z\in\C\setminus [\inf \Omega, \sup \Omega]\), the equations
\begin{align}\label{B.L.rn.A}
 \ran{\STOua{\sigma}{z}}
 &=\ran{\sigma(\Omega)}&
&\text{and}&
 \nul{\STOua{\sigma}{z}}
 &=\nul{\sigma(\Omega)}
\end{align}
 hold true.
 Furthermore,
\begin{align}\label{B.L.rn.B}
 \ran{\im \STOua{\sigma}{z}}
 &=\ran{\sigma(\Omega)}&
&\text{and}&
 \nul{\im \STOua{\sigma}{z}}
 &=\nul{\sigma(\Omega)}
\end{align}
 are valid for all \(z\in\C\setminus\R\) and
\begin{align}\label{B.L.rn.C}
 \ran{\re \STOua{\sigma}{w}}
 &=\ran{\sigma(\Omega)}&
&\text{and}&
 \nul{\re \STOua{\sigma}{w}}
 &=\nul{\sigma(\Omega)}
\end{align}
 are fulfilled for all \(w\in\C\) with \(\re w<\inf\Omega\) or \(\re w>\sup\Omega\).
\elem
\bproof
 For each \(\zeta\in\C\setminus\Omega\), let \(g_\zeta\colon\Omega\to\C\) be defined by \(g_\zeta(t)\defeq1/(t-\zeta)\).
 As already mentioned at the beginning of this section, we have \(g_\zeta\in\LnnH{\sigma}\) for each \(\zeta\in\C\setminus\Omega\).
 Because of \rremss{B.R.int*}{B.R.int+}, then the functions \(\re g_\zeta\) and \(\im g_\zeta\) both belong to \(\LnnH{\sigma}\) and, in view of \(\STOua{\sigma}{\zeta}=\int_\Omega g_\zeta\dif\sigma\), fulfill \(\re \STOua{\sigma}{\zeta}=\int_\Omega\re g_\zeta\dif\sigma\) and \(\im \STOua{\sigma}{\zeta}=\int_\Omega\im g_\zeta\dif\sigma\) for each \(\zeta\in\C\setminus\Omega\).
 Furthermore, we have \(\re g_\zeta(t)=(t-\re\zeta)/\abs{t-\zeta}^2\) and \(\im g_\zeta(t)=\im\zeta/\abs{t-\zeta}^2\) for each \(\zeta\in\C\setminus\Omega\). 
 Consider now an arbitrary \(z\in\C\setminus\R\).
 In the case \(\im z>0\), we have \(\im g_z>0\) and thus the application of \rlem{L0824} with \(f=g_z\), \(\eta=0\), and \(\theta=1\) yields then \eqref{B.L.rn.B} and \eqref{B.L.rn.A}.
 If \(\im z<0\), then \(-\im g_z>0\) and we can infer \eqref{B.L.rn.B} and \eqref{B.L.rn.A} in a similar way with \(\theta=-1\). 
 Now let \(w\in\C\) with \(\re w<\inf\Omega\) or \(\re w>\sup\Omega\).
 In the case \(\re w<\inf\Omega\), we have \(\re g_w>0\) and thus the application of \rlem{L0824} with \(f=g_w\), \(\eta=1\), and \(\theta=0\) yields then \eqref{B.L.rn.C} and \eqref{B.L.rn.A}.
 If \(\re w>\sup\Omega\), then \(-\re g_w>0\) and we can infer \eqref{B.L.rn.C} and \eqref{B.L.rn.A} in a similar way with \(\eta=-1\).
\eproof

 Recall that a matrix \(A\in\Cqq\) is called \emph{\tEP{matrix}} if \(\ran{A^\ad}=\ran{A}\) and that \(A\) is said to be \emph{\tAD} if each \(x\in\Cq\) with \(x^\ad Ax=0\) necessarily fulfills \(Ax=\Ouu{q}{1}\).
 The corresponding classes are denoted by \(\CEPq\) and \(\CADq\) (\tcf{}~\rdefn{D.EPAD}). 
 In view of \eqref{B.G.STO}, we obviously have \([\STOua{\sigma}{z}]^\ad=\STOua{\sigma}{\ko z}\) for all \(z\in\C\setminus\Omega\).
 Consequently, from \rlem{B.L.rn} we can infer \(\ran{[\STOua{\sigma}{z}]^\ad}=\ran{\STOua{\sigma}{z}}\), \tie{}, \(\STOua{\sigma}{z}\in\CEPq\) for all \(z\in\C\setminus[\inf\Omega,\sup\Omega]\).
 Regarding \rrem{L.AD<EP}, the values of \(\STOu{\sigma}\) satisfy a stronger condition:

\bleml{L1746}
 For all \(z\in\C\setminus[\inf\Omega,\sup\Omega]\), we have \(\STOua{\sigma}{z}\in\CADq\).
\elem
\bproof
 First consider an arbitrary \(z\in\C\setminus\R\).
 Let \(\eta\defeq\frac{1}{\im z}\).
 \rlem{B.L.STOJ} yields \(\im\rk{\eta\ek{\STOua{\sigma}{z}}}\in\Cggq\).
 From \rlem{B.L.rn} we can infer \(\rank\im \STOua{\sigma}{z}=\rank \STOua{\sigma}{z}\), implying \(\rank\im\rk{\eta\ek{\STOua{\sigma}{z}}}=\rank\rk{\eta\ek{\STOua{\sigma}{z}}}\).
 Thus, we can apply \rlem{L1305} to obtain \(\eta\ek{\STOua{\sigma}{z}}\in\CADq\).
 By virtue of \rrem{R1749}, then \(\STOua{\sigma}{z}\in\CADq\) follows. 
 Now consider an arbitrary \(x\in\R\setminus[\inf\Omega,\sup\Omega]\).
 Assume \(x>\sup\Omega\).
 From \rlem{B.L.STOJ} we can infer then \(-\STOua{\sigma}{x}=-\re \STOua{\sigma}{x}\in\Cggq\).
 Because of \rrem{L.AD<EP}, hence \(-\STOua{\sigma}{x}\in\CADq\).
 By virtue of \rrem{R1749}, thus \(\STOua{\sigma}{x}\in\CADq\).
 If \(x<\inf\Omega\), then \rlem{B.L.STOJ} yields \(\STOua{\sigma}{x}=\re \STOua{\sigma}{x}\in\Cggq\), implying \(\STOua{\sigma}{x}\in\CADq\) by \rrem{L.AD<EP}.
\eproof

\bleml{ab.L1400}
 For all \(z\in\C\setminus\R\), the following matrix inequalities hold true:
\begin{align*}
 \ek*{\STOua{\sigma}{z}}^\ad\ek*{\frac{1}{\im z}\im \STOua{\sigma}{z}}^\mpi\ek*{\STOua{\sigma}{z}}&\lleq\sigma(\Omega)&
&\text{and}&
 \ek*{\STOua{\sigma}{z}}\ek*{\sigma(\Omega)}^\mpi\ek*{\STOua{\sigma}{z}}^\ad&\lleq\frac{1}{\im z}\im \STOua{\sigma}{z}.
\end{align*}
\elem
\bproof
 Consider an arbitrary \(z\in\C\setminus\R\).
 Let \(g_z\colon\Omega\to\C\) be defined by \(g_z(t)\defeq1/(t-z)\) and let \(\Lambda,\Xi\colon\Omega\to\Cqq\) be defined by \(\Lambda(t)\defeq\abs{g_z(t)}^2\Iq\) and \(\Xi(t)\defeq\Iq\).
 As already mentioned at the beginning of this section, we have \(g_z\in\LnnH{\sigma}\).
 Because of \(\im g_z/\im z=\abs{g_z}^2\), we can conclude with \rremss{B.R.int+}{B.R.int*} then \(\abs{g_z}^2\in\LnnH{\sigma}\).
 According to \rrem{DFK2.2.5}, hence the pair \((\Lambda,\Xi)\) is \tli{\sigma}.
 Consequently, the matrix-valued function \(\Lambda\sigma_\trm'\Xi^\ad\) belongs to \(\ek{\Loa{1}{\trm}}^\x{q}\), where \(\trm\) denotes the trace measure of \(\sigma\) and \(\sigma_\trm'\) is the trace derivative of \(\sigma\).
 Let \(\Theta\colon\Omega\to\Cqq\) be defined by \(\Theta(t)\defeq g_z(t)\Iq\).
 Then \(\Theta\) is measurable and fulfills \(\Lambda\sigma_\trm'\Xi^\ad=\Theta\sigma_\trm'\Theta^\ad\).
 Thus, the pair \((\Theta,\Theta)\) is \tli{\sigma}, \tie{}, \(\Theta\in\qqLsqa{\sigma}\).
 By virtue of \rremss{DFK2.2.5}{B.R.STORJ}, we have
\[
 \int_\Omega\Theta\dif\sigma\Theta^\ad
 =\int_\Omega\Theta\sigma_\trm'\Theta^\ad\dif\trm
 =\int_\Omega\Lambda\sigma_\trm'\Xi^\ad\dif\trm
 =\int_\Omega\Lambda\dif\sigma\Xi^\ad
 =\int_\Omega\abs{g_z}^2\dif\sigma
 =\frac{1}{\im z}\im \STOua{\sigma}{z}.
\]
 Furthermore, \(\Xi\) belongs to \(\qqLsqa{\sigma}\) and fulfills
\[
 \int_\Omega\Xi\dif\sigma\Xi^\ad
 =\int_\Omega\Xi\sigma_\trm'\Xi^\ad\dif\trm
 =\int_\Omega\sigma_\trm'\dif\trm
 =\sigma(\Omega).
\]
 In view of \rrem{DFK2.2.3}, then the pairs \((\Theta,\Xi)\) and \((\Xi,\Theta)\) are both \tli{\sigma}.
 Using \rremss{DFK2.2.5}{B.R.L2HR}, we then can infer
\begin{align*}
 \int_\Omega\Theta\dif\sigma\Xi^\ad
 &=\int_\Omega g_z\dif\sigma
 =\STOua{\sigma}{z}&
&\text{and}&
 \int_\Omega\Xi\dif\sigma\Theta^\ad
 &=\rk*{\int_\Omega\Theta\dif\sigma\Xi^\ad}^\ad
 =\ek*{\STOua{\sigma}{z}}^\ad.
\end{align*}
 The proof is completed by applying \rlem{ab.L1459} once with \(\Phi=\Theta\) and \(\Psi=\Xi\) and a second time with \(\Phi=\Xi\) and \(\Psi=\Theta\).
\eproof

 In combination with \rlem{ab.L1400}, the following result reveals a certain minimality of the \tnnH{} matrix \(\sigma(\Omega)\) with respect to the L\"owner partial order:
 
\bleml{B.L.s<B}
 Let \(A\in\Cggq\) be such that \(\ran{\sigma(\Omega)}\subseteq\ran{A}\) and \(\ek{\STOua{\sigma}{z}}^\ad A^\mpi\ek{\STOua{\sigma}{z}}\lleq\frac{1}{\im z}\im \STOua{\sigma}{z} \) for all \(z\in\uhp\).
 Then \(\sigma(\Omega)\lleq A\).
\elem
\bproof
 Denote by \(F\) the restriction of \(\STOu{\sigma}\) onto \(\uhp\).
 Because of \rlem{B.L.STOR0}, then \(F\in\RFuq{0}\) and \(\sigmaFa{\R}=\sigma(\Omega)\).
 Taking additionally into account \rlem{B.L.rn} and the assumptions, we get, for all \(z\in\uhp\), then \(\ran{F(z)}\subseteq\ran{A}\) and furthermore \(\frac{1}{\im z}\im F(z)-\ek{F(z)}^\ad A^\mpi\ek{F(z)}\in\Cggq\).
 By virtue of \rlem{L.AEP}, hence the matrix \(\smat{A&F(z)\\\ek{F(z)}^\ad&\frac{1}{\im z}\im F(z)}\)is \tnnH{} for all \(z\in\uhp\).
 The application of \rlem{CRDFK06.L8-9} thus yields \(\sigmaFa{\R}\lleq A\), implying the assertion.
\eproof

\section{Particular pairs of matrices}\label{A.S.cp}
 Regular pairs of matrices considered in this section implement the extension of the set of complex matrices by corresponding points at infinity analogous to the transition from the affine to the projective space.
 In this sense, they can be thought of as homogeneous coordinates.
 
\bdefnl{A.D.cp}
 An ordered pair \(\copa{P}{Q}\) of complex \tpqa{matrices} \(P\) and \(Q\) is called \emph{\tcp{p}{q}}.
 Such a pair is said to be \emph{regular} if \(\rank\tmat{P\\ Q}=q\) and \emph{proper} if \(\rank Q=q\).
\edefn

 Each \tcp{p}{q} \(\copa{P}{Q}\) generates a linear relation \(R\defeq\setaca{(Qv,Pv)}{v\in\Cq}\) in the \(\C\)\nobreakdash-vector space \(\Cp\).
 In accordance with that, we associate to each \tcp{p}{q} \(\copa{P}{Q}\) the linear subspaces \(\gracp{P}{Q}\), \(\domcp{P}{Q}\), \(\rancp{P}{Q}\), \(\nulcp{P}{Q}\), and \(\mulcp{P}{Q}\).
 Obviously, we have \(\nulcp{P}{Q}\subseteq\domcp{P}{Q}\) and \(\mulcp{P}{Q}\subseteq\rancp{P}{Q}\).
 Consequently, we get the inequalities \(0\leq\rankcp{P}{Q}\leq q\).

\bleml{ab.L1532}
 Let \(\copa{P}{Q}\) be a \tcp{p}{q}.
 Then
\[
 \dim\Domcp{P}{Q}+\dim\rk*{\Mulcp{P}{Q}}
 =\dim\Gracp{P}{Q}
 =\dim\Rancp{P}{Q}+\dim\rk*{\Nulcp{P}{Q}}.
\]
\elem
\bproof
 Let \(U\defeq\nul{P}\) and \(V\defeq\nul{Q}\).
 The mappings \(\phi\colon U\to\Cp\) and \(\psi\colon V\to\Cp\) defined by \(\phi(u)\defeq Qu\) and \(\psi(v)\defeq Pv\), \tresp{}, are \(\C\)\nobreakdash-linear with \(\ker\phi=U\cap V\), \(\phi(U)=Q\rk{U}\), \(\ker\psi=V\cap U\), and \(\psi(V)=P\rk{V}\).
 Regarding \(\dim\ker\phi+\dim\phi(U)=\dim U\) and \(\dim\ker\psi+\dim\psi(V)=\dim V\), then \(\dim U=\dim(U\cap V)+\dim\rk{\nulcp{P}{Q}}\) and \(\dim V=\dim(U\cap V)+\dim\rk{\mulcp{P}{Q}}\)
 follow.
 The application of \rrem{A.R.RNT} yields \(\dim\Rancp{P}{Q}+\dim U=q\) and \(\dim\Domcp{P}{Q}+\dim V=q\).
 By virtue of \rrem{ab.L1044}, we see from \rrem{A.R.RNT} that
\(
 \dim\gracp{P}{Q}+\dim(U\cap V)
 =q
\).
 Taking all together, we obtain
\begin{align*}
 \dim\Domcp{P}{Q}+\dim\rk*{\Mulcp{P}{Q}}&=q-\dim(U\cap V)=\dim\Gracp{P}{Q}
\shortintertext{and}
 \dim\Rancp{P}{Q}+\dim\rk*{\Nulcp{P}{Q}} &=q-\dim(U\cap V)=\dim\Gracp{P}{Q}.\qedhere
\end{align*}
\eproof

 A generalization of \rrem{A.R.RNT} for \tcp{p}{q}s immediately follows from \rlem{ab.L1532}:

\breml{A.R.RNTcp}
 The equation \(\rankcp{P}{Q}+\dim\rk{\nulcp{P}{Q}}=\dim\domcp{P}{Q}\) holds true for each \tcp{p}{q} \(\copa{P}{Q}\).
\erem

\bleml{A.L.rg=q}
 Let \(\copa{P}{Q}\) be a \tcp{p}{q} with \(\rankcp{P}{Q}=q\).
 Then \(\rank P=q\), \(\rank Q=q\), \(\nulcp{P}{Q}=\set{\Ouu{p}{1}}\), and \(\mulcp{P}{Q}=\set{\Ouu{p}{1}}\).
 In particular, the pair \(\copa{P}{Q}\) is proper.
\elem
\bproof
 By assumption, we have
\(
 q
 =\rankcp{P}{Q}
 \leq\dim\Rancp{P}{Q}
 =\rank P
 \leq q
\).
 Consequently, \(\Mulcp{P}{Q}=\set{\Ouu{p}{1}}\) and \(\rank P=q\).
 Thus, \(\nul{P}=\set{\Ouu{q}{1}}\) and hence \(\nulcp{P}{Q}=\set{\Ouu{p}{1}}\).
 We infer from \(\Mulcp{P}{Q}=\set{\Ouu{p}{1}}\) furthermore \(\nul{Q}\subseteq\nul{P}\), implying \(\nul{Q}=\set{\Ouu{q}{1}}\).
 Therefore, \(\rank Q=q\), \tie{}, \(\copa{P}{Q}\) is proper.
\eproof

 Each proper \tcp{p}{q} \(\copa{P}{Q}\) is necessarily regular with \(\mulcp{P}{Q}=\set{\Ouu{p}{1}}\) and \(\rankcp{P}{Q}=\dim\rancp{P}{Q}\).
 Furthermore, using \rremsss{A.R.RNT}{ab.R1842x}{ab.L1044}, the following result is readily checked:

\breml{A.R.PQre}
 Let \(\copa{P}{Q}\) be a \tcp{p}{q}.
 Then \(\copa{P}{Q}\) is proper if and only if \(\det\rk{Q^\ad Q}\neq0\).
 Furthermore, the following statements are equivalent:
\baeqi{0}
 \il{A.R.PQre.i} \(\copa{P}{Q}\) is regular.
 \il{A.R.PQre.ii} \(\det\rk{P^\ad P+Q^\ad Q}\neq0\).
 \il{A.R.PQre.iv} \(\nul{P}\cap\nul{Q}=\set{\Ouu{q}{1}}\).
\eaeqi
\erem

 Using \rrem{A.R.rgcolreg}, we obtain furthermore:

\breml{A.R.PQV}
 Let \(\copa{P}{Q}\) be a \tcp{p}{q} and let \(V\in\Coo{q}{s}\).
 Let \(\phi\defeq PV\) and let \(\psi\defeq QV\).
 Then \(\copa{\phi}{\psi}\) is a \tcp{p}{s} fulfilling \(\psi^\ad\phi=V^\ad\rk{Q^\ad P}V\).
 If \(\rank V=s\) and \(\copa{P}{Q}\) is regular (\tresp{}, proper), then \(\copa{\phi}{\psi}\) is regular (\tresp{}, proper).
\erem

 It is readily checked that by the following definition an equivalence relation on the set of \tcp{p}{q}s is given:

\bdefnl{ab.N1142}
 Two \tcp{p}{q}s \(\copa{P}{Q}\) and \(\copa{S}{T}\) are said to be \emph{equivalent} if \(\gracp{P}{Q}=\gracp{S}{T}\).
 In this case, we write \(\copa{P}{Q}\cpaeq\copa{S}{T}\).
 Furthermore, denote by \(\cpcl{P}{Q}\) the corresponding equivalence class of a \tcp{p}{q} \(\copa{P}{Q}\).
\edefn

\breml{ab.R1457}
 \rremss{ab.L1418}{A.R.PQre} show that two regular \tcp{p}{q}s \(\copa{P}{Q}\) and \(\copa{S}{T}\) are equivalent if and only if there is an \(R\in\Cqq\) with \(\det R\neq0\) fulfilling \(S=PR\) and \(T=QR\).
\erem

\breml{A.R.cpesim}
 Each proper \tcp{q}{q} \(\copa{P}{Q}\) satisfies \(\det Q\neq0\) and \(\copa{P}{Q}\cpaeq\copa{PQ^\inv}{\Iq}\).
\erem

 Consequently, the set of equivalence classes of proper \tcp{q}{q}s can be identified with the set of complex \tqqa{matrices} by means of \(\cpcl{P}{Q}\mapsto A\defeq PQ^\inv\), where \(\domcp{P}{Q}=\Cq\), \(\rancp{P}{Q}=\ran{A}\), \(\nulcp{P}{Q}=\nul{A}\), and \(\mulcp{P}{Q}=\set{\Ouu{q}{1}}\).

 In the remaining part of this section, we are concerned with reducing certain \tcp{q}{q}s \(\copa{P}{Q}\), which satisfy a condition of the form \(\rancp{P}{Q}\subseteq\ran{M}\) with a given complex \tqpa{matrix} \(M\) of rank \(r\geq1\), to \tcp{r}{r}s \(\copa{\phi}{\psi}\) without loosing any information:

\bleml{A.L.detB}
 Let \(\theta\in\C\) with \(\abs{\theta}=1\) and let \(\copa{P}{Q}\) be a \tcp{q}{q}.
 Let \(A_\theta\defeq Q+\theta P\) and let \(B_\theta\defeq Q-\theta P\).
 Suppose that \(\det B_\theta\neq0\) and let \(K_\theta\defeq A_\theta B_\theta^\inv\).
 Then
 \begin{align}
  \ran{P}&=\ran{\Iq-K_\theta},&&&\nul{P}&=B_\theta^\inv\rk*{\nul{\Iq-K_\theta}},\label{A.L.detB.A}\\
  \ran{Q}&=\ran{\Iq+K_\theta},&&&\nul{Q}&=B_\theta^\inv\rk*{\nul{\Iq+K_\theta}},\label{A.L.detB.B}\\
  \Nulcp{P}{Q}&=\nul{\Iq-K_\theta},&&\text{and}&\Mulcp{P}{Q}&=\nul{\Iq+K_\theta}.\label{A.L.detB.D}
 \end{align}
 Furthermore,
\(
  \rankcp{P}{Q}=\rank\rk{\Iq-K_\theta}+\rank\rk{\Iq+K_\theta}-q
\)
 and \(\copa{P}{Q}\) is regular.
\elem
\bproof
We have \(K_\theta B_\theta=A_\theta\) and hence \(\rk{\Iq+K_\theta}B_\theta=B_\theta+A_\theta=2Q\) and \(\rk{\Iq-K_\theta}B_\theta=B_\theta-A_\theta=-2\theta P\).
 Thus, we can easily infer \eqref{A.L.detB.B} and \eqref{A.L.detB.A}, using \rrem{A.R.rnLAR}.
 From \eqref{A.L.detB.A} and \eqref{A.L.detB.B} we obtain
\begin{align*}
 \Nulcp{P}{Q}
 &=\rk{Q-\theta P}\rk*{\nul{P}}
 =B_\theta B_\theta^\inv\rk*{\nul{\Iq-K_\theta}}
 =\nul{\Iq-K_\theta}
\shortintertext{and}
 \Mulcp{P}{Q}
 &=\rk{Q-\theta P}\rk*{\nul{Q}}
 =B_\theta B_\theta^\inv\rk*{\nul{\Iq+K_\theta}}
 =\nul{\Iq+K_\theta},
 \end{align*}
 \tie{}, \eqref{A.L.detB.D}.
 Due to \rrem{A.R.RNT}, we have
\begin{align}\label{A.L.detB.2}
 \dim\ran{\Iq+K_\theta}+\dim\nul{\Iq+K_\theta}&=q,&
 \dim\ran{\Iq-K_\theta}+\dim\nul{\Iq-K_\theta}&=q.
\end{align}
 Taking into account \eqref{A.L.detB.A} and \eqref{A.L.detB.D}, we conclude from the first equation in \eqref{A.L.detB.2} then
\[\begin{split}
 \rankcp{P}{Q}
 &=\dim\ran{\Iq-K_\theta}-\dim\nul{\Iq+K_\theta}\\
 &=\rank\rk{\Iq-K_\theta}+\rank\rk{\Iq+K_\theta}-q.
\end{split}\]
 \rlem{ab.L1532} yields \(\dim\gracp{P}{Q}=\dim\rancp{P}{Q}+\dim\nulcp{P}{Q}\).
 Because of \eqref{A.L.detB.A}, \eqref{A.L.detB.D}, and the second equation in \eqref{A.L.detB.2}, we infer \(\rank\tmatp{P}{Q}=q\), \tie{}, \(\copa{P}{Q}\) is regular.
\eproof

 We think that the following result is well-known.
 However, we did not succeed in finding an available reference.
 
\blemnl{\tcf{}~\zitaa{Thi06}{\clem{1.6}{15}}}{ab.L1508}
 Let \(\copa{P}{Q}\) be a regular \tcp{q}{q} satisfying \(\im\rk{Q^\ad P}\in\Cggq\).
 Let \(A\defeq Q+\iu P\) and let \(B\defeq Q-\iu P\).
 Then \(\det B\neq0\) and the matrix \(K\defeq AB^\inv\) satisfies \(\normS{K}\leq1\).
\elem
\bproof
 We have
\[
 A^\ad A
 =\rk{Q^\ad-\iu P^\ad}\rk{Q+\iu P}
 =Q^\ad Q+\iu\rk{Q^\ad P-P^\ad Q}+P^\ad P
 =Q^\ad Q+P^\ad P-2\im\rk{Q^\ad P}
\]
 and
\[
 B^\ad B
 =\rk{Q^\ad+\iu P^\ad}\rk{Q-\iu P}
 =Q^\ad Q-\iu\rk{Q^\ad P-P^\ad Q}+P^\ad P
 =Q^\ad Q+P^\ad P+2\im\rk{Q^\ad P}.
\]
 In view of \rremss{A.R.kK}{A.R.XAX}, in particular \(B^\ad B\lgeq Q^\ad Q+P^\ad P\lgeq\Oqq\) follows.
 Using \rrem{ab.R1842x} and \rlem{A.R.rA<rB}, we infer then \(\nul{B}=\nul{B^\ad B}\subseteq\nul{Q^\ad Q+P^\ad P}\).
 From \rrem{A.R.PQre} we see furthermore \(\det\rk{P^\ad P+Q^\ad Q}\neq0\).
 Consequently, \(\nul{B}=\set{\NM}\), implying \(\det B\neq0\).
 Regarding \(KB=A\), we have moreover
\(
 B^\ad\rk{\Iq-K^\ad K}B
 =B^\ad B-A^\ad A
 =4\im\rk{Q^\ad P}
\).
 Taking into account \rremss{A.R.kK}{A.R.XAX}, we can conclude then \(\Iq-K^\ad K=4B^\invad\im\rk{Q^\ad P}B^\inv\lgeq\Oqq\).
 Thus, the application of \rrem{DFK.L1-1-12} yields \(\normS{K}\leq1\).
\eproof

\bleml{A.L.cpblo}
 Assume \(r\leq q\).
 Let \(U\in\Coo{q}{r}\) with \(U^\ad U=\Iu{r}\) and let \(\copa{\phi}{\psi}\) be an \tcp{r}{r}.
 Let \(P\defeq U\phi U^\ad\) and let \(Q\defeq U\psi U^\ad+\OPu{\ek{\ran{U}}^\orth}\).
 Then \(\copa{P}{Q}\) is a \tcp{q}{q} with \(\rancp{P}{Q}\subseteq\ran{U}\) fulfilling \(\det\rk{P^\ad P+Q^\ad Q}=\det\rk{\phi^\ad\phi+\psi^\ad\psi}\) and \(Q^\ad P=U\rk{\psi^\ad \phi}U^\ad\).
 In particular, \(\copa{P}{Q}\) is regular if and only if \(\copa{\phi}{\psi}\) is regular.
\elem
\bproof
 Observe that \(\rancp{P}{Q}=\ran{U\phi U^\ad}\subseteq\ran{U}\).
 Furthermore, we have
\(
 P^\ad P
 =U\phi^\ad U^\ad U\phi U^\ad
 =U\phi^\ad \phi U^\ad
\).
 Let \(N\defeq\OPu{\ek{\ran{U}}^\orth}\).
 By virtue of \rrem{R.P}, then \(\nul{N}=\ran{U}\), implying \(NU=\NM\).
 Consequently, we infer
\[%
 Q^\ad Q
 =U\psi^\ad U^\ad U\psi U^\ad+U\psi^\ad\rk{NU}^\ad+NU\psi U^\ad+N ^\ad N
 =U\psi^\ad \psi U^\ad+N
\]
 and
\[
 Q^\ad P
 =U\psi^\ad U^\ad U\phi U^\ad+N^\ad U\phi U^\ad
 =U\psi^\ad \phi U^\ad+NU\phi U^\ad
 =U\psi^\ad \phi U^\ad.
\]
 We are now going to show that
\beql{A.L.cpblo.0}
 \det\rk{P^\ad P+Q^\ad Q}
 =\det\rk{\phi^\ad\phi+\psi^\ad\psi}
\eeq
 holds true.
 Observe that, because of \rrem{A.R.PQre}, the asserted equivalence immediately follows from \eqref{A.L.cpblo.0}.
 Using the already shown identities, we get
\beql{A.L.cpblo.1}
 P^\ad P+Q^\ad Q
 =U\rk{\phi^\ad \phi+\psi^\ad \psi}U^\ad+N.
\eeq 
 First assume \(r=q\).
 Then the matrix \(U\) is unitary, implying \(N=\NM\).
 Thus, \eqref{A.L.cpblo.0} is a consequence of \eqref{A.L.cpblo.1}.
 
 Now we consider the case \(r<q\).
 Then there exists some \(V\in\Coo{q}{(q-r)}\) such that \(W\defeq\mat{U,V}\) is a unitary \tqqa{matrix}.
 In particular, we get
\begin{align}\label{A.L.cpblo.2}
 W^\ad W
 &=
 \bMat
  U^\ad U&U^\ad V\\
  V^\ad U&V^\ad V
 \eMat
 =
 \bMat
  \Iu{r}&\Ouu{r}{(q-r)}\\
  \Ouu{(q-r)}{r}&\Iu{q-r}
 \eMat&
&\text{and}&
 WW^\ad
 &=UU^\ad+VV^\ad
 =\Iq. 
\end{align}
 Because of \(U^\ad U=\Iu{r}\) and \rrem{A.R.P=UU*}, we have \(UU^\ad=\OPu{\ran{U}}\).
 In view of \rrem{R.P}, we obtain from the last equation in \eqref{A.L.cpblo.2} thus \(VV^\ad=\Iq-UU^\ad=N\).
 Taking additionally into account \eqref{A.L.cpblo.1} and \eqref{A.L.cpblo.2}, we hence can infer
\[
 W^\ad\rk{P^\ad P+Q^\ad Q}W
 =\matp{U^\ad}{V^\ad}\ek*{U\rk{\phi^\ad\phi +\psi^\ad\psi}U^\ad+VV^\ad}\mat{U,V}
 =
 \begin{pmat}[{|}]
  \phi^\ad\phi+\psi^\ad\psi&\Ouu{r}{(q-r)}\cr\-
  \Ouu{(q-r)}{r}&\Iu{q-r}\cr
 \end{pmat}.
\]
 In particular, \eqref{A.L.cpblo.0} holds true.
\eproof

\bpropnl{\tcf{}~\zitaa{MR1395706}{\clem{4.3}{270}}}{A.L.cpred}
 Let \(M\in\Cqp\) with rank \(r\geq1\), let \(u_1,u_2,\dotsc,u_r\) be an orthonormal basis of \(\ran{M}\), and let \(U\defeq\mat{u_1,u_2,\dotsc,u_r}\).
 Let \(\copa{P}{Q}\) be a regular \tcp{q}{q} fulfilling \(\im\rk{Q^\ad P}\in\Cggq\) and \(\rancp{P}{Q}\subseteq\ran{M}\) and let \(B\defeq Q-\iu P\).
 Then \(\det B\neq0\).
 Let \(\phi\defeq U^\ad PB^\inv U\) and let \(\psi\defeq U^\ad QB^\inv U\).
 Then \(\copa{\phi}{\psi}\) is a regular \tcp{r}{r} satisfying \(\psi^\ad \phi=\rk{B^\inv U}^\ad\rk{Q^\ad P}\rk{B^\inv U}\).
 Let \(S\defeq U\phi U^\ad\) and let \(T\defeq U\psi U^\ad+\OPu{\ek{\ran{M}}^\orth}\).
 Then \(\copa{S}{T}\) is a regular \tcp{q}{q} satisfying \(\det\rk{S^\ad S+T^\ad T}=\det\rk{\phi^\ad\phi+\psi^\ad\psi}\) and \(T^\ad S=B^\invad\rk{Q^\ad P}B^\inv\).
 Furthermore, \(\copa{P}{Q}\cpaeq\copa{S}{T}\) with \(S=PB^\inv\) and \(T=QB^\inv\).
\eprop
\bproof
 We only consider here the case \(r<q\).
 Then there exists some \(V\in\Coo{q}{(q-r)}\) such that \(W\defeq\mat{U,V}\) is a unitary \tqqa{matrix}.
 In particular, we get \eqref{A.L.cpblo.2}.
 Let \(A\defeq Q+\iu P\).
 From \rlem{ab.L1508} we see \(\det B\neq0\) and that \(K\defeq AB^\inv\) satisfies \(\normS{K}\leq1\).
 Consequently the matrix \(L\defeq W^\ad KW\) then satisfies \(\normS{L}\leq1\) as well.
 Furthermore, we have
\[
 L
 =\matp{U^\ad}{V^\ad}K\mat{U,V}
 =
 \bMat
  U^\ad KU&U^\ad KV\\
  V^\ad KU&V^\ad KV
 \eMat.
\]
 Obviously, \(B-A=-2\iu P\) and \(B+A=2Q\) hold true.
 Consequently, we obtain
\begin{align}\label{A.L.cpred.2}
 \frac{\iu}{2}\rk{\Iq-K}&=PB^\inv&
&\text{and}&
 \frac{1}{2}\rk{\Iq+K}&=QB^\inv.
\end{align}
 \rrem{A.R.P=UU*} yields \(UU^\ad=\OPu{\ran{M}}\).
 Hence, \(UU^\ad P=P\) follows.
 With \eqref{A.L.cpred.2} and \eqref{A.L.cpblo.2}, we can thus conclude
\[
 V^\ad\rk{\Iq-K}
 =-2\iu V^\ad PB^\inv
 =-2\iu V^\ad UU^\ad PB^\inv
 =\Ouu{(q-r)}{q},
\]
 implying \(V^\ad K=V^\ad\).
 Regarding \eqref{A.L.cpblo.2}, we infer for the lower blocks of \(L\) then \(V^\ad KU=V^\ad U=\Ouu{(q-r)}{r}\) and \(V^\ad KV=V^\ad V=\Iu{q-r}\).
 In particular, the lower right block \(V^\ad KV\) of \(L\) is unitary.
 Consequently, the application of \rrem{A.L.KAB0D} to \(L\) yields the \tbr{} 
\beql{A.L.cpred.4}
 L
 =
 \bMat
  U^\ad KU&\Ouu{r}{(q-r)}\\
  \Ouu{(q-r)}{r}&\Iu{q-r}
 \eMat.
\eeq
 First we verify the assertions for the pair \(\copa{S}{T}\):
 Using \eqref{A.L.cpblo.2} and \eqref{A.L.cpred.4}, we obtain
\begin{align*}
 \Iq&=UU^\ad+VV^\ad=UU^\ad UU^\ad+VV^\ad,&
 K&=WLW^\ad=\mat{U,V}L\matp{U^\ad}{V^\ad}=UU^\ad KUU^\ad+VV^\ad
\end{align*}
 and, consequently, 
\(
 \Iq-K=UU^\ad\rk{\Iq-K}UU^\ad\) as well as \(
 \Iq+K=UU^\ad\rk{\Iq+K}UU^\ad+2VV^\ad\).
 Because of \eqref{A.L.cpred.2}, then
\begin{align*}
 PB^\inv&=UU^\ad PB^\inv UU^\ad=U\phi U^\ad&
&\text{and}&
 QB^\inv&=UU^\ad QB^\inv UU^\ad+VV^\ad=U\psi U^\ad+VV^\ad
\end{align*}
 follow.
 In view of \eqref{A.L.cpblo.2} and \(UU^\ad=\OPu{\ran{M}}\), we infer from \rrem{R.P} furthermore \(VV^\ad=\Iq-UU^\ad=\OPu{\ek{\ran{M}}^\orth}\).
 Thus, we get \(PB^\inv=S\) and \(QB^\inv=T\).
 By virtue of \rrem{A.R.PQV}, hence \(\copa{S}{T}\) is a regular \tcp{q}{q} fulfilling \(T^\ad S=B^\invad \rk{Q^\ad P}B^\inv\).
 In addition, \rrem{ab.R1457} yields \(\copa{P}{Q}\cpaeq\copa{S}{T}\).
 It remains to show the assertions involving the pair \(\copa{\phi}{\psi}\):
 Regarding \(U^\ad U=\Iu{r}\) and \(\ran{U}=\ran{M}\), we can apply \rlem{A.L.cpblo} to the \tcp{r}{r} \(\copa{\phi}{\psi}\) to obtain \(\det\rk{S^\ad S+T^\ad T}=\det\rk{\phi^\ad\phi+\psi^\ad\psi}\) and to see that \(\copa{\phi}{\psi}\) is regular and that \(T^\ad S=U\rk{\psi^\ad\phi}U^\ad\) holds true.
 From the last equation we can infer then
\[
 \psi^\ad\phi
 =U^\ad U\psi^\ad\phi U^\ad U
 =U^\ad T^\ad S U
 =U^\ad\rk{QB^\inv}^\ad\rk{PB^\inv}U
 =\rk{B^\inv U}^\ad\rk{Q^\ad P}\rk{B^\inv U}.\qedhere
\]
\eproof

\section{Linear fractional transformations of matrices}\label{A.s1.lft}
 In this appendix, we consider a matricial generalization of the transformation \(z\mapsto\frac{az+b}{cz+d}\) of the extended complex plane.
 We thereby follow~\zitaa{MR1152328}{\csec{1.6}}, while restricting ourselves to the version \(Z\mapsto\rk{AZ+B}\rk{CZ+D}^\inv\) with denominator on the right side.
 Let \(\tmat{A & B\\ C & D}\) be the \tbr{} of a complex \taaa{(p+q)}{(p+q)}{matrix} \(M\) with \tppa{block} \(A\).
 If the set
\begin{align*}
 \lftfdb{C}{D}&\defeq\setaca*{Z\in\Cpq}{\det\rk{CZ+D}\neq0}\\
 \text{(\tresp{}, }\lftpdb{C}{D}&\defeq\setaca*{(P,Q)\in\xx{\Cpq}{\Cqq}}{\det\rk{CP+DQ}\neq0}\text{)}
\end{align*}
 is non-empty, then let the \emph{linear fractional transformation} \(\lftfpq{M}\colon\lftfdb{C}{D}\to\Cpq\) (\tresp{}, \(\lftppq{M}\colon\lftpdb{C}{D}\to\Cpq\)) be defined by
\begin{align*}
 \lftfpqa{M}{Z}&\defeq\rk{AZ+B}\rk{CZ+D}^\inv&
&\text{(\tresp{}, }&
 \lftppqA{M}{\copa{P}{Q}}&\defeq\rk{AP+BQ}\rk{CP+DQ}^\inv\text{)}.
\end{align*}
 In this context, the block matrix \(M=\tmat{A & B\\ C & D}\) is called the \emph{generating matrix} of the linear fractional transformation.
 For each matrix \(Z\in\lftfdb{C}{D}\), we obviously have \(\rk{Z,\Iq}\in\lftpdb{C}{D}\) and \(\lftfpqa{M}{Z}=\lftppqA{M}{\rk{Z,\Iq}}\).
 We first characterize the case, that the corresponding domain is non-empty:

\blemnl{\zitaa{MR3765778}{\clem{D.2}{104}}}{M.L.gltL161}%
 The following statements are equivalent:
 \baeqi{0}
  \il{M.L.gltL161.i} \(\lftfdb{C}{D}\neq\emptyset\).
  \il{M.L.gltL161.ii} \(\lftpdb{C}{D}\neq\emptyset\).
  \il{M.L.gltL161.iii} \(\rank\mat{C,D}=q\).
 \eaeqi
\elem

 The composition of two linear fractional transformations is again a linear fractional transformation with generating matrix \(M\) emerging from ordinary matrix multiplication:

\bpropnl{\tcf{}~\zitaa{MR3765778}{\cpropss{D.3}{104}{D.4}{105}}}{M.P.gltP163}%
 Let \(\tmat{A_1 & B_1\\ C_1 & D_1}\) and \(\tmat{A_2 & B_2\\ C_2 & D_2}\) be the \tbr{s} of two given complex \taaa{(p+q)}{(p+q)}{matrices} \(M_1\) and \(M_2\) with \tppa{block} \(A_1\) and \(A_2\), \tresp{}
 Let \(\tmat{A & B\\ C & D}\) be the \tbr{} of the product \(M\defeq M_2M_1\) with \tppa{block} \(A\).
\benui
 \il{M.P.gltP163.a} Suppose that the set
\(
 \mathcal{Q}
 \defeq\setaca{Z\in\lftfdb{C_1}{D_1}}{\lftfpqa{M_1}{Z}\in\lftfdb{C_2}{D_2}}
\)
 is non-empty.
 Then \(\mathcal{Q}\subseteq\lftfdb{C}{D}\) and \(\lftfpqa{M}{Z}=\lftfpqa{M_2}{\lftfpqa{M_1}{Z}}\) for all \(Z\in\mathcal{Q}\).
 \il{M.P.gltP163.b} Suppose that the set
\(
 \mathcal{PQ}
 \defeq\setaca{\copa{P}{Q}\in\lftpdb{C_1}{D_1}}{\lftppqA{M_1}{\copa{P}{Q}}\in\lftfdb{C_2}{D_2}}
\)
 is non-empty.
 Then \(\mathcal{PQ}\subseteq\lftpdb{C}{D}\) and \(\lftppqa{M}{\copa{P}{Q}}=\lftfpqa{M_2}{\lftppqa{M_1}{\copa{P}{Q}}}\) for all \(\copa{P}{Q}\in\mathcal{PQ}\).
 \eenui
\eprop

 In connection with the particular embedding of \tcp{r}{r}s into the class of \tcp{q}{q}s for \(r\leq q\) considered in \rlem{A.L.cpblo}, the following auxiliary result is of interest:

\bleml{A.L.red}
 Suppose \(q\geq2\) and let \(r\in\mn{1}{q-1}\).
 Let \(\mat{U,V}\) be the \tbr{} of a unitary \tqqa{matrix} \(W\) with \taaa{q}{r}{block} \(U\).
 Let \(\tmat{A & B\\ C & D}\) be the \tbr{} of a complex \taaa{2q}{2q}{matrix} \(M\) with \tqqa{block} \(A\) and let \(N\defeq\tmat{AW & BW\\ CW & DW}\).
\benui
 \il{A.L.red.a} Let \(f\in\Coo{r}{r}\) and let \(F\defeq\zdiag{f}{\Ouu{\rk{q-r}}{\rk{q-r}}}\).
 Then \(UfU^\ad\in\lftfdb{C}{D}\) if and only if \(F\in\lftfdb{CW}{DW}\).
 In this case, \(\lftfooa{q}{q}{M}{UfU^\ad}=\lftfooa{q}{q}{N}{F}\).
 \il{A.L.red.b} Let \(f,g\in\Coo{r}{r}\), let \(F\defeq\zdiag{f}{\Ouu{\rk{q-r}}{\rk{q-r}}}\), and let \(G\defeq\zdiag{g}{\Iu{q-r}}\).
 Then \(\rk{ UfU^\ad,UgU^\ad+\OPu{\ek{\ran{U}^\orth}}}\in\lftpdb{C}{D}\) if and only if \(\rk{F,G}\in\lftpdb{CW}{DW}\).
 In this case, \(\lftpooa{q}{q}{M}{\rk{UfU^\ad,UgU^\ad+\OPu{\ek{\ran{U}^\orth}}}}=\lftpooa{q}{q}{N}{\copa{F}{G}}\).
\eenui
\elem
\bproof
 \eqref{A.L.red.a} Because of \(WFW^\ad=UfU^\ad\) and \(W^\inv=W^\ad\), we have \(\ek{\rk{AW}F+\rk{BW}}W^\inv=A\rk{WFW^\ad}+B=A\rk{UfU^\ad}+B\) and similarly \(\ek{\rk{CW}F+\rk{DW}}W^\inv=C\rk{UfU^\ad}+D\).
 Consequently,~\eqref{A.L.red.a} follows.

 \eqref{A.L.red.b} As in the proof of \rlem{A.L.cpblo}, we have \eqref{A.L.cpblo.2} and we can conclude \(VV^\ad=\OPu{\ek{\ran{U}}^\orth}\).
 Beside \(WFW^\ad=UfU^\ad\), we get \(WGW^\ad=UgU^\ad+VV^\ad=UgU^\ad+\OPu{\ek{\ran{U}}^\orth}\).
 The equation \(\ek{\rk{AW}F+\rk{BW}G}W^\inv=A\rk{WFW^\ad}+B\rk{WGW^\ad}=A\rk{UfU^\ad}+B\rk{UgU^\ad+\OPu{\ek{\ran{U}}^\orth}}\) then follows from \(W^\inv=W^\ad\).
 Similarly, we obtain moreover \(\ek{\rk{CW}F+\rk{DW}G}W^\inv=C\rk{UfU^\ad}+D\rk{UgU^\ad+\OPu{\ek{\ran{U}}^\orth}}\).
 Consequently,~\eqref{A.L.red.b} follows.
\eproof

\section{Holomorphic matrix-valued functions}\label{B.S.hol}
 Let \(\dom\) be a \emph{domain}, \tie{}, an open, non-empty, and connected subset of \(\C\).
 A matrix-valued function \(F\colon\dom\to\Cpq\) is said to be \emph{holomorphic} if all entries \(f_{jk}\colon\dom\to\C\) of \(F=\mat{f_{jk}}_{\substack{j=1,\dotsc,p\\ k=1,\dotsc,q}}\) are holomorphic functions.
 In this case, the matrix-valued function \(F\) admits, for each \(z_0\in\dom\), a unique power series representation
\(
 F(z)
 =\sum_{n=0}^\infi\rk{z-z_0}^nA_n
\).
 The corresponding disk of convergence coincides with the largest open disk with center \(z_0\) lying entirely in \(\dom\).
 The coefficients \(A_n=\mat{a_{jk,n}}_{\substack{j=1,\dotsc,p\\ k=1,\dotsc,q}}\) are given by the Taylor series
\(
 f_{jk}(z)
 =\sum_{n=0}^\infi a_{jk,n}\rk{z-z_0}^n
\)
 at \(z_0\).
 Setting \(\der{F}{n}\) with the \(n\)\nobreakdash-th derivatives \(\der{f_{jk}}{n}\) of the infinitely differentiable functions \(f_{jk}\), we have
\(
 A_n
 =\frac{1}{n!}\dera{F}{n}{z_0}.
\)
 Basic results on holomorphic functions can be generalized to the matrix case considered here in an appropriate way:

\breml{B.R.Tay*}
 Let \(F\colon\dom\to\Cpq\) be holomorphic, let \(U\in\Coo{r}{p}\), and let \(V\in\Coo{q}{s}\).
 Then \(H\defeq UFV\) is holomorphic with \(\der{H}{n}=U\der{F}{n}V\) for all \(n\in\NO\).
\erem

 The \tCP{} for sequences of matrices determines the coefficients of the product of two matrix-valued power series:
 
\breml{B.R.cp}
 Let \(F\colon\dom\to\Cpq\) and \(G\colon\dom\to\Cqr\) be two holomorphic functions.
 Let \(z\in\dom\) and let the sequences \(\seq{A_n}{n}{0}{\infi}\) and \(\seq{B_n}{n}{0}{\infi}\) be given by \(A_n\defeq\frac{1}{n!}\dera{F}{n}{z}\) and \(B_n\defeq\frac{1}{n!}\dera{G}{n}{z}\), \tresp{}
 Then \(H\defeq FG\) is holomorphic and  the sequence \(\seq{C_n}{n}{0}{\infi}\) given by \(C_n\defeq\frac{1}{n!}\dera{H}{n}{z}\) coincides with the \tCPa{\(\seq{A_n}{n}{0}{\infi}\)}{\(\seq{B_n}{n}{0}{\infi}\)}.
\erem

 If, in the case \(p=q\), the values \(F(z)\) of the holomorphic matrix-valued function \(F\) are invertible matrices for all \(z\in\dom\), then the function \(G\colon\dom\to\Cqq\) defined by \(G(z)\defeq\ek{F(z)}^\inv\) is holomorphic as well.
 Now suppose that \(F\) satisfies only the weaker condition of having constant column space \(\ran{F(z)}\) and constant null space \(\nul{F(z)}\), independent of the argument \(z\in\dom\).
 Then, even in the case \(p\neq q\), the function \(G\colon\dom\to\Cqp\) defined by \(G(z)\defeq\ek{F(z)}^\mpi\) turns out to be holomorphic.
 Furthermore, the sequences of Taylor coefficients of \(G\) and \(F\) both belong to the class introduced in \rnota{M.N.Dpq} below and are mutually reciprocal in the sense of \rdefn{D1419}:

\bnotal{M.N.Dpq}
 Let \(\Dpqkappa\) be the set of all sequences \(\seqska\) of complex \tpqa{matrices} satisfying \(\bigcup_{j=0}^\kappa\ran{s_j}\subseteq\ran{s_0}\) and \(\nul{s_0}\subseteq\bigcap_{j=0}^\kappa\nul{s_j}\).
\enota

 The following is a specification of a result due to Campbell and Meyer~\zitaa{MR1105324}{\cthm{10.5.4}{227}}:

\bpropl{ab.L0850}
 Let \(F\colon\dom\to\Cpq\) be holomorphic.
 Then the following statements are equivalent:
\baeqi{0}
 \il{ab.L0850.i} The function \(G\colon\dom\to\Cqp\) defined by \(G(z)\defeq\ek{F(z)}^\mpi\) is holomorphic.
 \il{ab.L0850.ii} \(\ran{F(z)}=\ran{F(w)}\) and \(\nul{F(z)}=\nul{F(w)}\) for all \(z,w\in\dom\).
 \il{ab.L0850.iii} \(\seq{\frac{1}{n!}\dera{F}{n}{z}}{n}{0}{\infi}\in\Dpqinf\) for all \(z\in\dom\).
 \eaeqi
 If~\rstat{ab.L0850.i} is fulfilled and \(z_0\in\dom\), then \(\seq{\frac{1}{n!}\dera{G}{n}{z_0}}{n}{0}{\infi}\) is exactly the \trFa{\(\seq{\frac{1}{n!}\dera{F}{n}{z_0}}{n}{0}{\infi}\)}. 
\eprop
\bproof
 The equivalence of~\rstat{ab.L0850.i} and~\rstat{ab.L0850.ii} is an immediate consequence of~\zitaa{MR3014197}{\cprop{8.4}{43}}.
 Let~\rstat{ab.L0850.i} be fulfilled.
 Consider an arbitrary \(z_0\in\dom\).
 Because of~\zitaa{MR3014197}{\cthmss{8.9}{45,}{4.21}{26}}, the sequence \(\seq{\frac{1}{n!}\dera{F}{n}{z_0}}{n}{0}{\infi}\) belongs to \(\Dpqinf\) and \(\seq{\frac{1}{n!}\dera{G}{n}{z_0}}{n}{0}{\infi}\) is exactly the \trFa{\(\seq{\frac{1}{n!}\dera{F}{n}{z_0}}{n}{0}{\infi}\)}.
 In particular,~\rstat{ab.L0850.iii} holds true.
 Conversely, suppose that~\rstat{ab.L0850.iii} is fulfilled.
 From~\zitaa{MR3014197}{\cthm{8.9}{45}} we can then infer that the function \(G\) is holomorphic in all points \(z\in\dom\).
 Consequently,~\rstat{ab.L0850.i} holds true.
\eproof

 Next, we give analogous results for power series expansions at \(z_0=\infc\).
 To that end, let \(\rho\in(0,\infp)\) and suppose that the improper open annulus \(\diskc{\rho}\defeq\setaca{z\in\C}{\abs{z}>\rho}\) is entirely contained in \(\dom\).
 Furthermore, let a holomorphic matrix-valued function \(F\colon\dom\to\Cpq\) be given, admitting the series representation
\beql{B.G.Fsum}
 F(z)
 =\sum_{n=0}^\infi z^{-n}C_n
\eeq
 for all \(z\in\diskc{\rho}\) with certain complex \tpqa{matrices} \(C_0,C_1,C_2,\dotsc\)
 This is the matricial version of a special case of the general situation of a given complex-valued function \(f\) which is holomorphic in an annulus \(\mathcal{A}\defeq\setaca{z\in\C}{r<\abs{z-c}<R}\) centered at \(c\in\C\) with radii \(0\leq r<R\leq\infp\).
 As is well known, such a function \(f\) has a Laurent series \(f(z)=\sum_{\ell=-\infty}^\infi a_\ell\rk{z-c}^\ell\) at the point \(c\) converging on \(\mathcal{A}\) with uniquely determined coefficients \(a_\ell\in\C\).
 In the particular situation of interest considered here, we have \(c=0\), \(R=\infp\), and \(a_\ell=0\) for all \(\ell\in\N\).
 This case can be easily reduced to the ordinary power series expansion of holomorphic functions, discussed at the beginning of this section.
 By means of the substitution \(z\mapsto w\defeq1/z\), we can proceed to a holomorphic function \(\Phi\) defined on the open disk \(\disk{1/\rho}\defeq\setaca{w\in\C}{\abs{w}<1/\rho}\) with Taylor series \(\Phi(w)=\sum_{n=0}^\infi w^nC_n\) at the point \(w_0=0\):

\bleml{ab.L1106}
 Let \(F\colon\dom\to\Cpq\) be holomorphic, admitting for all \(z\in\diskc{\rho}\) the series representation \eqref{B.G.Fsum} with certain complex \tpqa{matrices} \(C_0,C_1,C_2,\dotsc\)
 Then \(\lim_{\zeta\to0}F(1/\zeta)=C_0\) and the matrix-valued function \(\Phi\colon\disk{1/\rho}\to\Cpq\) defined by \(\Phi(w)\defeq F(1/w)\) for \(w\neq0\) and by \(\Phi(0)\defeq\lim_{\zeta\to0}F(1/\zeta)\)
 is holomorphic with \(\frac{1}{n!}\dera{\Phi}{n}{0}=C_n\) for all \(n\in\NO\). 
\elem

 We continue with the analogue of \rrem{B.R.Tay*} for power series expansion at \(z_0=\infc\):

\breml{B.R.Lau*}
 Let \(F\colon\dom\to\Cpq\) be holomorphic, admitting the series representation \eqref{B.G.Fsum} for all \(z\in\diskc{\rho}\) with certain complex \tpqa{matrices} \(C_0,C_1,C_2,\dotsc\)
 Let \(U\in\Coo{r}{p}\) and let \(V\in\Coo{q}{s}\).
 Then \(H\defeq UFV\) is holomorphic and \(H(z)=\sum_{n=0}^\infi z^{-n}\rk{UC_nV}\) for all \(z\in\diskc{\rho}\).
\erem

 Likewise, \rrem{B.R.cp} can be modified in a well-known matter for power series expansion at \(z_0=\infc\):

\bleml{ab.L1324}
 Let \(F\colon\dom\to\Cpq\) and \(G\colon\dom\to\Cqr\) be holomorphic functions, admitting the series representations \(F(z)=\sum_{n=0}^\infi z^{-n}C_n\) and \(G(z)=\sum_{n=0}^\infi z^{-n}D_n\) for all \(z\in\diskc{\rho}\) with certain complex \tpqa{matrices} \(C_0,C_1,C_2,\dotsc\) and certain complex \taaa{q}{r}{matrices} \(D_0,D_1,D_2,\dotsc\), \tresp{}
 Let \(H\defeq FG\) and denote by \(\seq{E_n}{n}{0}{\infi}\) the \tCPa{\(\seq{C_n}{n}{0}{\infi}\)}{\(\seq{D_n}{n}{0}{\infi}\)}.
 Then \(H\) is holomorphic and \(H(z)=\sum_{n=0}^\infi z^{-n}E_n\) for all \(z\in\diskc{\rho}\).
\elem

 Using \rprop{ab.L0850}, we are able to expand the function \(z\mapsto\ek{F(z)}^\mpi\) under certain conditions at \(z_0=\infc\) into a series with coefficients given, according to \rdefn{D1419}, by the \trFa{the} sequence \(\seq{C_n}{n}{0}{\infi}\) from \eqref{B.G.Fsum}:

\bleml{ab.L0836}
 Let \(F\colon\dom\to\Cpq\) be holomorphic and let \(\seq{C_n}{n}{0}{\infi}\) be a sequence of complex \tpqa{matrices} such that \eqref{B.G.Fsum} and furthermore \(\ran{F(z)}=\ran{C_0}\) and \(\nul{F(z)}=\nul{C_0}\) hold true for all \(z\in\diskc{\rho}\).
 Let \(G\colon\diskc{\rho}\to\Cqp\) be defined by \(G(z)\defeq\ek{F(z)}^\mpi\) and denote by \(\seq{D_n}{n}{0}{\infi}\) the \trFa{\(\seq{C_n}{n}{0}{\infi}\)}.
 Then \(G\) is holomorphic and \(G(z)=\sum_{n=0}^\infi z^{-n}D_n\) for all \(z\in\diskc{\rho}\).
\elem
\bproof
 According to \rlem{ab.L1106}, we proceed to a holomorphic function \(\Phi\colon\disk{1/\rho}\to\Cpq\), which satisfies \(\frac{1}{n!}\dera{\Phi}{n}{0}=C_n\) for all \(n\in\NO\).
 Consider an arbitrary \(w\in\disk{1/\rho}\).
 If \(w=0\), then \(\Phi(w)=C_0\).
 In the case \(w\neq0\), we see that \(z\defeq1/w\) belongs to \(\diskc{\rho}\) and that \(\Phi(w)=F(z)\).
 Consequently, \(\ran{\Phi(w)}=\ran{C_0}\) and \(\nul{\Phi(w)}=\nul{C_0}\) for all \(w\in\disk{1/\rho}\).
 In particular, \(\ran{\Phi(w)}\) and \(\nul{\Phi(w)}\) are independent of \(w\in\disk{1/\rho}\).
 Let \(\Psi\colon\disk{1/\rho}\to\Cqp\) be defined by \(\Psi(w)\defeq\ek{\Phi(w)}^\mpi\).
 From \rprop{ab.L0850} we see then that \(\Psi\) is holomorphic and that \(\seq{\frac{1}{n!}\dera{\Psi}{n}{0}}{n}{0}{\infi}\) is exactly the \trFa{\(\seq{\frac{1}{n!}\dera{\Phi}{n}{0}}{n}{0}{\infi}\)}.
 Hence, we have \(\frac{1}{n!}\dera{\Psi}{n}{0}=D_n\) for all \(n\in\NO\) and thus \(\Psi(w)=\sum_{n=0}^\infi w^nD_n\) for all \(w\in\disk{1/\rho}\).
 Consider an arbitrary \(z\in\diskc{\rho}\).
 Then \(w\defeq1/z\) belongs to \(\disk{1/\rho}\setminus\set{0}\) and we have \(\Psi(w)=\ek{F(1/w)}^\mpi\), implying
 \[
  G(z)
  =\ek*{F(z)}^\mpi
  =\Psi(w)
  =\sum_{n=0}^\infi w^nD_n
  =\sum_{n=0}^\infi z^{-n}D_n.\qedhere
 \]
\eproof

 In the remaining part of this section, let \(\dom\) be again an arbitrary domain.
 Next, we consider the matricial generalization of a special class of holomorphic functions, which is well studied, especially in the generic case of \(\dom\) being the open unit disk:
 
\bnotal{A.N.SF}
 Denote by \(\SchFpq{\dom}\) the set of all functions \(S\colon\dom\to\Cpq\), which are holomorphic in \(\dom\) and satisfy \(\normS{S(z)}\leq1\) for all \(z\in\dom\).
\enota

 The matrix-valued functions belonging to \(\SchFpq{\dom}\) are called \emph{Schur functions (in \(\dom\))}.

\blemnl{\tcf{}~\zitaa{arXiv:1712.08358}{\clem{3.9}{7}}}{ab.P2028}
 Let \(S\in\SchFpq{\dom}\) and let \(U,V\in\Cpq\) with \(UU^\ad=\Ip\) and \(V^\ad V=\Iq\).
 For all \(z,w\in\dom\), then \(\ran{U+S(z)}=\ran{U+S(w)}\) and \(\nul{V+S(z)}=\nul{V+S(w)}\).
\elem

 We end this section with some remarks concerning meromorphic matrix-valued functions.
 A subset \(\mathcal{D}\) of \(\dom\) is said to be \emph{discrete in \(\dom\)} if \(\dom\) does not contain any accumulation point of \(\mathcal{D}\).
 So, according to the identity theorem for holomorphic functions, two holomorphic functions \(F,G\colon\dom\to\Cpq\) coincide if and only if the set \(\setaca{z\in\dom}{F(z)=G(z)}\) is not discrete in \(\dom\).
 Speaking in the following of a \emph{discrete subset of \(\dom\)}, we always mean a subset of \(\dom\), which is discrete \emph{in \(\dom\)}.
 For such a discrete subset \(\mathcal{D}\) of \(\dom\), the set \(\dom\setminus\mathcal{D}\) is a domain.
 
 A complex-valued function \(f\) is said to be \emph{meromorphic in \(\dom\)} if there exists a discrete subset \(\pol{f}\) of \(\dom\) such that \(f\) is a holomorphic function defined on the domain \(\dom\setminus\pol{f}\), which in each point from \(\pol{f}\) has a pole (of positive order).
 In particular, each holomorphic function \(f\colon\dom\to\C\) is meromorphic in \(\dom\) with \(\pol{f}=\emptyset\).
 We call a \(\Cpq\)\nobreakdash-valued function \(F\) \emph{meromorphic in \(\dom\)} if all entries \(f_{jk}\) of \(F=\mat{f_{jk}}_{\substack{j=1,\dotsc,p\\ k=1,\dotsc,q}}\) are complex-valued functions meromorphic in \(\dom\).
 In this case, the union \(\pol{F}\defeq\bigcup_{j=1}^p\bigcup_{k=1}^q\pol{f_{jk}}\) of the sets of poles of all entries \(f_{jk}\) is a discrete subset of \(\dom\).
 In particular, each holomorphic function \(F\colon\dom\to\Cpq\) is meromorphic in \(\dom\) with \(\pol{F}=\emptyset\).
 Since \(\dom\) is assumed to be connected, the set of complex-valued functions meromorphic in \(\dom\) has the algebraic structure of a field.
 Using the arithmetic of this field, the usual operations from matrix algebra can be formally carried over to matrix-valued functions, which are meromorphic in \(\dom\).
 Thus, corresponding sums and products of such matrix-valued functions are again meromorphic in \(\dom\).
 Furthermore, it is readily checked that the determinant \(\det F\) of a \(\Cqq\)\nobreakdash-valued function \(F\) meromorphic in \(\dom\) is a complex-valued function, which is meromorphic in \(\dom\).
 If \(\det F\) does not identically vanish, then the mapping \(F^\inv\) given by formal matrix inversion of \(F\) (seen as a matrix with entries in the field of complex-valued functions meromorphic in \(\dom\)) is again a \(\Cqq\)\nobreakdash-valued function, which is meromorphic in \(\dom\) with not identically vanishing determinant.

\bibliography{183arxiv}
\bibliographystyle{bababbrv}

\vfill\noindent
\begin{minipage}{0.5\textwidth}
 Universit\"at Leipzig\\
 Fakult\"at f\"ur Mathematik und Informatik\\
 PF~10~09~20\\
 D-04009~Leipzig
\end{minipage}
\begin{minipage}{0.49\textwidth}
 \begin{flushright}
  \texttt{
   fritzsche@math.uni-leipzig.de\\
   kirstein@math.uni-leipzig.de\\
   maedler@math.uni-leipzig.de
  } 
 \end{flushright}
\end{minipage}

\end{document}